\documentclass[11pt]{amsart}

\usepackage{amsmath, amssymb, amscd, mathrsfs, url, pinlabel,setspace,hyperref,verbatim,mathtools}
\usepackage[margin=1.25in,marginparwidth=1in,centering,letterpaper,dvips]{geometry}
\usepackage{color,dcpic,latexsym,graphicx,epstopdf,comment}
\usepackage[all]{xy}
\usepackage{tikz,pgfplots}
\usepackage{tikz-cd,circuitikz}
\usepackage{enumitem,upquote,tabularx,textcomp,setspace,color}
\usepackage{makecell}
\usepackage{caption}
\usepackage[T1]{fontenc}
\xyoption{rotate}
\DeclareMathAlphabet{\mathpzc}{OT1}{pzc}{m}{it}
\SetMathAlphabet{\mathpzc}{bold}{OT1}{pzc}{b}{it}

\newtheorem*{rep@theorem}{\rep@title}
\newcommand{\newreptheorem}[2]{%
\newenvironment{rep#1}[1]{%
 \def\rep@title{#2 \ref{##1}}%
 \begin{rep@theorem}}%
 {\end{rep@theorem}}}
\makeatother

\newreptheorem{proposition}{Proposition}
\newtheorem {theorem}{Theorem}
\newreptheorem{theorem}{Theorem}
\newtheorem {lemma}[theorem]{Lemma}
\newtheorem {proposition}[theorem]{Proposition}
\newtheorem {corollary}[theorem]{Corollary}

\numberwithin{equation}{section}
\numberwithin{theorem}{section}

\theoremstyle{definition}
\newtheorem{definition}[theorem]{Definition}
\newtheorem{construction}[theorem]{Construction}

\newtheorem{remark}[theorem]{Remark}

\newtheorem{example}[theorem]{Example}
\newtheorem*{ack}{Acknowledgement}
\newtheorem*{useai}{AI Declaration}

\newlist{pcases}{enumerate}{1}
\setlist[pcases]{
  label=\bf{Case~\arabic*:}\protect\thiscase.~,
  ref=\arabic*,
  align=left,
  labelsep=0pt,
  leftmargin=0pt,
  labelwidth=0pt,
  parsep=0pt
}
\newcommand{\case}[1][]{%
  \if\relax\detokenize{#1}\relax
    \def\thiscase{}%
  \else
    \def\thiscase{~#1}%
  \fi
  \item
}

\DeclareFontFamily{U}{mathx}{\hyphenchar\font45}
\DeclareFontShape{U}{mathx}{m}{n}{
      <5> <6> <7> <8> <9> <10>
      <10.95> <12> <14.4> <17.28> <20.74> <24.88>
      mathx10
      }{}
\DeclareSymbolFont{mathx}{U}{mathx}{m}{n}
\DeclareFontSubstitution{U}{mathx}{m}{n}
\DeclareMathAccent{\widecheck}{0}{mathx}{"71}

\usetikzlibrary{calc,intersections}
\tikzset{every picture/.style=thick}
\tikzset{link/.style = { white, double = black, line width = 1.75pt, double distance = 1.25pt, looseness=1.75 }}
\tikzset{crossing/.style = {draw, circle, dotted, minimum size=0.5cm, inner sep=0, outer sep=0}}
\pgfplotsset{compat=1.12}

\usepackage{extarrows}
\usepackage{amsbsy}
\usepackage{appendix}
\usepackage{float}
\usepackage{subfigure}
\usepackage[numbers,sort&compress]{natbib}
\usepackage{amsthm}
\usepackage{geometry}
\usepackage{fancyhdr}
\usepackage{overpic}
\usepackage{mathrsfs}
\usepackage{colonequals}
\usepackage{microtype}
\usepackage{diagbox}
\usepackage{tabularx}

\definecolor{lygreen}{HTML}{016646}

\newcommand{\Z}{\mathbb{Z}}
\newcommand{\N}{\mathbb{N}}
\newcommand{\R}{\mathbb{R}}
\newcommand{\C}{\mathbb{C}}
\newcommand{\B}{\mathrm{B}}
\newcommand{\BN}{\mathrm{B}_\bullet\mathbb{N}}

\newcommand{\bfP}{\mathbf{P}}
\newcommand{\bfJ}{\mathbf{J}}
\newcommand{\bfg}{\mathbf{g}}
\newcommand{\F}{\mathbb{F}}
\newcommand{\Q}{\mathbb{Q}}
\newcommand{\sfC}{\mathsf{C}}
\newcommand{\bA}{\mathbb{A}}
\newcommand{\bfa}{\boldsymbol{a}}
\newcommand{\bY}{\mathbb{Y}}
\newcommand{\bW}{\mathbb{W}}
\newcommand{\cL}{\mathcal{L}}
\newcommand{\cA}{\mathcal{A}}
\newcommand{\cB}{\mathcal{B}}
\newcommand{\cC}{\mathcal{C}}
\newcommand{\cD}{\mathcal{D}}
\newcommand{\sX}{\mathscr{X}}
\newcommand{\sY}{\mathscr{Y}}
\newcommand{\sW}{\mathscr{W}}
\newcommand{\sR}{\mathscr{R}}
\newcommand{\cG}{\mathcal{G}}
\newcommand{\cI}{\mathcal{I}}
\newcommand{\cM}{\mathcal{M}}
\newcommand{\cP}{\mathcal{P}}
\newcommand{\bI}{\mathbb{I}}
\newcommand{\bK}{\mathbb{K}}
\newcommand{\fB}{\mathfrak{B}}
\newcommand{\fC}{\mathfrak{C}}
\newcommand{\fD}{\mathfrak{D}}
\newcommand{\D}{\mathsf{D}}
\newcommand{\K}{\mathsf{K}}
\newcommand{\bfq}{\mathbf{q}}

\newcommand{\rN}{\mathrm{N}}
\newcommand{\adj}{\mathfrak{adj}}
\newcommand{\mnd}{\mathfrak{mnd}}
\newcommand{\Mod}{\mathrm{Mod}}
\newcommand{\End}{\mathrm{End}}
\newcommand{\Alg}{\mathrm{Alg}}
\newcommand{\Mnd}{\mathrm{Mnd}}
\newcommand{\Fil}{\mathsf{Fil}}
\newcommand{\sSpaceo}{\mathrm{sSpace}}
\newcommand{\sSpace}{\mathpzc{s}\mathcal{S}\!\mathpzc{pace}}
\newcommand{\Space}{\mathcal{S}\!\mathpzc{pace}}
\newcommand{\Set}{\mathrm{Set}}
\newcommand{\Cat}{\mathrm{Cat}}
\newcommand{\Cati}{\mathsf{Cat}_\infty}
\newcommand{\QCat}{\mathcal{QC}\!\mathpzc{at}}
\newcommand{\sCat}{\mathpzc{s}\mathcal{C}\!\mathpzc{at}}
\newcommand{\sSet}{\mathrm{sSet}}
\newcommand{\sCato}{\mathrm{sCat}}
\newcommand{\Ch}{\mathsf{Ch}}
\newcommand{\Int}{\mathrm{Int}}
\newcommand{\fInt}{\mathfrak{Int}}
\newcommand{\PB}{\mathrm{PBord}}
\newcommand{\fPB}{\mathfrak{PBord}}
\newcommand{\Bord}{\mathrm{Bord}}
\newcommand{\BI}{\mathrm{BI}}
\newcommand{\PBI}{\mathrm{PBI}}
\newcommand{\BIa}{\mathrm{BI}_a}
\newcommand{\Cob}{\mathrm{Cob}}
\newcommand{\CobI}{\mathrm{CobI}}
\newcommand{\Se}{\mathrm{Se}}
\newcommand{\CSe}{\mathrm{CSe}}
\newcommand{\Ndg}{\mathrm{N}^{\mathrm{dg}}}
\newcommand{\Nhc}{\mathrm{N}^{\mathrm{hc}}}
\newcommand{\Fun}{\mathsf{Fun}}
\newcommand{\fDiff}{\mathfrak{D}\mathrm{iff}}
\newcommand{\Man}{\mathrm{Man}}
\newcommand{\cAut}{\mathcal{A}ut}
\newcommand{\Top}{{\rm Top}}

\newcommand{\adm}{\mathrm{adm}}
\newcommand{\com}{\mathrm{com}}
\newcommand{\cadm}{\mathrm{cadm}}
\newcommand{\cyl}{\mathrm{cyl}}
\newcommand{\iin}{\mathrm{in}}
\newcommand{\out}{\mathrm{out}}
\newcommand{\fin}{\mathrm{fin}}
\newcommand{\h}{\mathrm{h}}
\newcommand{\e}{\mathrm{e}}
\newcommand{\proj}{\mathrm{proj}}
\newcommand{\free}{\mathrm{free}}
\newcommand{\sym}{\mathrm{Sym}}
\newcommand{\sd}{\mathrm{sd}}
\newcommand{\triv}{\mathrm{triv}}
\newcommand{\opp}{\operatorname{op}}
\newcommand{\loc}{\operatorname{loc}}
\newcommand{\inr}{\operatorname{int}}
\newcommand{\pt}{\mathrm{pt}}
\newcommand{\sm}{\mathrm{sm}}
\newcommand{\aff}{\mathrm{aff}}
\newcommand{\sub}{\mathrm{sub}}
\newcommand{\kih}{\mathrm{kih}}
\newcommand{\col}{\mathrm{col}}
\newcommand{\met}{\mathrm{met}}
\newcommand{\mmod}{\mathrm{mod}}
\newcommand{\sing}{\mathrm{sing}}

\newcommand{\bu}{\bullet}
\newcommand{\I}{\un{I}}
\newcommand{\ov}[1]{{\overline{#1}}}
\newcommand{\un}[1]{{\underline{#1}}}
\newcommand{\wti}[1]{{\widetilde{#1}}}
\newcommand{\what}[1]{{\widehat{#1}}}
\newcommand{\bbox}{\square}
\newcommand{\bs}{\backslash}
\newcommand{\al}{\alpha}
\newcommand{\be}{\beta}
\newcommand{\ga}{\gamma}
\newcommand{\Ga}{\Gamma}
\newcommand{\de}{\delta}
\newcommand{\ep}{\epsilon}
\newcommand{\epp}{{\epsilon^{\prime}}}
\newcommand{\pa}{\partial}
\newcommand{\ot}{\otimes}
\newcommand{\op}{\oplus}
\newcommand{\aand}{~\mathrm{and}~}
\newcommand{\rmif}{~\mathrm{if}~}
\newcommand{\rmfor}{~\mathrm{for}~}
\newcommand{\xra}{\xrightarrow}

\DeclareMathOperator{\disc}{disc}
\DeclareMathOperator{\Hom}{Hom}
\DeclareMathOperator{\Map}{Map}
\DeclareMathOperator{\id}{id}
\DeclareMathOperator{\fib}{fib}
\DeclareMathOperator{\cofib}{cofib}
\DeclareMathOperator{\Tcofib}{tot.cofiber}
\DeclareMathOperator{\glue}{glue}
\DeclareMathOperator*{\hocolim}{hocolim}
\DeclareMathOperator*{\colim}{colim}
\DeclareMathOperator{\ev}{ev}
\DeclareMathOperator{\coev}{coev}
\DeclareMathOperator{\Funo}{Fun}
\DeclareMathOperator{\Fg}{Fg}
\DeclareMathOperator{\diag}{diag}
\newcommand{\pr}{\operatorname{pr}}
\DeclareMathOperator{\Aut}{Aut}
\newcommand{\tr}{\operatorname{tr}}

\newcommand{\bpf}{\begin{proof}}
\newcommand{\epf}{\end{proof}}
\newcommand{\bthm}{\begin{theorem}}
\newcommand{\ethm}{\end{theorem}}
\newcommand{\bprop}{\begin{proposition}}
\newcommand{\eprop}{\end{proposition}}
\newcommand{\bcor}{\begin{corollary}}
\newcommand{\ecor}{\end{corollary}}
\newcommand{\blem}{\begin{lemma}}
\newcommand{\elem}{\end{lemma}}
\newcommand{\bdefn}{\begin{definition}}
\newcommand{\edefn}{\end{definition}}
\newcommand{\bcons}{\begin{construction}}
\newcommand{\econs}{\end{construction}}
\newcommand{\bexmp}{\begin{example}}
\newcommand{\eexmp}{\end{example}}
\newcommand{\brem}{\begin{remark}}
\newcommand{\erem}{\end{remark}}
\newcommand{\benu}{\begin{enumerate}}
\newcommand{\benum}{\begin{enumerate}[leftmargin=*]}
\newcommand{\eenu}{\end{enumerate}}

\newcommand{\Def}{Definition }
\newcommand{\Prop}{Proposition }
\newcommand{\Cor}{Corollary }
\newcommand{\Defs}{Definitions }

\title{An \texorpdfstring{$\infty$}{infinity}-categorical TQFT from instantons}
\author{Fan Ye}
\address{School of Mathematical Sciences\\Peking University}
\email{flyye@math.pku.edu.cn}

\makeatletter

\begin{document}
\begin{abstract}

In this paper, we lift the instanton TQFT from ordinary categories to the $\infty$-categorical setting by constructing a functor $CI:\mathrm{BI}^{\partial}\longrightarrow\mathsf{D}$ from an $\infty$-cobordism category $\mathrm{BI}^{\partial}$ for instantons to an $\infty$-derived category $\mathsf{D}$ of $2$-periodic chain complexes. The construction of $\mathrm{BI}^{\partial}$ is a modification of the $\infty$-cobordism category $\mathrm{Bord}_4$ introduced by Lurie and by Calaque--Scheimbauer using complete Segal spaces. The target $\infty$-category $\mathsf{D}$ is obtained via the dg nerve of the dg category of $2$-periodic chain complexes of finitely presented projective $\mathbb{Z}$-modules.

The information encoded by the functor $CI$ was already developed by Kronheimer--Mrowka using families of metrics on cobordisms. Our reinterpretation in terms of $\infty$-categories conceptually simplifies the construction of the hypercube complexes underlying the link spectral sequence. This viewpoint also clarifies hyperbox complexes and compression in Manolescu--Ozsv\'{a}th's work, and possible $A_\infty$-structures, especially their signs and higher coherence.

In addition, we lift the generalized cap product $\mu$-operators in instanton Floer homology to the chain level and construct explicit homotopies and higher homotopies establishing the commutativity of multiple $\mu$-operators in even degrees.

\end{abstract}
\maketitle
\tableofcontents
\section{Introduction}

Topological quantum field theory, or TQFT for short (also sometimes called topological field theory, or TFT), was introduced by Atiyah 
\cite{atiyah1988tqft} as a rigid symmetric monoidal functor from the cobordism category to the category of vector spaces. Recall that the cobordism category $\Cob_n$ (cf.\ Definition \ref{defn: cob cat}) consists of closed oriented smooth $(n-1)$-dimensional manifolds as objects and smooth oriented $n$-dimensional cobordisms as morphisms, for a fixed $n\in\mathbb{N}_+$. 

Throughout this paper, manifolds are smooth, Hausdorff, and second countable. Unless otherwise specified, they have no boundary; manifolds with boundary or corners are indicated explicitly.

The symmetric monoidal structure (cf.\ Definition \ref{defn: monoidal category}) can be viewed roughly as a tensor product on the category, which for $\Cob_n$ is induced by the disjoint union. The rigid structure 
(cf.\ Definition \ref{defn: rigidity}) means that every object is dualizable with respect to suitable evaluation and coevaluation morphisms. For $\Cob_n$, the dual of an object is given by orientation reversal, while the corresponding evaluation and coevaluation morphisms are given by trace and cotrace cobordisms, obtained by suitably modifying the ends of product cobordisms.

Frobenius algebras provide fundamental examples of TQFTs for 
$2$-dimensional cobordisms. For $4$-dimensional cobordisms, further examples of TQFTs arise from Floer theories, either through gauge theory, including the Yang--Mills and Seiberg--Witten equations (giving rise to instanton and monopole theories, respectively), or through symplectic geometry, such as Heegaard Floer theory. The image of an object under the corresponding TQFT functor is usually called its Floer homology, while the image of a morphism is called a cobordism map. For each Floer theory, the cobordism category requires suitable 
modifications; see Li--Ye \cite[\S 2.1]{LY2021} for more details.

We mainly focus on instanton Floer theory and follow the setup of 
Kronheimer--Mrowka \cite[\S 5.1]{kronheimer2011khovanov}, restricting to the case of empty singular submanifolds. We first describe the construction at the level of ordinary categories to motivate the corresponding $\infty$-categorical setting. The instanton cobordism category $\CobI$ consists of objects
\[
\bY=(Y,P,\phi,J)
\]
and morphisms
\[
\bW=(W,\mathbf{P},\Phi,\mathbf{J}):\bY_0\longrightarrow\bY_1,
\]
as described below. The details of the additional data on manifolds can be safely skipped on a first reading.
\begin{itemize}
    \item $Y$ is a (possibly disconnected and empty) closed oriented $3$-manifold;
    \item $P$ is a principal $SO(3)$-bundle over $Y$;
    \item $\phi$ is a subgroup of $H^1(Y;\Z/2)$;
    \item $J$ is an almost complex structure on $\R\times Y$ that induces the orientation of $Y$;
    \item $\bY$ should satisfy the \emph{non-integral} condition, i.e.\ there exists an embedded oriented surface $\Sigma\subset Y$ such that, for each component $Y_\al$ of $Y$, the second Stiefel--Whitney class $w_2(P)$ evaluates nontrivially on $\Sigma_\al=\Sigma\cap Y_\al$ and every element of $\phi$ restricts trivially to $\Sigma_\al$. We write these conditions as\[w_2(P)|_{\Sigma_\al}\neq 0\aand \phi|_{\Sigma_\al}=0~\text{for every }\al.\]Only the existence of $\Sigma$, rather than a specific choice, is required. By convention, the empty manifold satisfies the non-integral condition.
    \item $W:Y_0\to Y_1$ is an oriented $4$-dimensional cobordism;
    \item $\bfP$ is a principal $SO(3)$-bundle over $W$;
    \item $\Phi$ is a subgroup of $H^1(W;\Z/2)$;
    \item $\bfJ$ is an almost complex structure on $W$ that induces the orientation of $W$;
    \item The data $(\bfP,\Phi,\bfJ)$ restrict to those of $\bY_i$ for $i=0,1$. Two such cobordisms represent the same morphism if there is a diffeomorphism with a bundle isomorphism preserving $\Phi$ and the boundary identifications, under which their almost complex structures are homotopic relative to the boundary collars. For simplicity, we also call $\bW$ a \emph{cobordism}.
    \item The composition $\bW_1\circ \bW_0$ of $\bW_i=(W_i,\bfP_i,\Phi_i,\bfJ_i)$ for $i=0,1$ is obtained by composing $(W_i,\bfP_i,\bfJ_i)$ in the usual way and taking the largest subgroup $\Phi\subset H^1(W;\Z/2)$ whose restriction to $W_i$ is $\Phi_i$.
\end{itemize}
The principal bundles are used to formulate the instanton (anti-self-dual, or ASD) equations, and the subgroups specify the enlarged gauge groups. The almost complex structures determine the signs of the cobordism maps and the absolute $\Z/2$-gradings. Compatibility with the structures on the ends makes these maps grading-preserving. See \S \ref{subsect: Sketch of gauge theory construction} for the gauge-theoretic construction.

The symmetric monoidal structure on $\CobI$ is induced by disjoint union, with direct sums of the subgroups and unit the empty manifold. For $\bY=(Y,P,\phi,J)$, the pull-back subgroup on $[0,1]\times Y$ has boundary image $\diag\phi=\{(a,a):a\in\phi\}$, whereas the disjoint union $Y\sqcup (-Y)=\partial \big([0,1]\times Y)\big)$ carries $\phi\oplus\phi$. Thus the usual trace and cotrace cylinders define morphisms with the required boundary data only when $\phi=0$. Nevertheless, orientation reversal defines a duality functor, sending $\bY$ to $\bY^\vee=(-Y,P,\phi,J)$ and a cobordism $\bW:\bY_0\to\bY_1$ to the same decorated $4$-manifold regarded as $\bW^\vee:\bY_1^\vee\to\bY_0^\vee$. We therefore use categories with duality as in Definition \ref{defn: dual}, without requiring rigidity. For results about trace, cotrace, and self-gluing cobordisms, see Proposition \ref{prop: J trace self gluing}.

Note that a morphism between empty manifolds consists of a closed 
$4$-manifold $X$ equipped with additional data. The corresponding instanton invariant is a \emph{Donaldson invariant}. In the standard metric-independent setup, one assumes that $b_2^+(X)>1$ \cite{DK1992instanton,donaldson2002floer}. Since $b_2^+$ is not additive under composition in general, we restrict our attention to a subcategory of $\CobI$ that excludes closed components of cobordisms to obtain a well-defined functor.

The construction of the suitable subcategory is inspired by 
\cite[\S 6]{GTMW2009cob} (see also \cite[Definition Sketch 4.2.10]{lurie2009cob}). It has the same objects as $\CobI$ and consists of morphisms with \emph{positive boundaries}, namely, cobordisms whose connected components all have nonempty outgoing ends. Equivalently, the natural map $\pi_0(Y_1)\longrightarrow \pi_0(W)$ is surjective. Morphisms with positive boundaries are clearly closed under 
composition. We denote this subcategory by $\CobI^{\pa_+}$.

Similarly, we denote by $\CobI^{\pa_-}$ the subcategory consisting of 
morphisms with \emph{negative boundaries}, namely, those for which $\pi_0(Y_0)\longrightarrow \pi_0(W)$ is surjective. We further write $\CobI^{\pa}$ for the subcategory consisting of morphisms with both 
positive and negative boundaries, also called morphisms with \emph{full boundaries}. Note that cotrace cobordisms for $\phi=0$ are contained only in $\CobI^{\pa_+}$, whereas trace cobordisms for $\phi=0$ are contained only in $\CobI^{\pa_-}$. 

With the above construction, we can write the instanton functor as
\begin{equation}\label{eq: instanton TQFT}
\begin{aligned}
I:\CobI^{\pa}&\longrightarrow\Mod_2^{\fin}(\Z),\\
\bY&\longmapsto I(\bY),\\
\bW&\longmapsto I(\bW).
\end{aligned}
\end{equation}
Here $\Mod_2^{\fin}(\Z)$ is the category of finitely presented $\Z/2$-graded $\Z$-modules and grading-preserving homomorphisms\footnote{In the first arXiv version of this paper, we defined $\Mod_2^{\h,\fin}(\mathbb{Z})$ to include all not necessarily grading-preserving homomorphisms and required the image of $I$ to consist of homogeneous morphisms. This is unnecessary because cobordism maps are grading-preserving with respect to the $\Z/2$-grading from the almost complex structures.}, with the usual shift operator $[1]$. We can also replace the source by $\CobI^{\pa_\pm}$, meaning either $\CobI^{\pa_+}$ or $\CobI^{\pa_-}$. We use the same notation $I$ for all three functors. To make $I$ into a symmetric monoidal functor preserving the duality, we need to replace 
$\Mod_2^{\fin}(\mathbb{Z})$ by $\Mod_2^{\fin}(\mathbb{F})$ for a field 
$\mathbb{F}$, due to the Tor terms arising from the K\"{u}nneth formula.

The main result of this paper is to lift the above TQFT to an 
$\infty$-categorical framework. The idea of viewing TQFTs as functors between $\infty$-categories was introduced by Lurie \cite{lurie2009cob}. More generally, Lurie considered the $(\infty,n)$-category of cobordisms, which incorporates manifolds of lower dimensions, and studied functors between higher categories, known as \emph{extended} TQFTs. This perspective played a fundamental role in the formulation of the Baez--Dolan cobordism hypothesis \cite{BD1995tqft}. Here, we recall that the term $\infty$-category is an abbreviation for $(\infty,1)$-category. For readers who are unfamiliar with $\infty$-categories, standard references include \cite[Chapter 1]{HTTlurie2009}, \cite[Chapter 1]{HAlurie2017}, \cite[\S 1-2]{lurie2009cob}, \cite{Rezk2001completion,JT2007quasisegal}, and the survey articles \cite{Bergner2010survey1,Bergner2020survey2}.

Although these further developments are not needed in this paper, they provide important motivation for adopting the $\infty$-categorical viewpoint and may also be useful for future constructions of instanton analogues of bordered Floer theory \cite{LOT2018bordered} and cornered Floer theory \cite{DM2014cornered1,DLM2019cornered2}; see \S \ref{subsec: further direction} for more details. Lurie's work was subsequently refined and generalized by Calaque--Scheimbauer \cite{CS2019category}, Schommer-Pries \cite{SP2024tqft}, Steinebrunner \cite{Ste21tqft}, and Grady--Pavlov \cite{GP2023geopre,GP2022geocob}. See \cite[\S 1.1, \S 1.2]{GP2022geocob} for further surveys of related 
developments.

We will explain further motivations for introducing $\infty$-categories into instanton theory after stating the main theorem.
\bthm[{Proposition \ref{prop: homotopy cat BI} and Theorem \ref{thm: TQFT functor}}]\label{thm: main}
There is a duality-preserving symmetric monoidal functor between $\infty$-categories\begin{equation}\label{eq: functor CI}
    CI:\BI^{\pa}\to \D=\D_2(\Mod(\Z)^\fin),
\end{equation}such that the induced functor between the homotopy categories recovers \eqref{eq: instanton TQFT} after taking homology. Here $\BI^{\pa}$ is the $\infty$-cobordism category for instantons from Definition \ref{defn: BI}, and $\D$ is the $\infty$-derived category of finitely presented $2$-periodic chain complexes over $\Z$ from Definition \ref{defn: derived}. 

The same construction applies to $\BI^{\pa_\pm}$, whose homotopy categories are equivalent to $\CobI^{\pa_\pm}$, compatibly with the inclusions
\[
\BI^\pa\hookrightarrow\BI^{\pa_\pm},
\]with duality switching the signs in the sense of Definition \ref{defn: dual}.
\ethm
\brem
The notation $CI$ comes from the fact that the functor assigns to each object its instanton Floer chain complex. A primary motivation for using $\infty$-categories is to retain the chain-level information of instanton theory rather than passing to homology. A direct advantage of working at the chain level is that we can work over $\Z$ instead of a field $\F$ for the symmetric monoidal and dual structures.
\erem
\brem
We only use the setup of Kronheimer--Mrowka \cite{kronheimer2011khovanov} with empty singular submanifolds, while they studied instanton theory for closed $3$-manifolds with singular links and $4$-dimensional cobordisms with singular surfaces. We expect the functor $CI$ in Theorem \ref{thm: main} to extend to those cases, but the construction of the $\infty$-cobordism category with singular submanifolds needs more work. It is related to Baez--Dolan's tangle hypothesis \cite{BD1995tqft}; see also \cite[\S 4.4]{lurie2009cob} and \cite{SP2024tqft,Ste21tqft,AF2024tangle}.
\erem
\brem
It is natural to ask whether the $\infty$-categorical enhancement extends to other Floer theories of closed $3$-manifolds and $4$-dimensional cobordisms; see Li--Ye \cite[\S 2.1]{LY2021} for the ordinary categorical setup. 

As a gauge-theoretic theory, monopole Floer theory is expected to behave similarly to instanton Floer theory and to satisfy naturality of all orders. The main subtlety arises for cobordisms involving disconnected $3$-manifolds, where one must carefully choose the appropriate version of monopole Floer homology for each component and orient the relevant moduli spaces; see \cite[Remark 2.10]{QY2026triangle}. Once these issues are addressed, monopole Floer theory should also admit an $\infty$-categorical enhancement. Alternatively, one may restrict to the sub-$\infty$-category of connected closed $3$-manifolds and cobordisms between them. 

The situation for Heegaard Floer theory is more delicate. At present, only first-order naturality is available \cite{Juhasz2012,Gartner2023naturality,AM2025sign}, and cobordisms involving disconnected $3$-manifolds require additional ribbon-graph data \cite{Zemke2019}. Higher-order naturality may require further structure, such as a framing at the basepoint \cite{HHSZ23naturality}. Consequently, an $\infty$-categorical enhancement of Heegaard Floer theory appears less accessible at present.
\erem

We review some steps in the construction of the instanton functor $I$, which motivate the constructions of $\BI^{\pa}$ and $CI$. Here $\BI^\pa$ is the sub-$\infty$-category of $\BI$ consisting of morphisms with full boundaries\footnote{In the first arXiv version of this paper, we used $\BI$ to denote the $\infty$-category with full boundary condition.}. We first explain the ordinary categorical construction and then indicate how the additional homotopical information is incorporated in the $\infty$-categorical setting.

To construct $I$, we first consider an auxiliary cobordism category $\CobI_a^{\partial}$, whose objects and morphisms are of the form $(\bY,a)$ and $(\bW,\bfa)$, respectively, where $a$ and $\bfa$ denote auxiliary data including Riemannian metrics and holonomy perturbations and satisfy some admissibility conditions. We do need the boundary condition because the admissibility condition for closed $4$-manifolds is not for perturbations but for Riemannian metrics. These data are necessary for defining the Yang--Mills equations and constructing the instanton chain complexes. Following \cite{floer1988instanton,donaldson2002floer,kronheimer2011khovanov}, the functor $I$ is obtained from a functor
\begin{equation}\label{eq: Ia}
    I_a:\CobI_a^{\pa}\to \Mod_2^{\fin}(\mathbb{Z}).
\end{equation}

To eliminate the auxiliary data in $\CobI_a^{\pa}$, one considers the category $\Mod_2^{\fin}(\mathbb{Z})/\mathrm{Can}$ of modules up to canonical isomorphisms (cf.\ the discussion before \cite[Corollary 23.1.7]{kronheimer2007monopoles}), or, equivalently, the category of transitive systems of modules (cf.\ \cite[Definition 1.1]{Juhasz2012}). The functor $I_a$ in \eqref{eq: Ia} induces a functor
\[
I':\CobI^{\pa}\to \Mod_2^{\fin}(\mathbb{Z})/\mathrm{Can},
\]
and the functor $I$ in \eqref{eq: instanton TQFT} is obtained by composing $I'$ with the standard functor
\[
\Mod_2^{\fin}(\mathbb{Z})/\mathrm{Can}\longrightarrow\Mod_2^{\fin}(\mathbb{Z})
\]
given by taking colimits.

Recall that the composition law is given by
\begin{equation}\label{eq: composition law}
    I((\bW_1,\bfa_1)\circ(\bW_0,\bfa_0))
    =
    I(\bW_1,\bfa_1)\circ I(\bW_0,\bfa_0),
\end{equation}
where we omit the subscript $a$ in $I_a$ for simplicity. The fact that different choices of $a$ for a fixed $\bY$ are related by canonical isomorphisms follows from the special case of the composition law applied to product cobordisms with different auxiliary data on the incoming and outgoing ends.

To prove the composition law, one considers a $1$-parameter family of metrics on $W_1\circ W_0$ parametrized by $[0,\infty]$, obtained by stretching the metric in a neighborhood of the outgoing end $Y_1$ of $W_0$. In particular, the metric at $0$ is the one associated with $\bfa_1\circ\bfa_0$, while the metric at $\infty$ is the broken metric along $Y_1$ corresponding to $\bfa_1\sqcup\bfa_0$. An extension theorem for perturbations implies that the perturbations at $0$ and $\infty$ can also be extended to a family of perturbations over $[0,\infty]$. Through instanton theory, these families of metrics and perturbations induce a chain homotopy between the chain maps
\[
CI((\bW_1,\bfa_1)\circ(\bW_0,\bfa_0))
\quad\text{and}\quad
CI(\bW_1,\bfa_1)\circ CI(\bW_0,\bfa_0).
\]
Passing to homology yields the composition law \eqref{eq: composition law}. Note that only the \emph{existence} of the chain homotopy is needed; the homotopy itself is not recorded.

The main advantage of the $\infty$-categorical framework is that it records homotopies and higher homotopies. As in the ordinary categorical case, to construct the functor $CI$ in \eqref{eq: functor CI}, we construct an auxiliary $\infty$-cobordism category $\BI_a^{\pa}$, whose homotopy category is $\CobI_a^{\pa}$, together with a functor of $\infty$-categories
\begin{equation*}
    CI_a:\BI_a^{\pa}\to \D .
\end{equation*}

In particular, the $2$-morphisms in $\BI_a^{\pa}$ are given by two composable cobordisms together with a family of perturbations over $[0,\infty]$. We do not include the family of metrics explicitly, since it is canonically determined by the metrics on the two cobordisms through a standard procedure. The image of a $2$-morphism in $\BI_a^{\pa}$ under $CI_a$ is precisely the above chain homotopy.

More generally, the $k$-morphisms in $\BI_a^{\pa}$ are given by $k$ composable cobordisms together with a family of perturbations over $[0,\infty]^{k-1}$. This family corresponds to the family of metrics obtained by stretching neighborhoods of the outgoing ends of the first $k-1$ cobordisms. The images of these higher morphisms were used in the construction of the hypercube complexes associated to link surgeries and the link spectral sequence in \cite{kronheimer2011khovanov,scaduto2015instantons} (see also \cite{Bloom2011link} for analogous results in monopole theory). Consequently, the hypercube construction lifts naturally to $\BI_a^{\pa}$ or $\BI^{\pa}$. We state this formally in \S \ref{subsec: Cube complexes}.

The passage from $CI_a$ to $CI$ only involves the source $\infty$-category. More precisely, the $k$-morphisms in $\BI^{\pa}$ are given by $k$ composable cobordisms, and the extension theorem for perturbations implies that the forgetful functor
\[
\BI_a^{\pa}\longrightarrow\BI^{\pa}
\]
is an equivalence of $\infty$-categories (cf.\ Proposition \ref{prop: BIa=BI}). Therefore, we obtain $CI$ from $CI_a$ by precomposing with the inverse equivalence.

Although we have described the $k$-morphisms of $\BI_a^{\pa}$ and $\BI^{\pa}$, it remains nontrivial to verify that these data define $\infty$-categories. There are many equivalent models of $\infty$-categories. Following Lurie \cite{lurie2009cob} and Calaque--Scheimbauer \cite{CS2019category}, we use complete Segal spaces as concrete models for $\BI^{\pa}$ and $\BI_a^{\pa}$. Following Lurie's book \emph{Higher Algebra} \cite{HAlurie2017}, we use quasicategories as a concrete model for the $\infty$-derived category $\D$, which is obtained from the dg nerve of a suitable dg category. To construct $CI_a$, we regard $\D$ as a complete Segal space and construct a functor between complete Segal spaces. To describe the symmetric monoidal structure explicitly, we also use simplicial categories as an auxiliary model of $\infty$-categories \cite[\S 1.1.5]{HTTlurie2009}.
\begin{ack}
The author thanks Longke Tang for helpful discussions and patient explanations that motivated this work, as well as for writing Appendix \S \ref{sec: Triangle detection lemma}. The author is also grateful to Peter Kronheimer for discussions on the $\mu$-operators in \S \ref{sec: Generalized cap product}, to Michael Hopkins for introducing $\infty$-category theory while he was at Harvard, to Jacob Lurie for explaining the cobordism category during his visit to the Institute for Advanced Study, to Ian Zemke for conversations on hypercubes, to Yonghwan Kim for pointing out the reference \cite{Varolgunes2021cube}, and to Daren Chen, Hongjian Yang, and Sungkyung Kang for comments on an earlier version of the paper. The author further thanks John A. Baldwin, Zhenkun Li, and Steven Sivek for discussions regarding upcoming collaborative work, which partly motivated this work. The author was partially supported by Simons Collaboration Grant \#271133 through Peter Kronheimer while this work was in preparation. 
\end{ack}
\begin{useai}
The author used OpenAI Codex with the \texttt{gpt-5.6-sol}, \texttt{gpt-6-astra}, and \texttt{gpt-6.1-sol} models to expand details in the first arXiv version and to assist with proofreading and language editing. The author takes full responsibility for the manuscript.
\end{useai}

\section{Sketch of more constructions}
This section continues the introduction. We describe in more detail the instanton-theoretic constructions involving $\infty$-categories that are mentioned in the abstract. Finally, we outline the remaining sections.
\subsection{Cube complexes and exact triangles}\label{subsec: Cube complexes}
The main advantage of $\infty$-categories is that they retain the higher-homotopy data used in the construction of hypercube complexes in instanton and monopole theories \cite{kronheimer2011khovanov,scaduto2015instantons,Bloom2011link}. We first review the definition of hypercube complexes from \cite[\S 5]{MO2022link}, \cite[\S 2.2]{Zemke23bordered}, and \cite{BHL2019cancel}, and then describe their $\infty$-categorical analogue.
\bdefn\label{defn: hypercube}
Let $\bK$ be a commutative ring with identity and $2=0$. A \emph{hypercube complex} of dimension $n\in\N$ over $\bK$ has the form

\[
(C,D)=\left(\bigoplus_{\ep\in\{0,1\}^n}C_\ep[|\ep|-n],\sum_{\ep\le\epp\in\{0,1\}^n}D_{\ep,\epp}\right).
\]
where the $C_\ep$ are finitely presented projective $\bK$-modules carrying a $\Z$-homological grading or a $\Z/m$-homological grading for some $m\in\N_+$, and the $D_{\ep,\epp}:C_\ep\to C_\epp$ are module morphisms of degree $|\epp-\ep|-1$, where $|\cdot|$ denotes the number of $1$s. Here $C_\ep[|\ep|-n]_j=(C_\ep)_{j+|\ep|-n}$, so every summand of $D$ has total degree $-1$.

The maps satisfy the compatibility condition
\begin{equation}\label{eq: compatibility condition}
    \sum_{\ep\le \de\le \epp\in\{0,1\}^n}D_{\de,\epp}\circ D_{\ep,\de}=0,\text{ or equivalently } D\circ D=0.
\end{equation}
In particular, $(C,D)$ and $(C_\ep,D_{\ep,\ep})$ are chain complexes over $\bK$.
\edefn
The module morphism $D_{\ep,\epp}$ for $|\epp-\ep|=1$ is a chain map, while those for $|\epp-\ep|>1$ are chain homotopies and higher homotopies. In the $m$-periodic case, the latter are not encoded in the ordinary category $\Mod^{\fin}_m(\bK)$ but in the $\infty$-category $\D_m(\Mod^{\fin}(\bK))$. For general coefficients, the edge chain maps and the total differential are related by the signs in Construction \ref{cons: sign on hypercube}. For $m=2$ and $\bK=\Z$, we use the notation $\D$ introduced in Theorem \ref{thm: main}.

In this case, an $\infty$-categorical lift of a hypercube complex of dimension $n$ is a functor\[\sX:\bbox^n=(\Delta^1)^n\to \D,\]where $\Delta^1=\{0\to 1\}$ is the standard $1$-simplex. We call such a functor a \emph{cube complex}.
The product $(\Delta^1)^n$ has higher simplices, and suitable sums of their images recover the module morphisms $D_{\ep,\epp}$. We can further obtain the corresponding total chain complex via a total cofiber functor\[\Tcofib:\Fun(\bbox^n,\D)\to \D,\]which is defined by iterated cofibers on each $\Delta^1$ factor, where $\Fun$ denotes the functor $\infty$-category. The details of these constructions will be provided in \S \ref{subsec: Cubes of chain complexes}.
\brem
Conversely, a hypercube complex might be lifted to a cube complex using $\infty$-categorical techniques. We leave the details to future work.
\erem
\brem\label{rem: hyperbox}
As a generalization of a cube complex, an $\infty$-categorical version of the \emph{hyperbox complexes} of size $(d_1,\dots,d_n)$ introduced in \cite[\S 5]{MO2022link} and \cite[\S 2.2]{Zemke23bordered} is a functor
\[\prod_{i=1}^n\Delta^{d_i}\to\D.\]
Its compression is given by restriction
\[\bbox^n=\prod_{i=1}^n\Delta^1\xra{\prod_{i=1}^ns^{d_i}}\prod_{i=1}^n\Delta^{d_i}\to\D,\]
where $s^d:\Delta^1\to\Delta^d$ is induced by $(0<1)\mapsto(0<d)$. This formulation incorporates the signs through the dg nerve, whereas the hyperbox construction in \cite[\S 5]{MO2022link} is carried out in characteristic $2$.
\erem

By Theorem \ref{thm: main}, postcomposing a functor
\[
\sX:\bbox^n\to\BI^{\pa}
\]
with $CI$ gives a cube complex.

Since the $k$-morphisms in $\BI$ are given by $k$ composable cobordisms without auxiliary data, this construction is direct. In the case of link surgeries, the images of higher morphisms in $\bbox^n$ under $\sX$ are given by the collections of surgery cobordisms corresponding to different orders of composition. More generally, we have the following theorem.

\bthm[{Construction \ref{cons: cube of cobordisms}}]\label{thm: cube complex cob}
For $i\in \un{n}=\{1,\dots,n\}$, suppose that $\bW_i:\bY\to \bY_i$ are cobordisms supported in disjoint regions (cf.\ Definition \ref{defn: supported region}). Then there exists an $n$-cube of cobordisms\[\sX=\sX[\bY,\{\bW_i\}_{i\in\un{n}}]:\bbox^n\to \BI^{\pa}\]that sends $\{0\}^{i-1}\times \Delta^1\times \{0\}^{n-i}$ to $\bW_i$.
\ethm
Finally, the surgery exact triangle and the iteration procedure for link surgeries in \cite[\S 5-6]{scaduto2015instantons} also admit an $\infty$-categorical lift. Recall from \cite[Definition 1.1.1.9]{HAlurie2017} that, in a stable $\infty$-category, a fiber square is equivalent to a cofiber square, and vice versa. Moreover, by \cite[\S 1.1.2]{HAlurie2017}, a fiber square in a stable $\infty$-category induces a distinguished triangle in its homotopy category. In our setting, the $\infty$-category $\D$ is stable, and therefore $\Fun(\bbox^n,\D)$ is also stable. Using the triangle detection lemma for stable $\infty$-categories proved by Longke Tang in Appendix \S\ref{sec: Triangle detection lemma}, the surgery exact triangle and the iteration procedure become fiber squares in $\Fun(\bbox^n,\D)$. The details will be provided in \S\ref{subsec: Fiber squares for surgery cobordisms}.
\subsection{Enhanced cobordisms and \texorpdfstring{$\mu$}{mu}-operators}\label{subsec: enhanced cob}
For a closed $4$-manifold $X$ with a principal $SO(3)$ bundle $P$ satisfying suitable admissibility conditions to avoid reducible solutions, the Donaldson polynomial invariants introduced in \cite{Donaldson1990polynomial} (see also \cite[\S 5.2 and 9.2]{DK1992instanton}) define a linear map
\begin{equation}\label{eq: A(X) map}
    D_{(X,P)}:\bA(X)=\operatorname{Sym}^*H_{\operatorname{even}}(X;\Q)\otimes \Lambda^*H_{\operatorname{odd}}(X;\Q)\to \Q,
\end{equation}
where $\bA(X)$ is the graded commutative algebra generated by $H_*(X;\Q)$, with the subscripts denoting the even and odd degree parts, respectively.

The construction adapts to instanton cobordism maps and provides linear maps\begin{equation}\label{eq: A(W) map}
   I(\bW;-):I(\bY_0;\Q)\ot \bA(W)\to I(\bY_1;\Q) \end{equation}for a cobordism \[\bW=(W,\bfP,\Phi,\bfJ):\bY_0\to\bY_1\]with trivial subgroup data $\Phi=\phi_0=\phi_1=0$. Related references include \cite[\S 4.2]{BD1995fukaya}, \cite[\S 7.3]{donaldson2002floer}, \cite[\S II.(ii)]{KM1995structure}, \cite[\S 2.3]{Froyshov2002equivariant}, \cite[\S 7.2]{kronheimer2010knots}, and \cite{Donaldson1995mu,munoz1999ring}; see also \cite[\S 25.3 and 26.2]{kronheimer2007monopoles} for the analogous construction in monopole theory.

   For a fixed homogeneous element $z\in\bA(W)$, the corresponding map $I(\bW;z)$ will be called an \emph{enhanced cobordism map}. The identity element $1\in\bA(W)$ induces the original instanton cobordism map in \eqref{eq: instanton TQFT}. For the product cobordism of $\bY$, the enhanced cobordism map is an endomorphism of $I(\bY;\Q)$, called a \emph{generalized cap product} or a \emph{$\mu$-operator} and denoted by
\[
\mu(z):I(\bY;\Q)\to I(\bY;\Q)
\]
for $z\in\bA(I\times Y)=\bA(Y)$. For a homogeneous homology class $z$, the operator $\mu(z)$ has degree $4-\deg z$. These operators are independent of the auxiliary choices and satisfy
\begin{equation}\label{eq: additive z}
\mu(z_1+z_2)=\mu(z_1)+\mu(z_2),~\mu(cz_1)=c\mu(z_1)
\end{equation}
for $z_1,z_2\in\bA(Y)$ and $c\in\Q$. For homogeneous $z_1,z_2\in\bA(Y)$, we also have
\begin{equation}\label{eq: commute}
\mu(z_1)\circ\mu(z_2)=\mu(z_1z_2)=(-1)^{\deg z_1\deg z_2}\mu(z_2)\circ\mu(z_1).
\end{equation}
\brem
The constructions of \eqref{eq: A(X) map} and \eqref{eq: A(W) map}, including their comparison and composition arguments, work over $\Z[1/2]$, since the normalization factors only involve powers of $2$.
\erem
As in the construction of the functor $CI$ in Theorem \ref{thm: main}, we focus on the chain-level construction of \eqref{eq: A(W) map}. In addition to the auxiliary data consisting of Riemannian metrics and holonomy perturbations, we must choose additional data, which we call \emph{enhanced data}. Since it is difficult and technically involved to construct an $\infty$-cobordism category incorporating the enhanced data, we keep track of the dependence on the enhanced data explicitly in our results rather than packaging all the information into a single functor between $\infty$-categories.

For the concrete construction, we adapt the construction for closed $4$-manifolds in \cite[\S II.(ii)]{KM1995structure} to $4$-manifolds with ends. An important advantage of introducing the enhanced data is that it allows us to construct an explicit chain homotopy realizing \eqref{eq: commute} for homogeneous elements of even degree. We restrict our attention to even degrees for technical reasons. For simplicity, we write
\[
\bA_{\e}(X;1/2)=\operatorname{Sym}^*H_{\operatorname{even}}(X;\Z[1/2]).
\]
Our construction also extends to higher homotopies establishing the commutativity of three or more homogeneous elements of even degrees.

To state the result more explicitly, we introduce the following definition.

\bdefn\label{defn: BNr}
For $r\in\N$, let $\BN^r$ be the nerve of the category with one object and endomorphism monoid $\N^r$. Thus $(\BN^r)_n=(\N^r)^n$ for $n\in\N$. By \cite[\texttt{00JC}]{kerodon}, the inclusion $\Delta^1/\partial\Delta^1\hookrightarrow\BN$ induced by $1\in\N$ is inner anodyne. Hence its $r$-fold product is a categorical equivalence
\[(\Delta^1/\partial\Delta^1)^r\xra{\simeq}\BN^r.\]

As in the definition of a cube complex, we define an \emph{endomorphism cube complex} of dimension $r$ to be a functor\[\sX:\BN^r\to \D[1/2],\]where\[\D[1/2]=\D_2(\Mod(\Z[1/2])^\fin)\]is the $\infty$-derived category from Definition \ref{defn: derived}.
\edefn
\bthm[{Construction \ref{cons: mu maps} and Theorem \ref{thm: mu Y}}]\label{thm: mu Y main}
Suppose $\bY=(Y,P,\phi,J)$ is an object in $\BI$ with $\phi=0$ and let $\beta_p$, for $p\in\un{r}$, be homogeneous classes in $H_{\operatorname{even}}(Y;\Z[1/2])$. For a choice of enhanced data, we obtain an endomorphism cube complex
\[
\sX[\bY,\{\beta_p\}_{p\in\un{r}}]:\BN^r\to\D[1/2]
\]
that sends the unique object to $CI(\bY)$ and the $p$-th generating edge to a chain map inducing $\mu(\beta_p)=I(I\times\bY;\beta_p)$ on homology.
\ethm
\brem
For $\beta_p\in H_k(Y;\Z[1/2])$, the degree of $\mu(\beta_p)$ is $4-k$. Since we consider $2$-periodic chain complexes in $\D$, the homogeneous elements in $\bA_{\e}$ correspond to grading-preserving operators.
\erem
\brem
In Theorem \ref{thm: mu cube}, there is a more complicated construction that combines Theorems \ref{thm: cube complex cob} and \ref{thm: mu Y main} and provides a functor\[\BN^r \times \bbox^n\to \D[1/2].\]Given this construction, one might hope to upgrade the functor\[CI:\BI^{\pa}\to \D\]from Theorem \ref{thm: TQFT functor} to a symmetric monoidal functor\[CI^{\mu,r}:\BI^{\mu,r,\pa}\to \Fun(\BN^r,\D[1/2]).\]  One expects a symmetric monoidal structure on $\Fun(\BN^r,\D[1/2])$ whose effect on the generating operators is governed by the Leibniz rule.  More explicitly, for $C,D\in \D[1/2]$ and $S,T$ as operators on $C,D$, respectively, we define\[(C,S)\ot (D,T)=(C\ot D,S\ot \id+\id\ot T).\]Moreover, the $\infty$-category $\BI^{\mu,r,\pa}$ is an expected $\infty$-cobordism category for instantons with the enhanced data. Omitting the enhanced data, an object in $\BI^{\mu,r,\pa}$ is an object in $\BI^{\pa}$ together with $r$ (possibly disconnected and empty) embedded even-dimensional submanifolds that intersect each other transversely, and a morphism in $\BI^{\mu,r,\pa}$ is a morphism in $\BI^{\pa}$ together with cobordisms of those submanifolds that intersect each other transversely. We leave this construction to future work.
\erem
\brem\label{rem: special mu case}
In the following special case, the functor in Theorem \ref{thm: mu Y main} admits the alternative construction described in \cite[Remark 2.2]{GZ2021}, which ultimately originates from the proof of \cite[Corollary 7.2]{kronheimer2010knots}. Suppose that the homogeneous elements $\beta_p\in H_{\operatorname{even}}(Y;\Z[1/2])$ are represented by submanifolds $Z_p$ contained in embedded oriented surfaces $\Sigma_p\subset Y$ such that the second Stiefel--Whitney class $w_2(P)$ of the principal bundle $P$ over $Y$ is nontrivial when restricted to $\Sigma_p$. Here, $Z_p$ is allowed to be the whole surface $\Sigma_p$. We then consider the cobordism
\[
W_p:Y\sqcup S^1\times\Sigma_p\to Y
\]
obtained from $I\times Y$ by removing a neighborhood of $\Sigma_p$. Let $y_p=I(D^2\times\Sigma_p)(1)\in I(S^1\times\Sigma_p)$ be the relative invariant of the removed neighborhood. Then the $\mu$-operator is
\[I(W_p)(-\otimes\mu([Z_p])(y_p)):I(Y)\to I(Y),\]
where $\mu([Z_p])$ acts on $I(S^1\times\Sigma_p)$ and we suppress the bundles and auxiliary structures. If all the surfaces $\Sigma_p$ are pairwise disjoint in $Y$, then the functor in Theorem \ref{thm: mu Y main} can be recovered from the standard hypercube construction of chain complexes as in \cite[\S 4-7]{Bloom2011link}, \cite[\S 5-6]{scaduto2015instantons}, and \cite[\S 6-7]{kronheimer2011khovanov}. When the surfaces $\Sigma_p$ intersect, however, the above approach requires additional modifications, whereas the construction developed in this paper continues to apply. See also Remark \ref{rem: no bubble}. Finally, this construction also applies to elements in $H_1(Y;\mathbb Z[1/2])$ represented by submanifolds of the surfaces $\Sigma_p$. In this case, however, we need to use $\Delta^1$ instead of $\BN$, since the corresponding $\mu$-operators do not preserve the $\Z/2$-grading.
\erem

\subsection{Further direction: \texorpdfstring{$A_\infty$}{A-infinity}-structures}\label{subsec: further direction}

In this subsection, we describe a possible source of $A_\infty$-algebra and module structures in instanton Floer theory. Related structures in $\operatorname{Pin}(2)$- and $S^1$-equivariant monopole Floer theory were studied by Lin \cite{francesco2017connectedsum}, and their Morse-theoretic analogues were studied by Bloom \cite{Bloom2013morse} and Mescher \cite{Mescher2018Morse}. The categorical reformulation below was motivated by discussions with Longke Tang.

The equivalence between the $A_\infty$-, associative, and $E_1$-algebra descriptions follows from \cite[Proposition 4.1.3.19 and Example 5.1.0.7]{HAlurie2017}. For left modules, see \cite[Corollary 4.2.2.16]{HAlurie2017}.

The symmetric monoidal functor $CI:\BI^{\pa_+}\to\D$ of Theorem \ref{thm: main} carries $E_1$-algebras and their modules in $\BI^{\pa_+}$ to $E_1$-algebras and modules in $\D$. We use the positive boundary condition to allow unit cobordisms with empty incoming boundary.

For a monoidal $\infty$-category $\cC$, we write $\Alg(\cC)$ for the $\infty$-category of $E_1$-algebra objects in $\cC$ \cite[Notation 4.1.1.8]{HAlurie2017}. For $A,B\in\Alg(\cC)$, we write $\operatorname{LMod}_A(\cC)$, $\operatorname{RMod}_B(\cC)$, and ${}_{A}\operatorname{BMod}_{B}(\cC)$ for the $\infty$-categories of left $A$-modules, right $B$-modules, and $A$-$B$-bimodules, respectively; see \cite[Definition 4.2.1.13, Variant 4.2.1.36, and Definition 4.3.1.12]{HAlurie2017}. Under this convention, the operations have the form
\[
A^{\otimes(d-1)}\otimes M\longrightarrow M,
\quad
M\otimes B^{\otimes(d-1)}\longrightarrow M,
\quad
A^{\otimes p}\otimes M\otimes B^{\otimes q}\longrightarrow M.
\]
The mixed coherence between the two actions is part of the bimodule data; see \cite[Remark 4.3.1.19]{HAlurie2017}. The left- and right-module theories are exchanged by reversing the monoidal product; see \cite[Remark 4.2.1.37]{HAlurie2017}. If $\cC$ is symmetric monoidal, the reversal involution gives
\[
\operatorname{LMod}_A(\cC)\simeq\operatorname{RMod}_{A^{\mathrm{rev}}}(\cC);
\]
see \cite[Remark 4.6.3.2]{HAlurie2017}. Unless otherwise specified, module means left module below. Lin instead uses the right-module convention; see \cite[\S 1.2]{francesco2017connectedsum}.

Our approach to obtaining $E_1$-algebra and module objects in $\BI^{\pa_+}$ uses $(\infty,2)$-cobordism categories. Calaque--Scheimbauer construct a symmetric monoidal $(\infty,2)$-category $\Bord_4^{(\infty,2)}$ of surfaces, $3$-dimensional bordisms, and $4$-dimensional bordisms; see \cite[\S 6]{CS2019category}. We expect that the construction of $\BI$ in this paper extends to a symmetric monoidal $(\infty,2)$-category $\BI_2$. We further expect that the positive boundary condition on the $4$-dimensional bordisms can be incorporated into a symmetric monoidal $(\infty,2)$-category $\BI_2^{\pa_+}$ satisfying a monoidal equivalence
\[
\End_{\BI_2^{\pa_+}}(\varnothing)\simeq\BI^{\pa_+},
\]
under which horizontal composition corresponds to disjoint union. The construction of these $(\infty,2)$-cobordism categories for instantons is left to future work.

We first describe the expected geometric representatives of the algebra and module operations. Suppose $\bY=(Y,P,\phi,J)$ is an object of $\BI^{\pa_+}$ and that an oriented closed separating surface $\Sigma\subset Y$ determines a collared decomposition
\begin{equation}\label{eq:decomp}
Y=M\cup_\Sigma\bigl(\Sigma\times[-1,1]\bigr)\cup_\Sigma N.
\end{equation}
Assume that every component of $N$ meets $\Sigma$ and that the induced decorations make
\[
N:\varnothing\longrightarrow\Sigma,
\quad
M:\Sigma\longrightarrow\varnothing
\]
$1$-morphisms of $\BI_2^{\pa_+}$ and identify $\bY$ with $M\circ N$. Suppose that $N$ admits a coherent right adjoint $N^r$ in the sense of Definition \ref{defn: coherent-adjunction} below, which is geometrically represented by $-N$ with the orientation and cobordism direction reversed. After rounding corners, the composite $\bY_N=N^r\circ N$ is represented by
\[
Y_0=\partial(N\times I)=\bigl(\{\pm1\}\times N\bigr)\cup_{\{\pm1\}\times\Sigma}(I\times\Sigma).
\]
Thus
\[
Z_0=N\times I:\varnothing\longrightarrow Y_0
\]
is the expected geometric representative of the unit \[\eta:\id_{\varnothing}\Rightarrow N^r\circ N\]of the adjunction. The same product manifold, with outgoing face $\Sigma\times I$, represents the counit
\[
\epsilon:N\circ N^r\Rightarrow\id_\Sigma.
\]

Let $W_1=I\times Y$. For $d\geq2$, choose pairwise disjoint subintervals $I_1',\ldots,I_{d-1}'\subset\operatorname{int}(I)$ and remove collared copies of $I_i'\times N$ from $W_1$. After rounding corners, the resulting manifold is a cobordism
\begin{equation*}
W_d:Y\sqcup\bigsqcup_{i=1}^{d-1}Y_0\longrightarrow Y.
\end{equation*}
If $I''\subset I$ contains the intervals $I_i'$ and $\epsilon>0$ is sufficiently small, the portion of $W_d$ lying over
\[
I''\times\bigl(N\cup_\Sigma\Sigma\times[1-\epsilon,1]\bigr)
\]
gives a cobordism
\begin{equation*}
Z_{d-1}:\bigsqcup_{i=1}^{d-1}Y_0\longrightarrow Y_0.
\end{equation*}
With the convention that $Z_d$ has $d$ incoming copies of $Y_0$, the cobordisms $Z_d$ and $W_d$ should be viewed as candidate geometric representatives of the iterated product and right-module action
\[
\mu_d:\bY_N^{\circ d}\simeq N^r\circ(N\circ N^r)^{\circ(d-1)}\circ N
\xrightarrow{\,N^r\epsilon^{\circ(d-1)}N\,}N^r\circ N=\bY_N,
\quad d\geq1,
\]
and
\[
m_d:\bY\circ\bY_N^{\circ(d-1)}\simeq M\circ(N\circ N^r)^{\circ(d-1)}\circ N
\xrightarrow{\,M\epsilon^{\circ(d-1)}N\,}M\circ N=\bY,
\quad d\geq1,
\]
together with the unit $\eta$, geometrically represented by $Z_0$.

Lin \cite[\S 4]{francesco2017connectedsum} considered the special case $N=B^3$, for which $N\times I=B^4$. In the present setting, it is natural to expect associahedrally coherent families of metrics and perturbations on $Z_d$ and $W_d$ to give chain-level representatives of the product and right-module operations on $CI(\bY_N)$ and $CI(\bY)$. 

Lin works over $\Z/2$, where signs can be ignored \cite[\S 1]{francesco2017connectedsum}. Given $\BI_2^{\pa_+}$, a coherent adjunction $N\dashv N^r$, and the symmetric monoidal functor $CI$, the construction below requires no further choices of metrics, perturbations, or signs, since these are already encoded by $CI$.

We now work in an arbitrary $(\infty,2)$-category $\cC$, with no monoidal structure required. For background, see \cite[\S 2]{Haugseng2021monad}. 

\begin{definition}\label{defn: coherent-adjunction}
Let $\adj$ be the \emph{walking adjunction}, namely the free ordinary $2$-category generated by an adjunction; see \cite[Remark 3.3.8]{RV2016monad} and \cite[Notation 4.1]{Haugseng2021monad}. It has two objects $-$ and $+$, generating $1$-morphisms
\[
l:-\longrightarrow+,\quad r:+\longrightarrow-,
\]
and $2$-morphisms
\[
\eta:\id_-\Rightarrow r\circ l,\quad \epsilon:l\circ r\Rightarrow\id_+,
\]
subject to the triangle identities
\[
(\epsilon l)\circ(l\eta)=\id_l,
\quad
(r\epsilon)\circ(\eta r)=\id_r.
\]
Let $\mnd\subset\adj$ be the full sub-$2$-category on the object $-$; this is the \emph{walking monad}. Its generating endomorphism $t=r\circ l$ has unit $\eta:\id_-\Rightarrow t$ and multiplication
\[
\mu_t=r\epsilon l:t\circ t=r\circ l\circ r\circ l\Rightarrow r\circ l=t.
\]
We regard an ordinary $2$-category $\cA$ as a $\Cat_\infty$-enriched category with the same objects and mapping $\infty$-categories $\rN_\bu(\cA(x,y))$, the ordinary nerves of its hom-categories. This is the convention of \cite[Remark 2.14]{Haugseng2021monad}. See also \cite[Proposition 3.2.2, Observation 3.3.1, and Corollary 3.3.5]{RV2016monad}.

Then a \emph{monad} in an $(\infty,2)$-category $\cC$ is a functor $\mnd\to\cC$; see \cite[Notation 5.1 and Definition 5.2]{Haugseng2021monad}. Following \cite[Definition 4.2]{Haugseng2021monad}, a \emph{coherent adjunction} $N\dashv N^r$ in $\cC$ is a functor of $(\infty,2)$-categories
\[
\Phi:\adj\longrightarrow\cC
\]
such that
\[
\Phi(-)=\varnothing,~ \Phi(+)=\Sigma,~ \Phi(l)=N,~ \Phi(r)=N^r.
\]
Its low-dimensional data include the $2$-morphisms
\[
\eta:\id_{\varnothing}\Rightarrow N^r\circ N,
~
\epsilon:N\circ N^r\Rightarrow\id_\Sigma,
\]
and invertible $3$-morphisms representing the triangle identities
\[
(\epsilon N)\circ(N\eta)\simeq\id_N,
~
(N^r\epsilon)\circ(\eta N^r)\simeq\id_{N^r}.
\]
All further higher compatibility data are encoded by the functor $\Phi$.
\end{definition}
\brem\label{rem: contractible choice of adj}
It is enough to specify an adjunction $N\dashv N^r$ in the homotopy $2$-category $h_2\cC$. Here $h_2\cC$ denotes the $2$-category with hom-categories $h(\cC(x,y))$, where $\cC(x,y)$ is the mapping $\infty$-category of \cite[Notation 2.13(ii)]{Haugseng2021monad}. Indeed, \cite[Theorem 4.4]{Haugseng2021monad} states that evaluation at the generating left adjoint induces an equivalence
\[
\Map(\adj,\cC)\xrightarrow{\simeq}\Map(\Delta^1,\cC)_{\mathrm{ladj}}.
\]
The analogous result holds for right adjoints. In the simplicially enriched model of Riehl--Verity, existence follows from \cite[Theorem 4.3.11]{RV2016monad}; the trivial-fibration statement underlying evaluation follows from \cite[Propositions 4.4.12 and 4.4.17]{RV2016monad}, and the contractibility of the space of extensions of a fixed left adjoint is \cite[Theorem 4.4.18]{RV2016monad}. Consequently, if it is only specified that $N$ admits a right adjoint, the space of coherent adjunctions extending $N$ is contractible.
\erem

The main categorical result is as follows.
\begin{proposition}\label{prop:adjunction-algebra-module}
Let $\cC$ be an $(\infty,2)$-category, let $\varnothing,\Sigma,\Sigma'\in\cC$, and let $\mathcal V=\End_{\cC}(\varnothing)$. Consider $1$-morphisms
\[
\varnothing\xrightarrow{N}\Sigma\xrightarrow{Q}\Sigma'\xrightarrow{P}\varnothing.
\]
Suppose that $N\dashv N^r$ is a coherent adjunction, and set $\bY_N=N^r\circ N$ and $\bY=P\circ Q\circ N$. Then \[\bY_N\in\Alg(\mathcal V)\aand \bY\in\operatorname{RMod}_{\bY_N}(\mathcal V).\]If $P^\ell\dashv P$ is also a coherent adjunction and $\bY_P=P\circ P^\ell$, then \[\bY_P\in\Alg(\mathcal V)\aand\bY\in{}_{\bY_P}\operatorname{BMod}_{\bY_N}(\mathcal V).
\]
\end{proposition}
\begin{proof}
Let $\Phi_N:\adj\to\cC$ exhibit $N\dashv N^r$, carrying the generating $1$-morphisms $l:-\to+$ and $r:+\to-$ to $N$ and $N^r$, respectively. When $P^\ell\dashv P$ is given, let $\Phi_P:\adj\to\cC$ exhibit this adjunction, carrying $l$ and $r$ to $P^\ell$ and $P$.

The restrictions $\Phi_N|_{\mnd}$ and $\Phi_P|_{\mnd}$ are monads on $\varnothing$ whose underlying endomorphisms are $N^r\circ N$ and $P\circ P^\ell$, respectively. By \cite[Theorem 1.4 and Corollary 8.9]{Haugseng2021monad}, there is an equivalence
\[
\Mnd(\cC)_{\mathrm{colax},\varnothing}\simeq\Alg\bigl(\End_{\cC}(\varnothing)\bigr)
\]
compatible with the forgetful functors to endomorphisms, where we write $\Mnd(\cC)_{\mathrm{colax}}$ for the underlying $\infty$-category of the $(\infty,2)$-functor category $\operatorname{FUN}(\mnd,\cC)_{\mathrm{colax}}$ and $\Mnd(\cC)_{\mathrm{colax},\varnothing}$ for its fiber over $\varnothing$ under evaluation at the unique object of $\mnd$. Hence $\bY_N=N^r\circ N$ and $\bY_P=P\circ P^\ell$ carry the asserted $E_1$-algebra structures. Their units and multiplications are represented by
\[
\eta_N:\id_{\varnothing}\Rightarrow N^r\circ N,\quad N^r\epsilon_NN:(N^r\circ N)\circ(N^r\circ N)\Rightarrow N^r\circ N
\]
and
\[
\eta_P:\id_{\varnothing}\Rightarrow P\circ P^\ell,\quad P\epsilon_PP^\ell:(P\circ P^\ell)\circ(P\circ P^\ell)\Rightarrow P\circ P^\ell.
\]

We next make the module structures explicit. In $\adj$, set $t=r\circ l$, with unit $\eta:\id_-\Rightarrow t$, multiplication $\mu=r\epsilon l:t\circ t\Rightarrow t$, and counit $\epsilon:l\circ r\Rightarrow\id_+$. Composition gives $\adj(-,+)$ a right action of $\adj(-,-)$ and gives $\adj(+,-)$ a left action of $\adj(-,-)$, defined by $(F,T)\mapsto F\circ T$ and $(T,G)\mapsto T\circ G$. With respect to these actions, $l$ is a right module and $r$ is a left module over $t$, with action maps
\[
\rho_l=\epsilon l:l\circ t=l\circ r\circ l\Rightarrow l,~ \lambda_r=r\epsilon:t\circ r=r\circ l\circ r\Rightarrow r.
\]
The triangle identities and interchange give
\[
\rho_l\circ(l\eta)=\id_l,~ \rho_l\circ(\rho_lt)=\rho_l\circ(l\mu),~
\lambda_r\circ(\eta r)=\id_r,~ \lambda_r\circ(t\lambda_r)=\lambda_r\circ(\mu r).
\]
The off-diagonal hom-categories of $\adj$ and their ordinal-sum composition encode the corresponding iterated actions and all their coherences; see \cite[Remark 3.3.8]{RV2016monad}. Consequently, any $(\infty,2)$-functor $\Phi:\adj\to\cC$ carries these structures to coherent module objects in the appropriate mapping $\infty$-categories.

Applying $\Phi_N$ makes $N$ a right module over $\bY_N$ in the right $\mathcal V$-module $\infty$-category $\cC(\varnothing,\Sigma)$, with action
\[
\epsilon_NN:N\circ\bY_N=N\circ N^r\circ N\Rightarrow N.
\]
Postcomposition by $P\circ Q$ defines a right $\mathcal V$-linear functor
\[
P\circ Q\circ-:\cC(\varnothing,\Sigma)\longrightarrow\mathcal V,
\]
and therefore equips $\bY=P\circ Q\circ N$ with a right $\bY_N$-module structure whose action is represented by
\[
PQ\epsilon_NN:\bY\circ\bY_N=P\circ Q\circ N\circ N^r\circ N\Rightarrow P\circ Q\circ N=\bY.
\]

Similarly, applying $\Phi_P$ makes $P$ a left $\bY_P$-module in $\cC(\Sigma',\varnothing)$. The left $\mathcal V$-linear functor $-\circ Q\circ N:\cC(\Sigma',\varnothing)\to\mathcal V$ then gives $\bY$ the left action
\[
P\epsilon_PQN:\bY_P\circ\bY=P\circ P^\ell\circ P\circ Q\circ N\Rightarrow P\circ Q\circ N=\bY.
\]

The two adjunctions together with $Q$ determine a functor from the free $(\infty,2)$-category generated by two coherent adjunctions and a $1$-morphism between their target objects. In this free $(\infty,2)$-category, the cross-composite $r_P\circ Q\circ l_N$ carries the left action of $r_P\circ l_P$ and the right action of $r_N\circ l_N$, and composition and interchange supply their mixed coherences. Thus the right $\bY_N$-action on $\bY$ is a morphism of left $\bY_P$-modules, so that
\[
\bY\in\operatorname{RMod}_{\bY_N}\bigl(\operatorname{LMod}_{\bY_P}(\mathcal V)\bigr).
\]
The equivalence
\[
\operatorname{RMod}_{\bY_N}\bigl(\operatorname{LMod}_{\bY_P}(\mathcal V)\bigr)\simeq{}_{\bY_P}\operatorname{BMod}_{\bY_N}(\mathcal V)
\]
of \cite[Theorem 4.3.2.7 and Remark 4.3.3.8]{HAlurie2017} gives the asserted bimodule structure.
\end{proof}

\brem
Although Haugseng formulates the Gray-product construction under \cite[Assumption 3.5]{Haugseng2021monad}, condition (3) of that assumption is not used in \cite[Theorem 1.4 and Corollary 8.9]{Haugseng2021monad}. Their proofs pass through \cite[Theorem 7.3, Proposition 8.5, Definition 8.7, and Remark 8.8]{Haugseng2021monad} and require only conditions (1) and (2), which were established in \cite[Theorem 6.1]{Mae2021tensor}; see also \cite[Remark 3.6]{Haugseng2021monad}. Condition (3) enters instead in the comparison with the strict Gray products used by Zaganidis \cite{Zaganidis2017thesis} in \cite[Remark 4.12, Theorem 5.6, Corollary 5.8, and Theorem 1.2]{Haugseng2021monad}. Consequently, Proposition \ref{prop:adjunction-algebra-module} applies to an arbitrary $(\infty,2)$-category and requires no additional model hypothesis. If we only need the objectwise monad-algebra correspondence, it also follows directly from \cite[Proposition 8.5, Definition 8.7, and Remark 8.8]{Haugseng2021monad}.
\erem

Assume the existence of $\BI_2^{\pa_+}$ and suppose that the decorated cobordisms
\[
N:\varnothing\longrightarrow\Sigma,
~
M:\Sigma\longrightarrow\varnothing
\]
arising from \eqref{eq:decomp} define $1$-morphisms in $\BI_2^{\pa_+}$ and that $N\dashv N^r$ is a coherent adjunction. Set
\[
\bY_N=N^r\circ N,
~
\bY=M\circ N.
\]
Applying Proposition \ref{prop:adjunction-algebra-module} with $Q=\id_\Sigma$ and $P=M$, and using the monoidal equivalence
\[
\End_{\BI_2^{\pa_+}}(\varnothing)\simeq\BI^{\pa_+},
\]
gives
\[
\bY_N\in\Alg(\BI^{\pa_+}),
~
\bY\in\operatorname{RMod}_{\bY_N}(\BI^{\pa_+}).
\]
Thus the categorical construction recovers the right-module convention used in the geometric description above.

Since $CI$ is symmetric monoidal,
\[
\cA=CI(\bY_N)\in\Alg(\D),~\cM=CI(\bY)\in\operatorname{RMod}_{\cA}(\D).
\]
Finally, we outline a possible enhancement of $CI$ to an $(\infty,2)$-functor, although this is not needed for the preceding argument. A complete construction requires concrete models of $(\infty,2)$-categories, such as the complete $2$-fold Segal spaces in \cite{CS2019category}.

One first needs to construct a variant $\what{\D}$ of $\D$ to which $CI$ extends and which satisfies the ``nice'' hypothesis of \cite[Theorem 1.1]{Haugseng2017morita}. Assuming such a variant, one may seek an enhancement
\[
CI_2:\BI_2^{\pa_+}\longrightarrow\mathfrak{Alg}^{\mathrm{Mor}}_1(\what{\D}).
\]
Here $\mathfrak{Alg}^{\mathrm{Mor}}_1(\what{\D})$ denotes Haugseng's Morita $(\infty,2)$-category constructed in \cite[Theorem 1.1]{Haugseng2017morita}; see also \cite{Haugseng2023morita}. Its objects are $E_1$-algebras in $\what{\D}$, its $1$-morphisms are bimodules, and its $2$-morphisms are bimodule morphisms.

Following Baldwin--Bloom \cite{BBmonopole}, one may seek, for each oriented surface $\Sigma$, a finite collection $\{N_{\Sigma,i}\}_{i\in I_{\Sigma}}$ of adjointable $1$-morphisms that generate the relevant bordered $3$-manifolds under iterated surgery triangles; the corresponding undecorated statement is proved in \cite{Culler2014surgery}. This surgery-triad generation motivates an analogous Yoneda-type split-generation construction in an $A_\infty$-category; see \cite[\S 3.4]{Denis2014intro}. The heuristic assignments are as follows.

After choosing coherent right adjoints $N_{\Sigma,i}^r$, the assignment to an object $\Sigma$ is the matrix $E_1$-algebra
\[
\cA_\Sigma=\bigoplus_{i,j\in I_\Sigma}CI(N_{\Sigma,j}^r\circ N_{\Sigma,i}).
\]
The expected multiplication on matching summands is induced by the counits as in Proposition \ref{prop:adjunction-algebra-module}
\[
CI(N_{\Sigma,k}^r\circ N_{\Sigma,j})\otimes CI(N_{\Sigma,j}^r\circ N_{\Sigma,i})
\longrightarrow CI(N_{\Sigma,k}^r\circ N_{\Sigma,i}),
\]
and is expected to vanish on nonmatching summands. Equivalently, one expects this algebra to be the image under $CI$ of the matrix endomorphism algebra associated with the formal direct sum of the $N_{\Sigma,i}$ in the free semiadditive completion of $\BI_2^{\pa_+}$. See \cite[Proposition 2.8 and Theorem 2.13]{Glasman2018add} for the construction of the free semiadditive completion of a general $\infty$-category and \cite[\S (3k)]{seidelbook} for its $A_\infty$ category analogue.

For a $1$-morphism $Q:\Sigma_0\to\Sigma_1$, the assignment is the bimodule
\[
\cM_Q=\bigoplus_{\substack{i\in I_{\Sigma_0}\\j\in I_{\Sigma_1}}}
CI(N_{\Sigma_1,j}^r\circ Q\circ N_{\Sigma_0,i}).
\]
Applying Proposition \ref{prop:adjunction-algebra-module} componentwise and taking finite direct sums should equip $\cM_Q$ with an $(\cA_{\Sigma_1},\cA_{\Sigma_0})$-bimodule structure.

For a $2$-morphism $W:Q\Rightarrow Q'$, horizontal composition with the identity $2$-morphisms gives the candidate components
\[
W_{j,i}=\id_{N_{\Sigma_1,j}^r}\circ W\circ\id_{N_{\Sigma_0,i}}:
N_{\Sigma_1,j}^r\circ Q\circ N_{\Sigma_0,i}
\Longrightarrow
N_{\Sigma_1,j}^r\circ Q'\circ N_{\Sigma_0,i}.
\]
Geometrically, the identity $2$-morphisms are represented by the product cobordisms $I\times N_{\Sigma_1,j}^r$ and $I\times N_{\Sigma_0,i}$. Gluing them to the corresponding vertical boundary faces of the $4$-dimensional cobordism $W$ and rounding corners produces a $4$-dimensional cobordism $W_{j,i}$ between the closed $3$-manifolds displayed above. Applying $CI$ componentwise should give a bimodule morphism
\[
\cM_W=
\bigoplus_{\substack{i\in I_{\Sigma_0}\\j\in I_{\Sigma_1}}}
CI(W_{j,i}):
\cM_Q\longrightarrow\cM_{Q'}.
\]
The adjunction counits should give a comparison map
\[
\cM_{Q_2}\otimes_{\cA_{\Sigma_1}}\cM_{Q_1}
\longrightarrow\cM_{Q_2\circ Q_1}.
\]

The main remaining problem is to prove that this map is an equivalence. Here the relative tensor product is formed in $\what{\D}$ and computed by the geometric realization of the two-sided bar construction.

The intended strategy is to verify this map in elementary cases and then extend the equivalence using the surgery fiber squares developed in \S \ref{subsec: Fiber squares for surgery cobordisms}. This surgery-induction argument requires a generation result similar to \cite{Culler2014surgery} for cobordisms.

A toy model for this strategy in $\operatorname{Pin}(2)$-monopole Floer theory is provided by \cite[Theorem 4]{francesco2017connectedsum}, in which the $Q_i$ are $1$-morphisms from $3$-manifolds with spherical boundaries. For instanton theory, motivated by Mrowka's talk \cite{Mrowka2024prospects}, one might need to lift the $\Z/2$-grading to a $\Z/4$-grading for the verification of the elementary cases without explicitly calculating the higher multiplications.

Their naturality and higher coherences remain to be established in future work.

More speculatively, one may hope for a fully extended instanton TQFT given by a symmetric monoidal functor
\[
CI_4:\BI_4^{\pa_+}\longrightarrow\mathfrak{Alg}^{\mathrm{Mor}}_3(\what{\D}).
\]
Here $\BI_4^{\pa_+}$ would denote a suitable $(\infty,4)$-cobordism category for instanton theory, with the positive boundary condition imposed on $4$-dimensional bordisms, modeled on $\Bord_4^{(\infty,4)}$ constructed in \cite[\S 6]{CS2019category}. The target $\mathfrak{Alg}^{\mathrm{Mor}}_3(\what{\D})$ is Haugseng's Morita $(\infty,4)$-category constructed in \cite{Haugseng2017morita,Haugseng2023morita}. Its objects are $E_3$-algebras, its $1$-morphisms are $E_2$-algebras in bimodules, its $2$-morphisms are $E_1$-objects in bimodules in bimodules, and so forth.

\subsection{Organization of the remaining sections}
In \S\ref{sec: Preliminaries in infty-category}, we review the language of $\infty$-categories, with particular emphasis on the models given by quasicategories, complete Segal spaces, and simplicial categories, as well as their relationship. We also consider the symmetric monoidal structures and duality, and construct the $\infty$-derived category $\D$.

In \S\ref{sec: Preliminaries in cobordism category}, we review the construction of the $\infty$-cobordism category $\Bord_n$ by Calaque--Scheimbauer \cite{CS2019category}. We adapt only the results needed to construct the $\infty$-cobordism categories $\BI$ and $\BIa$, with particular emphasis on tangential structures. The verifications of the Segal condition and the symmetric monoidal condition in our setting are provided in Appendices \S\ref{sec: segal} and \S\ref{sec: symmetric monoidal}, respectively.

In \S\ref{sec: instanton cob cat}, we add subgroup data, the non-integral condition, and perturbations to the usual $\infty$-cobordism category $\Bord_4^\theta$ with a tangential structure $\theta$ in order to construct $\BI^\pa$ and $\BIa^\pa$, where $\BIa^\pa$ incorporates auxiliary data coming from Riemannian metrics and perturbations. We also show that the forgetful functor $\BIa^\pa\to\BI^\pa$ is an equivalence of $\infty$-categories.

In \S\ref{sec: Construction of the TQFT functor}, we reinterpret the instanton-theoretic information developed by Kronheimer--Mrowka \cite{kronheimer2011khovanov} as a TQFT functor $CI_a:\BIa^\pa\to\D$. We only sketch the gauge-theoretic construction and focus primarily on the $\infty$-categorical aspects.

In \S\ref{sec: Cube complexes and exact triangles}, we provide further details on the constructions in \S\ref{subsec: Cube complexes}. We upgrade the construction of the hypercube underlying the link spectral sequence to a cube complex and extend this construction to collections of cobordisms supported in disjoint regions. We also lift the surgery exact triangle and its link surgery generalizations to fiber squares in the stable $\infty$-category $\Fun(\bbox^n,\D)$. The $\infty$-categorical version of the triangle detection lemma is given in Appendix \S\ref{sec: Triangle detection lemma} by Longke Tang.

In \S\ref{sec: Generalized cap product}, we provide further details on the constructions in \S\ref{subsec: enhanced cob}. We lift the generalized cap product $\mu$-operators in instanton Floer homology to the chain level and study their interactions with cobordisms and cubes of chain complexes.

\section{Preliminaries on \texorpdfstring{$\infty$}{infinity}-categories}\label{sec: Preliminaries in infty-category}
In this section, we review necessary definitions and constructions for $\infty$-categories. By \emph{category} we mean an ordinary category. The term \emph{$\infty$-category} is an abbreviation for $(\infty,1)$-category. The main difference between $\infty$-categories and ordinary categories is that each hom-set is replaced by a hom-space, modeled here by a Kan complex. Although the symbol $\Map$ is common for hom-spaces, we retain $\Hom$. We distinguish ordinary categories, model categories, and $\infty$-categories notationally, writing, for example, $\Set$, $\Space$, and $\Cati$, respectively.

There are several equivalent models of $\infty$-categories; see \cite{Rezk2001completion,JT2007quasisegal}, \cite[Chapter 1]{HTTlurie2009}, and also the survey articles \cite{Bergner2010survey1, Bergner2020survey2}. The $\infty$-category $\Cati$ of small $\infty$-categories is also well-defined (cf.\ \cite[Definition 3.0.0.1]{HTTlurie2009} and Definition \ref{defn: homotopy coherent nerve}), so that the categorical notions, e.g.\ limits, colimits, and equivalences, can be used freely for $\infty$-categories.

The set-theoretic technicalities are discussed in \cite[\S 1.2.15]{HTTlurie2009}, where the most convenient way to resolve the issue is by assuming the existence of the strongly inaccessible cardinal $\kappa$ and using the
collection $\mathcal{U}(\kappa)$ of all sets having rank smaller than $\kappa$, which is a \emph{Grothendieck universe}.

We use several models of $\infty$-categories. The first is the \emph{quasicategory} model, formulated in terms of simplicial sets \cite{HTTlurie2009,HAlurie2017,kerodon}. The second is the \emph{complete Segal space} model used in the literature on cobordism categories \cite{lurie2009cob,CS2019category}. We distinguish this notion from the definitions in \cite{Rezk2001completion,JT2007quasisegal}, which we call \emph{Reedy fibrant complete Segal spaces}, following \cite[Warning 2.1.16]{lurie2009cob} and \cite[\S 1.4]{CS2019category}. We also use simplicial categories as an auxiliary model to describe the symmetric monoidal construction explicitly.

Joyal--Tierney \cite{JT2007quasisegal} construct two Quillen equivalences between the model categories of Reedy fibrant complete Segal spaces and quasicategories. We use the first, which we describe after introducing simplicial sets and simplicial spaces. The identity functor also gives a Quillen equivalence between the model categories of complete Segal spaces and Reedy fibrant complete Segal spaces. Simplicial categories are related to quasicategories by the homotopy coherent nerve \cite[\S 1.1.5 and Theorem 2.2.0.1]{HTTlurie2009}.
\subsection{Quasicategory}\label{subsec: quasicategories}
We start with the definition of simplicial sets.
\bdefn[Simplicial set]\label{defn: simplicial set}
Let $\Delta$ denote the \emph{simplex category}, whose objects are the finite ordered sets
\[
[m]=(0<\cdots<m),
\]
and whose morphisms are monotone maps. Let $\delta_m^j:[m-1]\to[m]$ and $\sigma_m^j:[m+1]\to[m]$ be the standard coface and codegeneracy maps, respectively. Every morphism in $\Delta$ is a composite of these maps. We omit the subscript $m$ when it is clear.

A \emph{simplicial set} is a functor $X_\bu:\Delta^{\opp}\to\Set$. It consists of sets $X_m$ for $m\in\N$, together with face maps
\[
d_j=X_\bu(\delta_m^j):X_m\longrightarrow X_{m-1}
\]
and degeneracy maps
\[
s_j=X_\bu(\sigma_m^j):X_m\longrightarrow X_{m+1}.
\]
We omit the subscript $\bu$ when no confusion can arise.

An element of $X_0$ is called a \emph{vertex}. An element $e\in X_1$ is called an \emph{edge}, with \emph{source} $d_1(e)$ and \emph{target} $d_0(e)$. More generally, an element of $X_m$ is called an \emph{$m$-simplex}. A simplex is \emph{degenerate} if it lies in the image of a degeneracy map and is \emph{nondegenerate} otherwise.

A \emph{simplicial subset} $Y_\bu$ of $X_\bu$ is a simplicial set such that $Y_m\subseteq X_m$ for every $m\in\N$, denoted by $Y_\bu\subseteq X_\bu$. If
\[
X_\bu=Y_\bu\sqcup Y_\bu',
\]
then $Y_\bu$ is called a \emph{summand} of $X_\bu$.

A simplicial set $X_\bu$ is \emph{connected} if it is nonempty and every summand is either empty or equal to $X_\bu$. A simplicial subset $Y_\bu\subseteq X_\bu$ is a \emph{connected component} if it is a connected summand. Let $\pi_0(X_\bu)$ denote the set of connected components of $X_\bu$.

Let
\[
\sSet=\Funo(\Delta^{\opp},\Set)
\]
denote the category of simplicial sets.

Products of simplicial sets are defined degreewise,
\[
(X_\bu\times Y_\bu)_n=X_n\times Y_n,
\]
with face and degeneracy maps acting componentwise.
\edefn
\bdefn[Standard simplex]
For fixed $n\in\N$, let $\Delta^n$ be the simplicial set that sends $[m]$ to $\Hom_{\Delta}([m],[n])$, called the \emph{standard $n$-simplex}. An $n$-simplex of $X_\bu$ can be identified with a morphism $\sigma:\Delta^n\to X_\bu$. Let $\partial \Delta^n$ be the simplicial subset of $\Delta^n$ that sends $[m]$ to \[\{\al\in\Hom_{\Delta}([m],[n])~\big|~\al\text{ is not surjective}\},\]which is called the \emph{boundary} of $\Delta^n$. Note that $\partial \Delta^n$ is the union of $(n+1)$ copies of standard $(n-1)$-simplex along some standard $(n-2)$-simplices in their boundaries. In particular, we know that $\partial \Delta^0=\varnothing$ is the trivial simplicial set that sends any $[m]$ to the empty set.

A \emph{simplicial operator} $\alpha:\Delta^r\to\Delta^q$ is the map of standard simplices induced by an order-preserving map $\alpha:[r]\to[q]$. By the Yoneda lemma, every map of simplicial sets between standard simplices arises uniquely in this way.

For $n\in\N_+$ and an integer $i\in[0,n]$, let $\Lambda^n_i$ be the simplicial subset of $\partial \Delta^n$ that sends $[m]$ to \[\{\al\in\Hom_{\Delta}([m],[n])~\big|~[n]\not\subseteq \al([m])\cup \{i\}\},\]which is called the \emph{$i$-th horn} in $\Delta^n$. It is called an \emph{inner horn} if $i\in(0,n)$ and an \emph{outer horn} if $i\in\{0,n\}$. Note that $\Lambda^n_i$ is the union of $n$ copies of standard $(n-1)$-simplex along some standard $(n-2)$-simplices in their boundaries.

The geometric realization of $\Delta^n$ is the standard geometric simplex
    \[
    |\Delta^n|=\big\{(x_0, \ldots, x_n) \in \R^{n+1}~\big|~\sum_i x_i = 1, \, x_i \geq 0 \big\}.
    \]
  
  We also consider the \emph{affine geometric simplex}
    \[
   |\Delta^n|_{\aff}= \big\{(x_0, \ldots, x_n) \in \R^{n+1}~\big|~\sum_i x_i = 1\big\}.
    \]
Moreover, we define\[|\partial \Delta^n|=\{(x_0,\dots,x_n)\in|\Delta^n|~\big|~ x_j=0\text{ for some }j\},\]\[|\Lambda^n_i|=\{(x_0,\dots,x_n)\in|\Delta^n|~\big|~ x_j=0\text{ for some }j\neq i\}.\]  Define their affine realizations $|\partial\Delta^n|_{\aff}$ and $|\Lambda_i^n|_{\aff}$ by replacing $|\Delta^n|$ with $|\Delta^n|_{\aff}$ in these formulas. Note that $|\partial \Delta^n|=\partial |\Delta^n|$, and $|\Lambda^n_i|$ is the union of all codimension-one faces of $|\Delta^n|$ except the face opposite its $i$-th vertex.
\edefn
\bdefn[Kan complex]\label{defn: space}
A map $p:X_\bu \to Y_\bu$ of simplicial sets is called a \emph{Kan fibration} if for each $n\in\N_+$ and $i\in[0,n]$, every strictly commutative diagram
\begin{equation*}
  \xymatrix{
\Lambda^n_i \ar[r] \ar[d] & X_\bu \ar[d]^{p} \\
\Delta^n \ar[r] \ar@{-->}[ur] & Y_\bu
}
\end{equation*}
admits a dotted arrow making the diagram commute. We call this an \emph{extension} of a map from $\Lambda^n_i$ to $\Delta^n$, or the \emph{right lifting property} for the inclusion $\Lambda^n_i\to \Delta^n$. We will use these notions for other maps.

A simplicial set $X_\bu$ is called a \emph{Kan complex} if the projection map $X_\bu\to \Delta^0$ is a Kan fibration. This condition is called the \emph{Kan condition} of the simplicial set. A Kan complex is a simplicial analogue of a topological space; see Example \ref{exmp: singular complex}.

Let $\Space$ be the category of simplicial sets with the Kan--Quillen model structure. Its cofibrations are the monomorphisms, its fibrations are the Kan fibrations, and its weak equivalences are the morphisms whose geometric realizations are weak homotopy equivalences.

Note that in \cite[\S 1, Conventions]{CS2019category}, a \emph{space} means a fibrant object in $\Space$, i.e.\, a Kan complex. In some literature, a space also denotes a vector space or a topological space, so we keep using the notion of Kan complex for clarity.
\edefn

\begin{example}\label{exmp: singular complex}
The \emph{singular complex} of a topological space provides the basic example of a Kan complex. Recall from \cite[\S 18]{Lee2013smooth} that a singular $n$-simplex in a topological space $X$ is a continuous map $\sigma:|\Delta^n|\to X$. This defines a functor
\[
S^{\sing}_\bu:\mathrm{Top}\to\sSet,
\]
which is right adjoint to the geometric realization functor $|\cdot|$. The Kan condition follows from the fact that each inclusion $|\Lambda_i^n|\to|\Delta^n|$ admits a continuous retraction.

Suppose first that $X$ is a smooth manifold without boundary. Let $S^{\sm}_\bu(X)$ denote the \emph{smooth complex} whose $n$-simplices are maps $\sigma:|\Delta^n|\to X$ that locally extend to smooth maps on open subsets of \ $|\Delta^n|_{\aff}$. By \cite[\S 3]{SP2019diff}, this agrees with the complex defined using the sub-diffeology on $|\Delta^n|$, and there is an \ identification
\[
S^{\sm}_\bu(X)=S^{\sub}_\bu(X).
\]
We also consider the \emph{affine complex} $S^{\aff}_\bu(X)$ whose $n$-simplices are smooth maps \ $|\Delta^n|_{\aff}\to X$. Restriction induces natural maps
\begin{equation}\label{eq: inclusion}
S^{\aff}_\bu(X)\longrightarrow S^{\sm}_\bu(X)\longrightarrow S^{\sing}_\bu(X).
\end{equation}
The first map is a weak homotopy equivalence by \cite[Theorem 1.1]{Kihara2024smooth}, while the composite is the weak homotopy equivalence of \cite[Theorem 4.15]{Bunk2022smooth}. Hence the second map is also a weak homotopy equivalence by the two-out-of-three property; see Lemma \ref{lem: homotopy properties}. A direct proof may also be obtained by adapting the proof of \cite[Theorem 18.7]{Lee2013smooth} to the simplicial setting. The first map in \eqref{eq: inclusion} need not be injective, since distinct smooth maps \ $|\Delta^n|_{\aff}\to X$ can have the same restriction to $|\Delta^n|$.

The Kan conditions for these complexes are more subtle. By \cite[Definition 4.3 and Corollaries 4.36 and 4.47]{CW2014diff}, $S^{\aff}_\bu(X)$ is Kan for every smooth manifold $X$ without boundary, whereas it is not Kan for any smooth manifold with nonempty boundary or corners. For a general diffeological space $X$, the complex $S^{\sub}_\bu(X)$ need not be Kan by \cite[Remark 3.2(2)]{Kihara2024smooth}; this result does not address the case of smooth manifolds. Kihara's singular complex $S^{\kih}_\bu(X)$, defined for a general diffeological space, is a Kan complex by \cite[Corollary 2.6(2)]{Kihara2024smooth}. The natural maps
\[
S^{\aff}_\bu(X)\longrightarrow S^{\sub}_\bu(X)\longrightarrow S^{\kih}_\bu(X)
\]
are weak homotopy equivalences by \cite[Theorem 1.1]{Kihara2024smooth}. Thus $S^{\kih}_\bu(X)$ provides a Kan replacement for both $S^{\aff}_\bu(X)$ and $S^{\sub}_\bu(X)$. We discuss the diffeological spaces and these complexes in \S \ref{subsec: diffeological spaces}.
\end{example}
\bdefn[Quasicategory]\label{bdefn: quasicategory}
A \emph{quasicategory}, historically also called a \emph{weak Kan complex}, is a simplicial set satisfying the \emph{weak Kan condition} (also called \emph{inner horn condition}), i.e., any map of simplicial sets from the inner horn $\Lambda^n_i$ for $i\in (0,n)$ has an extension to $\Delta^n$. A Kan complex is a special case of a quasicategory. From \cite[Theorem 2.2.5.1]{HTTlurie2009}, there is a model structure on $\sSet$, originally constructed by Joyal and denoted by $\QCat$, such that the weak equivalences are the categorical equivalences, the cofibrations are the monomorphisms, and fibrant objects are the quasicategories. 

Throughout this paper, a \emph{weak homotopy equivalence} means a weak equivalence in the Kan--Quillen model structure. A \emph{categorical equivalence}, equivalently an equivalence of quasicategories, means a weak equivalence in the Joyal model structure. We use the unqualified term ``weak equivalence'' only when the ambient model category is specified.

For an ordinary category $\cC$, its \emph{nerve} $\rN_\bu\cC=\rN_\bu(\cC)$, defined by $\rN_m(\cC)=\Funo([m],\cC)$, is a quasicategory. Moreover, by \cite[\texttt{002Z}]{kerodon}, the nerve provides a fully faithful functor\[\rN_\bu:\Cat\to \sSet\]for the category $\Cat$ of small categories. Hence quasicategories are generalizations of Kan complexes and ordinary categories.

Moreover, by \cite[{\texttt{004M}, \texttt{004N}}]{kerodon}, the nerve functor admits a left adjoint \[h:\sSet\to \Cat,\]called the homotopy category functor, which is defined as follows. \begin{itemize}
    \item The objects of $h(X_\bu)$ are the vertices of $X_\bu$.
    \item Every edge $e$ of $X_\bu$ determines a morphism $[e]$ in $h(X_\bu)$, with the same source and the target as $e$.
    \item The collection of morphisms in $h(X_\bu)$ is generated under composition by morphisms of the form $[e]$, subject only to the relations
\[[s_0(x)] = \id_x \text{ for }x \in  X_0 \aand [d_1(\sigma)] = [d _0(\sigma)]\circ [d_2(\sigma)]\text{ for }\sigma \in X_2.\]
\end{itemize} 

For a quasicategory $X_\bu$, an element $f:x\to y$ in $X_1$ is called an \emph{isomorphism} if it induces an isomorphism $[f]:x\to y$ in $h(X_\bu)$, namely there exists another morphism $[g]:y\to x$ in $h(X_\bu)$ such that $[g]\circ [f]=\id_x$ and $[f]\circ [g]=\id_y$.

By \cite[\texttt{01D0}]{kerodon}, given a quasicategory $X_\bu$, let $X^\simeq_\bu$ denote the simplicial subset of $X_\bu$ comprised of those simplices $\sigma :\Delta^n \to X_\bu$ which carry each edge of $\Delta^n$ to an isomorphism in $X_\bu$. We will refer to $X^\simeq_\bu$ as the \emph{core} of $X_\bu$, which is the maximal Kan complex contained in $X_\bu$ by \cite[\texttt{01D9}]{kerodon}.

For a quasicategory $X_\bu$, an element in $X_0$ is called an \emph{object} and an element in $X_1$ is called an \emph{morphism} or \emph{$1$-morphism} in some literature. An element in $X_m$ for $m>1$ is called a \emph{$m$-morphism} or simply a \emph{higher morphism}. In particular, an $m$-morphism of $\rN_\bu\cC$ for an ordinary category $\cC$ is a diagram\[C_0\xra{f_1}\cdots\xra{f_m}C_m\]for objects $C_i$ and morphisms $f_j$ in $\cC$ with $i\in[0,m]$ and $j\in [1,m]$.
\edefn
\bdefn[Functor simplicial set]\label{defn: functor simplicial set}
Let $X_\bu$ and $Y_\bu$ be simplicial sets. The \emph{functor simplicial set} $\Fun_\bu(X_\bu,Y_\bu)$ is defined by
\[
\Fun_m(X_\bu,Y_\bu)=\Hom_{\sSet}(\Delta^m\times X_\bu,Y_\bu)
\]
for every $m\geq 0$ \cite[\texttt{0060}]{kerodon}. The functor $\Fun_\bu(X_\bu,-)$ is right adjoint to $-\times X_\bu$, with natural identifications
\begin{equation}\label{eq: product functor adjunction}
\Hom_{\sSet}(Z_\bu\times X_\bu,Y_\bu)
\cong\Hom_{\sSet}(Z_\bu,\Fun_\bu(X_\bu,Y_\bu))
\end{equation}for any simplicial set $Z_\bu$.

For ordinary categories $\cC$ and $\cD$, there is a canonical isomorphism of simplicial sets
\[
\rN_\bu\bigl(\Funo(\cC,\cD)\bigr)\xrightarrow{\sim}\Fun_\bu\bigl(\rN_\bu(\cC),\rN_\bu(\cD)\bigr)
\]
\cite[\texttt{0062}]{kerodon}. Thus, the functor simplicial set generalizes the ordinary functor category. If $Y_\bu$ is a quasicategory, then $\Fun_\bu(X_\bu,Y_\bu)$ is also a quasicategory \cite[\texttt{0066}]{kerodon}. In particular, when both $X_\bu$ and $Y_\bu$ are quasicategories, $\Fun_\bu(X_\bu,Y_\bu)$ is called the \emph{functor quasicategory} from $X_\bu$ to $Y_\bu$ \cite[\texttt{0067}]{kerodon}. Its objects are called \emph{functors}, and its morphisms are called \emph{natural transformations} \cite[\texttt{0068}, \texttt{0069}]{kerodon}. A natural transformation is a \emph{natural equivalence} or \emph{natural isomorphism} if it is an isomorphism in $\Fun_\bu(X_\bu,Y_\bu)$, equivalently if every component $\eta_x:F(x)\to G(x)$ is an isomorphism in $Y_\bu$ \cite[\texttt{01DG}, \texttt{01DK}]{kerodon}.

It is common in the literature to identify an ordinary category with its nerve and thereby regard it as a quasicategory \cite[\texttt{006A}]{kerodon}. When no confusion is likely, we omit the subscript $\bu$ and simply write $\Fun$ in place of $\Fun_\bu$.

Let $X_\bu$ be a simplicial set and let $x,y\in X_0$ be two vertices. Following \cite[\texttt{01J5}]{kerodon}, define the \emph{morphism simplicial set} from $x$ to $y$ by
\[
\Hom_{X}(x,y)
=
\{x\}
\times_{\Fun(\{0\},X)}
\Fun(\Delta^1,X)
\times_{\Fun(\{1\},X)}
\{y\},
\]
where the two maps from $\Fun(\Delta^1,X)$ are induced by the inclusions $\{0\},\{1\}\hookrightarrow\Delta^1$. If $X_\bu$ is a quasicategory, then $\Hom_{X}(x,y)$ is a Kan complex by \cite[\texttt{01JC}]{kerodon}, called the \emph{morphism space} from $x$ to $y$.
\edefn
\bdefn[{Homotopy equivalence}]\label{defn: homotopy equivalence}
Two morphisms $F_0,F_1:X_\bu\to Y_\bu$ of simplicial sets are said to be \emph{homotopic} if they determine the
same element of $\pi_0\big(\Fun(X_\bu,Y_\bu)\big)$ \cite[\texttt{00TS}]{kerodon}. In particular, if there exists a morphism
\[
H:\Delta^1\times X_\bu\to Y_\bu
\]
satisfying
\[
F_i=H|_{\{i\}\times X_\bu},
\text{ for } i=0,1,
\]
then $F_0$ and $F_1$ are homotopic. Such a morphism $H$ is called a
\emph{homotopy} from $F_0$ to $F_1$
\cite[\texttt{00TT}]{kerodon}. Conversely by \cite[\texttt{00TU}]{kerodon}, $F_0$ and $F_1$ are homotopic
if and only if there exists a finite sequence
\[
F_0=G_0,G_1,\ldots,G_n=F_1
\]
such that, for each $1\leq j\leq n$, there exists either a homotopy from
$G_{j-1}$ to $G_j$ or a homotopy from $G_j$ to $G_{j-1}$. If $Y_\bu$ is a Kan complex, then $F_0$ and $F_1$ are homotopic if and only if there exists a homotopy from $F_0$ to $F_1$.

Given a morphism $F:X_\bu\to Y_\bu$ of simplicial sets, a \emph{simplicial homotopy inverse} of $F$ is a
morphism $G:Y_\bu\to X_\bu$ such that $G\circ F$ and $F\circ G$ are homotopic to the respective identity morphisms. In this case, $F$ is called a \emph{simplicial
homotopy equivalence}, and $X_\bu$ and $Y_\bu$ are said to be
\emph{simplicially homotopy equivalent}
\cite[\texttt{00U2}]{kerodon}. Some authors refer to this notion simply
as a homotopy equivalence. We retain the adjective ``simplicial'' and
reserve the term \emph{homotopy equivalence} for the following stronger
notion \cite[\texttt{02BT}]{kerodon}. By \cite[\texttt{00UB}]{kerodon}, simplicial homotopy equivalence is a weak homotopy equivalence, while the converse holds for Kan complexes.

Let $F:X_\bu\to Y_\bu$ be a functor between quasicategories. A \emph{homotopy inverse} of $F$ is a functor $G:Y_\bu\to X_\bu$ such that $G\circ F$ and $F\circ G$, regarded as objects of $\Fun(X_\bu,X_\bu)$ and $\Fun(Y_\bu,Y_\bu)$, respectively, are
isomorphic to the corresponding identity functors
\cite[\texttt{01DY}]{kerodon}. In this case, we call $F$ a
\emph{homotopy equivalence}, and we say that $X_\bu$ and $Y_\bu$ are
\emph{homotopy equivalent}. In the standard terminology, such an $F$ is
an equivalence of quasicategories.

By \cite[\texttt{01EA}]{kerodon}, a functor between quasicategories is a
homotopy equivalence if and only if it is a categorical equivalence.
Equivalently, it is a weak equivalence in the Joyal model structure
$\QCat$; see also \cite[Theorem 2.2.5.1]{HTTlurie2009}. Every homotopy
equivalence is a simplicial homotopy equivalence
\cite[\texttt{01E3}]{kerodon}, and the two notions coincide when both
$X_\bu$ and $Y_\bu$ are Kan complexes
\cite[\texttt{01E2}]{kerodon}. However, a simplicial homotopy
equivalence between quasicategories need not be a homotopy equivalence.
The simplest counterexample is the projection map $\Delta^1\to\Delta^0$, which admits the inclusion $\Delta^0=\{0\}\hookrightarrow\Delta^1$ as a simplicial homotopy inverse, but is not an equivalence of
quasicategories \cite[\texttt{01E3}]{kerodon}.

Let $Y_\bu\subseteq X_\bu$ be quasicategories. We call $Y_\bu$ a \emph{deformation retract} of $X_\bu$ if there is a morphism $H:\Delta^1\times X_\bu\to X_\bu$ from the identity to a map factoring through $Y_\bu$, whose restriction to $\Delta^1\times Y_\bu$ is the projection followed by the inclusion, and whose component at every vertex is an isomorphism. The resulting retraction is a homotopy inverse to the inclusion \cite[\texttt{03TF}]{kerodon}.
\edefn

\blem[{\cite[\texttt{00U7}, \texttt{0193}, \texttt{00UB}, \texttt{00UD}, \texttt{02L2}]{kerodon}}]\label{lem: homotopy properties}
Let $X_\bu\xra{f}Y_\bu\xra{g}Z_\bu$ be morphisms of simplicial sets. If any two of $f$, $g$, and $g\circ f$ are simplicial homotopy equivalences, then so is the third. Simplicial homotopy equivalences are preserved by finite products and are weak homotopy equivalences. Likewise, if any two of $f$, $g$, and $g\circ f$ are weak homotopy equivalences, then so is the third. Weak homotopy equivalences are also preserved by finite products. We refer to the first and third assertions as the \emph{two-out-of-three properties} for simplicial homotopy equivalences and weak homotopy equivalences, respectively.
\elem

Finally, we introduce the definition of trivial Kan fibration of simplicial sets and some basic properties.

\bdefn[{\cite[\texttt{006W}]{kerodon}}]\label{defn: trivial kan fibration}
A map $p:X_\bu \to Y_\bu$ of simplicial sets is called a \emph{trivial Kan fibration} \ if, for every $n\in\N$, it has the right lifting property for the inclusion $\partial \Delta^n \to \Delta^n$, i.e.\, for any strict commutative diagram
\begin{equation*}\label{eq: trivial kan}
  \xymatrix{
\partial \Delta^n \ar[r] \ar[d] & X_\bu \ar[d]^{p} \\
\Delta^n \ar[r] \ar@{-->}[ur] & Y_\bu
}
\end{equation*}
a dotted arrow exists making the diagram commute.

A map $p:X_\bu \to Y_\bu$ of simplicial sets is called a \emph{trivial fibration} if it has the right lifting property for any monomorphism. By \cite[\texttt{006Y}]{kerodon}, this definition is equivalent to a trivial Kan fibration, so we do not distinguish them.
\edefn
We will use the following lemmas for the trivial Kan fibration.

\blem[{\cite[\texttt{006Z}, \texttt{00U9}]{kerodon}}]\label{lem: trivial Kan fibration}
Let $p : X_\bu \to Y_\bu$ be a trivial Kan fibration of simplicial sets. Then we have the following.
\benum
    \item The map $p$ admits a section: that is, there is a map of simplicial sets $s : Y_\bu \to X_\bu$ such that the composition $p \circ s$ is the identity map $\id_{Y_\bu} : Y_\bu \to Y_\bu$.
    \item Let $s$ be any section of $p$. Then the composition $s \circ p : X_\bu \to X_\bu$ is fiberwise homotopic to the identity. That is, there exists a map of simplicial sets $h : \Delta^1 \times X_\bu \to X_\bu$, compatible with the projection to $Y_\bu$, such that $h|_{ \{0\} \times X_\bu } = s \circ p$ and $h|_{\{1\} \times X_\bu } = \id_{X_\bu}$.
\eenu
In particular, a trivial Kan fibration $p:X_{\bu}\to Y_{\bu}$ is a simplicial homotopy equivalence.
\elem
\blem[{\cite[\texttt{02L0}, \texttt{00XE}]{kerodon}}]\label{lem: kan filtered colimits}
The collection of trivial Kan fibrations is closed under filtered colimits, when regarded as a full subcategory of the functor category $\Funo([1],\sSet)$. Similarly, the collection of weak homotopy equivalences is closed under filtered colimits.
\elem

\subsection{Complete Segal space}\label{subsec: complete segal space}
Then we consider simplicial spaces.
\bdefn[Simplicial space]\label{defn: simplicial space}

Let
\[
\sSpaceo=\Funo(\Delta^{\opp},\sSet)=\Funo(\Delta^{\opp}\times\Delta^{\opp},\Set)
\]
be the category of \emph{simplicial spaces}, also called \emph{bisimplicial sets} in some literature. Let $\sSpace_c$ and $\sSpace_f$ denote the category $\sSpaceo$ equipped with the injective and projective model structures, respectively \cite[Proposition A.3.3.2]{HTTlurie2009}. These model categories are Quillen equivalent via the identity functor, where the one from $\sSpace_f$ to $\sSpace_c$ is the left adjoint. The fibrant objects in the former are exactly those satisfying the Reedy fibrancy condition (called \emph{v-fibrant} in \cite{JT2007quasisegal}), while the fibrant objects in the latter are the levelwise fibrant ones, i.e.\ those for which each level is a Kan complex. Note that the Reedy fibrancy condition implies the levelwise Kan condition, but the converse is not true in general. By \cite[Lemma 7.4]{ED2014Reedy}, there is always a Reedy fibrant replacement, though we only need a simpler version from \cite[Proposition 4.10]{JT2007quasisegal} for simplicial spaces from quasicategories; see the end of this subsection. 

A simplicial space will be denoted by either $X_{\bu,\bu}$, where the first subscript denotes the simplicial level and the second subscript denotes the spatial level, or simply $X_\bu$, where each $X_k$ is a simplicial set. Sometimes we also omit all subscripts. By \emph{level} we mean simplicial level of a simplicial space.

Products of simplicial spaces are computed levelwise,
\[
(X\times Y)_{k,l}=X_{k,l}\times Y_{k,l},
\]
with all face and degeneracy maps acting componentwise.
\edefn

\bdefn[Homotopy fiber product]\label{defn: homotopy fiber product}
Let $f:X\to Z$ and $g:Y\to Z$ be maps of simplicial sets. Choose a Kan fibrant replacement $j:Z\to Z^{\mathrm{Kan}}$ and define the \emph{homotopy fiber product}
\begin{equation}\label{eq: path model}
    X\times_Z^hY
=
X\times_{\Fun(\{0\},Z^{\mathrm{Kan}})}
\Fun(\Delta^1,Z^{\mathrm{Kan}})
\times_{\Fun(\{1\},Z^{\mathrm{Kan}})}Y
\end{equation}
using $j\circ f$ and $j\circ g$. Its weak homotopy type is independent of $j$ by \cite[\texttt{010B}, \texttt{0111}, \texttt{032E}]{kerodon}. If $Z$ is not Kan, the analogous construction in \eqref{eq: path model} formed directly over $Z$ is only the \emph{oriented fiber product}; see \cite[\texttt{0327}, \texttt{01KF}]{kerodon}. We also refer \eqref{eq: path model} as the \emph{path model} of the homotopy fiber product.

A commutative square
\[
\xymatrix{
W\ar[r]\ar[d]&X\ar[d]^f\\
Y\ar[r]_g&Z
}
\]
is a \emph{homotopy pullback square}, or a \emph{homotopy Cartesian square}, if the induced map
\[
W\longrightarrow X\times_Z^hY
\]
is a weak homotopy equivalence. Equivalently, choose a factorization
\[
X\xra{w}X'\xra{p}Z,
\]
where $w$ is a weak homotopy equivalence and $p$ is a Kan fibration. Such a factorization exists by \cite[\texttt{00UV}]{kerodon}, and the square is homotopy Cartesian if and only if
\[
W\longrightarrow X'\times_ZY
\]
is a weak homotopy equivalence; see \cite[\texttt{010L}, \texttt{010N}]{kerodon}. This agrees with \cite[\S 2.3]{Rezk1998fibrations}.

If the square is a strict pullback, so that $W=X\times_ZY$, it is homotopy Cartesian precisely when
\[
X\times_ZY\longrightarrow X\times_Z^hY
\]
is a weak homotopy equivalence. In particular, this holds if $f$ is a Kan fibration.
\edefn

\blem[{\cite[\texttt{0111}]{kerodon}}]\label{lem: homotopy invariance fiber product}
Suppose that
\[
\xymatrix{
X\ar[r]\ar[d]^{w_X}&Z\ar[d]^{w_Z}&Y\ar[l]\ar[d]^{w_Y}\\
X'\ar[r]&Z'&Y'\ar[l]
}
\]
is commutative and $w_X,w_Y,w_Z$ are weak homotopy equivalences. After choosing compatible models as in Definition \ref{defn: homotopy fiber product}, the induced map
\[
X\times_Z^hY\longrightarrow X'\times_{Z'}^hY'
\]
is a weak homotopy equivalence.
\elem

\bdefn[{Segal space \cite[\S 1.3-1.4]{CS2019category}}]\label{defn: segal space}
A \emph{Segal space} is a fibrant object in $\sSpace_f$, i.e.\ a levelwise Kan complex, that further satisfies the \ \emph{Segal condition}. For any $n,m\in\N$, the commutative square
\begin{equation}\label{eq: Segal}
\xymatrix{
X_{m+n}\ar[r]\ar[d] & X_m\ar[d]\\
X_n\ar[r] & X_0
}
\end{equation}
induced by the maps
\[
[m]\to[m+n],\quad (0<\cdots<m)\mapsto(0<\cdots<m),
\]
\[
[n]\to[m+n],\quad (0<\cdots<n)\mapsto(m<\cdots<m+n),
\]
\[
[0]\to[m],\quad (0)\mapsto(m),
\]
and
\[
[0]\to[n],\quad (0)\mapsto(0),
\]
is a homotopy pullback square.

There is also an alternative description of the Segal condition. For $m\in\N_+$ and $i\in[1,m]$, define
\[
g_i:[1]\to[m],\quad (0<1)\mapsto(i-1<i).
\]
Then the map
\[
X_m\to X_1\times^h_{X_0}\cdots\times^h_{X_0}X_1
\]
induced by $g_1,\dots,g_m$ is a \ weak homotopy equivalence.

Although a simplicial space need not be levelwise Kan, the Segal condition is still well-defined, since homotopy fiber products can be formed in the category of simplicial sets. We call a simplicial space satisfying the Segal condition a \emph{Segal simplicial set}.

Following \cite[Warning 2.1.16]{lurie2009cob} and \cite[Remark 1.6]{CS2019category}, we also define a \emph{Reedy fibrant Segal space} to be a fibrant object in $\sSpace_c$ satisfying the Segal condition. The Reedy fibrancy condition implies that the canonical map
\[
X_m\times_{X_0}X_n\to X_m\times^h_{X_0}X_n
\]
is a weak homotopy equivalence, so the homotopy fiber product in the Segal condition can be replaced by the strict fiber product.

Moreover, there are model structures $\sSpace_c^{\Se}$ and $\sSpace_f^{\Se}$ on $\sSpaceo$, obtained from $\sSpace_c$ and $\sSpace_f$ by left Bousfield localizations, such that their fibrant objects are Reedy fibrant Segal spaces and Segal spaces, respectively. In particular, every levelwise weak homotopy equivalence is a weak equivalence.  These model categories are Quillen
equivalent via the identity functor, where the one from $\sSpace_f^{\Se}$ to $\sSpace_c^{\Se}$ is the left adjoint.  By \emph{equivalence} of (Reedy fibrant) Segal spaces we mean weak equivalence in those model categories.
\edefn
\bdefn[{Morphism space, \cite[Definition 1.9]{CS2019category}}]\label{defn: morphism space segal}
For a Segal space $X=X_{\bu,\bu}$, elements in $X_{0,0}$ are called \emph{objects}. For two objects $x,y\in X_{0,0}$, the morphism space is\[\Hom_X(x,y)=\{x\}\times_{X_0}^h X_1\times_{X_0}^h\{y\},\]which is a Kan complex.

The \emph{homotopy category} $h(X)$ of a Segal space $X=X_{\bu,\bu}$ is the ordinary category whose objects are the elements of $X_{0,0}$ and whose morphisms between objects $x,y$ are
\[
\Hom_{h(X)}(x,y)=\pi_0(\Hom_X(x,y)).
\]
For objects $x,y,z$, the following diagram induces the composition of morphisms, as a weak homotopy equivalence between Kan complexes is a simplicial homotopy equivalence and hence induces a bijection on $\pi_0$:
\begin{multline*}
\bigl(\{x\}\times_{X_0}^h X_1\times_{X_0}^h\{y\}\bigr)\times\bigl(\{y\}\times_{X_0}^h X_1\times_{X_0}^h\{z\}\bigr)\\
\longrightarrow\{x\}\times_{X_0}^h X_1\times_{X_0}^h X_1\times_{X_0}^h\{z\}
\xleftarrow{\simeq}\{x\}\times_{X_0}^h X_2\times_{X_0}^h\{z\}
\longrightarrow\{x\}\times_{X_0}^h X_1\times_{X_0}^h\{z\}.
\end{multline*}

A map $f:X_{\bu,\bu}\to Y_{\bu,\bu}$ of Segal spaces is a \emph{DK (Dwyer--Kan) equivalence} if the induced map on homotopy categories is essentially surjective, and for each pair of objects $x,y\in X_{0,0}$ the induced map
\[\Hom_{X}(x,y)\to\Hom_{Y}(f(x),f(y))\]is a \ weak homotopy equivalence.
\edefn
\brem
The notion of a simplicial space is taken from \cite{Rezk2001completion}, where a \emph{space} is simply a simplicial set. The space in \cite{CS2019category} denotes a Kan complex. Although the fibrancy condition is not explicitly mentioned in \cite[Definition 1.4]{CS2019category}, it is implicit in the statement that the fibrant objects in $\sSpace_f^{\Se}$ are Segal spaces. In \cite[Definition 2.1.8]{lurie2009cob}, a simplicial space is defined to be a functor $X_\bu:\Delta^{\opp}\to \mathrm{Top}$. From Example \ref{exmp: singular complex}, Kan complexes are simplicial analogues of topological spaces.
\erem
\bdefn[Complete Segal space]\label{defn: complete segal space}

For a Segal space $X_\bu$, let $X_1^{\simeq}\subset X_1$ be the union of those connected components whose points represent isomorphisms in the homotopy category $h(X_\bu)$. Since $s_0(x)$ represents $\id_x$, the degeneracy map factors through $X_1^{\simeq}$. The Segal space $X_\bu$ is \emph{complete} if $s_0:X_0\to X_1^{\simeq}$ is a weak homotopy equivalence; see \cite[Definition 1.12]{CS2019category}. A Reedy fibrant complete Segal space is a Reedy fibrant Segal space satisfying this condition. For every Reedy fibrant Segal space $X$, Rezk \cite[\S 14]{Rezk2001completion} constructed functorially a Reedy fibrant complete Segal space $\widehat{X}$ together with a DK equivalence $i_X:X\to \widehat{X}$, where $\widehat{X}$ is called the \emph{(Rezk) completion} of $X$. Here the completion is obtained by first taking some explicit simplicial space $\widetilde{X}$ and then taking its Reedy fibrant replacement. 

For a (non-Reedy fibrant) Segal space or Segal simplicial set, we can also take its Reedy fibrant replacement and then apply the completion construction. To be compatible with \cite[\Def 1.18]{CS2019category}, we still call this operation on a Segal space \emph{completion}.

There are model structures $\sSpace_c^{\CSe}$ and $\sSpace_f^{\CSe}$ on $\sSpaceo$ obtained from $\sSpace_c^{\Se}$ and $\sSpace_f^{\Se}$ by some further left Bousfield localizations such that fibrant objects are Reedy fibrant complete Segal spaces and complete Segal spaces, respectively. Those model categories are again Quillen
equivalent via the identity functor, where the one from $\sSpace_f^{\CSe}$ to $\sSpace_c^{\CSe}$ is the left adjoint. The completion of a Segal space also provides the fibrant replacement in $\sSpace_f^{\CSe}$. From \cite[Theorem 1.15]{CS2019category}, a morphism between Segal spaces is a weak equivalence in $\sSpace_f^{\CSe}$ if and only if it is a DK equivalence. On a Segal simplicial set, we denote the fibrant replacement in $\sSpace_c^{\CSe}$ by $\widehat{R}$, which is the composition of the Reedy fibrant replacement and the completion as above.
\edefn
\bdefn\label{defn: discrete and constant}
There are two standard ways to regard a simplicial set $X_\bu$ as a simplicial space. These constructions are often left implicit in \cite{Rezk2001completion,CS2019category}; here we distinguish them notationally.
\begin{itemize}
    \item The \emph{discrete simplicial space} $\disc X$ is defined by
    \[
    (\disc X)_{k,l}=X_k.
    \]
    It agrees with $p_1^*X_\bu$ from \eqref{eq: i1 p1} and regards $X_\bu$ as a simplicial space with discrete levels. By \cite[Example 1.8]{CS2019category}, it satisfies the Segal condition if and only if $X_\bu$ is the nerve of an ordinary category. It is always Reedy fibrant, and hence levelwise Kan, but it does not in general satisfy the completeness condition.
    \item The \emph{constant simplicial space} $cX$ is defined by
    \[
    (cX)_{k,l}=X_l.
    \]
    This construction is commonly used when $X_\bu$ is a Kan complex, for example the singular complex of a topological space. The Segal maps formed using ordinary fiber products, as well as the completeness map, are isomorphisms. Hence $cX$ is a complete Segal space whenever $X_\bu$ is a Kan complex. It need not be Reedy fibrant.
\end{itemize}
For a quasicategory $X_\bu$, we introduce a third associated simplicial space from \cite[p.301]{JT2007quasisegal}\footnote{The reference used the notation $\Ga$ instead of $\Psi$. We reserve $\Ga$ for $\Ga$-object in Definition \ref{defn: sym str Ga object}.}
\begin{equation}\label{eq: Reedy fiber replacement}
    \Psi(X_\bu)_{k,\bu}=\Fun(\Delta^k,X_\bu)^\simeq,
\end{equation}which is a Reedy fibrant complete Segal space and, by \cite[Proposition 4.10]{JT2007quasisegal}, a fibrant replacement of $\disc X_\bu$ in $\sSpace_c^{\CSe}$. This construction is functorial in $X_\bu$, with $\Psi(F)$ for a functor $F$ induced levelwise by postcomposition.

For simplicial sets $K$ and $L$, define their \emph{external product} $K\boxtimes L=\disc K\times cL$, so that
\begin{equation}\label{eq: external product}
(K\boxtimes L)_{m,n}=K_m\times L_n.
\end{equation}
For a simplicial space $X$, there is a natural identification
\[
\Hom_{\sSpaceo}(\Delta^k\boxtimes L,X)
\cong\Hom_{\sSet}(L,X_{k,\bu}).
\]
In particular, by \eqref{eq: Reedy fiber replacement}, a morphism $\Delta^k\boxtimes\Delta^l\to\Psi(\D)$ is equivalently a simplicial map
\[
\Delta^l\longrightarrow\Fun(\Delta^k,\D)^\simeq.
\]
\edefn
\bexmp\label{exmp: concrete completion}
A useful example illustrating both the distinction and the relation between the constructions in Definition \ref{defn: discrete and constant} is provided by \cite[Remark 14.1 and Lemma 14.2]{Rezk2001completion}. Let $G$ be a group, regarded as a category $\B G$ with one object with morphisms in $G$, and let $\B_\bu G=\rN_\bu(\B G)$ denote its nerve. Then $\disc\B_\bu G$ is a Reedy fibrant Segal space. There are natural maps
\[
\disc\B_\bu G\longrightarrow\Psi(\B_\bu G)\longleftarrow c\B_\bu G,
\]
where the left-hand map is a Rezk weak equivalence (i.e.\ in $\sSpace_c^{\CSe}$) and the right-hand map is a Reedy weak equivalence (i.e.\ in $\sSpace_c$, and hence also in $\sSpace_c^{\CSe}$). Thus, $\Psi(\B_\bu G)$ simultaneously provides a completion of $\disc\B_\bu G$ and a Reedy fibrant replacement of $c\B_\bu G$. Since $c\B_\bu G$ is a complete Segal space, we can also regard it as the fibrant model of $\disc\B_\bu G$ in $\sSpace_f^{\CSe}$.
\eexmp
As for quasicategories, natural transformations between morphisms of Segal spaces can be defined using the discrete simplicial interval $\disc\Delta^1$.
\bdefn[Natural transformations]\label{defn: natural transformation segal}
Let $F_0,F_1:\cC\to\cD$ be morphisms of Segal spaces. A \emph{natural transformation} $\eta:F_0\Rightarrow F_1$ is a morphism
\[
\eta:\cC\times\disc\Delta^1\longrightarrow\cD
\]
whose restrictions to $\cC\times\{i\}$ are $F_i$ for $i=0,1$. It is a \emph{natural equivalence}, or \emph{natural isomorphism}, if, for every object $X$, the component
\[
\eta_X:F_0(X)\longrightarrow F_1(X)
\]
represents an isomorphism in $h(\cD)$.
\edefn
Finally, we state the Quillen
equivalence between $\QCat$ and $\sSpace_c^{\CSe}$ from \cite{JT2007quasisegal}. We consider the functors
\[
\begin{aligned}
i_1:\Delta &\longrightarrow \Delta\times\Delta,
&
p_1:\Delta\times\Delta &\longrightarrow \Delta,\\[4pt]
[k] &\longmapsto ([k],[0]),
&
([k],[l]) &\longmapsto [k],
\end{aligned}
\]
which induce pullback functors\begin{equation}\label{eq: i1 p1}\begin{aligned}
    i_1^*:\sSpaceo &\longrightarrow \sSet,
&
p_1^*:\sSet &\longrightarrow \sSpaceo,\\[4pt]
X_{\bu,\bu} &\longmapsto X_{\bu,0},
&
X_\bu &\longmapsto \disc X_\bu.
\end{aligned}
\end{equation}
\bthm[{\cite[Theorem 4.11]{JT2007quasisegal}}]\label{thm: quillen equ}
The adjunction\[p_1^*:\sSet\rightleftarrows\sSpaceo:i_1^*\]is a Quillen equivalence between the model categories $\QCat$ and $\sSpace_c^{\CSe}$, with $p_1^*$ the left adjoint.
\ethm

In particular, the pair $(p_1^*,i_1^*)$ is a Quillen adjunction, which means $p_1^*$ preserves cofibrations and $i_1^*$ preserves fibrations. However, the functor $p_1^*$ does not necessarily preserve the fibrations, which means $p_1^*X_\bu$ for a quasicategory $X_\bu$ may not be a Reedy fibrant complete Segal space\footnote{In the first arXiv version of this paper, we incorrectly regarded $p_1^*X_\bu$ as a complete Segal space and used it as the target of the $\infty$-categorical TQFT constructed in \S\ref{sec: Construction of the TQFT functor}. To correct this, we replace the target as described below and accordingly modify the construction on higher spatial simplices in \S\ref{subsec: Morphisms and higher morphisms}.}. 

To resolve this issue, we use the Reedy fibrant complete Segal space $\Psi(X_\bu)$ from \eqref{eq: Reedy fiber replacement} for a quasicategory $X_\bu$, which is a fibrant replacement of $p_1^*X_\bu$ in $\sSpace_c^{\CSe}$. Note that by construction,
\begin{equation}\label{eq: Ga 0}
    (i_1^*\Psi(X_\bu))_k=\Psi(X_\bu)_{k,0}=X_k.
\end{equation}
Thus $i_1^*\Psi(X_\bu)\cong X_\bu$.
\brem\label{rem: construct functor}
The Quillen equivalence in Theorem \ref{thm: quillen equ} implies that, for a Reedy fibrant complete Segal space $X_{\bu,\bu}$, all its homotopical information is encoded in its first row $X_{\bu,0}$, which is a quasicategory. Moreover, by \cite[Corollary 3.6]{JT2007quasisegal}, every row $X_{\bu,l}$ of a Reedy fibrant Segal space is a quasicategory. In practice, however, to construct a functor from the $\infty$-category represented by a Segal space $X_{\bu,\bu}$ to a quasicategory $Y_\bu$, it is often easier to construct a map
\[
X_{\bu,\bu}\longrightarrow\Psi(Y_\bu)
\]
and then use the universal property of the completion to obtain a map between complete Segal spaces
\[
\widehat{X}\longrightarrow\Psi(Y_\bu),
\]where $\widehat{X}$ denotes the completion of a Reedy fibrant replacement of $X$. Such an induced map is unique up to contractible space of choices, following \cite[Proposition 9.1(3)]{Rezk2001completion} and the fact that the completion map is a weak equivalence in $\sSpace_c^{\CSe}$.

Constructing the corresponding functor directly at the level of quasicategories is less convenient, since one needs to consider the concrete construction of $\widehat{X}_{\bu,0}$. However, the Reedy fibrant replacement from \cite[Lemma 7.4]{ED2014Reedy} is not a closed formula.
\erem

\subsection{Simplicial category}\label{subsec: simplicial category}
In this subsection, we introduce simplicial categories as a third model of $\infty$-categories for describing symmetric monoidal structures.

\bdefn[Simplicial category]\label{defn: simplicial category}
A \emph{simplicial category} is a category enriched in simplicial sets. A \emph{simplicial functor} assigns objects and maps of simplicial sets between mapping spaces, preserving composition and identities. Write $\sCato$ for the category of small simplicial categories and simplicial functors. Its Bergner model structure from \cite[Proposition A.3.2.4 and Theorem A.3.2.24]{HTTlurie2009} is denoted by $\sCat$; the fibrant objects are the \emph{locally Kan} simplicial categories, whose mapping spaces are Kan complexes.

We define the simplicial category $\Cati^\Delta$ with small quasicategories as objects and
\[
\Hom_{\Cati^\Delta}(X_\bu,Y_\bu)=\Fun_\bu(X_\bu,Y_\bu)^\simeq,
\]
where $\Fun_\bu$ is the functor simplicial set from Definition \ref{defn: functor simplicial set} and $(-)^\simeq$ is its core from Definition \ref{bdefn: quasicategory}. Thus $\Cati^\Delta$ is locally Kan. 
\edefn

\bdefn[Homotopy coherent realization]\label{defn: simplicial resolution}\label{defn: homotopy coherent realization}
For a simplicial set $X_\bu$, its \emph{homotopy coherent realization} is
\[
\fC[X_\bu]=\colim_{\Delta^r\to X_\bu}\fC[\Delta^r],
\]
where the colimit is taken in $\sCato$. Here $\fC[\Delta^r]$ has objects $[r]$ and, for $i\leq j$,
\[
\Hom_{\fC[\Delta^r]}(i,j)=\rN_\bu\{T\subseteq\{i,\ldots,j\}\mid i,j\in T\},
\]
with the subsets ordered by inclusion; the mapping space is empty for $i>j$. Composition is induced by union. A monotone map $\nu:[r]\to[s]$ acts by $i\mapsto\nu(i)$ and $T\mapsto\nu(T)$; see \cite[Definitions 1.1.5.1 and 1.1.5.3]{HTTlurie2009}.

For a small ordinary category $\cC$, write $\fC[\cC]=\fC[\rN_\bu\cC]$, using the ordinary nerve of Definition \ref{bdefn: quasicategory}. Its objects are those of $\cC$. A string $\rho=(x_0\xra{f_1}\cdots\xra{f_r}x_r)$, with $r\geq1$, gives a \emph{cubical chart}
\[
\Hom_{\fC[\Delta^r]}(0,r)=(\Delta^1)^{r-1}\longrightarrow\Hom_{\fC[\cC]}(x_0,x_r).
\]
Its vertices are the \emph{cut sets} $C=T\setminus\{0,r\}\subseteq\{1,\ldots,r-1\}$, and its $u$-simplices are chains $C_0\subseteq\cdots\subseteq C_u$, allowing repetitions. Each cut set gives the consecutive block composites, with identities omitted. Composition concatenates strings and inserts the joining index into every cut set. The faces $t_i=0,1$ compose adjacent arrows and the two shorter charts, respectively. Deleting an internal identity replaces its adjacent cut coordinates by $\max(t_{i-1},t_i)$; at an endpoint it drops the adjacent coordinate.

Regard $\cC$ as a simplicial category with constant mapping simplicial sets. The \emph{augmentation}
\begin{equation}\label{eq: augmentation}
    \epsilon_\cC:\fC[\cC]\longrightarrow\cC
\end{equation}
is the identity on objects and sends every simplex in the chart of $\rho$ to the constant simplex at $f_r\cdots f_1$; the empty string gives the identity.

\edefn

\bdefn[Homotopy coherent nerve]\label{defn: homotopy coherent nerve}
For a simplicial category $\cD$, its \emph{homotopy coherent nerve} is the simplicial set
\[
\Nhc_\bu(\cD)_r=\Hom_{\sCato}(\fC[\Delta^r],\cD),
\]
with simplicial operators induced by precomposition with $\fC[\nu]$ \cite[Definition 1.1.5.5]{HTTlurie2009}. In particular, the homotopy coherent nerve of an ordinary category $\cC$ is the ordinary nerve $\rN_\bu \cC$. It is right adjoint to $\fC[-]$. If $\cD$ is locally Kan, then $\Nhc_\bu(\cD)$ is a quasicategory \cite[Proposition 1.1.5.10]{HTTlurie2009}. In particular, we define $\Cati=\Nhc_\bu(\Cati^\Delta)$ \cite[Definition 3.0.0.1]{HTTlurie2009}.

\edefn

The following theorem gives simplicial categories as another model of $\infty$-categories.
\bthm[{\cite[Theorem 2.2.5.1]{HTTlurie2009}}]\label{thm: coherent nerve Quillen equivalence}
The adjunction
\[
\fC[-]:\sSet\rightleftarrows\sCato:\Nhc_\bu
\]
is a Quillen equivalence between the model categories $\QCat$ and $\sCat$, with $\fC[-]$ the left adjoint.
\ethm

\bcor\label{cor: augmentation of coherent realization}
Let $\cC$ be a small ordinary category. The augmentation $\epsilon_\cC:\fC[\cC]\to\cC$ is a weak equivalence in $\sCat$.
\ecor
\bpf
Since $\cC$ is fibrant in $\sCat$ and its homotopy coherent nerve is cofibrant in $\QCat$, Theorem \ref{thm: coherent nerve Quillen equivalence} implies that $\epsilon_\cC$ is a weak equivalence in $\sCat$, as its adjoint is the identity of that nerve. An alternative proof is from \cite[Proposition 1.2.7 and Theorem 2.2.3]{riehl2023coherentstructures}.
\epf
We also use the following rectification lemma for symmetric monoidal structures.
\blem[Rectification]\label{lem: rectification of coherent diagrams}
Let $\cC$ be a small ordinary category and let $\fB:\fC[\cC]\to\Cat_\infty^\Delta$ be a simplicial functor. There exists a functor $B:\cC\to\QCat$ with quasicategory values and a zigzag of objectwise weak equivalences in $\QCat$
\[
\epsilon_\cC^*B\simeq\fB.
\]

\elem
\bpf
From Corollary \ref{cor: augmentation of coherent realization}, the augmentation $\epsilon_\cC:\fC[\cC]\to\cC$ is a weak equivalence in $\sCat$. Applying the rectification argument of \cite[Corollary 3.1.4.4 and Proposition A.3.3.8]{HTTlurie2009} to quasicategory-valued diagrams gives a strict functor $B:\cC\to\QCat$ with quasicategory values and the required zigzag.
\epf
\subsection{Symmetric monoidal structure}

In this subsection, we introduce symmetric monoidal structures. We use the notion of a \emph{$\Ga$-object}, also called a \emph{commutative monoid object}, as in \cite[Remarks 2.4.2.2 and 2.4.2.3]{HAlurie2017}.

\bdefn[Symmetric monoidal structure]\label{defn: monoidal category}
An (ordinary) category $\cC$ is called a \emph{monoidal category} if there is a tensor product operation $\otimes : \cC \times \cC \to \cC$ that is unital and associative (but not necessarily commutative) up to coherent isomorphism. If the operation is also commutative up to coherent isomorphism, then $\cC$ is called a \emph{symmetric monoidal category}. A monoidal category $\cC$ is called \emph{closed} if for each object $Y$ of $\cC$, the functor $-\ot Y:\cC\to \cC$ has a right adjoint functor $[Y,-]:\cC\to \cC$, called the \emph{internal-hom functor}, i.e.\, there is a natural isomorphism of the hom functors\[\Hom_{\cC}(-\ot Y,-)\simeq \Hom_{\cC}(-,[Y,-]):\cC^{\opp}\times \cC\to \Set\]

A \emph{(symmetric) monoidal dg category} is defined similarly, with the tensor product a dg functor $\sfC\otimes_{\mathrm{dg}}\sfC\to\sfC$ and the structural isomorphisms closed degree-zero dg-natural transformations. It is \emph{closed} if each $-\otimes Y$ has a dg right adjoint $[Y,-]$, so that the adjunction gives an isomorphism of morphism complexes.
\edefn
\bdefn[{$\Ga$-object, \cite[\Defs 3.1, 3.3, and 3.6]{CS2019category}}]\label{defn: sym str Ga object}
Let $\Ga$ be \emph{Segal's category} \cite{Segal1974}\footnote{Our convention is opposite to Segal's original one; see \cite[Remark 3.2]{CS2019category}.}, whose objects are $\langle m\rangle=\{0,\ldots,m\}$ for $m\in\N$, pointed at $0$, and whose morphisms are basepoint-preserving functions. Thus $\Ga$ is a skeleton of $\mathrm{Fin}_*$. The \emph{Segal morphisms} are
\begin{equation}\label{eq: segal morphism}
    \ga_\beta:\langle m\rangle\to\langle1\rangle,~j\mapsto\delta_{\beta j},~1\leq\beta\leq m.
\end{equation}

Let $\cM$ be a model category. A \emph{$\Ga$-diagram} in $\cM$ is a functor $A:\Ga\to\cM$. It is called a \emph{$\Ga$-object} if
\begin{equation}\label{eq: A condition}
\prod_{\beta=1}^m A(\ga_\beta):A\langle m\rangle\longrightarrow(A\langle1\rangle)^m
\end{equation}
is a weak equivalence in $\cM$ for every $m\geq0$. We call this the \emph{$\Ga$-condition}. The maps
\[
A(\theta):A\langle m\rangle\longrightarrow A\langle n\rangle,~\theta:\langle m\rangle\to\langle n\rangle,
\]
are called the \emph{$\Ga$-structure maps} of $A$. We call $A\langle1\rangle$ its \emph{underlying object}, and also call $A\langle1\rangle$ \emph{symmetric monoidal} when $A$ is understood. The case $m=0$ gives $A\langle0\rangle\simeq *$, and the map induced by $\langle0\rangle\to\langle1\rangle$ determines the unit.

A \emph{symmetric monoidal functor} is a morphism of $\Ga$-objects, namely a natural transformation $A\to B$, which therefore commutes with all $\Ga$-structure maps. We also call it a \emph{$\Ga$-morphism}.
\edefn
\brem\label{rem: equivalence of ga object}
From \cite[Corollary 2.13]{Ste21tqft}, if the values of a $\Ga$-object $A$ are all quasicategories or all complete Segal spaces, then it induces a symmetric monoidal structure on the homotopy category of its underlying object. A $\Ga$-morphism induces a symmetric monoidal functor between their homotopy categories. Moreover, from \cite[Theorem 2.12]{Ste21tqft}, the equivalence of models in Theorem \ref{thm: quillen equ} identifies the corresponding $\infty$-categories of $\Ga$-objects in quasicategories and complete Segal spaces.
\erem
In instanton construction, strict compatibility with the $\Ga$-structure maps is difficult to achieve. We use the following way to construct a symmetric monoidal functor via a homotopy coherent natural transformation.

\bdefn[{\cite[Definition 1.3.5]{riehl2023coherentstructures}}]\label{defn: homotopy coherent Gamma transformation}
The category $\sSpaceo$ of Definition \ref{defn: simplicial space} is simplicially enriched in the spatial level by
\[
\Hom_{\sSpaceo}(X,Y)_r=\Hom_{\sSpaceo}(\Delta^r\otimes X,Y),~(Z\otimes X)_{k,l}=Z_l\times X_{k,l}.
\]
Composition uses the diagonal of $\Delta^r$. For a small ordinary category $\cC$, a \emph{homotopy coherent diagram} is a simplicial functor $\fC[\cC]\to\sSpaceo$, with $\fC[\cC]$ as in Definition \ref{defn: homotopy coherent realization}. A \emph{homotopy coherent natural transformation} between two such diagrams $A,B$ is a simplicial functor
\[
H:\fC[\cC\times[1]]\longrightarrow\sSpaceo
\]
restricting to $A$ and $B$ along the endpoint inclusions $\fC[\cC]\to\fC[\cC\times[1]]$, denoted by $H:A\Rightarrow B$.

An \emph{ordinary diagram} $A:\cC\to\sSpaceo$ is a functor of ordinary categories, regarded as the coherent diagram $A\epsilon_\cC$ via the augmentation in \eqref{eq: augmentation}. Taking $\cC=\Ga$ applies this definition to the $\Ga$-diagrams of Definition \ref{defn: sym str Ga object}.
\edefn

\blem\label{lem: coherent Gamma morphism}
Let $A,B:\Ga\to\sSpace_f$ be $\Ga$-objects whose values satisfy the Segal condition. A homotopy coherent natural transformation $H:A\Rightarrow B$ in Definition \ref{defn: homotopy coherent Gamma transformation} induces a $\Ga$-morphism
\[
\widehat R A\longrightarrow\widehat R B,
\]
hence a symmetric monoidal functor in the sense of Definition \ref{defn: sym str Ga object}. Here $\widehat R$ is a functorial fibrant replacement in the injective model structure on $(\sSpace_c^{\CSe})^\Ga_c$ \cite[Proposition A.3.3.2]{HTTlurie2009}, taken at diagram level; its values are complete Segal space replacements as in Definition \ref{defn: complete segal space}. The underlying functor represents $H(\id_{\langle1\rangle},0\to1)$.
\elem
\bpf
Apply the rectification argument in the proof of Lemma \ref{lem: rectification of coherent diagrams} in $\sSpace_f$, with $\cC=\Ga\times[1]$. This gives a strict $\Ga$-morphism $A'\to B'$ with $A'\simeq A$ and $B'\simeq B$ through zigzags of objectwise levelwise weak equivalences. These preserve the $\Ga$- and Segal conditions. Applying $\widehat R$ gives the assertion, since completion preserves the $\Ga$-condition by \cite[Lemma 3.7]{CS2019category}. Evaluation at $(\id_{\langle1\rangle},0\to1)$ gives the stated underlying functor after completion and transport along the endpoint equivalences.
\epf

\subsection{Rigidity and duality}
In this subsection, we introduce the definitions about the rigidity and the duality, and study the interaction with symmetric monoidal structures.

\bdefn[Rigidity]\label{defn: rigidity}
Let $\cC$ be a symmetric monoidal category. An object $X\in\cC$ is \emph{dualizable} if there exist an object $X^\vee$ and morphisms
\[
\ev_X:X\otimes X^\vee\to\mathbf{1},~
\coev_X:\mathbf{1}\to X^\vee\otimes X
\]
satisfying the triangle identities
\begin{equation}\label{eq: dual comp}
\begin{aligned}
(\ev_X\otimes\id_X)\circ(\id_X\otimes\coev_X)&=\id_X,\\
(\id_{X^\vee}\otimes\ev_X)\circ(\coev_X\otimes\id_{X^\vee})&=\id_{X^\vee}.
\end{aligned}
\end{equation}
We call $X^\vee$ a \emph{dual} of $X$. The category $\cC$ is \emph{rigid} if every object is dualizable.

For a symmetric monoidal (complete) Segal space or quasicategory $\cC$, we call $\cC$ \emph{rigid} if $h(\cC)$ is rigid, following \cite[Definition 2.3.5]{lurie2009cob}.
\edefn
\brem
Lurie also defines dualizability directly for objects of a monoidal quasicategory \cite[Definition 4.6.1.7]{HAlurie2017}. In the symmetric monoidal case, this agrees with dualizability in its homotopy category by \cite[Lemma 4.6.1.6 and Remark 4.6.1.12]{HAlurie2017}. For a fixed dualizable object, the space of choices of a dual together with its duality data is contractible by \cite[Lemma 4.6.1.10]{HAlurie2017}.
\erem

\bdefn[Duality]\label{defn: dual}
A \emph{category with duality} is a category $\cC$ equipped with a functor
\[
(-)^\vee:\cC^{\opp}\longrightarrow\cC
\]
and a natural isomorphism $\eta:\id_\cC\to(-)^{\vee\vee}$ satisfying
\begin{equation}\label{eq: bidual compatibility}
\eta_X^\vee\circ\eta_{X^\vee}=\id_{X^\vee}.
\end{equation}
We call $(-)^\vee$ the \emph{duality functor} and $\eta$ the \emph{bidual isomorphism}. Evaluation and coevaluation are not part of this structure. The duality is \emph{strict} if $\eta=\id$, so in particular $(-)^{\vee\vee}=\id_\cC$.

A \emph{duality-preserving functor} $F:\cC\to\cD$ is equipped with a natural isomorphism
\[
\xi_X:F(X^\vee)\longrightarrow F(X)^\vee
\]
satisfying
\begin{equation}\label{eq: duality-preserving functor}
\xi_{X^\vee}\circ F(\eta_X^\cC)=\xi_X^\vee\circ\eta_{F(X)}^\cD.
\end{equation}

A \emph{dg category with duality} is defined similarly, using the graded dg opposite $\sfC^{\opp}$, a dg functor $(-)^\vee$, and a closed degree-zero dg-natural isomorphism $\eta$ satisfying \eqref{eq: bidual compatibility} strictly. A \emph{duality-preserving dg functor} consists of a dg functor $F$ and a closed degree-zero dg-natural isomorphism $\xi$ satisfying \eqref{eq: duality-preserving functor} strictly. Again, the duality is \emph{strict} if $\eta=\id$.

A \emph{symmetric monoidal (dg) category with duality} is a symmetric monoidal (dg) category whose duality functor is symmetric monoidal and whose bidual isomorphism $\eta$ is monoidal. A duality-preserving (dg) functor between such categories is \emph{symmetric monoidal} if $F$ is symmetric monoidal and $\xi$ is monoidal.

For a complete Segal space or quasicategory $\cC$, a \emph{duality} is understood in the homotopy coherent sense of \cite[Definition 2.3]{HLS2025}. We denote its underlying contravariant equivalence by \[(-)^\vee:\cC^{\opp}\to\cC\]and its bidual natural equivalence by \[\eta:\id_\cC\Rightarrow(-)^{\vee\vee}.\]Equation \eqref{eq: bidual compatibility} holds with its higher coherences. A functor \emph{preserves duality} if it is a morphism of these structures. Its underlying comparison is a natural equivalence \[\xi:F\circ(-)^\vee\Rightarrow(-)^\vee\circ F^{\opp}\]satisfying \eqref{eq: duality-preserving functor} coherently. For complete Segal spaces, natural equivalence has the meaning of Definition \ref{defn: natural transformation segal}. In the symmetric monoidal case we use \cite[Definition 2.7]{HLS2025}.

For a pair of categories $(\cC_+,\cC_-)$, a duality means a duality on $\cC_+\sqcup\cC_-$ whose duality functor exchanges the two summands. Thus \[(-)^\vee:\cC_\pm^{\opp}\to\cC_\mp\aand\eta_\pm:\id_{\cC_\pm}\to(-)^{\vee\vee}\]satisfy \eqref{eq: bidual compatibility} for both signs. Duality-preserving pairs of functors are defined in the same way. This convention also applies to complete Segal spaces and quasicategories, with all data homotopy coherent. If $\cC_\pm$ are symmetric monoidal, we require the structures on both components to be symmetric monoidal.
\edefn
\brem\label{rem: rigidity imply duality}
For an ordinary rigid symmetric monoidal category, taking duals defines a duality; see \cite[\S 2.10, after Example 2.10.14]{EGNO2015}. For a rigid symmetric monoidal quasicategory, the corresponding homotopy coherent symmetric monoidal duality is supplied by \cite[Theorem 5.8]{HLS2025}.
\erem
From \cite[\S3.1]{Dotto2016Equivariant} and \cite[\S3.3]{DMPR2021THR}, we have the following construction of a (dg) category with strict duality from a (dg) category with duality.
\bcons[Strictification of duality]\label{cons: strict dual}
Let $\cC$ be a category with duality. Define $\cC^{\sd}$ to have objects $(X,Y,u)$, where $u:Y\xra{\cong}X^\vee$. A morphism $(X,Y,u)\to(X',Y',u')$ is a pair $(f:X\to X',g:Y'\to Y)$ satisfying $u\circ g=f^\vee\circ u'$. Set
\[
(X,Y,u)^{\vee_{\sd}}=(Y,X,u^\vee\circ\eta_X),~
(f,g)^{\vee_{\sd}}=(g,f).
\]
Then $(-)^{\vee_{\sd}\vee_{\sd}}=\id$ strictly. The first-coordinate projection $P:\cC^{\sd}\to\cC$ is an equivalence preserving duality with comparison $\xi_{(X,Y,u)}=u$.

The same construction applies to a dg category $\sfC$ with duality. Here $u$ is a closed degree-zero isomorphism, and the mapping complexes and the dual of a homogeneous map $f$ are
\[
\Map_{\sfC^{\sd}}\bigl((X,Y,u),(X',Y',u')\bigr)=\Map_{\sfC}(X,X'),~
f^{\vee_{\sd}}=u^{-1}\circ f^\vee\circ u'.
\]
Thus $\sfC^{\sd}$ has strict dg duality, and $P:\sfC^{\sd}\to\sfC$ is a duality-preserving dg equivalence. 

If the original (dg) category is symmetric monoidal with duality in the sense of Definition \ref{defn: dual}, then its strictification inherits a symmetric monoidal structure by componentwise tensor product, which commutes strictly with duality. In this case, $P$ is a duality-preserving symmetric monoidal equivalence.
\econs
The duality on complete Segal spaces can be obtained via the following construction for simplicial spaces with strict duality.

\bcons[Duality on simplicial spaces]\label{cons: special dual}
For a simplicial space $X$, let $X^{\opp,\opp}$ have the same simplices as $X$, with the operator for $\alpha\times\beta$ replaced by that for $\alpha^\vee\times\beta^\vee$, where
\begin{equation}\label{eq: alpha beta dual}
\alpha^\vee(z)=l-\alpha(l'-z),~\beta^\vee(w)=k-\beta(k'-w)
\end{equation}
for $\alpha:[l']\to[l]$ and $\beta:[k']\to[k]$. Thus both the simplicial and spatial directions are reversed.

A \emph{strict duality} on $X$ is an involution $(-)^\vee:X^{\opp,\opp}\to X$. Write $\sSpaceo^{\sd}$ for the category of such objects and maps $f:X\to Y$ satisfying
\[
f\circ(-)^\vee=(-)^\vee\circ f^{\opp,\opp}.
\]

Double reversal preserves $\sSpace_c$, $\sSpace_f$, and their Segal and complete-Segal localizations. Pairing generating (trivial) cofibrations with their reflected copies, exchanged by duality, gives corresponding model structures on $\sSpaceo^{\sd}$, denoted by adding $\sd$ to the superscript, such as $\sSpace_c^{\CSe,\sd}$. Fibrations and weak equivalences are detected in the underlying model structures, and the forgetful functors preserve cofibrations. The small object argument gives functorial replacements. Objects with strict duality need not be cofibrant. We denote by $\widehat R^{\sd}$ the functorial composite of cofibrant replacement and fibrant replacement through a trivial cofibration, also in the relevant diagram model categories.

The opposite in Definition \ref{defn: dual} reverses only the simplicial direction. Spatial reversal is naturally related to the identity by weak homotopy equivalences, compatibly with reversal and products, by \cite[{\texttt{03JL},\texttt{048L}}]{kerodon}. Hence, after the above replacement $\widehat{R}^{\sd}$ in $\sSpace_c^{\CSe,\sd}$, strict duality induces a duality in the sense of that definition.
\econs
We extend Lemma \ref{lem: rectification of coherent diagrams} to diagrams with strict duality.

\bdefn[Diagrams with strict duality]\label{defn: diagrams with strict duality}
Write $\sCato^{\sd}$ for the category of small simplicial categories $\cD$ with involutive reflections
\[
(-)^\vee:\Hom_{\cD}(x,y)^{\opp}\longrightarrow\Hom_{\cD}(x,y)
\]
preserving composition and identities, and simplicial functors preserving them. Reflection reverses mapping simplices, not arrows; ordinary categories carry the identity reflection.

The category $\sSpaceo^{\sd}$ is enriched in simplicial sets with reflection, using $\Hom_{\sSpaceo}(X,Y)$ with reflection induced by parameter reversal and the dualities of $X$ and $Y$.

Put $\Delta^1_{\sd}=\rN_\bu(I[1])$, where $I[1]$ is the walking isomorphism with objects $0,1$. For a small ordinary category $\cC$, define $\fC^{\sd}[\cC]$ as in Definition \ref{defn: homotopy coherent realization}, replacing the charts by
\[
C_\rho^{\sd}=(\Delta^1_{\sd})^{r-1}
\]
and using entrywise $\max$ when deleting internal identity arrows. Its chart simplices are arbitrary lists of cut sets, and list reversal defines the reflection. Inclusion of the chains $C_0\subseteq\cdots\subseteq C_u$ and composition of strings give
\[
\fC[\cC]\longrightarrow\fC^{\sd}[\cC]\xra{\epsilon_\cC^{\sd}}\cC,
\]
whose composite is $\epsilon_\cC$.

A \emph{homotopy coherent diagram with strict duality} on $\cC$ is an enriched functor $\fB:\fC^{\sd}[\cC]\to\sSpaceo^{\sd}$. An ordinary diagram $A:\cC\to\sSpaceo^{\sd}$ determines the coherent diagram $A\epsilon_\cC^{\sd}$.
\edefn

\blem[Rectification with strict duality]\label{lem: rectification with strict duality}
Let $\cC$ be a small ordinary category and let $\fB:\fC^{\sd}[\cC]\to\sSpaceo^{\sd}$ be an enriched functor. There exists a functor $B:\cC\to\sSpace_c^{\CSe,\sd}$ with fibrant values and a zigzag of objectwise weak equivalences in $\sSpace_c^{\CSe,\sd}$
\[
(\epsilon_\cC^{\sd})^*B\simeq\fB.
\]
The construction is functorial for enriched natural transformations of $\fB$. Together with the equivalence represented by this zigzag, $B$ is determined up to a contractible space of choices. When $\cC=\cC'\times[1]$, the same holds with cofibrant-fibrant diagrams $B(-,0)$ and $B(-,1)$ and the restrictions of that equivalence to the two ends fixed.

In particular, for $\cC=\Ga\times[1]$, the conclusion of Lemma \ref{lem: coherent Gamma morphism} holds with strict duality under the same hypotheses on the endpoint diagrams, using $\widehat R^{\sd}$ from Construction \ref{cons: special dual} in $(\sSpace_c^{\CSe,\sd})^\Ga_c$.
\elem
\bpf
The entrywise operation $t_i\mapsto\min(s,t_i)$, for $s\in\Delta^1_{\sd}$, contracts each augmentation fiber and respects the chart identifications since $\min$ distributes over $\max$. Hence $\epsilon_\cC^{\sd}$ is a weak homotopy equivalence on mapping spaces.

As in Lemma \ref{lem: rectification of coherent diagrams}, extension and restriction along $\epsilon_\cC^{\sd}$ give a Quillen equivalence on underlying projective diagram categories by \cite[Proposition A.3.3.8(1)]{HTTlurie2009}. These operations commute with the reflections of Definition \ref{defn: diagrams with strict duality}, so Construction \ref{cons: special dual} gives the corresponding Quillen equivalence for diagrams with strict duality. Applying the functorial replacements $\widehat R^{\sd}$ gives $B$ and the displayed zigzag. The same Quillen equivalence gives the asserted uniqueness of $B$ together with the equivalence represented by the zigzag, also with the specified data at both ends fixed.

For $\cC=\Ga\times[1]$, the functor $B$ gives a map $B(-,0)\to B(-,1)$. Applying the preceding Quillen equivalence with $\cC=\Ga$ to the restrictions of the displayed zigzag gives zigzags of weak equivalences between $B(-,i)$ and the $\Ga$-diagram prescribed at that endpoint, for $i=0,1$. Apply $\widehat R^{\sd}$ to this map and these zigzags in $(\sSpace_c^{\CSe,\sd})^\Ga_c$. Each map in the resulting zigzags is a homotopy equivalence between cofibrant-fibrant diagrams. Composing along the zigzags, using homotopy inverses where necessary, gives the required map preserving strict duality. The remaining argument is that of Lemma \ref{lem: coherent Gamma morphism}.
\epf

\subsection{Periodic chain complexes}
In this subsection, we modify the construction in \cite[\S 1.3]{HAlurie2017} to obtain the $\infty$-homotopy category and the $\infty$-derived category of periodic chain complexes, i.e.\, chain complexes with a $\Z/m$-homological grading for $m\in\N_+$ instead of a $\Z$-homological grading. If $m$ is odd, we assume that the additive category under consideration has characteristic $2$, i.e.\ $2f=0$ for every morphism $f$. Under this assumption, every expression $(-1)^pf$ with $p\in\Z/m$ is independent of the choice of a representative of $p$. See \cite{Saito2023cyclic} for a discussion of the ordinary homotopy and derived categories of periodic complexes.
\bdefn[{\cite[\Def 1.3.2.1]{HAlurie2017} and \cite[Definition 3.1]{Saito2023cyclic}}]\label{defn: periodic chain complexes}
    Let $\cA$ be an additive category and let $m\in\N_+$. If $m$ is odd, assume that $\cA$ has characteristic $2$. Let $\sfC=\Ch_m(\cA)$ be the dg category of \emph{$m$-periodic chain complexes} over $\cA$, defined as follows.
     \begin{itemize}
        \item An object $X\in\sfC$ is of the form $(X_p,d_p)_{p\in\Z/m}$, where $X_p\in\cA$, $d_p\colon X_p\to X_{p-1}$, such that $d_{p-1}\circ d_p=0$ for all $p\in\Z/m$. We may simply write $d_p$ as $d$ or $\partial$.
        \item For $X=(X_*,d_{*}^X),Y=(Y_*,d_*^Y)\in\sfC$, let $\Map_{\sfC}(X,Y)_*$ be the chain complex of abelian groups defined by
       \[\Map_{\sfC}(X,Y)_p=\prod_{q\in\Z/m}\Hom_\cA(X_q,Y_{q+p}),\quad (d_pf)(x)=d_{q+p}^Y(f(x))-(-1)^pf(d^X_qx).\]
       \item For $f\in\Map_{\sfC}(X,Y)_p$ and $g\in\Map_{\sfC}(Y,Z)_q$, composition gives an associative map\[\Map_{\sfC}(X,Y)_p\times \Map_{\sfC}(Y,Z)_q\to \Map_{\sfC}(X,Z)_{p+q},\quad(f,g)\mapsto g\circ f,\]which satisfies the Leibniz rule\[d_{p+q}(g\circ f)=d_qg\circ f+(-1)^qg\circ d_pf.\]The signs are well-defined by the assumption on $m$ and $\cA$.
       \item For $X\in \sfC$, the identity morphism $\id_X\in \Map_{\sfC}(X,X)_0$ is the identity map of chain complexes.
    \end{itemize}
    \edefn

 \bdefn[{\cite[Construction 1.3.1.6, Proposition 1.3.1.10]{HAlurie2017}}]\label{defn: dg nerve}
Let $\cA$ be an additive category and let $m\in\N_+$. If $m$ is odd, assume that $\cA$ has characteristic $2$. Let
\[
\K_m(\cA)=\Ndg_\bu(\Ch_m(\cA))
\]
be the dg nerve of $\Ch_m(\cA)$, which is a quasicategory, called the \emph{$\infty$-homotopy category} of $m$-periodic chain complexes over $\cA$. More explicitly, $\Ndg_n(\sfC)$ for a dg category $\sfC$ consists of all pairs
\[
\bigl((X_i)_{i\in[0,n]},(f_I)\bigr)
\]
with the following properties.
\begin{itemize}
\item For each $i\in[0,n]$, the object $X_i=X_{i,*}$ belongs to $\sfC$.
\item For every subset
\[
I=\{i_-<i_p<i_{p-1}<\cdots<i_1<i_+\}\subseteq[n],~p\in\N,
\]
where there are no intermediate indices when $p=0$, there is an element
\[
f_I\in\Map_{\sfC}(X_{i_-},X_{i_+})_p
\]
satisfying
\begin{equation}\label{eq: df formula}
df_I=\sum_{j=1}^p(-1)^j\left(f_{I\setminus\{i_j\}}-f_{\{i_j<\cdots<i_+\}}\circ f_{\{i_-<\cdots<i_j\}}\right).
\end{equation}
\end{itemize}
Here the strict inequalities in the notation for $I$ only record the order of distinct indices in $[n]$ and impose no nondegeneracy condition on the corresponding simplex.

For a monotone map $\alpha:[l]\to[n]$, the induced map is
\[
\begin{aligned}
\alpha^*:\Ndg_n(\Ch_m(\cA))&\longrightarrow\Ndg_l(\sfC),\\
\bigl((X_i)_{0\leq i\leq n},(f_I)\bigr)&\longmapsto
\bigl((X_{\alpha(j)})_{0\leq j\leq l},(g_J)\bigr),
\end{aligned}
\]
where
\[
g_J=
\begin{cases}
f_{\alpha(J)}&\text{if $\alpha|_J$ is injective},\\
\id_{X_i}&\text{if $J=\{j,j'\}$ and $\alpha(j)=\alpha(j')=i$},\\
0&\text{otherwise}.
\end{cases}
\]
\edefn

\brem
In $\K_m(\cA)$, objects are $m$-periodic chain complexes, morphisms are chain maps, and $2$-simplices encode chain homotopies.
\erem

\brem
Let $\Ch(\cA)$ be the dg category of $\Z$-graded chain complexes over $\cA$. There is a functor
\[
\iota:\Ch_m(\cA)\longrightarrow\Ch(\cA)
\]
defined by
\[
\iota(X)=(\wti X_p,\wti d_p)_{p\in\Z},\qquad
(\wti X_p,\wti d_p)=(X_{\ov p},d_{\ov p})
\quad\text{for }\ov p\equiv p\pmod m.
\]
The induced map on each mapping complex is degreewise injective but need not be surjective. Applying the dg nerve gives a functor
\[
\K_m(\cA)\longrightarrow\K(\cA)=\Ndg_\bu(\Ch(\cA)),
\]
which we continue to denote by $\iota$.
\erem

For later use, we record the following immediate consequence of Definition \ref{defn: dg nerve}.

\blem[{\cite[\texttt{00PV}]{kerodon}}]\label{lem: dg nerve of a poset}
Let $\sfC$ be a dg category and let $P$ be a finite poset. A functor
\[
F:\rN_\bu(P)\longrightarrow\Ndg_\bu(\sfC)
\]
consists of objects $X_i\in\sfC$ for $i\in P$ and components
\[
f_I\in\Map_{\sfC}(X_{i_-},X_{i_+})_p,~I=\{i_-<i_p<\cdots<i_1<i_+\}\subseteq P,~p\in\N,
\]
satisfying \eqref{eq: df formula}. Its value on a strict chain $I$ consists of the objects $(X_i)_{i\in I}$ and components $(f_J)_{J\subseteq I,\,|J|\geq2}$. The components for weak chains are determined by degeneracies,
\begin{equation}\label{eq: degeneracy condition}
    f_{\{i\le i\}}=\id_{X_i},~f_I=0\text{ if $I$ has a repetition and $|I|\geq3$}.
\end{equation}
\elem
\begin{proof}
Apply \cite[\texttt{00PV}]{kerodon} to $\rN_\bu(P)$, multiplying each degree-$p$ component in the reference by $(-1)^p$ to match the sign convention of Definition \ref{defn: dg nerve} from \cite[Construction 1.3.1.6]{HAlurie2017}.
\end{proof}

\brem\label{rem: dg nerve of simplex and product}
For $P=[k]$, we have $\rN_\bu(P)=\Delta^k$, and every subset of $[k]$ is a strict chain. Lemma \ref{lem: dg nerve of a poset} therefore recovers the description of a $k$-simplex of $\Ndg_\bu(\sfC)$ in Definition \ref{defn: dg nerve}. Since the nerve preserves products, we have
\[
\rN_\bu([l]\times[k])
\cong\rN_\bu([l])\times\rN_\bu([k])
=\Delta^l\times\Delta^k.
\]
We will apply Lemma \ref{lem: dg nerve of a poset} in \S\ref{sec: Construction of the TQFT functor} with $P=[l]\times[k]$.
\erem

For $\Z$-graded chain complexes, the derived $\infty$-category is obtained by inverting quasi-isomorphisms. Under appropriate boundedness and projectivity hypotheses, it can also be modeled by the dg nerve of a category of projective complexes \cite[Definition 1.3.2.7 and Theorem 1.3.4.4]{HAlurie2017}. We use the analogous projective dg nerve for periodic complexes, since this model is explicit and our applications concern free modules. See Remark \ref{rem: free module}.
\bdefn[{\cite[\Def 1.3.2.4]{HAlurie2017}}]
Let $\cA$ be an abelian category. An object $P$ in $\cA$ is said to be \emph{projective} if, for every epimorphism $M\to N$ in $\cA$, the induced map $\Hom_{\cA}(P,M) \to \Hom_{\cA}(P,N)$ is surjective. Let $\cA_{\proj}$ be the full subcategory of $\cA$ spanned by the projective objects. We say that $\cA$ has \emph{enough projective objects} if, for every object $M$ in $\cA$, there exists an epimorphism $P\to M$, where $P$ is a projective object of $\cA$.
\edefn
\brem\label{rem: free module}
If $\cA=\Mod(\bK)$ is the category of modules over some commutative ring $\bK$ with identity, then the full subcategory $\cA_{\free}$ of $\cA$ spanned by free modules is a full subcategory of $\cA_{\proj}$. If $\bK$ is a PID, then $\cA_{\free}=\cA_{\proj}$. Also note that $\Mod(\bK)$ has enough projective objects for any $\bK$.
\erem
\bdefn\label{defn: derived}
Let $\cA$ be an abelian category with enough projective objects and let $m\in\N_+$. If $m$ is odd, assume that $\cA$ has characteristic $2$. Let $\D_m(\cA)=\Ndg_\bu(\Ch_m(\cA_{\proj}))$, called the \emph{$\infty$-derived category} of $m$-periodic chain complexes over $\cA$. Identifying $\D_m(\cA)$ with the localization of $\K_m(\cA)$ at quasi-isomorphisms requires additional hypotheses in general, since a degreewise projective periodic complex need not be homotopically projective.
\edefn
As a direct generalization of \cite[\Prop 1.3.2.10 and \Cor 1.3.2.18]{HAlurie2017}, we obtain the following result, for which we provide a short proof.
\bprop
The $\infty$-categories $\K_m(\cA)$ and $\D_m(\cA)$, when they are defined, are stable in the sense of \cite[Definition 1.1.1.9]{HAlurie2017}.
\eprop
\bpf Regard the mapping complexes of $\Ch_m(\cA)$ as periodic $\Z$-graded complexes via $\Z\to\Z/m$. The proof of \cite[\Prop 1.3.2.10]{HAlurie2017} then carries over verbatim. Indeed, mapping cylinders and mapping cones preserve $m$-periodicity, so the same argument shows that $\K_m(\cA)$ admits finite colimits. Its zero complex is a zero object, and its suspension is induced by the shift $X\mapsto X[-1]$, which is an equivalence. Hence $\K_m(\cA)$ is stable by \cite[\Cor 1.4.2.27]{HAlurie2017}. Moreover, the full dg-subcategory $\Ch_m(\cA_{\proj})$ is closed under shifts and mapping cones, since projective objects are closed under finite direct sums. Therefore $\D_m(\cA)$ is stable as well. \epf

We next define the symmetric monoidal structure on $\D_m(\cA)$.

\bdefn\label{defn: symm on D}
Let $m\in\N_+$ and let $\cA$ be an abelian category with enough projectives and a closed symmetric monoidal structure as in \Def \ref{defn: monoidal category}. If $m$ is odd, assume that $\cA$ has characteristic $2$.

Following \cite[\texttt{00P1}--\texttt{00P4}]{kerodon}, we equip $\Ch_m(\cA)$ with the Koszul tensor product. For $X=(X_*,d^X_*)$, $Y=(Y_*,d^Y_*)$, and homogeneous $x\in X_i$, $y\in Y_j$, set
\[
X\ot Y=\Bigl(\bigoplus_{a+b=p}X_a\ot Y_b\Bigr)_{p\in\Z/m},~
d^\ot(x\ot y)=d^X_i(x)\ot y+(-1)^i x\ot d^Y_j(y).
\]
The tensor unit is the unit of $\cA$ supported in degree zero. For homogeneous maps $f:X\to X'$ and $g:Y\to Y'$, set
\[
(f\ot g)(x\ot y)=(-1)^{|g||x|}f(x)\ot g(y).
\]
The symmetry is the chain isomorphism
\[
\sigma:X\ot Y\longrightarrow Y\ot X,~
\sigma(x\ot y)=(-1)^{|x||y|}y\ot x.
\]
Assume that the unit of $\cA$ is projective and that $\cA_{\proj}$ is closed under tensor products. Then $\Ch_m(\cA_{\proj})$ is a symmetric monoidal dg category. 

\edefn
\bcons[Symmetric monoidal structure on the dg nerve]\label{cons: symm on D}
Let $\sfC$ be a small symmetric monoidal dg category over a commutative ring and put $\D=\Ndg_\bu(\sfC)$. We construct a simplicial functor
\begin{equation}\label{eq: resolved target}
\fB:\fC[\Ga]\longrightarrow\Cat_\infty^\Delta,~\fB\langle n\rangle=\D^n,
\end{equation}
and rectify it by Lemma \ref{lem: rectification of coherent diagrams} to a strict functor
\begin{equation}\label{eq: strict Gamma target}
B:\Ga\longrightarrow\QCat,~B\langle n\rangle\simeq\D^n.
\end{equation}
The Segal maps of $\fB$ are projections, so $B$ satisfies the $\Ga$-condition.

A pointed map $\theta:\langle m\rangle\to\langle n\rangle$ acts on objects by
\[
(X^\alpha)_{\alpha\in [1,m]}\longmapsto\Bigl(\bigotimes_{\theta(\alpha)=\beta}X^\alpha\Bigr)_{\beta\in [1,n]},
\]
with each fiber increasingly ordered and tensors parenthesized from the left; singleton tensors are unchanged and empty tensors give the unit. For a string
\[
\rho=(\langle n_0\rangle\xra{\theta_1}\cdots\xra{\theta_r}\langle n_r\rangle),
\]
with $r\geq1$, recall the chart of Definition \ref{defn: homotopy coherent realization}
\[
C_\rho=(\Delta^1)^{r-1}\longrightarrow\fC[\Ga](\langle n_0\rangle,\langle n_r\rangle).
\]

We define $\fB_\rho:C_\rho\times\D^{n_0}\to\D^{n_r}$ simplexwise. For each $\Delta^u\to C_\rho$ and $f:\Delta^k\to\Delta^u\times\D^{n_0}$, we specify the output objects and components of Definition \ref{defn: dg nerve}, then verify the simplicial and chart relations.

A vertex $C\subseteq\{1,\ldots,r-1\}$ cuts $\rho$ after $\theta_j$ for each $j\in C$. Compose within each block, then apply the block tensor operations from left to right, retaining parentheses and units. Fix $\beta\ne0$ in $\langle n_r\rangle $ and put \[S_{\rho,\beta}=(\theta_r\cdots\theta_1)^{-1}(\beta)=\{\alpha_1<\cdots<\alpha_s\}.\] Tensor products indexed by $S_{\rho,\beta}$ are taken in this order and parenthesized from the left. For a $u$-simplex $C_0\subseteq\cdots\subseteq C_u$ of $C_\rho$, let $o^v$ be the left-to-right order of these labels in the tensor expression at $C_v$, for $0\leq v\leq u$.

We use the example $\rho=(\langle3\rangle\xra{\theta_1}\langle2\rangle\xra{\theta_2}\langle1\rangle)$ throughout, with $\theta_1^{-1}(1)=\{1,3\}$, $\theta_1^{-1}(2)=\{2\}$, and $\theta_2(1)=\theta_2(2)=1$. Here $C_\rho=\Delta^1$, $S_{\rho,1}=\{1,2,3\}$, and its two vertices give
\[
\begin{array}{c|c|c}
\text{cut}&\text{tensor expression}&\text{label order}\\ \hline
C_0=\varnothing&(X^1\ot X^2)\ot X^3&o^0=(1,2,3)\\
C_1=\{1\}&(X^1\ot X^3)\ot X^2&o^1=(1,3,2).
\end{array}
\]

Fix $f$ over $C_0\subseteq\cdots\subseteq C_u$. Write $c_i\in[u]$ for the projection of vertex $i$ to $\Delta^u$ and $f^\alpha:\Delta^k\to\D$ for the input projections, with $X_i^\alpha=f^\alpha(i)$. Let $X_i^\otimes$ be the tensor expression at $C_{c_i}$ and
\[
\kappa_i:\bigotimes_{\alpha\in S_{\rho,\beta}}X_i^\alpha\xra{\cong}X_i^\otimes
\]
its structural isomorphism with the label order $o^{c_i}$, determined by symmetric monoidal coherence. 

Recall from Definition \ref{defn: dg nerve} that a $k$-simplex of $\D$ has a degree-$p$ component indexed by each $I=\{i_-<i_p<\cdots<i_1<i_+\}\subseteq[k]$, where $p=|I|-2\geq0$. For $p=0$, set
\[
(\fB_\rho f)^\beta_I=\kappa_{i_+}\circ\Bigl(\bigotimes_{\alpha\in S_{\rho,\beta}}f_I^\alpha\Bigr)\circ\kappa_{i_-}^{-1}.
\]
For an edge in $C_\rho\times\{X\}$ with a fixed input tuple $X$, all $f_I^\alpha$ are identities, so its image is $\kappa_{i_+}\circ\kappa_{i_-}^{-1}$.

In our example, $\kappa_i=\id$ for $c_i=0$. For $c_i=1$, it is
\[
\begin{aligned}
(X_i^1\ot X_i^2)\ot X_i^3&\xra{\cong}X_i^1\ot(X_i^2\ot X_i^3)\\
&\xra{\id\ot\sigma}X_i^1\ot(X_i^3\ot X_i^2)\xra{\cong}(X_i^1\ot X_i^3)\ot X_i^2,
\end{aligned}
\]
where the outer arrows are associativity isomorphisms and $\sigma$ is the symmetry.

For $p\geq1$, distribute the $p$ intermediate vertices of $I$ among the tensor factors. Such a choice is a map $\nu:\{1,\ldots,p\}\to S_{\rho,\beta}$, where $\nu(j)$ is the factor assigned to $i_j$. For each factor $\alpha$, delete $i_j$ if $\alpha$ precedes $\nu(j)$ in $o^{c_{i_j}}$, split at $i_j$ if $\alpha$ follows $\nu(j)$, and retain it without splitting if $\alpha=\nu(j)$. Keep the endpoints $i_-,i_+$, and include each splitting vertex in both adjacent subchains. Let $f^\alpha_{I,\nu}$ be the composite of the resulting subchain components, with later subchains on the left. It has degree $|\nu^{-1}(\alpha)|$, and these degrees sum to $p$.

Put $\epsilon(\nu)=\#\{j<k\mid\nu(j)>\nu(k)\}$, the sign exponent for grouping the indices $j$ by increasing factor label. Taking the tensor factors in that label order, set
\begin{equation}\label{eq: cubical target component}
(\fB_\rho f)^\beta_I=\kappa_{i_+}\circ\Bigl(\sum_\nu(-1)^{\epsilon(\nu)}\bigotimes_{\alpha\in S_{\rho,\beta}}f^\alpha_{I,\nu}\Bigr)\circ\kappa_{i_-}^{-1}.
\end{equation}
This also recovers the edge formula for $p=0$. For $S_{\rho,\beta}=\varnothing$, the parenthesized expression is the identity of the unit when $p=0$ and zero otherwise.

For a $2$-simplex in our example at $C_0=\varnothing$, all $\kappa_i$ are identities. Keeping, deleting, or splitting at the only intermediate vertex gives $f^\alpha_{012}$, $f^\alpha_{02}$, or $f^\alpha_{12}\circ f^\alpha_{01}$, respectively. The three assignments $\nu(1)=3,2,1$ give
\[
\begin{aligned}
(\fB_\rho f)^1_{012}={}&f^1_{02}\ot f^2_{02}\ot f^3_{012}\\
&+f^1_{02}\ot f^2_{012}\ot(f^3_{12}\circ f^3_{01})\\
&+f^1_{012}\ot(f^2_{12}\circ f^2_{01})\ot(f^3_{12}\circ f^3_{01}).
\end{aligned}
\]

For verification, use the cubical description in \cite[Definition 4.3]{rivera2020colimit} and the proof of Proposition 4.5 of that reference. Put $Q_I=[0,\infty]^p$, with ordered coordinates $q_1,\ldots,q_p$ indexed by $i_1,\ldots,i_p$ and the product orientation from $0$ to $\infty$. Each $f^\alpha$ gives a chain map on its cellular chains, taking the top cell to $f_I^\alpha$. Deletion, splitting, and retention at $i_j$ correspond to $q_j=0$, $q_j=\infty$, and leaving $q_j$ free, respectively. The sum in \eqref{eq: cubical target component} applies the iterated Serre diagonal in each coordinate $q_j$ with order $o^{c_{i_j}}$ \cite[Example 2.2 and Proposition 4.6]{rivera2020colimit}. These diagonals are chain maps and send each endpoint to that endpoint in every factor. The remaining coordinates retain their orders on each face. The cubical boundary formula, graded interchange, and cancellation of the intermediate $\kappa_{i_j}$ therefore give the dg nerve differential identity.

An interior repeated vertex gives two coordinates with the same order. The Serre diagonal commutes with combining them by $\max$, and with forgetting a coordinate at an endpoint repetition. This gives the degeneracy identities.

Let $t_i$ be the cut coordinate of $C_\rho$ after $\theta_i$. At $t_i=0$, the prescriptions agree with those for $d_i\rho$, obtained by composing $\theta_i$ and $\theta_{i+1}$. At $t_i=1$, assigning each intermediate vertex first to a block and then to a factor gives the same summands as \eqref{eq: cubical target component}. Graded interchange gives the same signs, and symmetric monoidal naturality and coherence identify the endpoint maps. Empty intermediate fibers contribute only the unit in degree zero. Thus
\begin{equation}\label{eq: cubical target faces}
\fB_\rho|_{t_i=0}=\fB_{d_i\rho},~
\fB_\rho|_{t_i=1}=\fB_{\rho_+}\circ \fB_{\rho_-}
\end{equation}
where $\rho_-,\rho_+$ are the prefix and suffix at $i$. Deleting an identity arrow leaves the tensor expressions unchanged and combines the adjacent cut coordinates by $\max$, or deletes one at an endpoint. Hence $\fB$ is a simplicial functor.

\econs

\brem
Hinich gives an alternative construction via his lax symmetric monoidal dg nerve \cite[\S 3.2.3, (33)--(34)]{Hinich2015Ndg}. In the non-periodic case, Lurie constructs a symmetric monoidal structure on the derived $\infty$-category of modules over a commutative ring \cite[Remark 7.1.2.12]{HAlurie2017}. We leave comparison with Construction \ref{cons: symm on D} to future work.
\erem
Finally, we consider the duality for $\D_m(\cA)$.
\bdefn\label{defn: dual on D}
Following the notation in Definition \ref{defn: symm on D}, suppose further that every projective object of $\cA$ is dualizable and that its dual is projective. Put $\sfC=\Ch_m(\cA_{\proj})$. For $X=(X_*,d_*)\in\sfC$, define its \emph{graded dual} by
\[
X^\vee=\bigl(X^\vee_p=(X_{-p})^\vee,~
d^\vee_p=(-1)^{p+1}(d_{-p+1})^\vee\bigr)_{p\in\Z/m}.
\]
For a homogeneous map $f:X\to Y$ of degree $n$, define
\[
(f^\vee)_p=(-1)^{np}(f_{-p-n})^\vee.
\]
The evaluation and coevaluation maps are induced by those in $\cA$, with a factor $(-1)^p$ on the summands corresponding to $X_p$. Thus $\D_m(\cA)$ is rigid.

The double dual satisfies
\[
(X^{\vee\vee})_p=X_p^{\vee\vee},~
d_p^{\vee\vee}=-(d_p)^{\vee\vee}.
\]
Hence the \emph{graded bidual isomorphism} $\eta_X:X\to X^{\vee\vee}$ has components
\[
(\eta_X)_p=(-1)^p\eta_{X_p},
\]
where $\eta_{X_p}:X_p\to X_p^{\vee\vee}$ is the bidual isomorphism in $\cA$. The identities $d(f^\vee)=(df)^\vee$ and $(g\circ f)^\vee=(-1)^{|f||g|}f^\vee\circ g^\vee$ make $(-)^\vee:\Ch_m(\cA_{\proj})^{\opp}\to\Ch_m(\cA_{\proj})$ a dg functor, where $\opp$ denotes the graded dg opposite. The maps $\eta_X$ form a closed degree-zero dg-natural isomorphism satisfying \eqref{eq: bidual compatibility} strictly. These data make $\sfC$ a dg category with duality and induce a duality on $\D_m(\cA)=\Ndg_\bu(\sfC)$, as also follows from its rigidity by Remark \ref{rem: rigidity imply duality}.

To obtain a strict-duality model, put $\D^{\sd}=\Ndg_\bu(\sfC^{\sd})$, using Construction \ref{cons: strict dual}. Duality reverses a chain $I$ and sends its degree-$p$ component, where $p=|I|-2$, to
\begin{equation}\label{eq: sign reversal fI}
f_I\longmapsto(-1)^{p(p-1)/2}f_I^{\vee_{\sd}}.
\end{equation}
The sign follows from the dg nerve differential identity and the graded dg opposite. Thus $\Psi(\D^{\sd})\in\sSpaceo^{\sd}$, and projection induces a duality-preserving equivalence $\D^{\sd}\simeq\Ndg_\bu(\sfC)$.
\edefn
We now incorporate duality into the tensor diagram of Construction \ref{cons: symm on D}.

\bcons[Tensor diagrams with strict duality]\label{cons: strict tensor diagram}
Let $\sfC$ be a symmetric monoidal dg category with duality and put $\D^{\sd}=\Ndg_\bu(\sfC^{\sd})$. Apply Construction \ref{cons: symm on D} to the symmetric monoidal strictification of Construction \ref{cons: strict dual}. Its formulas extend to the charts $C_\rho^{\sd}$ of Definition \ref{defn: diagrams with strict duality}, since their verification uses only the tensor order at each vertex. This gives
\[
\fB^{\sd}:\fC^{\sd}[\Ga]\longrightarrow\mathrm{Cat}_\infty^\Delta,~\fB^{\sd}\langle n\rangle=(\D^{\sd})^n.
\]
The reflection formula \eqref{eq: sign reversal fI}, applied to the cubical components, and tensor-duality compatibility at the endpoints give
\begin{equation}\label{eq: strict tensor reflection}
(\fB^{\sd}_\rho f^\vee)_{I^\vee}^{\beta}=(-1)^{p(p-1)/2}\bigl((\fB^{\sd}_\rho f)_I^\beta\bigr)^{\vee_{\sd}}.
\end{equation}
Here $p=|I|-2$, and the left side uses the reversed chart list, i.e.\ $(C_u,\ldots,C_0)$ in place of $(C_0,\ldots,C_u)$, and the reversed input simplex. Thus $\Psi(\fB^{\sd})=\Psi\circ\fB^{\sd}$ is a homotopy coherent diagram with strict duality in the sense of Definition \ref{defn: diagrams with strict duality}. Lemma \ref{lem: rectification with strict duality} gives a fibrant $\Ga$-diagram $B^{\sd}$ in $\sSpace_c^{\CSe,\sd}$. After forgetting duality, the projection $P$ of Construction \ref{cons: strict dual} induces $B^{\sd}\simeq\Psi(B)$.
\econs
\bexmp\label{exmp: finite module}
Let $\bK$ be a commutative Noetherian ring and put $\cA=\Mod(\bK)^\fin$. This is an abelian category with enough projectives, and $\cA_{\proj}$ consists of finitely generated projective modules and is closed under duality. Definition \ref{defn: dual on D} gives the rigid quasicategory $\D_m(\cA)$ a strict-duality model
\[
\D_m^{\sd}(\cA)=\Ndg_\bu\bigl(\Ch_m(\cA_{\proj})^{\sd}\bigr)\simeq\D_m(\cA).
\]
Construction \ref{cons: strict tensor diagram} represents its symmetric monoidal structure by a $\Ga$-diagram of complete Segal spaces with strict duality.
\eexmp

\section{Preliminaries on cobordism category}\label{sec: Preliminaries in cobordism category}
In this section, we review the construction by Calaque--Scheimbauer \cite{CS2019category} about $\infty$-version of the cobordism category. In the reference, the authors used the $(\infty,n)$-category for $n\in\N$ to construct the fully extended cobordism category, while we only need the special case of $(\infty, 1)$-category (i.e.\, $\infty$-category) that extends $n$-dimensional cobordisms only down to $(n-1)$-dimensional objects. Moreover, we also consider the symmetric monoidal structure, tangential structure (e.g.\ orientations, Riemannian metrics, almost complex structures), and the dual structure on the $\infty$-cobordism category. For simplicity, we only adapt necessary constructions and results in \cite{CS2019category} and omit details. For a quick review, see \cite[\S 4.1]{CHS2025tqft}.

\subsection{Diffeological spaces}\label{subsec: diffeological spaces}
In this subsection, we recall the notions and properties of diffeological spaces used below; see \cite{CW2014diff,Kihara2024smooth}.
\bdefn[{\cite[\S 2]{CW2014diff} and \cite[Definition 2.1]{Kihara2024smooth}}]\label{defn: diffeological space}
Let $X$ be a set. A \emph{diffeology} on $X$ is a collection $\fD_X$ of maps $p:U\to X$, called \emph{plots}, for each open $U\subset\R^d$ and each $d\in\N$, satisfying the following conditions.
\benum
\item Every constant map $p:U\to X$ belongs to $\fD_X$.
\item If $p:U\to X$ is a map and there is an open cover $\{U_i\}$ of $U$ such that $p|_{U_i}\in\fD_X$ for every $i$, then $p\in\fD_X$.
\item If $p:U\to X$ belongs to $\fD_X$ and $F:V\to U$ is smooth, where $V\subset\R^m$ is open for some $m\in\N$, then $p\circ F\in\fD_X$.
\eenu
A set equipped with a diffeology is a \emph{diffeological space}. Given a diffeological space $X$, the final topology induced by its plots, namely a subset $V\subset X$ is open if and only if $p^{-1}(V)$ is open for any plot $p$ on $X$, is called the \emph{$\fD$-topology} on $X$.

For a diffeological space $X$ and a subset $A\subset X$, the largest diffeology on $A$ making the inclusion   $A\hookrightarrow X$  smooth is called the \emph{sub-diffeology}. It consists of   the maps $U\to A$ whose composites with $A\hookrightarrow X$ are plots  of $X$. For a diffeological space $X$ with an equivalence relation $\simeq$, the smallest diffeology on   $X/\simeq$  making the quotient map   $X\twoheadrightarrow X/\simeq$  smooth is called the \emph{quotient diffeology}.   Its plots locally lift  through the quotient map  to plots of $X$ . The smallest diffeology on $X$ containing a   collection $A=\{U_i\to X\}_{i\in I}$ is  the diffeology \emph{generated} by $A$.   Its plots  locally either factor through the given maps   by smooth maps  or are constant.

A map $f:X\to Y$ between diffeological spaces is \emph{smooth} if $f\circ p\in\fD_Y$ for every $p\in\fD_X$. We write $C^\infty(X,Y)$ for the set of smooth maps. The collection of all diffeological spaces and smooth maps forms a category, which we denote $\fDiff$. An isomorphism in $\fDiff$ will be called a \emph{diffeomorphism}. 

A smooth map $X\to Y$ is continuous when both $X$ and $Y$ are equipped with
the $\fD$-topology, and so this defines a functor $\fD : \fDiff \to \Top$. Every topological space $X$ has a natural diffeology, called the \emph{continuous diffeology}, whose plots $U\to X$ are the continuous maps. A continuous map $X\to Y$ is smooth when both $X$ and $Y$ are equipped with the continous diffeology, and so this defines a functor $\fC : \Top \to  \fDiff$. The functor $\fC$ is the right adjoint of $\fD$.

  The \emph{functional diffeology}  on $C^\infty(X,Y)$   consists of the maps $p:U\to C^\infty(X,Y)$ for which the adjoint map $U\times X\to Y$ is smooth. With this diffeology, $-\times Y$ is left adjoint to $C^\infty(Y,-)$; hence  $\fDiff$ is cartesian closed  \cite[Proposition 2.2(2)]{Kihara2024smooth} .

\edefn
\bexmp
For a smooth manifold $X$ without boundary, the standard diffeology   consists of the smooth maps $U\to X$  in the ordinary sense. The corresponding $\fD$-topology is the same as the original topology of $X$. For  manifolds with corners, we   use  the standard diffeology   defined by corner charts  whose local models carry the sub-diffeology.   By \cite[\S 3]{SP2019diff} , diffeologically smooth maps between manifolds with corners are precisely smooth maps in the   usual sense.
 \eexmp

 The main diffeological space used in this paper is the following one and its variants.
 \bdefn[{\cite[\S 2.3]{SP2024tqft}}]\label{defn: embedded manifold diffeology}
Let $N$ be a smooth manifold, possibly with boundary and corners. Let $\psi_n(N)$ be the set of smooth $n$-dimensional submanifolds $W\subset N$ without boundary that are closed in $N$ and disjoint from $\partial N$. A map $p:U\to\psi_n(N)$ is a plot if
\[
\Gamma_p=\{(u,x)\in U\times N\mid x\in p(u)\}
\]
is a closed smooth $(n+\dim U)$-dimensional submanifold of $U\times N$ and the projection $\Gamma_p\to U$ is a submersion. These plots define the \emph{embedded-manifold diffeology} on $\psi_n(N)$.
\edefn
The following result records the two functorialities of this diffeology that will be used for cutting, rescaling, and stabilization.
\blem[{\cite[\S 2.3 and Lemma 2.10]{SP2024tqft}}]\label{lem: embedded manifold functoriality}
Let $N$ and $N'$ be smooth manifolds of the same dimension, and assume that the maps below preserve the interiors when boundaries are present. A submersion $f:N\to N'$ induces a smooth pullback map
\[
f^*:\psi_n(N')\longrightarrow\psi_n(N),~W\longmapsto W\times_{N'}N.
\]
A closed embedding $g:N\hookrightarrow N''$ induces a smooth pushforward map
\[
g_*:\psi_n(N)\longrightarrow\psi_n(N''),~W\longmapsto g(W).
\]
More generally, let $f_y:N\hookrightarrow N'$ be a smooth family of open embeddings and let $g_y:N\hookrightarrow N''$ be a smooth family such that $(y,x)\mapsto(y,g_y(x))$ is a closed embedding $Y\times N\hookrightarrow Y\times N''$. These families induce smooth maps
\[
Y\times\psi_n(N')\longrightarrow\psi_n(N),~(y,W)\longmapsto f_y^{-1}(W),
\]
and
\[
Y\times\psi_n(N)\longrightarrow\psi_n(N''),~(y,W)\longmapsto g_y(W).
\]
Thus $\psi_n(-)$ is contravariant for submersions between ambient manifolds of the same dimension and covariant for closed embeddings.
\elem
\bpf
The fixed maps send plots to plots by \cite[\S 2.3]{SP2024tqft}. Although \cite[Lemma 2.10]{SP2024tqft} is stated in terms of continuity for the plot topology, its proof verifies that the parametrized maps send plots to plots and hence proves the stated diffeological smoothness.
\epf
 \bdefn[{\cite[\S 2.2]{Kihara2024smooth}}]\label{defn: smooth homotopy equivalence}
Let $I=[0,1]$ carry the standard diffeology. A \emph{smooth homotopy} from smooth maps $f_0,f_1:X\to Y$ is a smooth map $H:X\times I\to Y$ such that $H|_{X\times\{i\}}=f_i$ for $i=0,1$. A smooth map $f:X\to Y$ is a \emph{smooth homotopy equivalence} if there is a smooth map $g:Y\to X$ such that $g\circ f$ and $f\circ g$ are smoothly homotopic to the respective identity maps.
\edefn

For convenience, we retain the notation for the singular functors $S^\aff_\bu$ and $S^\sub_\bu$ introduced in Example \ref{exmp: singular complex}, defined using $|\Delta^l|_\aff$ and $|\Delta^l|$ with their standard diffeologies, respectively. The Kihara complex $S^\kih_\bu$ is defined using Kihara's standard simplices \cite[Definition 2.4]{Kihara2024smooth}; in degrees at least two, their diffeologies differ from the sub-diffeologies on $|\Delta^l|$. For every diffeological space $X$, there is a natural inclusion
\[
S^\sub_\bu(X)\hookrightarrow S^\kih_\bu(X).
\]

Following \cite[Definition 5.6]{CS2019category}, we formulate smooth families in the affine model and restrict them to $|\Delta^l|$ for the analytic construction of the instanton functor. We do not construct these families directly in the $S^\sub_\bu$ model, since \cite[Lemma 2.17]{GR2010cob} does not directly give the required global description over the manifold with corners $|\Delta^l|$. The Kihara model provides an explicit Segal space model, but is not used to parametrize instanton moduli spaces. Fibrant replacement and completion may instead be applied directly to the affine model, as explained in Remark \ref{rem: affine construction}. We distinguish these models by their superscripts.
\blem\label{lem: affine Kihara comparison}
The functors $S^\aff_\bu$, $S^\sub_\bu$, and $S^\kih_\bu$ are right adjoints and therefore preserve limits, especially pullbacks. For every diffeological space $X$, the comparison morphisms
\[
S^\aff_\bu(X)\longrightarrow S^\sub_\bu(X)\hookrightarrow S^\kih_\bu(X)
\] 
are natural in $X$, compatible with the canonical product isomorphisms, and weak homotopy equivalences. Moreover,
\[
S^{\aff}_0(X)=S^{\sub}_0(X)=S^{\kih}_0(X)=X,
\]
the comparisons are the identity in degree zero, $S^\kih_\bu(X)$ is a Kan complex, and all three functors send smoothly homotopic maps to homotopic morphisms of simplicial sets in the sense of Definition \ref{defn: homotopy equivalence}. In particular, they send smooth homotopy equivalences to simplicial homotopy equivalences.
 \elem
\bpf
The adjointness and preservation of limits are \cite[Remark 3.2(1)]{Kihara2024smooth}. The comparison morphisms are constructed in \cite[Proposition 3.4]{Kihara2024smooth} and are weak homotopy equivalences by \cite[Theorem 1.1]{Kihara2024smooth}. Their compatibility with products and their degree-zero description follow from the degreewise construction. The Kan property follows from \cite[Corollary 2.6(2)]{Kihara2024smooth}, and preservation of smooth homotopies is \cite[Lemma 3.3]{Kihara2024smooth}.
\epf

\bdefn\label{defn: simplicial diffeological space}
A \emph{simplicial diffeological space} is a functor
\[
X_\bu:\Delta^{\opp}\longrightarrow\fDiff.
\]
This is the diffeological analogue of a simplicial space in Definition \ref{defn: simplicial space}. For $\mmod\in\{\aff,\sub,\kih\}$, applying $S_\bu^\mmod$ levelwise gives a simplicial space whose $(k,l)$-simplices are
\[
S_l^\mmod(X_k).
\]
\edefn
Consequently, applying the comparison morphism in Lemma \ref{lem: affine Kihara comparison} levelwise to a simplicial diffeological space, or objectwise to a $\Ga$-object, gives a natural levelwise weak homotopy equivalence compatible with all $\Ga$-structure maps. Homotopy invariance of homotopy fiber products and two-out-of-three then transfer the Segal and $\Ga$-conditions from the affine model to the Kihara model. The required affine statements are proved in Appendices \S \ref{sec: segal} and \S \ref{sec: symmetric monoidal}.

\subsection{Interval category}
We first define a simplicial space $\Int_\bu^{\kih}$ of intervals in $\R$ and then prove that it is a complete Segal space. The level $k$ for $\Int_k^{\kih}$ is called the \emph{simplicial level}, while the level in each fixed $\Int_k^{\kih}$ is called the \emph{spatial level}. The main results in this subsection follow from \cite[\S 4]{CS2019category}.
\bdefn\label{defn: Int}
For $k\in\N$, let
\begin{equation}\label{eq: intk}
    \begin{aligned}
        \fInt_k= \{ (\un{a}, \un{b}) = &(a_0, \ldots, a_k, b_0, \ldots, b_k) ~\big|~ a_j < b_j \text{ for } 0 \leq j \leq k \\\text{ and } & a_{j-1} \leq a_j \text{ and } b_{j-1} \leq b_j \text{ for } 1 \leq j \leq k \} \subset \R^{2k+2}
    \end{aligned}
\end{equation}
with the sub-diffeology. For $(\un a,\un b)\in\fInt_k$, set
  \[
I_j^c=[a_j,b_j],~B(\I)=(a_0,b_k),~I_j=I_j^c\cap B(\I),
\] 
and write $\I=(I_0\le\cdots\le I_k)$. Here $B$ is called the \emph{boxing map}. Thus $I_i\le I_j$ means $a_i\le a_j$ and $b_i\le b_j$. For $t\in[0,1]$, define
  \begin{equation}\label{eq: middle}
      m_t(I_j)=(1-t)a_j+tb_j,~m(I_j)=m_{1/2}(I_j).
  \end{equation}
For $\mmod\in\{\aff,\kih\}$, set
\[
\Int_k^{\mmod}=S^{\mmod}_\bu(\fInt_k).
\]
More explicitly, an $l$-simplex of the affine model $\Int_k^\aff$ consists of smooth maps
  \[
|\Delta^l|_{\aff} \to \R,~s \mapsto a_j(s), \, b_j(s),
\] 
for the \emph{affine geometric simplex}
    \[
    |\Delta^l|_{\aff}=\big\{(x_0, \ldots, x_l) \in \R^{l+1}~\big|~\sum_i x_i = 1\big\}\] and $j = 0, \ldots, k$ such that for every $s \in |\Delta^l|_{\aff}$, the following inequalities hold.
  \[
a_j(s) < b_j(s) \text{ for } j = 0, \ldots, k,
\] 
  \[
a_{j-1}(s) \leq a_j(s) \text{ and } b_{j-1}(s) \leq b_j(s) \text{ for } j = 1, \ldots, k.
\] 

For such a family, we use the pointwise notation $I_j^c(s)$, $B(\I(s))$,   $I_j(s)$, and $m_t(I_j(s))$ . We denote the simplex by
\[
(I_0(s)\le\cdots\le I_k(s))_{s\in|\Delta^l|_\aff}\text{ or }\I(s)_{s\in|\Delta^l|_\aff}.
\]
We call
\[
B(\I(s))_{s\in|\Delta^l|_\aff}=\bigcup_{s\in|\Delta^l|_\aff}B(\I(s))\times\{s\}\subset\R\times|\Delta^l|_\aff
\]
its total space.

   From Lemma \ref{lem: affine Kihara comparison}, the two models have the same $0$-simplices.

For a morphism $g:[m]\to[k]$ in $\Delta$, define the smooth map
\begin{equation}\label{eq: maps g}
g^*:\fInt_k\longrightarrow\fInt_m,~(I_0\le\cdots\le I_k)\longmapsto(I_{g(0)}\le\cdots\le I_{g(m)}),
\end{equation}

where the intervals on the right are intersected with the new boxing interval $(a_{g(0)},b_{g(m)})$. These maps are strictly functorial. Applying $S^{\mmod}_\bu$ levelwise makes $\Int_\bu^{\mmod}$ a simplicial space for every model. In particular, the $j$-th simplicial face and degeneracy maps delete and repeat the $j$-th interval, respectively. We call $\Int_\bu^\kih$ the \emph{interval category}.
\edefn
\brem
In \cite[\S 4]{CS2019category}, the superscript $c$ denotes the presentation by ordered closed intervals   $I_j^c$ , while the presentation without $c$ uses   $I_j=I_j^c\cap B(\I)$. Since $I_j^c$  is the closure of   $I_j$  in $\R$, the two presentations are isomorphic.
\erem

\brem
As stated with affine complexes, however, \cite[Proposition 4.7]{CS2019category} does not have the required levelwise fibrancy. Indeed, for $k\ge1$, $\fInt_k$ is a manifold with nonempty corners, so $\Int_k^\aff$ is not a Kan complex by \cite[Corollary 4.47]{CW2014diff}. The weak-contractibility argument remains valid, and passing to the Kihara model supplies the Kan property.
\erem
\bprop\label{prop: Int constant}
For every $\mmod\in\{\aff,\kih\}$, each $\Int_k^{\mmod}$ is weakly contractible and $*\to\Int_\bu^{\mmod}$ is a levelwise weak homotopy equivalence. Moreover, every $\Int_k^\kih$ is Kan and $\Int_\bu^\kih$ is a complete Segal space.
\eprop
\bpf
The tuples $a_j=0$, $b_j=1$ define $*\to\Int_\bu^{\mmod}$. Linear interpolation, with a smooth cutoff in the homotopy parameter, contracts each convex space $\fInt_k$ to this tuple. Lemma \ref{lem: affine Kihara comparison} therefore gives levelwise weak contractibility and the Kan property for $\Int_k^\kih$. The Segal maps are weak equivalences since their sources and homotopy-fiber-product targets are weakly contractible. Finally, the unique component of $\Int_1^\kih$ contains degenerate edges, so all its edges are equivalences. Thus $s_0:\Int_0^\kih\to(\Int_1^\kih)^\simeq=\Int_1^\kih$ is a weak equivalence, proving completeness.
\epf
\subsection{\texorpdfstring{$\infty$}{Infinity}-cobordism category}
In this subsection, we review the $\infty$-cobordism category following \cite[\S 5]{CS2019category}. Each simplicial level $(\PB_n^{V,\mmod})_k$, where $\mmod\in\{\aff,\kih\}$, is obtained by applying a singular functor to a diffeological space $(\fPB_n^V)_k$. In $(\PB_n^{V,\mmod})_{k,l}$, the first index is the simplicial level and the second is the spatial level.

Since $|\Delta^l|_{\aff}\cong\R^l$ has no boundary, \cite[Lemma 2.17]{GR2010cob} identifies an affine family, after box rescaling, with a closed smooth submanifold of $|\Delta^l|_{\aff}\times V\times(0,1)$ whose projection to $|\Delta^l|_{\aff}$ is a submersion.

\bdefn\label{defn: pbord}
Let $V$ be a finite-dimensional $\R$-vector space,   identified  with some $\R^r$ , and let $n\in\N_+$. For $k\in\N$, let $(\fPB_n^V)_k$ be the set of pairs $(M,\I)$, where $\I\in\fInt_k$ and $M\subset V\times B(\I)$ is a smooth $n$-dimensional submanifold without boundary that is closed as a subset and bounded in the $V$-direction, such that the projection $\pi:M\to B(\I)$ is proper and is a submersion over every $I_j$ .

 For $\I=(I_0\le\cdots\le I_k)$, the box rescaling
\[
B(\I)\longrightarrow(0,1),~t\longmapsto\frac{t-a_0}{b_k-a_0}
\]
sends $M$ to a submanifold of $V\times(0,1)$ and thereby defines the inclusion
  \[
(\fPB_n^V)_k\hookrightarrow\psi_n(V\times(0,1))\times\fInt_k.
\] 
A map $p:U\to(\fPB_n^V)_k$, written $p(u)=(M_u,\I(u))$, is a plot if its composite with the above inclusion is smooth and the projection
\[
\{(u,x,t)\mid(x,t)\in M_u\}\longrightarrow\{(u,t)\mid t\in B(\I(u))\},~(u,x,t)\longmapsto(u,t),
\]
is proper. These conditions are local on $U$ and preserved by smooth base change, so they define a diffeology. Box rescaling, restriction to a smaller boxing interval, and stabilization under a linear inclusion $V\hookrightarrow V'$ are smooth by Lemma \ref{lem: embedded manifold functoriality} and preservation of properness under these operations.

For $\mmod\in\{\aff,\kih\}$, set
\[
(\PB_n^{V,\mmod})_k=S_\bu^{\mmod}((\fPB_n^V)_k).
\]
An $l$-simplex of $(\PB_n^{V,\aff})_k$ is equivalently a tuple
\[
(M,\I(s)_{s\in|\Delta^l|_{\aff}}),~\I(s)=(I_0(s)\leq\cdots\leq I_k(s)),
\]
satisfying the following conditions.

\benu
\item\label{p1} $\I(s)_{s\in|\Delta^l|_{\aff}}$ is an $l$-simplex in $\Int_k^\aff$.
\item\label{p2} $M$ is a closed $(n+l)$-dimensional submanifold of \[V \times B(\I(s))_{s\in|\Delta^l|_{\aff}}\subset V\times \R\times |\Delta^l|_{\aff}\]  whose fibers $M_s$ are bounded in the $V$-direction and which satisfies
 \benum
    \item the composition $\pi: M \hookrightarrow V \times B(\I(s))_{s\in|\Delta^l|_{\aff}} \twoheadrightarrow B(\I(s))_{s\in|\Delta^l|_{\aff}}$ is a proper map;
    \item the composition $\pi^l: M \hookrightarrow V \times B(\I(s))_{s\in|\Delta^l|_{\aff}} \twoheadrightarrow |\Delta^l|_{\aff}$ is a submersion.
\eenu
\item\label{p3} For $0 \leq j \leq k$, at every
\[
x \in \pi^{-1}(\bigcup_{s\in|\Delta^l|_{\aff}}I_j(s)\times\{s\}),
\]
the map $\pi$ is submersive.
\eenu
\edefn
 We  provide a series of remarks to   clarify the  definition.
\brem
Our conventions differ from those of \cite{CS2019category} in the following respects.
\begin{itemize}
\item The affine presentation in \cite[Definition 5.6]{CS2019category} additionally imposes a prescribed triviality condition outside $|\Delta^l|\subset|\Delta^l|_\aff$. We omit this condition, since it is not intrinsic to the definition of the affine simplex.
\item In \cite[Definition 6.3]{CS2019category}, Calaque--Scheimbauer use $(\PB_n^{(\infty,d),V})_k$ for their general $d$-fold construction; our case corresponds to $d=1$. Since $\Bord_n$ always denotes the complete Segal space model and the ordinary cobordism category is denoted by $\Cob_n$, we omit the superscript $\infty$ from both $\PB_n^{V,\mmod}$ and $\Bord_n$.
\item In \cite[\S\S 5.2.1--5.2.2]{CS2019category}, $(\fPB_n^V)_k$ is treated as a topological space, and smooth families are understood in the sense of \cite[Definition 2.16]{GR2010cob}. By \cite[Theorem A.5]{SP2024tqft}, the $\fD$-topology on the ambient space $\psi_n(V\times(0,1))$ agrees with the topology of Galatius--Randal-Williams \cite[\S 2.1]{GR2010cob}. Our diffeological formulation records the smooth families used to define the affine simplices, information that is not determined by the underlying topology alone.
\end{itemize}
\erem
\brem[{\cite[Remark 5.3 and Proposition 8.10]{CS2019category}}]\label{rem: 0-simplices}
An element
  \[
(M,\I)\in(\fPB_n^V)_k=(\PB_n^{V,\kih})_{k,0}
\] 
can be regarded as a collection of $k$ composable cobordisms as follows.
\begin{itemize}
\item Condition \eqref{p3} implies that any $t_j\in I_j$ is a regular value of $\pi$. Hence $\pi^{-1}(t_j)$ is a closed $(n-1)$-manifold, where compactness follows from the properness of $\pi$. The Morse lemma implies that the preimages $\pi^{-1}(t_j)$ of different choices of $t_j$ in the same $I_j$ are diffeomorphic, with a choice-independent isotopy class. In particular, for $i\in[0,k]$ we take
\[
Y_i=\pi^{-1}(m(I_i))
\]
for $m(I_i)$ defined in \eqref{eq: middle}.

\item For $(M,\I)\in(\PB_n^\aff)_{k,0}$ and $j\in[0,k-1]$, set
\[
M_j=\pi^{-1}\big([m_{1/3}(I_j),m_{2/3}(I_{j+1})]\big),
\]
which is a cobordism from $Y_j$ to $Y_{j+1}$. We write $(M,(M_j)_{j\in[0,k-1]})$ for the corresponding $k$ composable cobordisms with collars.
\end{itemize}
\erem
\brem[{\cite[Remarks 5.8 and 5.9]{CS2019category}}]\label{rem: smooth fiber}
Given an affine $l$-simplex   $(M,\un I(s)_{s\in|\Delta^l|_\aff})$  as in Definition \ref{defn: pbord}, the fiber $M_s$ of $\pi^l$ at   $s\in|\Delta^l|_\aff$, together with $\un I(s)$, defines an element of $(\fPB_n^V)_k$. We write  $\pi_s:M_s\to B(\I(s))$   for the induced map. Choose smooth functions $c_0<c_1$  and $\ep>0$   such that  $c_0(s)\in I_0(s)$, $c_1(s)\in I_k(s)$, and
  \[
(c_0(s)-\ep(s),c_1(s)+\ep(s))=B(\un I(s)).
\] 
The   argument following \cite[Definition 2.6]{GTMW2009cob}, applied with $(X,W,\pi,(\pi,f))$ in their notation corresponding to $(|\Delta^l|_\aff,M,\pi^l,\pi)$ in our notation, shows that $\pi^l$ is a smooth fiber bundle . Since   $|\Delta^l|_\aff$  is contractible,   this bundle admits a noncanonical trivialization $M\cong M_{s_0}\times|\Delta^l|_\aff$ for every $s_0\in|\Delta^l|_\aff$.
 \erem

  On $(\PB_n^{V,\mmod})_{k,l}$, the spatial maps are those from singular complexes, and the simplicial maps are generalized from \eqref{eq: maps g}. We record both explicitly in the  affine model.

\bdefn[{\cite[\Defs 5.10 and 5.15]{CS2019category}}]\label{defn: f and g}
We define the spatial and simplicial maps on   $(\PB_n^{V,\aff})_{k,l}$ .

We first define spatial maps. Let $f:[p]\to [l]$ be a morphism in the simplex category $\Delta$, which could be either $\delta^j$, $\sigma^j$, namely the standard face and degeneracy maps in the simplex category $\Delta$, or their compositions. Let $|f|:|\Delta^p|_{\aff}\to |\Delta^l|_{\aff}$ be the induced affine map. The induced map on the spatial level is  \[\begin{aligned}
    f^\Delta: (\PB_n^{V,\aff})_{k,l}&\to (\PB_n^{V,\aff})_{k,p}\\
            (M, \I(s)_{s\in|\Delta^l|_{\aff}})&\mapsto (\bigcup_{t\in|\Delta^p|_{\aff}}M_{|f|(t)}\times \{t\}, \I(|f|(t))_{t\in|\Delta^p|_{\aff}}),
\end{aligned}\] where $M_{|f|(t)}$ is the fiber of $M\to |\Delta^l|_{\aff}$ at $|f|(t)$ as in Remark \ref{rem: smooth fiber}. We also write $f^\Delta(M)$ and $f^\Delta(\un{I})$ for the maps on the corresponding entries.

Then we define simplicial maps. Let $g:[m]\to[k]$ be a morphism in $\Delta$. The induced map on the simplicial level is  \[\begin{aligned}
    g^*: (\PB_n^{V,\aff})_{k,l}&\to (\PB_n^{V,\aff})_{m,l}\\
            (M, \I(s)_{s\in|\Delta^l|_{\aff}})&\mapsto (\pi^{-1}\big(B(g^*(\I(s)))_{s\in|\Delta^l|_{\aff}}\big), g^*\I(s)_{s\in|\Delta^l|_{\aff}}),
\end{aligned}\] where \[g^*\I(s)=(I_{g(0)}(s)\le \cdots \le I_{g(m)}(s))\text{ for }\I(s)=(I_0(s)\le\cdots \le I_k(s)).\]
Similar to $f^\Delta$, we also write $g^*$ for maps on entries. Note that for the fiber $M_s$ at $s\in|\Delta^l|_{\aff}$, we have $(g^*M)_s=g^*M_s$.

As special cases of $f^\Delta$ and $g^*$, we still write the spatial face and degeneracy maps by $d_j^\Delta,~s_j^\Delta$ and simplicial face and degeneracy maps by $d_j,~s_j$.  The maps $g^*$ are smooth by Lemma \ref{lem: embedded manifold functoriality}, since restriction to the smaller boxing interval also preserves the properness condition in Definition \ref{defn: pbord}. Thus   $(\fPB_n^V)_\bu$  is a simplicial diffeological space. Applying $S^{\mmod}_\bu$ levelwise gives a simplicial space $(\PB_n^{V,\mmod})_{\bu,\bu}$.
\edefn
\bprop\label{prop: segal space}
Both $(\PB_n^{V,\aff})_{\bu,\bu}$ and $(\PB_n^{V,\kih})_{\bu,\bu}$ satisfy the Segal condition. Thus $(\PB_n^{V,\aff})_{\bu,\bu}$ is a Segal simplicial set in the sense of Definition \ref{defn: segal space}, while  $(\PB_n^{V,\kih})_{\bu,\bu}$ is  a Segal space.
\eprop

\brem\label{rem: special proof}
We prove Proposition \ref{prop: segal space} in Appendix \S \ref{sec: segal}. Since the affine model is not levelwise Kan, the affine path model used in the proof of \cite[Proposition 5.19]{CS2019category} does not automatically compute a homotopy pullback. The appendix verifies Rezk's pullback criterion for the relevant maps and thereby shows that these models compute the homotopy fiber products of Definition \ref{defn: homotopy fiber product}. It then transfers the Segal condition to the Kihara model along the natural weak homotopy equivalence of Lemma \ref{lem: affine Kihara comparison}.
\erem
Based on Proposition \ref{prop: segal space}, we define the $\infty$-cobordism category   $\Bord_n$  via homotopy colimit and completion.
\bdefn
Let $\R^\infty$ be countably infinite-dimensional Euclidean space.   For $\mmod\in\{\aff,\kih\}$, define
\[
\PB_n^\mmod=\hocolim_{V\subset\R^\infty}(\PB_n^{V,\mmod})_{\bu,\bu}.
\] 
  The  weak homotopy equivalence from Lemma \ref{lem: affine Kihara comparison}  induces a levelwise weak homotopy equivalence
  \begin{equation}\label{eq: kappa}
\kappa:\PB_n^\aff\longrightarrow\PB_n^\kih.
\end{equation} 
  This follows from  homotopy invariance of homotopy colimits.

Also note that each manifold $M$ can be embedded in some $V$ by Whitney's embedding theorem\footnote{Indeed, we need some modified version of Whitney's embedding theorem for cobordisms; see \cite[\S 8]{CS2019category}}. We will omit $\bu$ for simplicity. 

Let   $\Bord_n$  be the completion of   $\PB_n^\kih$ as in Definition \ref{defn: complete segal space}. It is a complete Segal space, called the \emph{$n$-dimensional $\infty$-cobordism category}. We omit the model superscript only after completion.
\edefn
\brem
Since the model structure $\sSpace_f^{\Se}$ for Segal spaces is obtained from the projective model structure of simplicial spaces $\sSpace_f$ by left Bousfield localization, the homotopy colimit can also be computed in the model structure in $\sSpace_f$.
\erem

We have an explicit description of each simplicial level.
\bprop[{\cite[Proposition 8.17]{CS2019category}}]\label{prop: classifying spaces}
For every $k\in\N$, there is a zigzag of weak homotopy equivalences of topological spaces
\[
|(\PB_n^\kih)_k|\xleftarrow[\simeq]{\kappa_k}|(\PB_n^\aff)_k|
\xrightarrow{\simeq}\coprod_{[(M,(M_j)_{j\in[0,k-1]})]}B\operatorname{Diff}(M,(M_j)_{j\in[0,k-1]}),
\]
where the coproduct is indexed by diffeomorphism classes of $k$ composable cobordisms with collars, and $\operatorname{Diff}(M,(M_j)_{j\in[0,k-1]})$ is the group of diffeomorphisms of $M$ intertwining the composed cobordisms and their collars.
\eprop
Finally, we review the symmetric monoidal structure on the $\infty$-cobordism category   $\Bord_n$  from \cite{CS2019category}.

\bdefn[{\cite[Definition 7.1]{CS2019category}}]\label{defn: symmetric monoidal str}
Let $V$ be a finite-dimensional vector space, let $n\in\N_+$, and let $m,k\in\N$. For $m=0$, set $(\fPB_n^V\langle0\rangle)_k=\fInt_k$. For $m>0$, let $(\fPB_n^V\langle m\rangle)_k$ be the sub-diffeological space of $((\fPB_n^V)_k)^m$ consisting of $m$-tuples with common interval data and pairwise disjoint underlying manifolds. We write its elements as
\[
(M_1,\ldots,M_m,\I).
\]
For $\mmod\in\{\aff,\kih\}$, set
\[
(\PB_n^{V,\mmod}\langle m\rangle)_k=S^\mmod_\bu((\fPB_n^V\langle m\rangle)_k),~
\PB_n^\mmod\langle m\rangle=\hocolim_{V\subset\R^\infty}\PB_n^{V,\mmod}\langle m\rangle.
\]
Completing $\PB_n^\kih\langle m\rangle$ gives a complete Segal space $\Bord_n\langle m\rangle$. Note that $\Bord_n\langle0\rangle\simeq\Int_\bu^\kih$, which is equivalent to the constant complete Segal space by Proposition \ref{prop: Int constant}.

  For $\mmod\in\{\aff,\kih\}$, define the functor
  \[
A^{\mmod}:\Ga\longrightarrow\sSpace,~\langle m\rangle\longmapsto\PB_n^\mmod\langle m\rangle,
\]  
where to a morphism $f: \langle m \rangle \to \langle k \rangle$ in $\Ga$, we assign the map
  \[
\begin{aligned}
    \PB_n^\mmod \langle m \rangle &\to \PB_n^\mmod \langle k \rangle,\\(M_1, \ldots, M_m, \I) &\mapsto \left( \coprod_{\al \in f^{-1}(1)} M_\al, \ldots, \coprod_{\al \in f^{-1}(k)} M_\al, \I \right).
\end{aligned}\] 
 By Proposition \ref{prop: segal space}, $A^\kih\langle1\rangle$ is a Segal space.
 \edefn
\bprop\label{prop: symmetric condition}
The   stabilized Kihara $\Ga$-object $A^\kih$ satisfies the $\Ga$-condition in the sense of Definition \ref{defn: sym str Ga object}. Its objectwise completion  endows $\Bord_n$ with a symmetric monoidal structure   whose unit is represented by  $(\varnothing,\I)$.
\eprop
\brem
We prove Proposition \ref{prop: symmetric condition} in Appendix \S \ref{sec: symmetric monoidal},  adapted from \cite[Proposition 7.2]{CS2019category}.  
 \erem
\brem\label{rem: affine construction}
Compatible morphisms $F^{V,\aff}:\PB_n^{V,\aff}\to\Psi(\D)$, or a homotopy coherent cocone, induce $F^\aff:\PB_n^\aff\to\Psi(\D)$. Since $\kappa$ is a levelwise weak homotopy equivalence, $F^\aff$ determines a functor $\Bord_n\to\Psi(\D)$ up to a contractible space of choices, as in Remark \ref{rem: construct functor}. Equivalently, one may apply a fibrant replacement functor $\what R$ in $\sSpace_c^{\CSe}$ directly to the affine model.

The affine Segal condition in Proposition \ref{prop: segal space} and Lemma \ref{lem: homotopy invariance fiber product} show that a levelwise Kan replacement $(\PB_n^\aff)^{\mathrm{Kan}}$ is a Segal space. Lemma \ref{lem: smooth mapping path} computes its mapping spaces by the affine path model. Since completion is a DK equivalence,
\[
h(\Bord_n)\simeq h\bigl((\PB_n^\aff)^{\mathrm{Kan}}\bigr),
\]
with composition given by geometric gluing.

The symmetric monoidal version is obtained similarly, provided that the morphisms $F^{V,\aff}$ form a morphism of the corresponding $\Ga$-objects. In the formal approach, $\what{R}$ is applied at the level of $\Ga$-objects. This produces a symmetric monoidal functor by Definition \ref{defn: sym str Ga object}.
\erem

\subsection{Relation to ordinary cobordism category}\label{subsec: Relation to ordinary cobordism category}
In this subsection, we compare the $\infty$-cobordism category   $\Bord_n$  with the ordinary cobordism category $\Cob_n$. The main results in this subsection follow from \cite[\S 8.3]{CS2019category}.
\bdefn[{\cite[\Def 8.18 and the proof of \Prop 8.20]{CS2019category}}]\label{defn: cob cat}
Given $n\in\N_+$, the \emph{$n$-dimensional (ordinary) cobordism category} $\Cob_n$ is defined as follows.

\begin{itemize}
    \item Objects are (possibly disconnected) closed $(n-1)$-dimensional smooth manifolds.
    \item A morphism from $M$ to $N$ is a diffeomorphism class of $n$-dimensional cobordisms from $M$ to $N$. Here an $n$-dimensional cobordism is a compact smooth manifold $\Sigma$ with boundary $\partial \Sigma=\partial_{\iin}\Sigma\sqcup \partial_{\out}\Sigma$, together with specified diffeomorphisms \[\partial_{\iin}\Sigma \cong M \aand \partial_{\out}\Sigma\cong N.\]A diffeomorphism of cobordisms intertwines the specified boundary identifications. We write $\Sigma:M\to N$, call $M$ and $N$ the \emph{incoming end} (or \emph{source}) and the \emph{outgoing end} (or \emph{target}) of $\Sigma$, respectively, and call $M\sqcup N$ the \emph{end} of $\Sigma$.
    \item Composition of morphisms $\Sigma_1: M_0 \to M_1$ and $\Sigma_2: M_1 \to M_2$ is given by the diffeomorphism class of the gluing $\Sigma_1 \sqcup_{M_1} \Sigma_2$\footnote{To make the composition become a smooth manifold, one needs to pick collars of boundaries, though the choice does not affect the diffeomorphism class; see also \cite[Remark 1.1.2]{lurie2009cob}.}.
    \item The identity morphism on $M$ is the diffeomorphism class of the \emph{product cobordism} $[0,1]\times M$ viewed as a morphism from $M$ to $M$.
\end{itemize}
Moreover,   disjoint union endows  $\Cob_n$  with  a symmetric monoidal structure .
\edefn

By Remark \ref{rem: affine construction}, the affine path model computes the morphisms in $h(\Bord_n)$. Up to isomorphism, choose objects $x,y$ represented at a common finite-dimensional stage $V$. Lemma \ref{lem: smooth mapping path} identifies the corresponding homotopy fiber product, up to weak homotopy equivalence, with the affine path model. Its $0$-simplices consist of $(N,\I)\in(\PB_n^{V,\aff})_{1,0}$ together with smooth paths joining $d_1(N,\I)$ to $x$ and $d_0(N,\I)$ to $y$. After stabilization, the connected components of these models compute $\Hom_{h(\Bord_n)}(x,y)$, with composition induced by cobordism gluing. This gives the following proposition.

\bprop[{\cite[Proposition 8.20]{CS2019category}}]\label{prop: homotopy cat cob}
The construction in Remark \ref{rem: 0-simplices} induces a symmetric monoidal equivalence
\[
h(\Bord_n)\simeq\Cob_n.
\]
\eprop

\subsection{Tangential structures}\label{subsec: Tangential structures}
In this subsection, we review the construction of the generalization of the $\infty$-cobordism category that includes manifolds with additional (tangential) structures. We also mention some concrete tangential structures that are used in instanton theory. The main results in this subsection follow from \cite[\S 9]{CS2019category}, though we modify some notation and definitions to be compatible with \cite{GTMW2009cob,GR2010cob}.

\bdefn[{\cite[\Defs 9.1 and 9.2]{CS2019category}, \cite[\S 5]{GTMW2009cob}, and \cite[Definition 2.12]{GR2010cob}}]\label{defn: tangential str} Let $X$ be a topological space and let $\theta: X \to BO(n)$ be a continuous map, where $BO(n)$ is the classifying space of $O(n)$. An \emph{$(X,\theta)$-structure}, or shortly a \emph{$\theta$-structure}, on an $n$-dimensional vector bundle $V\to M$ is a bundle map $\varphi: V\to \theta^*\ga$, where $\ga$ is the canonical $n$-dimensional vector bundle over $BO(n)$ and a bundle map is a continuous map of the total space of the vector bundles that restricts to a linear isomorphism on each fiber. In other words, we have $\varphi=(f,\triv)$, where
\benum
    \item $f:M \to X$ is a continuous map, and
    \item $\triv:V \cong (\theta\circ f)^*\ga$ is an isomorphism of vector bundles over $M$.
\eenu
A $\theta$-structure on a smooth manifold $M$ is a $\theta$-structure on $TM$, which is usually called a \emph{tangential structure}. More generally, for $k\in\N$, a $\theta_{n-k}$-structure on $M$ is a $\theta$-structure on $\R^k\op TM$. The set of $\theta_{n-k}$-structured $(n-k)$-dimensional smooth manifolds is denoted by $\Man_{n-k}^{\theta_{n-k}}$.

A special case of the tangential structure is obtained from a topological group $G$ and a representation $\rho:G\to GL(n,\R)$. Let \[\theta_{\rho}: BG \xra{B\rho} BGL(n,\R)\xra{\simeq}BO(n)\]be the induced map, where the second map is a chosen homotopy inverse to the map induced by $O(n)\hookrightarrow GL(n,\R)$. A $\theta_\rho$-structure is also called a $(G,\rho)$-structure. If $\rho$ is some standard representation, then we also call it $G$-structure. The set of $(G,\rho)$-structured $n$-dimensional smooth manifolds is denoted by $\Man_n^{(G,\rho)}$, or $\Man_n^G$ if $\rho$ is standard.
\edefn
\begin{example}\label{exmp: tangential str}
    Examples of $G$-structures provide usual structures in the following cases.
    \benu
    \item A $GL^+(n,\R)$-structure is an orientation.
    \item An $O(n)$-structure is a Riemannian metric.
    \item Note that $SO(n)=GL^+(n,\R)\cap O(n)$. Then an $SO(n)$-structure is an orientation together with a Riemannian metric.
    \item When $n=2m$ is even, a $GL(m,\C)$-structure is an almost complex structure. Note that $GL(m,\C)\subset GL^+(2m,\R)$, so that an almost complex structure induces an orientation.
    \item When $n=2m$ is even, note that $U(m)=GL(m,\C)\cap O(2m)=GL(m,\C)\cap SO(2m)$. Then a $U(m)$-structure is an almost complex structure together with a compatible Riemannian metric.
    \item When $G=\operatorname{Spin}(n)$ and $\rho:\operatorname{Spin}(n)\to SO(n)\to GL(n,\R)$, then an $\operatorname{Spin}(n)$-structure is a spin structure. We can similarly consider the spin$^c$ structure.
    \eenu
\end{example}
\begin{example}\label{exmp: tangential str 2}
    We also provide examples of $\theta$-structures that are not $G$-structures.

Let $Y$ be a topological space with an $O(n)$-action, set $X=EO(n)\times_{O(n)}Y$, and let $\theta:X\to BO(n)$ be the projection. If $Y$ has a trivial $O(n)$-action, then a $\theta$-structure on $M$ is a map from $M$ to $Y$. If $Y=O(n)/SO(n)\times Z$ with the trivial action on $Z$, or equivalently $X=BSO(n)\times Z$, then a $\theta$-structure on $M$ is a map from $M$ to $Z$ together with an orientation and a Riemannian metric on $M$. We may also take $X=BG\times Z$ for another topological group $G$. In particular, when $Z=BG'$ and $\theta:BG\times BG'\to BO(n)$ is trivial on the second factor, a $\theta$-structure is a $G$-structure together with a principal $G'$-bundle.
\end{example}
    As in the undecorated case, we   define the $\infty$-cobordism category $\Bord_n^\theta$ of $\theta$-structured manifolds by the following procedure.
\benu
    \item Given a fixed finite-dimensional vector space $V$ and $k\in\N$, define $(\fPB_n^{\theta,V})_k$.
     
    \item Equip $(\fPB_n^{\theta,V})_k$ with the tangential diffeology defined below and apply $S^\aff_\bu$ and $S^\kih_\bu$. The affine simplices admit the explicit description below; compare \cite[Definition 2.16 and Lemma 2.17]{GR2010cob}.
    \item Introduce simplicial maps on $(\fPB_n^{\theta,V})_k$ to make $(\PB_n^{\theta,V,\kih})_\bu$ a simplicial space and prove the Segal condition.
    \item Take the homotopy colimit over all subspaces $V\subset\R^\infty$ to obtain $\PB_n^{\theta,\kih}$ and take the completion to obtain the complete Segal space $\Bord_n^\theta$.
    \item Define $\PB_n^{\theta,V,\kih}\langle m\rangle$ for all $m\in\N$ and endow $\Bord_n^\theta$ with a symmetric monoidal structure.
\eenu
For simplicity, we give only the explicit affine description of $(\PB_n^{\theta,V,\aff})_{k,l}$; the remaining constructions are induced from the corresponding diffeological spaces.
\bdefn[{\cite[Definitions 9.8 and 9.10]{CS2019category} and \cite[\S 2.4]{SP2024tqft}}]\label{defn: pbord tangential structure}
For a smooth manifold $N$, let $\psi_n^\theta(N)$ be the set of pairs $(M,\varphi)$, where $M\in\psi_n(N)$ and $\varphi$ is a $\theta$-structure on $M$. More explicitly, given an open set $U\subset\R^r$ and a set-theoretic map
\[
p:U\longrightarrow\psi_n^\theta(N),~u\longmapsto(M_u,\varphi_u),
\]
set
\[
\Gamma_p=\{(u,x)\in U\times N\mid x\in M_u\},~\pi_p:\Gamma_p\longrightarrow U.
\]
Writing $\varphi_u=(f_u,\triv_u)$, we declare $p$ to be a plot if $u\mapsto M_u$ is a plot in $\psi_n(N)$ and there are a continuous map $f_p:\Gamma_p\to X$ and an isomorphism of vector bundles
\[
\triv_p:\ker(D\pi_p)\xrightarrow{\cong}(\theta\circ f_p)^*\ga
\]
such that, for every $u\in U$,
\[
f_p\circ\iota_u=f_u,~\iota_u^*\triv_p=\triv_u,
\]
under the canonical identification
\[
\iota_u^*\ker(D\pi_p)\cong TM_u,
\]
where $\iota_u:M_u\hookrightarrow\Gamma_p$ is the fiber inclusion. These plots define the \emph{tangential diffeology} on $\psi_n^\theta(N)$. This is the specialization of \cite[\S 2.4 and Example 2.11]{SP2024tqft} obtained by taking $B=X$ and $E=\theta^*\ga$.

Let $(\fPB_n^{\theta,V})_k$ consist of tuples
\[
(M, \varphi, \un{I}=(I_0 \leq \cdots \leq I_{k})),
\]
where $(M,\un I)\in(\fPB_n^V)_k$ and $(M,\varphi)\in\Man_n^\theta$. Box rescaling of $(M,\varphi)$ and forgetting $\varphi$ give an inclusion
\[
(\fPB_n^{\theta,V})_k\hookrightarrow\psi_n^\theta(V\times(0,1))\times(\fPB_n^V)_k.
\]
We equip $(\fPB_n^{\theta,V})_k$ with the resulting sub-diffeology.

For $\mmod\in\{\aff,\kih\}$, define
\[
(\PB_n^{\theta,V,\mmod})_k=S^{\mmod}_\bu((\fPB_n^{\theta,V})_k).
\]
  
An $l$-simplex of the affine model is a smooth map
\[
|\Delta^l|_\aff\longrightarrow(\fPB_n^{\theta,V})_k.
\]
Its image under the forgetful map is an element of $(\PB_n^{V,\aff})_{k,l}$, while its tangential data are encoded by the tangential diffeology on $\psi_n^\theta(V\times(0,1))$.

  The forgetful map is smooth and induces
  \[
\Fg:(\PB_n^{\theta,V,\mmod})_k\to(\PB_n^{V,\mmod})_k.
\]
\edefn
As a direct generalization of Lemma \ref{lem: embedded manifold functoriality} and the results in \cite[\S 2.3--2.4, Lemma 2.10, and Example 2.11]{SP2024tqft}, we obtain the following lemma.
   \blem\label{lem: tangential diffeology functoriality}
The tangential diffeology is preserved by the pullback and pushforward families in Lemma \ref{lem: embedded manifold functoriality}. In particular, every morphism $g:[m]\to[k]$ in $\Delta$ induces a smooth map
\[
g^*:(\fPB_n^{\theta,V})_k\longrightarrow(\fPB_n^{\theta,V})_m,
\]
and every linear inclusion $V\hookrightarrow V'$ induces a smooth stabilization map
\[
(\fPB_n^{\theta,V})_k\longrightarrow(\fPB_n^{\theta,V'})_k.
\]
The finite disjoint-union maps used in the $\Ga$-object construction are also smooth. Consequently, $(\fPB_n^{\theta,V})_\bu$ is a simplicial diffeological space, and the constructions in Definition \ref{defn: pbord tangential structure} are simplicial spaces compatible with stabilization and the symmetric monoidal construction.
\elem
\brem\label{rem: smooth geometric models}
Definition \ref{defn: pbord tangential structure} only requires the additional structures to vary continuously. For the structures in Example \ref{exmp: tangential str}, possibly with additional principal Lie group bundles as in Example \ref{exmp: tangential str 2}, we may instead use a smooth model. In this model, the bundles, their identifications and the structures on $\ker(D\pi_p)$ are smooth over the total space $\Gamma_p$ of each plot. After extending fixed smooth face data to compatible collars, we apply \cite[Theorem II.11]{MW2009} to the principal bundles and \cite[Theorem 11]{Wockel2009} to the sections defining the tangential structures and bundle identifications, over parameter disks and homotopies. This gives a levelwise weak equivalence between the smooth and continuous models for fixed $V$.
\erem
  \bdefn
The symmetric monoidal category $\Cob_n^\theta$ has closed $(n-1)$-manifolds  $M$   equipped with  $\theta_{n-1}$  -structures  $\varphi_M$   as objects. A morphism $(M,\varphi_M)\to(N,\varphi_N)$  is an equivalence class of  compact  $\theta$-structured cobordisms   $(\Sigma,\varphi):M\to N$ with specified identifications
\[
\varphi|_{\partial_\iin\Sigma}\cong\varphi_M,~
\varphi|_{\partial_\out\Sigma}\cong\varphi_N,
\] 
  under which the distinguished $\R$-directions  of the boundary   structures correspond  to the inward and outward   normal directions , respectively. The equivalence relation is generated by isomorphisms of $\theta$-structured cobordisms preserving the boundary identifications and homotopies of the $\theta$-structure relative to boundary collars.  Composition is defined by collar gluing, identities are product cobordisms, and the  symmetric monoidal structure   is  disjoint union.
 \edefn
Motivated by Definition \ref{defn: rigidity}, if we can regard a product cobordism of $(M,\varphi_M)$ as a cobordism from $\varnothing$ to $(M,\varphi_M)\sqcup (M,\ov{\varphi}_M)$, then $(M,\ov{\varphi}_M)$ is the dual of $(M,\varphi_M)$. In other words, if there is an automorphism on the set of $\theta$-structures on $\R\times TM$ that reverses the direction of $\R$, then we obtain a dual of $(M,\varphi_M)$. The dual $(\Sigma,\ov{\varphi})$ of a $\theta$-structured cobordism $(\Sigma,\varphi):(M,\varphi_M)\to(N,\varphi_N)$ has the same underlying $\theta$-structured manifold but different specified boundary identifications, so that it defines a morphism $(\Sigma,\ov{\varphi}):(N,\ov{\varphi}_N)\to(M,\ov{\varphi}_M)$.
\begin{example}\label{exmp: dual}
We provide the dual of $\theta_{n-1}$-structured manifold $(M,\varphi_M)$ for those in Example \ref{exmp: tangential str}.
\begin{itemize}
    \item The dual of an orientation on $\R\times TM$ is obtained by reversing both orientations on $\R$ and $TM$, i.e.\ the dual of an oriented manifold $M$ is its orientation reversal $-M$. Note that the dual of an oriented cobordism has the same orientation rather than the opposite one.
    \item The dual of a Riemannian metric is itself, which does not affect the metric on $TM$.
    \item The dual of an almost complex structure $J$ is itself, but regarded as the one on $(-\R)\times T(-M)=\R\times TM$. For simplicity, we still write $J$ for the dual.
    \item The dual of a spin or spin$^c$ structure is the one induced by $-M$, i.e.\ adding a minus sign to the Clifford multiplication.
\end{itemize}
\end{example}
As a direct generalization of \Prop \ref{prop: homotopy cat cob}, we obtain the following.
\bprop\label{prop: homotopy cat cob theta}
There is an equivalence of categories   \[h(\Bord_n^\theta)\simeq \Cob_n^{\theta}.\] Moreover, this equivalence respects the symmetric monoidal structures. If $\theta$ is obtained from $X=BG\times Z$ in Example \ref{exmp: tangential str 2} for $G$ from Example \ref{exmp: tangential str}, then both   $\Bord_n^\theta$  and $\Cob_n^{\theta}$ are rigid and the equivalence is a rigid functor.
\eprop

Finally, we show that the Riemannian metric does not affect the $\infty$-cobordism category. More precisely, we have the following proposition. Compare \cite[Example 2.4.22]{lurie2009cob}.

\bprop\label{prop: metric}
When $G=O(n)$ in Example \ref{exmp: tangential str}, the forgetful functor  \[\Fg:\Bord_n^{O(n)}\to \Bord_n\] is an equivalence of complete Segal spaces.
\eprop
\bpf
View an $O(n)$-structure as a Riemannian metric. At every simplicial level and for every $V$, restriction of the ambient Euclidean metric on $V\times\R$ defines a smooth section of the forgetful map. Convex interpolation from a given metric to the induced metric is a smooth homotopy from the identity to the composite of this section with the forgetful map, while the opposite composite is the identity. This construction is compatible with the spatial and simplicial operators. Lemma \ref{lem: affine Kihara comparison} therefore gives a levelwise weak homotopy equivalence
\[
\Fg^V:(\PB_n^{O(n),V,\kih})_{\bu,\bu}\longrightarrow(\PB_n^{V,\kih})_{\bu,\bu}.
\]
Filtered homotopy colimits preserve this equivalence, and completion preserves the resulting DK equivalence.
\epf

Based on Proposition \ref{prop: metric}, to construct a functor \[F:\Bord_n\to\cC\] for a complete Segal space $\cC$, we can first construct \[F':\Bord_n^{O(n)}\to\cC\] and then apply the equivalence of complete Segal spaces. Thus we may add a Riemannian metric when constructing a TQFT. We will use this observation for the instanton TQFT.

\subsection{Cylindrical condition}
In gauge theory, the tangential data are usually required to be translation invariant near the ends. We use a fixed-collar strengthening of the cylindricality condition in \cite[Definition 5.6]{SP2024tqft}.

\bdefn\label{defn: cylindrical}
A tuple
\[
(M,\varphi,\un I=(I_0\leq\cdots\leq I_k))\in(\fPB_n^{\theta,V})_k
\]
is \emph{cylindrical on collars} if, for every $j\in[0,k]$, setting
\[
Y_j=\pi^{-1}(m(I_j)),
\]
as in Remark \ref{rem: 0-simplices}, we have
\[
M\cap\bigl(V\times[m_{1/3}(I_j),m_{2/3}(I_j)]\bigr)=Y_j\times[m_{1/3}(I_j),m_{2/3}(I_j)],
\]
and the restriction of $\varphi$ to this product is induced by a $\theta_{n-1}$-structure on $Y_j$ under $T(Y_j\times[m_{1/3}(I_j),m_{2/3}(I_j)])\cong TY_j\oplus\R$. If $k=1$, we also say that the tuple is \emph{cylindrical on the ends}.

Let $(\fPB_{n,\cyl}^{\theta,V})_k$ be the resulting sub-diffeological space of $(\fPB_n^{\theta,V})_k$. For $\mmod\in\{\aff,\kih\}$, set
\[
(\PB_{n,\cyl}^{\theta,V,\mmod})_k=S_\bu^{\mmod}((\fPB_{n,\cyl}^{\theta,V})_k).
\]
Deleting or repeating intervals preserves the cylindrical condition, so these spaces form a simplicial diffeological space.
\edefn

Let $\Bord_{n,\cyl}^{\theta}$ denote the stabilization and completion of the cylindrical model. Then we have the following proposition.
\bprop[{\cite[Lemma 3.4]{GR2010cob} and \cite[Lemmas 6.4, 6.5, and 6.7]{SP2024tqft}}]\label{prop: cylindrical}
For every $k\in\N$, the inclusion induces a weak homotopy equivalence
\[
(\PB_{n,\cyl}^{\theta,V,\kih})_{k,\bu}\longrightarrow(\PB_n^{\theta,V,\kih})_{k,\bu}.
\]
Consequently, $(\PB_{n,\cyl}^{\theta,V,\kih})_{\bu,\bu}$ satisfies the Segal condition and
\[
\Bord_{n,\cyl}^{\theta}\simeq\Bord_n^{\theta}.
\]
\eprop

\bcor\label{cor: cylindrical representative}
Every morphism in $h(\Bord_n^\theta)$ is isomorphic to one represented by a cobordism cylindrical on the ends.
\ecor
\brem\label{rem: cylindrical product metrics}
The collar adjustments in the proof of \cite[Lemma 3.4]{GR2010cob} preserve smooth families, so Proposition \ref{prop: cylindrical} also holds in the smooth model. Together with Remark \ref{rem: smooth geometric models}, this gives a levelwise weak equivalence between the smooth and continuous cylindrical models.

For the $U(m)$-structures in Example \ref{exmp: tangential str}, we may require the compatible metric on each cylindrical collar $Y\times I$ to be $dt^2+g$, where $t$ is the interval coordinate and $g$ is a Riemannian metric on $Y$.

Indeed, let $J$ be the almost complex structure and $\bfg$ its compatible Riemannian metric on a collar, both independent of $t$. Set
\[
\bfg_s=(1-s)\bfg+s(dt^2+\bfg|_{TY}),~0\le s\le1,
\]
and let $J_s$ be obtained from $J$ by transport along the positive bundle isometry from $\bfg$ to $\bfg_s$. The pair $(J_s,\bfg_s)$ remains cylindrical and compatible, and is unchanged if $\bfg$ is already a product metric. The constructions agree on overlapping collars. Relative extension of the metrics, followed by the same transport of $J$, extends this deformation over the cobordism in parameter families relative to fixed product data. Hence the inclusions of the spaces satisfying the product condition induce levelwise weak homotopy equivalences. The same argument applies with an additional principal bundle as in Example \ref{exmp: tangential str 2}, which is kept fixed throughout the deformation.
\erem

\section{\texorpdfstring{$\infty$}{Infinity}-cobordism category for instantons}\label{sec: instanton cob cat}

In this section, we modify the construction of \S\ref{sec: Preliminaries in cobordism category} to obtain $\infty$-cobordism categories for instantons, with and without auxiliary data.
\subsection{Basic setup}\label{subsec: basic setup}
Let $n=2m=4$, let $X_a=BU(2)\times BSO(3)$ be as in Example \ref{exmp: tangential str 2}, and let $\theta_a$ be the corresponding map, where the subscript $a$ indicates ``auxiliary data''. Using the smooth model of Remark \ref{rem: smooth geometric models}, with product metrics on collars as in Remark \ref{rem: cylindrical product metrics}, we describe $\Cob_n^{\theta_a}$ by objects $(Y,P,J,g)$, where
\begin{itemize}
\item $Y$ is a (possibly disconnected and empty) closed oriented $3$-manifold;
\item $P$ is a principal $SO(3)$ bundle over $Y$;
\item $J$ is a translation-invariant almost complex structure on $\R\times Y$ that induces the orientation of $Y$;
\item $g$ is a Riemannian metric on $Y$ such that the product metric $dt^2+g$ is compatible with $J$.
\end{itemize}
Its morphisms are equivalence classes of $(W,\bfP,\bfJ,\bfg)$, where
\begin{itemize}
\item $W:Y_-\to Y_+$ is an oriented $4$-dimensional cobordism;
\item $\bfP$ is a principal $SO(3)$ bundle over $W$;
\item $\bfJ$ is an almost complex structure on $W$ that induces the orientation of $W$;
\item $\bfg$ is a Riemannian metric on $W$ compatible with $\bfJ$.
\end{itemize}
Two such tuples represent the same morphism if there is a diffeomorphism preserving the boundary identifications under which their $\theta_a$-structures are homotopic relative to boundary collars.

By Corollary \ref{cor: cylindrical representative}, we may choose a representative that is cylindrical on the ends. Following the convention $\R\oplus TY$ in Definition \ref{defn: tangential str}, we use interval-first cylindrical coordinates and consider the open manifold
\begin{equation}\label{eq: open end}
W^*=\big((-\infty,0]\times Y_-\big)\cup W\cup\big([1,\infty)\times Y_+\big).
\end{equation}
Here $\partial W=(-Y_-)\sqcup Y_+$, and the cylinders carry the interval-first product orientation. The interval-last notation of Definition \ref{defn: cylindrical} is understood through the factor-switching map $(y,t)\mapsto(t,y)$, with all tangential data transported accordingly. We call $W^*\setminus W$ the \emph{open ends} of $W$ and extend the tangential structures to them as translation-invariant product structures.

Let $X$ be obtained from $X_a$ by replacing $BU(2)$ with $BGL(2,\C)$, i.e.\ by forgetting the Riemannian metric, and let $\theta$ be the corresponding map from Example \ref{exmp: tangential str 2}. In this section, we construct two modified $\infty$-cobordism categories $\BI$ and $\BIa$ for instantons from $\Bord_{n,\cyl}^{\theta,\infty}$ and $\Bord_{n,\cyl}^{\theta_a,\infty}$, respectively. Motivated by the ordinary cobordism categories for instantons in \cite[\S 4--5]{kronheimer2011khovanov}, we incorporate the following information. See also \S \ref{subsect: Sketch of gauge theory construction}.

\begin{itemize}
\item \textbf{Subgroup}, cf.\ \cite[\S 5.2]{kronheimer2011khovanov}.
\end{itemize}

Objects of $\BI$ include a subgroup $\phi\subset H^1(Y;\Z/2)$, and morphisms include a subgroup $\Phi\subset H^1(W;\Z/2)$. For $(W,\Phi):(Y_-,\phi_-)\to(Y_+,\phi_+)$, the image of $\Phi$ under restriction to $Y_\pm$ is $\phi_\pm$. If $W=W_1\circ W_0$ and $\rho$ denotes the restriction map, composition is defined by
\[
\Phi=\rho^{-1}(\Phi_0\oplus\Phi_1)\subset H^1(W;\Z/2).
\]
Thus $\Phi$ is the maximal subgroup whose restrictions lie in $\Phi_0$ and $\Phi_1$, as required in the discussion following \cite[Proposition 5.4]{kronheimer2011khovanov}. By the Mayer--Vietoris sequence, $\Phi$ contains the image of
\[
H^0(W_0\cap W_1;\Z/2)\longrightarrow H^1(W;\Z/2),
\]
represented by the Poincar\'e duals of the components of the common boundary. For a smooth family parametrized by a smooth manifold with corners, we require that, locally on the parameter space, any smooth trivialization identifies the subgroups in all fibers with one fixed subgroup.

From Example \ref{exmp: dual}, the dual of an object $(Y,P,J,g)$ is $(-Y,P,J,g)$, where the same $J$ is viewed as an almost complex structure under the identification $T((-\R)\times(-Y))\cong T(\R\times Y)$. We use the same notation for the corresponding subgroup $\phi\subset H^1(-Y;\Z/2)$.

\begin{itemize}
\item \textbf{Non-integrality}, cf.\ \cite[Definitions 3.1 and 5.1]{kronheimer2011khovanov}.
\end{itemize}

Objects of $\BI$ and $\BIa$ are required to satisfy the \emph{non-integral} condition, i.e., there exists an embedded oriented surface $\Sigma\subset Y$ such that, for every component $Y_\al$ of $Y$, writing $\Sigma_\al=\Sigma\cap Y_\al$, we have
\begin{equation}\label{eq: witness}
   \langle w_2(P),[\Sigma_\al]\rangle\neq0~\text{and}~\phi|_{\Sigma_\al}=0 
\end{equation}
We call $\Sigma$ a \emph{witness surface} for the non-integral condition. Our convention is that the empty manifold satisfies the non-integral condition.

\begin{itemize}
\item \textbf{Exclusion of closed components}
\end{itemize}

As mentioned in the introduction, we need to exclude closed components of morphisms. As in the ordinary cobordism category, we introduce the sub-$\infty$-categories
\[
\BI^{\pa_+},~\BI^{\pa_-},~\BI^{\pa}\subset\BI
\]
in which the morphisms and higher morphisms have positive, negative, and full boundaries, respectively. We similarly define the corresponding sub-$\infty$-categories of $\BIa$.

\begin{itemize}
\item \textbf{Perturbation}, cf.\ \cite[\S 3.4, 3.5, 3.7, 3.9]{kronheimer2011khovanov}.
\end{itemize}

Objects of $\BIa$ also include an \emph{admissible perturbation} $\pi$. A perturbation lies in a Banach space $\cP(Y)$ of real sequences $\pi=\{\pi_i\}_{i\in Q}$ indexed by a countable set $Q$ and satisfying
\[
\|\pi\|_{\cP(Y)}=\sum_{i\in Q}C_i|\pi_i|<\infty
\]
for fixed positive weights $C_i$. The set $Q$ and the weights depend on $(Y,P)$, but not on $g$, $J$, or the orientation of $Y$. An element of $Q$ corresponds to a collection of immersions $\bfq=(q_1,\dots,q_l)$, where $q_j:S^1\times D^2\to Y$ and the maps agree on $\{p\}\times D^2$ for a fixed $p\in S^1$. Such a collection induces a gauge-invariant function $f_{\bfq}$ on the space of connections in $P$, and
\[
f_\pi=\sum_{i\in Q}\pi_i f_{\bfq_i}.
\]
The choice of the weights ensures convergence and the analytic bounds required in gauge theory. The function $f_\pi$ is called a \emph{holonomy perturbation}. For a connection $A=B+c\,dt$ on $\R\times P$, it induces a self-dual $2$-form $\widehat V_{\pi,g}(A)$ with respect to the product metric $dt^2+g$.

When $\phi\subset H^1(Y;\Z/2)$ is nontrivial, the holonomy perturbation must be $\phi$-invariant as in \cite[p. 153]{kronheimer2011khovanov}. We suppress this dependence from the notation.

In instanton theory, we consider the solution space of the perturbed ASD equation
\begin{equation}\label{eq: perturbed ASD}
F_A^{+_g}+\widehat V_{\pi,g}(A)=0,
\end{equation}
where
\[
F_A^{+_g}=\frac{1}{2}(F_A+*_gF_A).
\]
Here $*_g$ is the Hodge star for the product metric $dt^2+g$ on $\R\times Y$. Thus both the self-dual part $F_A^{+_g}$ and the perturbation term $\widehat V_{\pi,g}(A)$ depend on $g$. Let $\cM(Y,P,g,J,\pi)$ denote the moduli space of solutions modulo the determinant-$1$ gauge group, where $J$ orients the moduli space. If $\phi$ is nontrivial, quotienting by the corresponding larger gauge group gives
\[
\cM(Y,P,g,J,\phi,\pi)=\cM(Y,P,g,J,\pi)/\phi.
\]
Let $\breve \cM(Y,P,g,J,\phi,\pi)$ denote the quotient of the nonconstant trajectories by the $\R$-translation. The differential and the relation $\partial^2=0$ are obtained from the compactified moduli spaces of dimensions $0$ and $1$, respectively.

The admissibility condition is modeled on \cite[Hypothesis 3.6]{kronheimer2011khovanov}: the perturbed critical points are irreducible and nondegenerate, and the relevant trajectory moduli spaces are regular. Let
\begin{equation}\label{eq: adm Y}
    \cP^{\adm}(Y)\subset\cP(Y)
\end{equation}
be the subset of admissible perturbations, where we suppress the dependence on the remaining data. The almost complex structure determines the absolute $\Z/2$-grading of the chain complex.

By \cite[Proposition 3.3, Lemma 3.4, and Proposition 3.5]{kronheimer2011khovanov}, admissible perturbations exist. For the dual structure, there is a canonical identification of unoriented moduli spaces
\begin{equation}\label{eq: moduli dual}
    \breve{\cM}(-Y,P,g,J,\phi,-\pi)=\breve{\cM}(Y,P,g,J,\phi,\pi),
\end{equation}
The orientation comparison in the monopole setting is given in \cite[Corollary 22.5.5]{kronheimer2007monopoles}. The chain complex on the left is identified with the cochain complex on the right, with signs in the dual differential depending on degree.

A morphism with auxiliary data further includes a pair $\Pi=(\pi,\pi')$, where $\pi\in\cP^{\adm}(\partial W)$ is the perturbation on the ends and
\[
\pi'=(\pi'_-,\pi'_+)\in\cP(W)=\cP(-Y_-)\times\cP(Y_+)
\]
is a secondary perturbation satisfying the regularity condition of \cite[Proposition 3.8]{kronheimer2011khovanov}. Together they induce the twice-perturbed ASD equation
\begin{equation}\label{eq: perturbed ASD 2}
F_A^{+_{\bfg}}+\eta(t)\widehat V_{\pi,\bfg}(A)+\psi(t)\widehat V_{\pi',\bfg}(A)=0
\end{equation}
on the open ends of $W^*$, where $\eta$ is a fixed cut-off function and $\psi$ is a fixed bump function. Since the secondary term is supported away from infinity, this equation has the same asymptotic critical points as \eqref{eq: perturbed ASD}. Away from the open ends the equation is unperturbed. The secondary perturbation is admissible if
\[
\cM(W,\bfP,\bfg,\bfJ,\Phi,\Pi)
\]
is regular. Its admissibility depends on $W$ and $\pi$, even though $\pi'$ is supported near $\partial W$. Given $\pi\in\cP^{\adm}(\partial W)$, let
\[
\cP^{\adm}(W,\pi)\subset\cP(W)
\]
denote the admissible secondary perturbations, suppressing the remaining data, and set
\[
\cP^{\adm}(W)=\bigcup_{\pi\in\cP^{\adm}(\partial W)}\{\pi\}\times\cP^{\adm}(W,\pi).
\]
The $0$- and $1$-dimensional moduli spaces define the instanton cobordism map and verify the chain-map relation, respectively.

For a higher morphism, let $G$ be a finite-dimensional oriented smooth manifold with corners parametrizing a smooth family of cobordisms with fixed ends, called a \emph{parameter space}. Here a codimension-$r$ \emph{face} means the closure of a connected component of the locus where exactly $r$ local boundary coordinates vanish. We assume that all faces are embedded and admit compatible collars. We write
\[
(W_q,\bfP_q,\bfg_q,\bfJ_q,\Phi_q,\Pi_q)
\]
for the fiber over $q\in G$. A smooth family $\Pi_G=(\pi,\pi'_G)$ determines the parametrized moduli space
\begin{equation}\label{eq: parametrized moduli space}
\cM_G(W,\bfP,\bfg,\bfJ,\Phi,\Pi_G)
=\bigsqcup_{q\in G}\{q\}\times\cM(W_q,\bfP_q,\bfg_q,\bfJ_q,\Phi_q,\Pi_q).
\end{equation}
For each stratum $F\subset G$, write $\cM_F(W,\bfP,\bfg,\bfJ,\Phi,\Pi_G)$ for the restriction of \eqref{eq: parametrized moduli space} to $F$. We call $\Pi_G$ \emph{admissible} if $\cM_F(W,\bfP,\bfg,\bfJ,\Phi,\Pi_G)$ is regular as a moduli space of pairs $(q,A)$ for every stratum $F$ of $G$. Thus the parametrized equations are transverse on every stratum. Together with the fiberwise orientations determined by $\bfJ$, the orientation of $G$ orients the parametrized moduli space over $\operatorname{int}(G)$. We orient the faces of $G$ by the outward-normal-first convention, iterated in the specified order on higher-codimensional strata.

We write
\[
\cP_G^{\adm}(W,\pi)\subset C^\infty(G,\cP(W))
\]
for the admissible smooth families of secondary perturbations and set
\begin{equation}\label{eq: adm W}
    \cP_G^{\adm}(W)
=
\bigcup_{\pi\in\cP^{\adm}(\partial W)}
\{\pi\}\times\cP_G^{\adm}(W,\pi).
\end{equation}
For a closed union of faces $G'\subset G$, let $\cP_{G'}^{\adm}(W)$ consist of the pairs $(\pi,\pi'_{G'})$, where $\pi\in\cP^{\adm}(\partial W)$, such that $\pi'_{G'}$ is smooth on every face contained in $G'$, the restrictions from intersecting faces agree on their intersection, and the corresponding parametrized moduli space over every stratum of $G'$ is regular.

We have the following extension property for admissible families of secondary perturbations.

\bprop\label{prop: extension of perturbation}
Let $G$ be a compact finite-dimensional oriented smooth manifold with corners, and let $G'\subset G$ be a closed union of faces. Assume that every connected component of each $W_q$ has nonempty boundary. If $(\pi,\pi'_{G'})\in\cP_{G'}^{\adm}(W)$, then there exists $(\pi,\pi'_G)\in\cP_G^{\adm}(W)$ such that
\[
\pi'_G|_{G'}=\pi'_{G'}.
\]
In particular, taking $G'=\varnothing$ shows that $\cP_G^{\adm}(W,\pi)$ is nonempty for every $\pi\in\cP^{\adm}(\partial W)$.
\eprop
\bpf
For a parameter manifold without boundary and without given boundary data, this follows from \cite[Proposition 3.10]{kronheimer2011khovanov} and the subsequent discussion of families of cobordisms, using the argument of \cite[Proposition 24.4.10]{kronheimer2007monopoles}. For the relative statement, first extend the given family smoothly from $G'$ using compatible collars. We then apply the Sard--Smale argument of \cite[Proposition 24.4.10]{kronheimer2007monopoles} inductively over the strata of $G$, beginning with the strata of greatest codimension and using perturbations supported away from $G'$ and the strata already treated. This gives an extension regular on every stratum, as required. This stratumwise regularity is the convention used for cornered parameter spaces in \cite[\S 3.9]{kronheimer2011khovanov}.
\epf
\subsection{Constructions without auxiliary data}
In this subsection, we incorporate the subgroup, non-integrality, and boundary conditions into the diffeological construction of \S \ref{subsec: Tangential structures}. We use the tangential cobordism notation of Definition \ref{defn: pbord tangential structure}, its cylindrical variant from Definition \ref{defn: cylindrical}, and the singular functors recalled in Example \ref{exmp: singular complex} and \S \ref{subsec: diffeological spaces}.

\bdefn\label{defn: PBIV k0}
Fix a finite-dimensional vector space $V$ and $k\in\N$. Let $(\mathfrak{PBI}^V)_k$ be the set of tuples
\[
(M,\varphi,\Phi,\un I=(I_0\leq\cdots\leq I_k)),
\]
where $(M,\varphi,\un I)\in(\fPB_{n,\cyl}^{\theta,V})_k$ and $\Phi\subset H^1(M;\Z/2)$ is a subgroup. Let
\[
Y_i=\pi^{-1}(m(I_i))~\text{for }i\in[0,k]
\]
and, for $k>0$, let
\[
M_j=\pi^{-1}([m_{1/3}(I_j),m_{2/3}(I_{j+1})])~\text{for }j\in[0,k-1],
\]
as in Remark \ref{rem: 0-simplices}. Write $\rho_j:H^1(M;\Z/2)\to H^1(M_j;\Z/2)$ for restriction and require
\begin{equation}\label{eq: maximal subgroup}
\Phi=\rho_{\un I}^{-1}\left(\bigoplus_{j=0}^{k-1}\rho_j(\Phi)\right),~
\rho_{\un I}=\bigoplus_{j=0}^{k-1}\rho_j.
\end{equation}
For $k=0$, no maximality condition is imposed. Thus, for $k>0$, $\Phi$ is the largest subgroup with the given restrictions to the cobordisms $M_j$, as in the discussion following \cite[Proposition 5.4]{kronheimer2011khovanov}. In particular, the Mayer--Vietoris sequence implies that $\Phi$ contains the Poincar\'{e} duals of the components of $Y_i$ for $i\in[1,k-1]$.

Let $P_i$ be the restriction of the $SO(3)$-bundle to $Y_i$ and let
\[
\phi_i=\operatorname{im}(\Phi\to H^1(Y_i;\Z/2)).
\]
We require $(Y_i,P_i,\phi_i)$ to be non-integral for every $i$. Explicitly, for every component $Y_{i,\alpha}$, there is a closed oriented surface $\Sigma_{i,\alpha}\subset Y_{i,\alpha}$ such that
\[
\langle w_2(P_i),[\Sigma_{i,\alpha}]\rangle\neq0~\text{and}~
\phi_i|_{\Sigma_{i,\alpha}}=0.
\]
The empty manifold is non-integral by convention.

We equip $(\mathfrak{PBI}^V)_k$ with the following diffeology. A map
\[
p:U\longrightarrow(\mathfrak{PBI}^V)_k
\]
is a plot if its underlying family is a plot of $(\fPB_{n,\cyl}^{\theta,V})_k$ and, locally on $U$, any smooth trivialization of the family identifies all subgroups $\Phi_u\subset H^1(M_u;\Z/2)$ with one fixed subgroup. This condition is independent of the trivialization because a smooth family of fiber diffeomorphisms acts locally constantly on cohomology.

Let $(\mathfrak{PBI}^{V,\pa_+})_k$, $(\mathfrak{PBI}^{V,\pa_-})_k$, and $(\mathfrak{PBI}^{V,\pa})_k$ be the sub-diffeological spaces obtained by requiring every $M_j$ to have positive, negative, and full boundary, respectively. Here, for $j\in[0,k-1]$, positive boundary means that $\pi_0(Y_{j+1})\to\pi_0(M_j)$ is surjective, negative boundary means that $\pi_0(Y_j)\to\pi_0(M_j)$ is surjective, and full boundary means both.
\edefn

\brem
By the Mayer--Vietoris sequence, requiring $\Phi$ to contain the Poincar\'{e} duals of the intermediate boundaries gives only
\[
\Phi=\rho_{\un I}^{-1}(\rho_{\un I}(\Phi)).
\]
This suffices for the composition formula in \cite[Proposition 5.4]{kronheimer2011khovanov}. In \eqref{eq: maximal subgroup}, we take $\Phi$ to be the largest subgroup with the given restrictions to the cobordisms $M_j$, as in \cite[p. 158]{kronheimer2011khovanov}. Thus the subgroup on a composite is determined by the subgroup data on its pieces, as required by the Segal condition. For example, the diagonal subgroup $\langle(1,1)\rangle\subset(\Z/2)\oplus(\Z/2)$ projects onto both factors, but the largest subgroup with these projections is $(\Z/2)\oplus(\Z/2)$.
\erem

\bdefn\label{defn: PBIV kl}
For $\mmod\in\{\aff,\kih\}$, define
\[
(\PBI^{V,\mmod})_k=S^{\mmod}_\bu((\mathfrak{PBI}^V)_k),
\]
and define the variants with boundary conditions similarly. The affine $l$-simplices are smooth families over $|\Delta^l|_{\aff}$ whose subgroup data satisfy the local triviality condition in Definition \ref{defn: PBIV k0}.

The spatial maps are induced by smooth base change. For $g:[m]\to[k]$, set $M^g=\pi^{-1}(B(g^*\un I))$ and let $\Phi|_{M^g}$ denote the image of $\Phi$ under restriction. For $m>0$, let $(M^g_\nu)_{\nu\in[0,m-1]}$ be the composable cobordisms determined by $(M^g,g^*\un I)$ as in Remark \ref{rem: 0-simplices}, and let
\[
\rho^g_\nu:H^1(M^g;\Z/2)\longrightarrow H^1(M^g_\nu;\Z/2)
\]
be restriction. Define
\[
\Phi^g=(\rho^g)^{-1}\left(\bigoplus_{\nu=0}^{m-1}\rho^g_\nu(\Phi|_{M^g})\right),~\rho^g=\bigoplus_{\nu=0}^{m-1}\rho^g_\nu.
\]
For $m=0$, set $\Phi^g=\Phi|_{M^g}$. If $h:[l]\to[m]$, each $M^{g\circ h}_\mu$ is either a union of consecutive $M^g_\nu$ or a cylindrical identity cobordism arising from a repeated interval. Successive applications of the Mayer--Vietoris sequence in the first case, and the canonical cohomology identification in the second, give $(\Phi^g)^h=\Phi^{g\circ h}$. Hence $(g\circ h)^*=h^*g^*$. These maps preserve plots and define simplicial diffeological spaces $(\mathfrak{PBI}^V)_\bu$ and their boundary variants. Applying $S^{\mmod}_\bu$ levelwise gives simplicial spaces $(\PBI^{V,\mmod})_{\bu,\bu}$. Each simplicial level $(\PBI^{V,\kih})_{k,\bu}$ is a Kan complex by Lemma \ref{lem: affine Kihara comparison}.
\edefn

\bprop\label{prop: segal space PBIV}
Both $(\PBI^{V,\aff})_{\bu,\bu}$ and $(\PBI^{V,\kih})_{\bu,\bu}$ satisfy the Segal condition. Thus the former is a Segal simplicial set and the latter is a Segal space. The same holds for the variants with boundary conditions.
\eprop
\bpf
Fix $m,l\in\N$ and set $k=m+l$. If $m=0$ or $l=0$, one of the two maps to $(\PBI^{V,\aff})_0$ is the identity and hence a Kan fibration. Thus the homotopy fiber product is computed by the strict pullback, and the Segal map is a weak homotopy equivalence. We may therefore assume $m,l>0$.

By Lemma \ref{lem: smooth mapping path}, the target of the affine Segal map is weakly equivalent to the affine path model in Appendix \S \ref{sec: segal}. The cutting, gluing, and rescaling constructions there apply to the underlying cobordisms. We verify that they lift to the subgroup data.

Suppose a matching path identifies the outgoing end of $(M^-,\Phi^-)$ with the incoming end of $(M^+,\Phi^+)$, and write
\[
M=M^-\cup_YM^+,~\rho:H^1(M;\Z/2)\longrightarrow H^1(M^-;\Z/2)\oplus H^1(M^+;\Z/2).
\]
Since the path is a family of objects with locally constant subgroup data, it identifies the images of $\Phi^-$ and $\Phi^+$ in $H^1(Y;\Z/2)$. Define
\begin{equation}\label{eq: gluing subgroup}
\Phi=\rho^{-1}(\Phi^-\oplus\Phi^+).
\end{equation}
The Mayer--Vietoris sequence shows that the restrictions $\Phi\to\Phi^-$ and $\Phi\to\Phi^+$ are surjective. Indeed, for $a^-\in\Phi^-$, choose $a^+\in\Phi^+$ with the same restriction to $Y$. Then $(a^-,a^+)$ lifts to $H^1(M;\Z/2)$, and every such lift lies in $\Phi$. The other restriction is treated similarly.

Let $(M_j,\Phi_j)_{j\in[0,k-1]}$ be the composable cobordisms and their subgroup data. Since $\Phi^-$ and $\Phi^+$ satisfy \eqref{eq: maximal subgroup}, formula \eqref{eq: gluing subgroup} is equivalently
\[
\Phi=\{a\in H^1(M;\Z/2)~\mid~a|_{M_j}\in\Phi_j\text{ for every }j\in[0,k-1]\}.
\]
Thus $\Phi$ satisfies \eqref{eq: maximal subgroup}. On subgroup data, cutting and gluing are inverse, iterated gluing is associative, and the construction commutes with the face and degeneracy maps. In smooth local trivializations the cohomology groups, restriction maps, and subgroups are constant, so the construction preserves plots and commutes with smooth base change. The gluing and rescaling homotopies of Appendix \S \ref{sec: segal} therefore lift to the subgroup-decorated affine model. Non-integrality and the boundary conditions are preserved as well. Hence the affine path-model Segal map is a simplicial homotopy equivalence, and Lemma \ref{lem: smooth mapping path} gives the affine Segal condition.

Finally, naturality of the comparison in Lemma \ref{lem: affine Kihara comparison} gives a commutative square between the affine and Kihara Segal maps. The induced map on their sources is a weak homotopy equivalence by that lemma, while the induced map on their targets is a weak homotopy equivalence by Lemma \ref{lem: homotopy invariance fiber product}. The two-out-of-three property from Lemma \ref{lem: homotopy properties} gives the Kihara Segal condition. Its levels are Kan by Lemma \ref{lem: affine Kihara comparison}, so it is a Segal space. The same argument applies to the variants with boundary conditions.
\epf

\bdefn\label{defn: BI}
For $\mmod\in\{\aff,\kih\}$, let
\[
\PBI^{\mmod}=\hocolim_{V\subset\R^\infty}(\PBI^{V,\mmod})_{\bu,\bu}.
\]
Let $\BI$ be the completion of $\PBI^{\kih}$, and define $\BI^{\pa_+}$, $\BI^{\pa_-}$, and $\BI^{\pa}$ similarly. Their symmetric monoidal structures are induced by disjoint union as in Definition \ref{defn: symmetric monoidal str}. The dual of $(Y,P,J,\phi)$ is $(-Y,P,J,\phi)$.
\edefn
Combining Propositions \ref{prop: homotopy cat cob theta} and \ref{prop: cylindrical} with the construction above gives the following.

\bprop\label{prop: homotopy cat BI}
There is a duality-preserving symmetric monoidal equivalence
\[
h(\BI)\simeq\CobI,
\]
where $\CobI$ is the instanton cobordism category defined in the introduction.
\eprop

\subsection{Abstract description of auxiliary data}\label{subsec: abstract auxiliary data}
In this subsection, we organize auxiliary data abstractly, following the family cobordism maps in \cite[\S3.9, Formula (23)]{kronheimer2011khovanov} and the description of functors to a dg nerve in Lemma \ref{lem: dg nerve of a poset} and Remark \ref{rem: dg nerve of simplex and product}.

We use the metric refinement $\theta_a$ of $\theta$ from \S\ref{subsec: basic setup}. Let $\PBI_a^{V,\aff}$ denote the simplicial space of Definitions \ref{defn: PBIV k0} and \ref{defn: PBIV kl}, using the smooth model of Remark \ref{rem: smooth geometric models}, with $\theta$ replaced by $\theta_a$ and product metrics on collars as in Remark \ref{rem: cylindrical product metrics}. For a fixed cylindrical almost complex structure, the space of all compatible cylindrical metrics is nonempty and convex. The argument of Proposition \ref{prop: metric}, together with Remarks \ref{rem: smooth geometric models} and \ref{rem: cylindrical product metrics}, Proposition \ref{prop: cylindrical}, and Lemma \ref{lem: affine Kihara comparison}, therefore shows that the forgetful map
\begin{equation}\label{eq: Fg V}
\Fg^V:\PBI_a^{V,\aff}\longrightarrow\PBI^{V,\aff}
\end{equation}
is a levelwise weak homotopy equivalence. Thus this refinement presents the same $\infty$-cobordism category after taking the Kihara model and the completion.

 For each $\sigma\in(\PBI_a^{V,\aff})_{k,l}$, we seek a collection of auxiliary data $\bfa$, called an \emph{auxiliary collection}. If the collection satisfies some admissibility and nondegeneracy condition, the Floer complexes and family cobordism maps determine an element of
\[
\Psi(\D)_{k,l}=\bigl(\Fun(\Delta^k,\D)^\simeq\bigr)_l\subseteq\operatorname{Hom}_{\sSet}(\Delta^l\times\Delta^k,\D).
\]
For the instanton construction, the target is
\[
\D=\D_2(\Mod(\Z)^\fin)=\Ndg_\bu\bigl(\Ch_2((\Mod(\Z)^\fin)_\proj)\bigr),
\]
as in Definition \ref{defn: derived}. Thus its objects are $2$-periodic chain complexes of finitely presented projective $\Z$-modules.

Unlike Definition \ref{defn: PBIV kl}, we define auxiliary collections and their simplicial and spatial maps directly. The admissibility conditions are given in \S\ref{subsec: perturbation}.

Set $P_{l,k}=[l]\times[k]$ with the product partial order. Its nerve is $\rN_\bu(P_{l,k})=\Delta^l\times\Delta^k$. We call the elements of $P_{l,k}$ vertices and write $i=(z,w)$, where $z\in[0,l]$ and $w\in[0,k]$. Let $\cI_{l,k}$ be the set of finite \emph{weak chains} in $P_{l,k}$ with at least two entries, written as \begin{equation}\label{eq: notation}
I=\{i_-\leq i_p\leq\cdots\leq i_1\leq i_+\}\text{ or }(i_-,i_p,\dots,i_1,i_+),~p\geq0.
\end{equation}The number $|I|=p+2$ counts entries, including repetitions. 

By Lemma \ref{lem: dg nerve of a poset}, a simplicial map $\Delta^l\times\Delta^k\to\D$ is specified by objects $X_i=X_{(z,w)}$ and elements $f_I$ of degree $|I|-2$ for strict chains $I$, satisfying the differential identities in that lemma. These identities involve deleting an intermediate entry and composing at that entry. Following Construction \ref{cons: symm on D}, we assign a cube $Q_I$ of dimension $|I|-2$ to every weak chain, with zero faces corresponding to deletion and infinite faces to composition. Keeping repetitions makes the auxiliary collections closed under simplicial and spatial pullbacks.

For a simplicial map, the components on weak chains are determined by degeneracies, as in \eqref{eq: degeneracy condition}. The instanton construction need not satisfy these strict identities. We therefore first construct semisimplicial maps and use \cite{tanaka2018functors} to recover simplicial maps up to coherent equivalence. For repeated edges, we require only that $f_{(i,i)}$ be an isomorphism in $\D$, rather than the identity; see Proposition \ref{prop: abstract auxiliary theory}.

Let $e_z\in|\Delta^l|$ be the $z$-th standard vertex. For $\sigma\in(\PBI_a^{V,\aff})_{k,l}$, write
\begin{equation}\label{eq: sigma}
    \sigma(e_z)=(M_{e_z},\varphi_{e_z},\Phi_{e_z},(I_0(e_z)\leq\cdots\leq I_k(e_z))).
\end{equation}
Following Definition \ref{defn: PBIV k0}, the decorated $3$-manifold $\bY_{(z,w)}$ associated with $(z,w)$ has underlying manifold and subgroup
\[
Y_{(z,w)}=\pi_{e_z}^{-1}(m(I_w(e_z))),~\phi_{(z,w)}=\operatorname{im}\bigl(\Phi_{e_z}\to H^1(Y_{(z,w)};\Z/2)\bigr),
\]
where $\pi_{e_z}$ is the projection to the interval coordinate and $m$ denotes the midpoint of an interval, as in \eqref{eq: middle}. Its bundle, cylindrical almost complex structure, and metric $g_{(z,w)}$ are inherited from the $\theta_a$-structure $\varphi_{e_z}$. We also write $\bY_i=\bY_{(z,w)}$ when $i=(z,w)$.

\bdefn\label{defn: cubical auxiliary data}
Let $\sigma\in(\PBI_a^{V,\aff})_{k,l}$. For $I\in\cI_{l,k}$ as in \eqref{eq: notation}, set
\[
Q_I=Q_{p+1}=[0,\infty]^p,~p=|I|-2.
\]
The coordinate $q_j$ corresponds to the entry indexed by $j$, including when entries coincide. Write $I\backslash\{j\}$ for deletion of that entry. Equip $[0,\infty]$ with the smooth structure given by $x=q/(1+q)$ and orient $Q_I$ by $dx_1\wedge\cdots\wedge dx_p$. In particular, $Q_I=Q_1=*$ when $|I|=2$. Set
\[
I_j^-=(i_-,i_p,\dots,i_j),~I_j^+=(i_j,i_{j-1},\dots,i_+).
\]
The coordinates $q_j$ give canonical identifications
\begin{equation}\label{eq: abstract elementary faces}
\{q_j=0\}\cong Q_{I\backslash\{j\}},~\{q_j=\infty\}\cong Q_{I_j^+}\times Q_{I_j^-}.
\end{equation}
More generally, for disjoint $A,B\subseteq\{1,\ldots,p\}$, let
\[
F_{A,B}=\{q\in Q_I\mid q_j=0\text{ for }j\in A,~q_j=\infty\text{ for }j\in B\}.
\]
Delete the entries indexed by $A$ and split the remaining chain at the entries indexed by $B$. Write $I_0,\ldots,I_b$, where $b=|B|$, for the resulting consecutive chains, ordered from $i_-$ to $i_+$. Thus $I_0$ begins at $i_-$, $I_b$ ends at $i_+$, and the last entry of $I_t$ is the first entry of $I_{t+1}$. Then
\begin{equation}\label{eq: abstract general face}
F_{A,B}\cong Q_{I_b}\times\cdots\times Q_{I_0}.
\end{equation}
The factors are ordered from the outgoing end to the incoming end. These identifications commute with further deletion and splitting.

An \emph{auxiliary collection} over $\sigma\in(\PBI_a^{V,\aff})_{k,l}$ is a tuple
\[
\bfa=\bigl((a_i)_{i\in P_{l,k}},(\bfa_I)_{I\in\cI_{l,k}}\bigr)
\]
with the following properties.
\benum
\item For each $i=(z,w)\in P_{l,k}$, the datum $a_i$ consists of auxiliary data on $\bY_{(z,w)}$. In the instanton construction, write $a_i=(g_{(z,w)},\pi_{(z,w)})$, where $g_{(z,w)}$ is the metric inherited from $\sigma$ and $\pi_{(z,w)}$ is a primary perturbation.
\item For each $I\in\cI_{l,k}$, the datum $\bfa_I$ consists of a $Q_I$-family of decorated cobordisms $\bW_I$ from $\bY_{i_-}$ to $\bY_{i_+}$, together with auxiliary data having fixed end data $a_{i_-}$ and $a_{i_+}$. These families are compatible with the cobordism and spatial data of $\sigma$ along $I$. For a constant chain with value $i$, the underlying family is the specified ordered family of product cylinders over $\bY_i$, with its sections and end identifications; for $I=(i,i)$, it is $\R\times\bY_i$. Its auxiliary data satisfy the same end and face conditions as for other chains. We also use $\bfa_I$ to denote the auxiliary data alone. In the instanton construction, these include a smoothing function $T_I$ determining the metric family $\bfg_I(T_I)$ and a secondary perturbation family $\Pi_I$. The primary perturbations at the two ends are fixed by $a_{i_-}$ and $a_{i_+}$, and both primary and secondary perturbations are used there as in \eqref{eq: perturbed ASD 2}, including for constant chains.
\item On every face \eqref{eq: abstract general face}, the data satisfy
\begin{equation}\label{eq: abstract auxiliary face}
\bfa_I|_{F_{A,B}}=\bigl(\pr_b^*\bfa_{I_b},\ldots,\pr_0^*\bfa_{I_0}\bigr).
\end{equation}
For $B=\varnothing$, the right-hand side is a single family. Otherwise, it is a broken family whose adjacent pieces have matching data at their common end and depend on independent parameters. All identifications agree on intersections of faces.
\eenu
Let $\mathcal A(\sigma)$ denote the set of auxiliary collections over $\sigma$.
\edefn

We next assemble the sets $\mathcal A(\sigma)$ into a simplicial space by defining the simplicial and spatial maps.

\bdefn\label{defn: abstract auxiliary model}
For $k,l\geq0$, set
\begin{equation}\label{eq: abstract auxiliary total space}
(\PBI_a^{V,\mathrm{col}})_{k,l}=\coprod_{\sigma\in(\PBI_a^{V,\aff})_{k,l}}\{\sigma\}\times\mathcal A(\sigma).
\end{equation}

An \emph{auxiliary-data model} consists of these sets together with maps
\begin{equation}\label{eq: be al}
\begin{aligned}
\beta^*:\mathcal A(\sigma)&\longrightarrow\mathcal A(\beta^*\sigma),&&\beta:[k']\longrightarrow[k],\\
\alpha^*:\mathcal A(\sigma)&\longrightarrow\mathcal A(\alpha^*\sigma),&&\alpha:[l']\longrightarrow[l],
\end{aligned}
\end{equation}
for all monotone maps, called the \emph{simplicial} and \emph{spatial maps}, respectively. We require the identity and composition laws, and require the simplicial maps to commute with the spatial maps.

For $\rho=\alpha\times\beta:P_{l',k'}\to P_{l,k}$, write $\rho^*\bfa=\alpha^*\beta^*\bfa$. For $J\in\cI_{l',k'}$, let $\rho J$ denote its ordered image chain, retaining repetitions. Identify $Q_J$ with $Q_{\rho J}$ by matching coordinates with the same index. Require
\begin{equation}\label{eq: weak chain pullback}
(\rho^*\bfa)_j=a_{\rho(j)},~(\rho^*\bfa)_J=\bfa_{\rho J}.
\end{equation}
The specified identifications preserve the end data, commute with deletion and splitting, and satisfy the identity and composition laws.
These pullbacks make \eqref{eq: abstract auxiliary total space} a simplicial space, with a forgetful morphism
\[
p_a^V:\PBI_a^{V,\mathrm{col}}\longrightarrow\PBI_a^{V,\aff}.
\]

Assume that auxiliary collections admit disjoint unions over common weak chain parameters, compatible with \eqref{eq: weak chain pullback} and \eqref{eq: abstract auxiliary face}. Let $\sigma\in(\PBI_a^{V,\aff}\langle m\rangle)_{k,l}$, with labelled parts $\sigma^1,\ldots,\sigma^m$ as in Definition \ref{defn: symmetric monoidal str}. A lift
\[
(\sigma,(\bfa^\alpha)_{\alpha=1}^m)\in(\PBI_a^{V,\col}\langle m\rangle)_{k,l}
\]
consists of an auxiliary collection $\bfa^\alpha\in\mathcal A(\sigma^\alpha)$ for each $\alpha=1,\ldots,m$, as in Definition \ref{defn: cubical auxiliary data}. A pointed map $\theta:\langle m\rangle\to\langle n\rangle$ acts by disjoint union over its nonzero fibers,
\begin{equation}\label{eq: common parameter}
a_i^\be=\bigsqcup_{\theta(\alpha)=\be}a_i^\alpha,~
\bfa_I^\be(q)=\bigsqcup_{\theta(\alpha)=\be}\bfa_I^\alpha(q),~1\leq\be\leq n.
\end{equation}
An empty fiber gives the empty collection. These operations define the $\Ga$-diagram of auxiliary collections.
\edefn
Let $\PBI_a^{\partial,V,\mathrm{col}}$ denote the simplicial subspace satisfying the full boundary condition of Definition \ref{defn: PBIV k0}. In \S\ref{subsec: perturbation}, we impose admissibility and prove that
\[
p_a^{\partial,V,\adm}:\PBI_a^{\partial,V,\adm}\longrightarrow\PBI_a^{\partial,V,\aff}
\]
is a levelwise trivial Kan fibration. Together with \eqref{eq: Fg V}, this allows stabilization and completion as in Remark \ref{rem: affine construction}, also for the $\partial_\pm$ variants. The following proposition gives the required properties of the family maps without imposing strict compatibility with degeneracies.

\bprop\label{prop: abstract auxiliary theory}
Let $\sfC$ be a dg category and set $\D=\Ndg_\bu(\sfC)$. The following assignments yield a morphism of simplicial spaces
\[
X:\PBI_a^{\partial,V,\adm}\longrightarrow\Psi(\D)
\]
through a contractible space of choices. The construction is natural, up to coherent homotopy, under morphisms of simplicial spaces preserving the assignments.
\benum
\item To each vertex $(\bY,a)$, assign an object $X(\bY,a)\in\sfC$. For each oriented family $(\bW,\bfa)$ over a parameter manifold $G$ occurring in an admissible collection, its faces, or its broken products, assign
\[
m_G(\bfa)\in\Map_{\sfC}(X_-,X_+)_n,~n=\dim G,~X_\pm=X(\bY_\pm,a_\pm).
\]
These assignments are invariant under the specified orientation-preserving identifications, including the simplicial and spatial pullbacks, change sign on reversing the orientation of $G$, and satisfy
\begin{equation}\label{eq: abstract family identities}
\begin{aligned}
d\,m_G(\bfa)&=m_{\partial G}(\bfa|_{\partial G}),\\
m_{G_+\times G_-}(\bfa_+,\bfa_-)&=m_{G_+}(\bfa_+)\circ m_{G_-}(\bfa_-).
\end{aligned}
\end{equation}
Here $m_{\partial G}$ is the sum over codimension-one faces with their boundary orientations. Broken products have independent parameters and matching end data, with factors ordered from the outgoing end to the incoming end. For chain complexes, the mapping-complex differential is
\begin{equation}\label{eq: d m G}
d\,m_G=\partial_+\circ m_G-(-1)^n m_G\circ\partial_-.
\end{equation}
Only the $Q_I$-families in Definition \ref{defn: cubical auxiliary data} and their faces and broken products are needed; the notation $G$ is used to state the identities uniformly.
\item For $(\sigma,\bfa)\in(\PBI_a^{\partial,V,\adm})_{k,l}$, put
\[
X_i=X(\bY_i,a_i),~f_I=m_{Q_I}(\bfa_I),~I\in\cI_{l,k}.
\]
Require $f_I$ to be an isomorphism in $\D$ whenever $I=((z,w),(z',w))$ with $z\leq z'$. No strict identity or vanishing condition as in \eqref{eq: degeneracy condition} is imposed on repeated chains.
\eenu
 
\eprop
\bpf
We use Tanaka's construction \cite{tanaka2018functors} to restore compatibility with degeneracy maps. Write $K_{\mathrm{inj}}$ for the restriction of a simplicial set $K$ to $\Delta_{\mathrm{inj}}^{\opp}$. Freely adjoining degeneracies gives a simplicial set $K^+$ and an augmentation
\[
(K^+)_n=\coprod_{s:[n]\twoheadrightarrow[m]}K_m,~
\epsilon:K^+\longrightarrow K,~\epsilon(s,u)=s^*u.
\]

Put $\mathcal X=\PBI_a^{\partial,V,\adm}$. For every weak chain $I$, the cube orientation and \eqref{eq: abstract auxiliary face} give
\[
df_I=\sum_{j=1}^p(-1)^j\bigl(f_{I\backslash\{j\}}-f_{I_j^+}\circ f_{I_j^-}\bigr).
\]
For $x=(\sigma,\bfa)\in\mathcal X_{k,l}$, these identities on all subchains, keeping their positions, define a semisimplicial map $F_x:K_{\mathrm{inj}}\to\D_{\mathrm{inj}}$, where $K=\Delta^l\times\Delta^k$, and hence a simplicial map $F_x^+:K^+\to\D$.

Write
\[
\Psi^+(\D)_{k,l}=\{F:K^+\to\D\mid F(e)\in\D^\simeq_1\},
\]
where $e$ ranges over the copies in $K^+$ of the edges $(z,w)\to(z',w)$ of $K$, with $z\leq z'$. Pullback defines the simplicial operators. By item (2) and \eqref{eq: weak chain pullback}, the maps $F_x^+$ assemble into $X^+:\mathcal X\to\Psi^+(\D)$. Precomposition with $\epsilon$ gives
\begin{equation}\label{eq: auxiliary rectification}
\epsilon^*:\Psi(\D)\longrightarrow\Psi^+(\D).
\end{equation}

We verify that $\epsilon^*$ is a levelwise weak homotopy equivalence as follows. Fix $k$ and regard $\Delta^\bu\times\Delta^k$ and $(\Delta^\bu\times\Delta^k)^+$ as cosimplicial objects in the model structure on marked simplicial sets over a point \cite[Proposition 3.1.3.7]{HTTlurie2009}. Mark the edges in $\Delta^l\times\{w\}$ for every vertex $w$ of $\Delta^k$, and their copies in $(\Delta^l\times\Delta^k)^+$. Mapping these marked objects into $\D$, marked by its equivalences, gives precisely $\Psi(\D)_{k,\bu}$ and $\Psi^+(\D)_{k,\bu}$, respectively.

By \cite[Theorem 1.1]{tanaka2018functors}, followed by localization at the remaining marked edges, $\epsilon$ is a weak equivalence in each degree. Projection to $\Delta^k$, with only degenerate edges marked, is also a weak equivalence, since inverting all edges of $\Delta^l$ makes it equivalent to a point. The boundary inclusions in the $\Delta^l$-direction, before and after applying $(-)^+$, are monomorphisms. Thus $\Delta^\bu\times\Delta^k$ and $(\Delta^\bu\times\Delta^k)^+$, with the above markings, are cosimplicial resolutions of $\Delta^k$ in the sense of \cite[Definition 16.1.2(1)]{Hirschhorn2003model}. Since $\D$ with these markings is fibrant, \cite[Corollary 16.5.5(1)]{Hirschhorn2003model} shows that $\epsilon^*$ is a weak homotopy equivalence for each $k$.

Let $j:\Psi^+(\D)\to\what{R}(\Psi^+(\D))$ be the fibrant replacement in $\sSpace_c^{\CSe}$. Since $\Psi(\D)$ is fibrant, the weak equivalence $j\epsilon^*$ gives the required $X:\mathcal X\to\Psi(\D)$, together with $j\epsilon^*X\simeq jX^+$, through a contractible space of choices. Naturality gives compatibility under morphisms preserving the assigned data, up to coherent homotopy.
\epf

To extend Proposition \ref{prop: abstract auxiliary theory} to the symmetric monoidal setting, we compare common-parameter disjoint unions with the independent-parameter tensor operations of Construction \ref{cons: symm on D}. We use the enlarged charts of Definition \ref{defn: diagrams with strict duality} so that the comparison data are closed under reflection.

\bcons[Comparison families]\label{cons: comparison families}
Fix an auxiliary collection $x=(\sigma,(\bfa^\alpha)_{\alpha=1}^m)\in(\PBI_a^{V,\col}\langle m\rangle)_{k,l}$. Consider a \emph{mixed string} of $r\geq1$ arrows in $\Ga\times[1]$
\[
\rho=((\langle n_0\rangle,\epsilon_0)\xra{\Theta_1}\cdots\xra{\Theta_r}(\langle n_r\rangle,\epsilon_r)),~n_0=m,~0=\epsilon_0\leq\cdots\leq\epsilon_r=1.
\]
Here levels refer to the $[1]$-coordinate. An arrow or block is called \emph{level-$0$} or \emph{level-$1$} if both endpoints have coordinate $0$ or $1$, respectively, and \emph{mixed} if its endpoint coordinates are $0$ and $1$. A cut after $\Theta_i$ is at level $\epsilon_i$. Use the chart $C_\rho^{\sd}$ of Definition \ref{defn: diagrams with strict duality}, whose simplices are arbitrary cut-set lists $(C_0,\ldots,C_u)$.

For a pointed map $\theta:\langle m\rangle\to\langle n\rangle$ and $\beta\ne0$ in $\langle n\rangle$, write $S_{\theta,\beta}=\theta^{-1}(\beta)=\{\alpha_1<\cdots<\alpha_s\}$. On a nonempty fiber, define the piecewise-smooth \emph{staircase path} $h^\otimes:[0,\infty]\to[0,\infty]^{S_{\theta,\beta}}$, with coordinates in the displayed order, by
\begin{equation}\label{eq: standard comparison path}
h^\otimes(q)=\begin{cases}
\bigl(\underbrace{0,\ldots,0}_{s-j-1},\dfrac{q-j}{j+1-q},\underbrace{\infty,\ldots,\infty}_{j}\bigr),&j\leq q<j+1,~0\leq j<s-1,\\[3pt]
(q-s+1,\infty,\ldots,\infty),&s-1\leq q\leq\infty.
\end{cases}
\end{equation}
It increases one coordinate at a time, from $\alpha_s$ to $\alpha_1$; for $s=1$ it is the identity. An empty fiber gives the empty family. For $\Theta=(\theta,\epsilon\leq\epsilon')$ and $\theta(\alpha)\ne0$, let $h_\Theta^\alpha$ be the $\alpha$-coordinate of this staircase on $S_{\theta,\theta(\alpha)}$ when $\epsilon=\epsilon'=1$, and set $h_\Theta^\alpha(q)=q$ otherwise. Thus level-$1$ arrows use staircases, and all other arrows use diagonals.

Write $\Theta_{a,b}=\Theta_b\cdots\Theta_{a+1}$ for $0\leq a\leq b\leq r$, and let $\theta_{a,b}$ be its underlying pointed map, with empty composites equal to identities. Fix $\beta\ne0$ in $\langle n_r\rangle$ and put $S_{\rho,\beta}=\theta_{0,r}^{-1}(\beta)$. For a cut set $C=\{a_1<\cdots<a_{t-1}\}\subseteq\{1,\ldots,r-1\}$, put $a_0=0$, $a_t=r$, and $\gamma_j=\theta_{0,a_j}(\alpha)$ for $\alpha\in S_{\rho,\beta}$. Define
\begin{equation}\label{eq: comparison stage paths}
h_{\rho,C}^\alpha=h_{\Theta_{a_0,a_1}}^{\gamma_0}\circ\cdots\circ h_{\Theta_{a_{t-1},a_t}}^{\gamma_{t-1}},~\alpha\in S_{\rho,\beta}.
\end{equation}
Level-$1$ blocks use the label orders and parenthesizations of Construction \ref{cons: symm on D}. See Figure \ref{fig: comparison staircase} for examples.

For a $u$-simplex $(C_0,\ldots,C_u)$ of $C_\rho^{\sd}$, put $h_v=(h_{\rho,C_v}^\alpha)_{\alpha\in S_{\rho,\beta}}$. For a weak chain
\[
J=((c_-,i_-),(c_p,i_p),\ldots,(c_1,i_1),(c_+,i_+))\text{ in }[u]\times([l]\times[k]),
\]
let $I$ be its projection to $[l]\times[k]$, retaining all repetitions. Set $Q_J=Q_I=[0,\infty]^p$, with coordinates identified by position and the deletion, splitting, and orientation conventions of Definition \ref{defn: cubical auxiliary data}. Define the family by
\begin{equation}\label{eq: comparison family pullback}
q_j^\alpha=h_{c_j}^\alpha(q_j),~
\bfa_J(q)=\bigsqcup_{\alpha\in S_{\rho,\beta}}\bfa_I^\alpha(q_1^\alpha,\ldots,q_p^\alpha).
\end{equation}
Subdivide the $q_j$-interval at the breakpoints of $h_{c_j}^\alpha$; the resulting products of intervals are called \emph{cells}. On each cell the family is smooth in the compactified coordinates, with the prescribed broken geometry at infinite coordinates. Since all paths fix both endpoints, these families satisfy
\begin{equation}\label{eq: comparison chain boundary}
\bfa_J|_{q_j=0}=\bfa_{J\backslash\{j\}},~
\bfa_J|_{q_j=\infty}=(\pr_+^*\bfa_{J_j^+},\pr_-^*\bfa_{J_j^-}).
\end{equation}

\begin{figure}[!htbp]
\centering
{\small $\rho=\bigl((\langle3\rangle,0)\xra{(\id,0\to1)}(\langle3\rangle,1)\xra{(\theta_1,1\to1)}(\langle2\rangle,1)\xra{(\theta_2,1\to1)}(\langle1\rangle,1)\bigr).$\par}
\smallskip
\begin{tikzpicture}[font=\small,>=Stealth,line cap=round,line join=round]
\definecolor{modelblue}{RGB}{42,88,119}
\definecolor{modelrust}{RGB}{153,89,64}
\definecolor{modelguide}{RGB}{186,190,194}
\begin{scope}[x=2.5cm,y=2.5cm]
\node at (.5,1.42) {(a) Two components};
\draw[modelguide,line width=.45pt] (0,0) rectangle (1,1);
\draw[line width=.45pt,->] (0,0)--(1.16,0) node[right] {$q_j^1$};
\draw[line width=.45pt,->] (0,0)--(0,1.16) node[left] {$q_j^2$};
\node[below left] at (0,0) {$0$};
\node[below] at (1,0) {$\infty$};
\node[left] at (0,1) {$\infty$};
\draw[modelblue,line width=1.05pt,densely dashed] (0,0)--(1,1);
\draw[modelblue,line width=1.05pt,->] (.36,.36)--(.55,.55);
\node[modelblue,above=3pt] at (.64,.64) {$h_0$};
\draw[modelrust,line width=1.05pt] (0,0)--(0,1)--(1,1);
\draw[modelrust,line width=1.05pt,->] (0,.34)--(0,.56);
\draw[modelrust,line width=1.05pt,->] (.38,1)--(.62,1);
\node[modelrust,above=3pt] at (.64,1) {$h_1=h^\otimes$};
\node[modelrust,anchor=west] at (.07,.88) {$q_j=1$};
\foreach \u/\v in {0/0,0/1,1/1} {\fill (\u,\v) circle[radius=1.2pt];}
\end{scope}
\begin{scope}[shift={(6cm,0)},x=2.5cm,y=2.5cm]
\node at (.52,1.42) {(b) Three components};
\coordinate (a) at (0,0); \coordinate (b) at (.73,0);
\coordinate (c) at (.30,.24); \coordinate (d) at (1.03,.24);
\coordinate (e) at (0,.77); \coordinate (f) at (.73,.77);
\coordinate (g) at (.30,1.01); \coordinate (h) at (1.03,1.01);
\draw[modelguide,line width=.45pt] (a)--(b)--(d)--(c)--cycle (e)--(f)--(h)--(g)--cycle (a)--(e) (b)--(f) (c)--(g) (d)--(h);
\draw[modelblue,line width=1.05pt,densely dashed] (a)--(h);
\draw[modelblue,line width=1.05pt,->] (.36,.353)--(.55,.539);
\node[modelblue,below right=2pt] at (.44,.53) {$h_0$};
\draw[modelguide!40!black,line width=1.05pt,densely dotted,->] (e)--(h);
\node[modelguide!40!black,below=3pt] at (.55,.90) {$h_1$};
\draw[modelrust,line width=1.05pt,->] (a)--(e);
\draw[modelrust,line width=1.05pt,->] (e)--(g);
\draw[modelrust,line width=1.05pt,->] (g)--(h);
\node[modelrust,above=3pt] at (.68,1.01) {$h_2=h^\otimes$};
\foreach \u in {a,e,g,h} {\fill (\u) circle[radius=1.2pt];}
\node[below left] at (a) {$000$};
\node[left] at (e) {$001$};
\node[above left] at (g) {$011$};
\node[right] at (h) {$111$};
\node[below=3pt] at (.50,0) {$q_j^1$};
\node[right=3pt] at (.98,.16) {$q_j^2$};
\node[left=3pt] at (0,.35) {$q_j^3$};
\end{scope}
\begin{scope}[shift={(0,-4.65cm)},x=2.5cm,y=2.5cm]
\node at (.52,1.48) {(c) A chart simplex};
\coordinate (v0) at (.03,.12); \coordinate (v1) at (.54,.89); \coordinate (v2) at (1.05,.12);
\fill[modelblue!6] (v0)--(v1)--(v2)--cycle;
\draw[modelblue,line width=.7pt,->] (v0)--(v1) node[midway,left] {$d_2$};
\draw[modelblue,line width=.7pt,->] (v1)--(v2) node[midway,right] {$d_0$};
\draw[modelblue,line width=.7pt,->] (v0)--(v2) node[midway,below] {$d_1$};
\node at (.54,.40) {$\Delta^2$};
\foreach \u in {v0,v1,v2} {\fill (\u) circle[radius=1.2pt];}
\node[below left,align=center] at (v0) {$C_0=\varnothing$\\$h_0$};
\node[above,align=center] at (v1) {$C_1=\{2\}$\\$h_1$};
\node[below right,align=center] at (v2) {$C_2=\{1,2\}$\\$h_2$};
\node at (.52,-.40) {$\Delta^2\longrightarrow C_\rho=(\Delta^1)^2$};
\end{scope}
\begin{scope}[shift={(6cm,-4.65cm)},x=2.5cm,y=2.5cm]
\node at (.5,1.42) {(d) Cells in $Q_J=[0,\infty]^2$};
\fill[modelrust!12] (0,.5) rectangle (.5,1);
\draw[modelguide,line width=.45pt] (0,0) rectangle (1,1);
\draw[modelguide,line width=.45pt] (.5,0)--(.5,1) (0,.5)--(1,.5);
\draw[line width=.45pt,->] (0,0)--(1.16,0) node[right] {$q_1$};
\draw[line width=.45pt,->] (0,0)--(0,1.16) node[left] {$q_2$};
\node[modelrust] at (.25,.75) {$G$};
\node[below left] at (0,0) {$0$};
\node[below] at (.5,0) {$1$};
\node[below] at (1,0) {$\infty$};
\node[left] at (0,.5) {$1$};
\node[left] at (0,1) {$\infty$};
\node at (.5,-.34) {$G=[0,1]\times[1,\infty]$};
\end{scope}
\end{tikzpicture}
\caption{(a) The diagonal $h_0$ and staircase $h_1$. (b) and (c) For the displayed string, $\theta_1(1)=\theta_1(2)=1$, $\theta_1(3)=2$, and $\theta_2(1)=\theta_2(2)=1$. The paths give the vertices of the chart simplex in subfigure (c), which lies in $C_\rho$, not $Q_J$; $h_1$ and $h_2$ share the first edge. The same simplex and its reversed list belong to $C_\rho^{\sd}$. (d) Using subfigure (a) with $p=2$ and $c_1=c_2=1$, on $G$ the parameters $q_1,q_2$ vary components $2,1$, respectively. Axes use $q/(1+q)$; the digits in subfigure (b) denote the endpoints $0,\infty$.}
\label{fig: comparison staircase}
\end{figure}

A \emph{comparison collection} $\wti{\bfa}$ consists of full auxiliary data $\wti{\bfa}_{J,G}$ on these families, for all mixed strings, chart simplices, and cells $G\subset Q_J$, satisfying the following conditions compatibly on intersections.
\benum
\item The data have the prescribed ends, agree on common cell faces, and satisfy \eqref{eq: comparison chain boundary} with $\bfa$ replaced by $\wti{\bfa}$. They commute with simplicial and spatial pullbacks and the chart identifications of Definition \ref{defn: diagrams with strict duality}. For a single mixed arrow, they equal the original data.
\item On $t_i=0$, use the data for $d_i\rho$; on $t_i=1$, compose the data for the two strings cut at $i$. Level-$0$ cuts use common-parameter unions, while level-$1$ cuts use independent products of the prescribed cells $G_\alpha$ of the shorter mixed families,
\begin{equation}\label{eq: comparison auxiliary prescriptions}
\wti{\bfa}_{\prod_\alpha G_\alpha}=\bigsqcup_\alpha\pr_\alpha^*\wti{\bfa}_{G_\alpha}^\alpha.
\end{equation}
Substitution pulls back whole factor cells without further subdividing fixed data. Consecutive initial level-$0$ arrows use the pullback from their composite.
\item An interval factor of $G$ is \emph{unused} if it is nondegenerate and every component parameter in \eqref{eq: comparison family pullback} is independent of it. On such a factor, the full data are pulled back from the remaining factors.
\eenu

Let $(\PBI_a^{V,\com}\langle m\rangle)_{k,l}$ consist of the pairs $(x,\wti{\bfa})$ above, with simplicial and spatial maps induced by pullback. For $\theta:\langle m\rangle\to\langle m'\rangle$, act on $x$ as in Definition \ref{defn: abstract auxiliary model}. If $\rho$ starts at $(\langle m'\rangle,0)$, prepend $(\theta,\id_0)$ and restrict to the face $t_1=1$, identified with $C_\rho^{\sd}$. Thus
\begin{equation}\label{eq: comparison source restriction}
\theta_*(x,\wti{\bfa})=(\theta_*x,\theta_*\wti{\bfa}),~(\theta_*\wti{\bfa})_\rho=\wti{\bfa}_{((\theta,\id_0),\rho)}|_{t_1=1}.
\end{equation}
The initial level-$0$ and chart conditions give the composition and identity laws and compatibility with pullbacks. Thus these operations define a $\Ga$-diagram with a forgetful morphism to $\PBI_a^{V,\col}$.
\econs

Impose the full boundary condition of Definition \ref{defn: PBIV k0} and admissibility for the original and comparison families as in \S\ref{subsec: perturbation}. The resulting $\Ga$-diagram $A$ has underlying simplicial space $\PBI_a^{\partial,V,\cadm}$ and strict duality from Construction \ref{cons: reflected model families}. The same applies to the $\partial_\pm$ variants.

\bprop\label{prop: symmetric monoidal auxiliary theory}
Let $A$ be the $\Ga$-diagram above, let $\sfC$ be a symmetric monoidal dg category, and use $\D$ and $\fB$ from Construction \ref{cons: symm on D}. Suppose the assignments $X$ and $m_G$ satisfy the hypotheses of Proposition \ref{prop: abstract auxiliary theory} and the following conditions.
\benum
\item The family maps extend to all comparison cells, their faces, and prescribed independent or broken products, satisfying Proposition \ref{prop: abstract auxiliary theory}(1). For every prescribed subdivision $G=\bigcup_\nu G_\nu$ with induced orientations, require
\[
m_G(\wti{\bfa})=\sum_{\dim G_\nu=\dim G}m_{G_\nu}(\wti{\bfa}|_{G_\nu}).
\]
\item There are associative, unital, and symmetric identifications
\[
\iota:X\Bigl(\bigsqcup_\alpha(\bY^\alpha,a^\alpha)\Bigr)\xra{\cong}\bigotimes_\alpha X(\bY^\alpha,a^\alpha)
\]
satisfying, for independent parameter spaces with the displayed product orientation,
\begin{equation}\label{eq: independent parameter tensor}
\iota_+\circ m_{\prod_\alpha G_\alpha}\Bigl(\bigsqcup_\alpha\pr_\alpha^*\bfa^\alpha\Bigr)=\Bigl(\bigotimes_\alpha m_{G_\alpha}(\bfa^\alpha)\Bigr)\circ\iota_-.
\end{equation}
\item The family maps vanish on pullbacks along coordinate projections $\pr:G\to G'$ omitting one interval factor, with $\dim G=\dim G'+1$.
\item For every two-entry chain $J$ with constant $[k]$-coordinate, require
\[
m_{Q_J}(\wti{\bfa}_{J,Q_J}):X_-\xra{\simeq}X_+\text{ in }\D.
\]
\eenu
Then the maps induced by Proposition \ref{prop: abstract auxiliary theory} at each $\langle m\rangle\in\Ga$, pulled back to $A$, extend to a homotopy coherent natural transformation
\[
H:A\epsilon_\Ga\Longrightarrow\Psi(\fB),
\]
where $\epsilon_\Ga$ is the augmentation in \eqref{eq: augmentation}.

Suppose moreover that $\sfC$ has symmetric monoidal duality and $(X,m_G,\iota)$ lift to $\sfC^{\sd}$, denoted by the same symbols, with $\iota$ compatible with duality and
\begin{equation}\label{eq: comparison family duality}
\begin{aligned}
X(v^\vee)&=X(v)^{\vee_{\sd}},\\
m_{G^\vee}(\bfa^\vee)&=(-1)^{p(p-1)/2}m_G(\bfa)^{\vee_{\sd}},~p=\dim G,
\end{aligned}
\end{equation}
for every vertex $v$ and every original cube or top-dimensional comparison cell $G$, with coordinate orientations and the reflection of Construction \ref{cons: reflected model families}. Then the morphism represented by $H$ lifts to a map
\[
\widehat R^{\sd}A\longrightarrow\widehat R^{\sd}B^{\sd}
\]
preserving strict duality, where $B^{\sd}$ is from Construction \ref{cons: strict tensor diagram} and $\widehat R^{\sd}$ is the replacement of Construction \ref{cons: special dual}. Both constructions are coherently natural under morphisms preserving the assigned data and, when present, duality.
\eprop
\bpf
Use $\Psi^+$ and $\epsilon^*$ from Proposition \ref{prop: abstract auxiliary theory}. For each chart simplex $C_0\subseteq\cdots\subseteq C_u$, apply the block operations at each $C_v$ from left to right, using disjoint unions on level-$0$ blocks, disjoint unions followed by $X$ on the unique mixed block, and tensor products on level-$1$ blocks, over the corresponding fibers. Let $\iota_{C_v,i}$ be the identification from $X$ of the total disjoint union to the resulting object, given by item (2) and symmetric monoidal coherence. Set
\begin{equation}\label{eq: comparison map component}
(H_\rho^+)^\beta_J=\iota_{C_{c_+},i_+}\circ\Bigl(\sum_{G\subset Q_J}m_G(\wti{\bfa}_{J,G})\Bigr)\circ\iota_{C_{c_-},i_-}^{-1},
\end{equation}
where the sum is over top-dimensional cells with their subdivision orientations. Common cell faces cancel, so items (1) and (4) allow Proposition \ref{prop: abstract auxiliary theory} to apply on each chart. The tensor formula on level-$1$ cuts follows from item (2). Zero and level-$0$ faces use Construction \ref{cons: comparison families}, and unused parameters vanish by item (3). Thus these components define $H^+:A\epsilon_\Ga\Longrightarrow\Psi^+(\fB)$.

By the proof of Proposition \ref{prop: abstract auxiliary theory}, $\epsilon^*:\Psi(\fB)\to\Psi^+(\fB)$ is an objectwise levelwise weak equivalence. After fibrant replacement of $\Psi^+(\fB)$, construct the chart maps
\[
H_\rho:C_\rho\otimes A\langle n_0\rangle\longrightarrow\Psi(\D^{n_r})
\]
and homotopies comparing $\epsilon^*H_\rho$ with $H_\rho^+$ inductively on string length, fixing the endpoint diagrams and the maps and homotopies at each $\langle m\rangle\in\Ga$. The face rules of Definition \ref{defn: homotopy coherent realization} give compatible boundary data. Since
\[
\partial C_\rho\otimes A\langle n_0\rangle\hookrightarrow C_\rho\otimes A\langle n_0\rangle
\]
is a levelwise monomorphism, restriction to the boundary gives a fibration of mapping spaces into fibrant targets. The weak equivalence $\epsilon^*$ and the maps $H_\rho^+$ therefore allow the same mapping-space argument to extend the maps and homotopies relative to the boundary. Charts containing identity arrows are determined by the identity-deletion rules. This gives $H$.

For duality, the same formula on arbitrary chart lists defines $H^{\sd,+}$ on $\fC^{\sd}[\Ga\times[1]]$, with endpoints $A\epsilon_\Ga^{\sd}$ and $\Psi^+(\fB^{\sd})$. Equations \eqref{eq: comparison family duality} and \eqref{eq: strict tensor reflection}, together with duality compatibility of the tensor identifications, make this diagram compatible with strict duality. Weak-chain reversal commutes with $(-)^+$ and $\epsilon$, so $\epsilon^*$ also preserves reflection.

Apply Lemma \ref{lem: rectification with strict duality} to $H^{\sd,+}$ over $\Ga\times[1]$. The zigzag defining $B^{\sd}$ and $\epsilon^*$ identify the two resulting endpoint diagrams with $A$ and $B^{\sd}$ up to weak equivalence. Applying $\widehat R^{\sd}$ and composing along these zigzags, using homotopy inverses where necessary as in that lemma, gives the asserted map. Restricting $H^{\sd,+}$ to the original charts and applying $P:\sfC^{\sd}\to\sfC$ from Construction \ref{cons: strict dual} recovers $H^+$, hence the morphism represented by $H$. Coherent naturality follows from Proposition \ref{prop: abstract auxiliary theory} and Lemma \ref{lem: rectification with strict duality}.
\epf
\subsection{Constructions with auxiliary data}\label{subsec: perturbation}
In this subsection, we construct auxiliary collections for instantons, including model metric families and perturbation families. We then impose admissibility conditions on the perturbations and prove that the resulting forgetful maps are levelwise trivial Kan fibrations.

Fix $\sigma\in(\PBI_a^{V,\aff})_{k,l}$ and use the notation of \S\ref{subsec: abstract auxiliary data}. For a weak chain $I$ in \eqref{eq: notation}, write $i_j=(z_j,w_j)$, with $i_{p+1}=i_-$ and $i_0=i_+$. We use the open ends and interval-first orientation of \eqref{eq: open end}. We construct model $Q_I$-families of metrics, generalizing \cite[\S 3.9]{kronheimer2011khovanov}, first for fixed $z$ or fixed $w$, then for general chains.

\bcons[Ordered cylinder model]\label{cons: ordered cylinder model}
For a chain $I$ with fixed $z$, work with $\sigma(e_z)$ and the midpoint sections $Y_w=\pi^{-1}(m(I_w))$. We first describe the insertion of a product cylinder $[0,q]\times Y$ using the collars of Definition \ref{defn: cylindrical}. On a collar with metric $\epsilon^2dr^2+g_Y$, fix an even function $\chi\in C_c^\infty((-1,1),[0,1])$ equal to $1$ near $0$. Following \cite[\S 3.9]{kronheimer2011khovanov}, set
\begin{equation}\label{eq: model broken metric}
f_u(r)=1+\chi(r)\left(\frac{1+u^2}{1+u^2r^2}-1\right),~\bfg_q=\epsilon^2f_u(r)dr^2+g_Y,
\end{equation}
where $u\geq0$ is determined by the added length
\begin{equation}\label{eq: inserted length}
q=\epsilon\int_{-1}^1\bigl(\sqrt{f_u(r)}-1\bigr)\,dr.
\end{equation}
The added length is smooth and strictly increasing in $u^2$, with range $[0,\infty)$. The \emph{arclength coordinates}
\[
\ell(r)=-\epsilon+\epsilon\int_{-1}^r\sqrt{f_u(v)}\,dv,~\ell^-=1+\epsilon+\ell,~\ell^+=\ell-\epsilon-q
\]
identify the finite collar with $[-\epsilon,\epsilon+q]\times Y$ and its two limiting ends with $[1,\infty)\times Y$ and $(-\infty,0]\times Y$, respectively. The metric is product in these coordinates and unchanged at $q=0$. Near $q=\infty$, the function $u^{-2}$ extends smoothly in $(1+q)^{-1}$ and vanishes to infinite order at $0$, giving smoothness at the broken collar. Transport the remaining tangential data by the same identifications.

Add open ends to the cobordism between the endpoint sections of $I$ as in \eqref{eq: open end}. Insert cylinders of lengths $q_p,\ldots,q_1$ at the sections indexed by $w_p,\ldots,w_1$, in this order. If the endpoint sections coincide, start with $\R\times Y$. At a repeated section, insert the next cylinder at the right endpoint of the preceding one, as in Figure \ref{fig: successive cylinder insertions}. We call this the \emph{ordered cylinder model}. Near a breaking face, use the arclength coordinate $\ell\in[0,q_j]$ on each inserted cylinder and the translated coordinate $\ell-q_j$ near its right endpoint. The limiting metrics are product on the cylindrical ends, and the single-cylinder calculation gives smoothness at all infinite faces.

\begin{figure}[htbp]
\centering
\resizebox{.93\linewidth}{!}{%
\begin{tikzpicture}[x=1cm,y=1cm,font=\small,>=Stealth,line cap=round,line join=round]
\definecolor{modelblue}{RGB}{42,88,119}
\tikzset{band/.style={draw=black!70,line width=.5pt},junction/.style={draw=modelblue,line width=.8pt}}
\fill[black!5] (0,-.17) rectangle (1.8,.17); \fill[modelblue!14] (1.8,-.17) rectangle (5,.17); \fill[modelblue!7] (5,-.17) rectangle (8.2,.17);
\fill[modelblue!14] (8.2,-.17) rectangle (11.4,.17); \fill[black!5] (11.4,-.17) rectangle (14,.17); \draw[band] (0,-.17)--(14,-.17) (0,.17)--(14,.17);
\foreach \a in {1.8,5,8.2,11.4} {\draw[junction] (\a,-.24)--(\a,.24);}
\node at (3.4,0) {$q_3$}; \node at (6.6,0) {$q_2$}; \node at (9.8,0) {$q_1$}; \node[below] at (1.8,-.25) {$0$}; \node[below] at (5,-.25) {$q_3$};
\node[below] at (8.2,-.25) {$q_3+q_2$}; \node[below] at (11.4,-.25) {$q_3+q_2+q_1$}; \draw[->,black!55] (7,-.76)--(7,-1.2);
\node[anchor=west] at (7.25,-.98) {$q_3,q_2,q_1\to\infty$}; \fill[black!5] (0,-1.77) rectangle (.95,-1.43); \fill[modelblue!14] (.95,-1.77) rectangle (2.7,-1.43);
\draw[band,->] (0,-1.77)--(2.7,-1.77); \draw[band,->] (0,-1.43)--(2.7,-1.43); \fill[modelblue!14] (3.5,-1.77) rectangle (4.85,-1.43);
\fill[modelblue!7] (4.85,-1.77) rectangle (6.2,-1.43); \draw[band,<->] (3.5,-1.77)--(6.2,-1.77); \draw[band,<->] (3.5,-1.43)--(6.2,-1.43);
\draw[junction] (4.85,-1.85)--(4.85,-1.35); \fill[modelblue!7] (7,-1.77) rectangle (8.35,-1.43); \fill[modelblue!14] (8.35,-1.77) rectangle (9.7,-1.43);
\draw[band,<->] (7,-1.77)--(9.7,-1.77); \draw[band,<->] (7,-1.43)--(9.7,-1.43); \draw[junction] (8.35,-1.85)--(8.35,-1.35);
\fill[modelblue!14] (10.5,-1.77) rectangle (12.85,-1.43); \fill[black!5] (12.85,-1.77) rectangle (14,-1.43); \draw[band,<-] (10.5,-1.77)--(14,-1.77);
\draw[band,<-] (10.5,-1.43)--(14,-1.43); \node[below] at (1.35,-1.94) {incoming piece}; \node[below] at (4.85,-1.94) {$\mathbb R\times Y$};
\node[below] at (8.35,-1.94) {$\mathbb R\times Y$}; \node[below] at (12.25,-1.94) {outgoing piece};
\end{tikzpicture}%
}
\caption{Three successive insertions at the same section in the ordered cylinder model. Each band is a product with $Y$. At complete breaking, the two common endpoint sections give two product-cylinder components.}\label{fig: successive cylinder insertions}
\end{figure}
\econs
\brem\label{rem: generalization of broken}
The preceding construction extends to a Riemannian manifold $(M,g)$ with disjoint isometrically embedded product collars $[-1,1]\times Y_i$, $i=1,\ldots,k$. Independently stretching these collars gives a \emph{neck-stretching family} over $[0,\infty]^k$, unchanged outside the collars and smooth at the infinite faces in the arclength coordinates above. At $q_i=\infty$, the manifold is cut along $\{0\}\times Y_i$ and completed by two cylindrical ends. The sections need not separate $M$, and no decomposition into composable cobordisms is required.
\erem
\bcons[Continuation cobordism model]\label{cons: continuation cobordism model}
For a chain $I$ with fixed $w$, Remark \ref{rem: smooth fiber} gives a smooth fiber bundle $\mathcal Y_w\to|\Delta^l|$ with fibers $Y_w(s)=\pi_s^{-1}(m(I_w(s)))$ and metrics $g_s$. Fix once and for all a smooth nondecreasing function $\eta:\R\to[0,1]$ such that
\[
\eta(t)=0\text{ for }t\leq-1,~\eta(t)=1\text{ for }t\geq1,~\eta(-t)=1-\eta(t).
\]
For finite $q\in Q_I$, put $\tau_j=\sum_{h=j+1}^p q_h$ for $j\in[0,p]$ and set
\begin{equation}\label{eq: model continuation paths}
\begin{aligned}
\gamma_I(q,t)={}&(1-\eta(t))e_{z_-}+\eta(t-\tau_0)e_{z_+}\\
&+\sum_{j=1}^p\bigl(\eta(t-\tau_j)-\eta(t-\tau_{j-1})\bigr)e_{z_j}.
\end{aligned}
\end{equation}
The coefficients are nonnegative and sum to $1$; see Figure \ref{fig: model continuation paths} for an example. Define the \emph{continuation cobordism model} by
\[
\gamma_I(q,\cdot)^*\mathcal Y_w=\{(t,y)\mid y\in Y_w(\gamma_I(q,t))\},~\bfg_I(q)=dt^2+g_{\gamma_I(q,t)}.
\]
This pullback is noncanonically diffeomorphic to $\R\times Y_{i_-}$. The splitting used in the metric is determined by the unique tangent lift $(1,v)$ of $\partial_t$ with $v$ orthogonal to $T_yY_w(\gamma_I(q,t))$ in $V$ with its standard metric. Use this splitting for the remaining tangential data. The external ends are product for $t\leq-1$ and $t\geq\tau_0+1$; translate their coordinates as in \eqref{eq: open end}.

At $q_j=0$, the coefficient of $e_{z_j}$ vanishes, deleting the corresponding entry. For $q_j>2$, the path is constant at $e_{z_j}$ on $[\tau_j+1,\tau_{j-1}-1]$. As $q_j\to\infty$, the component paths in coordinates $t$ and $t-\tau_{j-1}$ are independent of $q_j$ on each compact subset once $q_j$ is sufficiently large. Thus the family extends smoothly in $(1+q_j)^{-1}$ to the breaking face. For a constant path, the pullback over $[\tau_j,\tau_{j-1}]$ is a product cylinder of length $q_j$, as in the ordered cylinder model.

\begin{figure}[htbp]
\centering
\resizebox{.93\linewidth}{!}{%
\begin{tikzpicture}[font=\small,>=Stealth,line cap=round,line join=round]
\definecolor{modelzero}{RGB}{93,98,105}
\definecolor{modelblue}{RGB}{42,88,119}
\definecolor{modelrust}{RGB}{153,89,64}
\definecolor{modelguide}{RGB}{186,190,194}
\pgfmathdeclarefunction{modelspatialeta}{1}{%
  \begingroup
  \pgfmathparse{#1}\let\modelspatialx\pgfmathresult
  \ifdim\modelspatialx pt>0pt
    \ifdim\modelspatialx pt<1pt
      \ifdim\modelspatialx pt<.5pt
        \pgfmathparse{exp(1/(1-\modelspatialx)-1/\modelspatialx)}%
        \pgfmathparse{\pgfmathresult/(1+\pgfmathresult)}%
      \else
        \pgfmathparse{exp(1/\modelspatialx-1/(1-\modelspatialx))}%
        \pgfmathparse{1/(1+\pgfmathresult)}%
      \fi
    \else
      \def\pgfmathresult{1}%
    \fi
  \else
    \def\pgfmathresult{0}%
  \fi
  \pgfmathsmuggle\pgfmathresult
  \endgroup
}
\begin{scope}[shift={(.65cm,0)},x=3.8cm,y=2.65cm]
\node at (.5,1.18) {(a) Path in $|\Delta^2|$}; \draw[modelguide,line width=.4pt] (0,0)--(.5,.866025)--(1,0)--cycle; \draw[modelrust,line width=1.05pt] (0,0)--(.5,.866025)--(1,0);
\draw[modelblue,line width=1.05pt,densely dashed] plot[domain=0:1.5,samples=161] ({.5*(modelspatialeta(\x)-modelspatialeta(\x-.5))+modelspatialeta(\x-.5)},{.866025*(modelspatialeta(\x)-modelspatialeta(\x-.5))});
\draw[modelzero,line width=.85pt,dash pattern=on 4pt off 2pt] (0,0)--(1,0); \draw[modelzero,line width=.85pt,->] (.36,0)--(.46,0);
\draw[modelrust,line width=1.05pt,->] (.17,.294449)--(.24,.415692); \draw[modelrust,line width=1.05pt,->] (.74,.450333)--(.81,.32909);
\foreach \a/\b in {0/0,.5/.866025,1/0} {\fill (\a,\b) circle[radius=1.25pt];}
\node[below left] at (0,0) {$e_{z_-}$}; \node[above] at (.5,.866025) {$e_{z_1}$}; \node[below right] at (1,0) {$e_{z_+}$};
\end{scope}
\begin{scope}[shift={(7.825cm,0)},x=.825cm,y=1.85cm]
\node at (2.3,1.69) {(b) $\eta(t)-\eta(t-q_1)$}; \draw[modelguide,line width=.4pt,densely dotted] (-1,1)--(6,1);
\draw[line width=.4pt,->] (-1.3,-.045)--(6.7,-.045) node[right] {$t$}; \draw[line width=.4pt,->] (-1,-.045)--(-1,1.23);
\foreach \a/\b in {-1/{-1},0/0,1/1,4/4,6/6} {\draw[line width=.4pt] (\a,-.045)--(\a,-.095);\node[below] at (\a,-.095) {$\b$};}
\node[left] at (-1.05,1) {$1$}; \draw[modelrust,line width=1.05pt] plot[domain=-.12:3.7,samples=281] ({2*\x-1},{modelspatialeta(\x)-modelspatialeta(\x-2.5)});
\draw[modelblue,line width=1.05pt,densely dashed] plot[domain=-.12:3.7,samples=281] ({2*\x-1},{modelspatialeta(\x)-modelspatialeta(\x-.5)});
\draw[modelzero,line width=.85pt,dash pattern=on 4pt off 2pt] (-1.24,0)--(6.4,0); \draw[modelguide,line width=.45pt,densely dotted] (1,0)--(1,1.12) (4,0)--(4,1.12);
\draw[line width=.45pt,<->] (1,1.13)--(4,1.13); \node[above] at (2.5,1.13) {$q_1-2$};
\end{scope}
\begin{scope}[shift={(1.1cm,-.9cm)}]
\draw[modelzero,line width=.85pt,dash pattern=on 4pt off 2pt] (0,0)--(.5,0); \node[right] at (.58,0) {$q_1=0$}; \draw[modelblue,line width=1.05pt,densely dashed] (3.9,0)--(4.4,0);
\node[right] at (4.48,0) {$q_1=1$}; \draw[modelrust,line width=1.05pt] (8.1,0)--(8.6,0); \node[right] at (8.68,0) {$q_1=5$};
\end{scope}
\begin{scope}[shift={(2.25cm,-4.65cm)},x=1.6cm,y=2.3cm]
\node at (.5,1.29) {(c) Coefficients for $q_1=1$};
\path[fill=modelzero!27] (-1.2,0)
  -- plot[domain=-.1:1.6,samples=201] ({2*\x-1},{1-modelspatialeta(\x)})
  --(2.2,0)--cycle;
\path[fill=modelblue!35]
  plot[domain=-.1:1.6,samples=201] ({2*\x-1},{1-modelspatialeta(\x)})
  -- plot[domain=1.6:-.1,samples=201] ({2*\x-1},{1-modelspatialeta(\x-.5)})
  --cycle;
\path[fill=modelrust!32]
  plot[domain=-.1:1.6,samples=201] ({2*\x-1},{1-modelspatialeta(\x-.5)})
  --(2.2,1)--(-1.2,1)--cycle;
\draw[modelzero!75,line width=.6pt] plot[domain=-.1:1.6,samples=201] ({2*\x-1},{1-modelspatialeta(\x)});
\draw[modelblue!85,line width=.6pt] plot[domain=-.1:1.6,samples=201] ({2*\x-1},{1-modelspatialeta(\x-.5)}); \draw[black!65,line width=.65pt] (-1.2,1)--(2.2,1);
\node at (-.7,.48) {$e_{z_-}$}; \node at (.5,.48) {$e_{z_1}$}; \node at (1.7,.48) {$e_{z_+}$}; \draw[line width=.4pt,->] (-1.2,0)--(2.46,0) node[right] {$t$};
\draw[line width=.4pt,->] (-1,-.04)--(-1,1.12);
\foreach \a/\b in {-1/{-1},0/0,1/1,2/2} {
\draw[line width=.4pt] (\a,0)--(\a,-.025); \node[below] at (\a,-.025) {$\b$};
}
\node[left] at (-1.03,1) {$1$};
\end{scope}
\begin{scope}[shift={(7.825cm,-4.65cm)},x=.825cm,y=2.3cm]
\node at (2.5,1.29) {(d) Coefficients for $q_1=5$};
\path[fill=modelzero!27] (-1.24,0)
  -- plot[domain=-.12:3.7,samples=201] ({2*\x-1},{1-modelspatialeta(\x)})
  --(6.4,0)--cycle;
\path[fill=modelblue!35]
  plot[domain=-.12:3.7,samples=201] ({2*\x-1},{1-modelspatialeta(\x)})
  -- plot[domain=3.7:-.12,samples=201] ({2*\x-1},{1-modelspatialeta(\x-2.5)})
  --cycle;
\path[fill=modelrust!32]
  plot[domain=-.12:3.7,samples=201] ({2*\x-1},{1-modelspatialeta(\x-2.5)})
  --(6.4,1)--(-1.24,1)--cycle;
\draw[modelzero!75,line width=.6pt] plot[domain=-.12:3.7,samples=201] ({2*\x-1},{1-modelspatialeta(\x)});
\draw[modelblue!85,line width=.6pt] plot[domain=-.12:3.7,samples=201] ({2*\x-1},{1-modelspatialeta(\x-2.5)}); \draw[black!65,line width=.65pt] (-1.24,1)--(6.4,1);
\node at (-.5,.48) {$e_{z_-}$}; \node at (2.5,.48) {$e_{z_1}$}; \node at (5.5,.48) {$e_{z_+}$}; \draw[line width=.4pt,->] (-1.24,0)--(6.8,0) node[right] {$t$};
\draw[line width=.4pt,->] (-1,-.04)--(-1,1.12);
\foreach \a/\b in {-1/{-1},0/0,1/1,4/4,6/6} {
\draw[line width=.4pt] (\a,0)--(\a,-.025); \node[below] at (\a,-.025) {$\b$};
}
\node[left] at (-1.03,1) {$1$};
\end{scope}
\begin{scope}[shift={(.5cm,-5.55cm)}]
\fill[modelzero!27] (0,-.07) rectangle (.4,.07); \node[anchor=west] at (.53,0) {$1-\eta(t)$}; \fill[modelblue!35] (4.15,-.07) rectangle (4.55,.07);
\node[anchor=west] at (4.68,0) {$\eta(t)-\eta(t-q_1)$}; \fill[modelrust!32] (9.5,-.07) rectangle (9.9,.07); \node[anchor=west] at (10.03,0) {$\eta(t-q_1)$};
\end{scope}
\end{tikzpicture}
}
\caption{The paths $\gamma_I(q,t)$ for $I=\{(z_-,w)<(z_1,w)<(z_+,w)\}$, with $\tau_1=0$, $\tau_0=q_1$, and cutoff transition interval $[-1,1]$. Subfigures (a) and (b) show the paths and middle coefficients for $q_1=0,1,5$. Subfigures (c) and (d) show the stacked coefficients for $q_1=1,5$. At $q_1=0$, the middle entry is deleted. For $q_1>2$, the constant segment at $e_{z_1}$ on $[1,q_1-1]$ gives a product cylinder of length $q_1-2$.}\label{fig: model continuation paths}
\end{figure}

\econs

\bcons[Model metric family]\label{cons: model metric family}
For a general chain $I$, we combine Constructions \ref{cons: ordered cylinder model} and \ref{cons: continuation cobordism model}. Write $\bW_{a,b}(s)$ for the decorated cobordism in $M_s$ between $Y_a(s)$ and $Y_b(s)$. For finite $q\in Q_I$, put $s_j=\gamma_I(q,\tau_j)$, where $\tau_j=\sum_{h=j+1}^p q_h$. We first glue the cobordisms $\bW_{w_{j+1},w_j}(s_j)$ in the order $j=p,\ldots,0$, inserting
\[
\left.\gamma_I(q,\cdot)^*\mathcal Y_{w_j}\right|_{[\tau_j,\tau_{j-1}]},~j\in[1,p],
\]
between consecutive cobordisms. The boundary fibers match, and when $q_j=0$ we identify the common boundary directly. Attach the incoming and outgoing continuation ends given by the pullbacks of $\mathcal Y_{w_-}$ over $(-\infty,0]$ and $\mathcal Y_{w_+}$ over $[\tau_0,\infty)$, respectively. The inherited metrics and tangential data are smooth on each piece but need not glue smoothly along their common boundaries.

For convenience, write $[a_j,b_j]$ for the interval occupied by $\bW_{w_{j+1},w_j}(s_j)$ in a common interval coordinate $\ell$, without rescaling. These endpoints are determined by
\begin{equation}\label{eq: mixed portion intervals}
\begin{aligned}
a_p&=0,\\
b_j&=a_j+m(I_{w_j}(s_j))-m(I_{w_{j+1}}(s_j)),&&j=p,\ldots,0,\\
a_{j-1}&=b_j+q_j,&&j=p,\ldots,1.
\end{aligned}
\end{equation}
No compact cobordism is inserted between coincident sections, but their entries and parameters are kept. Before smoothing, the parameter $t$ of $\gamma_I(q,t)$ equals $\tau_j$ on $[a_j,b_j]$ and $\tau_j+\ell-b_j$ on $[b_j,a_{j-1}]$.

Choose a smooth convex function $H$ equal to $\max(t,0)$ outside a sufficiently small neighborhood of $0$ and satisfying $H(-t)=H(t)-t$, and set
\begin{equation}\label{eq: mixed smoothed time}
T_I(q,\ell)=\ell-\sum_{j=0}^p\bigl(H(\ell-a_j)-H(\ell-b_j)\bigr).
\end{equation}
Use the same $H$ for all chains over $\sigma$. Choose the smoothing small enough that the changes occur only in the product regions of Definition \ref{defn: cylindrical}. A common choice is possible by compactness of $|\Delta^l|$, the finite number of cuts, and the bound on $\partial_t\gamma_I$ independent of the chain length. Then $T_I=\tau_j$ away from these product regions, while on the open ends
\begin{equation}\label{eq: T_I open end}
   T_I(q,\ell)=\ell~\text{for }\ell\ll0,~
T_I(q,\ell)=\tau_0+\ell-b_0~\text{for }\ell\gg b_0. 
\end{equation}
Figure \ref{fig: smoothing continuation parameter} illustrates this smoothing.

On each product region associated with $w$, pull back $\mathcal Y_w$ along $\ell\mapsto\gamma_I(q,T_I(q,\ell))$ and use
\begin{equation}\label{eq: mixed cutoff metric}
\bfg_I(T_I)(q)=d\ell^2+g_{Y_w(\gamma_I(q,T_I(q,\ell)))}.
\end{equation}
Use the splitting and tangential pullbacks of Construction \ref{cons: continuation cobordism model}, and keep the original data elsewhere. For disconnected cobordisms, the construction is understood componentwise.

At $q_j=0$, adjacent terms in \eqref{eq: mixed smoothed time} cancel. Near $q_j=\infty$, use
\[
(\ell,T_I)~\text{and}~(\ell-a_{j-1},T_I-\tau_{j-1})
\]
on the incoming and outgoing components, respectively. On compact subsets of the incoming and outgoing components, these formulas agree with the corresponding component formulas for sufficiently large $q_j$. Outside bounded neighborhoods of the two ends of a sufficiently long inserted cylinder, $T_I=\tau_j+\ell-b_j$ and $\gamma_I(q,T_I)=e_{z_j}$. Thus the family is smooth at infinity and satisfies
\begin{equation}\label{eq: model metric faces}
\begin{aligned}
\bfg_I(T_I)|_{q_j=0}&=\bfg_{I\backslash\{j\}}(T_{I\backslash\{j\}}),\\
\bfg_I(T_I)|_{q_j=\infty}
&=\bigl(\bfg_{I_j^+}(T_{I_j^+}),\bfg_{I_j^-}(T_{I_j^-})\bigr).
\end{aligned}
\end{equation}
Cancellation and the translations commute, so iterated face restrictions agree. For fixed $z$, the metric recovers the ordered cylinder model. For fixed $w$, all $a_j=b_j$, so $T_I=\ell$ and the construction recovers the continuation cobordism model.
\econs

\begin{figure}[htbp]
\centering
\resizebox{\linewidth}{!}{%
\begin{tikzpicture}[font=\small,>=Stealth,line cap=round,line join=round]
\definecolor{smoothblue}{RGB}{42,88,119}
\definecolor{oldgray}{RGB}{110,113,118}
\begin{scope}[x=1.48cm,y=1.35cm]
\node at (2.05,3.17) {(a) The reparametrization $\ell\mapsto t$};
\foreach \c in {0,1,3,4} {
\fill[smoothblue!9] ({\c-.25},-.6) rectangle ({\c+.25},2.82);
}
\draw[black!16,densely dotted] (-.55,0)--(4.7,0); \draw[black!16,densely dotted] (-.55,2)--(4.7,2);
\draw[smoothblue,line width=1.2pt] plot coordinates {
(-0.550000,-0.550000) (-0.535000,-0.535000) (-0.520000,-0.520000) (-0.505000,-0.505000) (-0.490000,-0.490000) (-0.475000,-0.475000) (-0.460000,-0.460000) (-0.445000,-0.445000)
(-0.430000,-0.430000) (-0.415000,-0.415000) (-0.400000,-0.400000) (-0.385000,-0.385000) (-0.370000,-0.370000) (-0.355000,-0.355000) (-0.340000,-0.340000) (-0.325000,-0.325000)
(-0.310000,-0.310000) (-0.295000,-0.295000) (-0.280000,-0.280000) (-0.265000,-0.265000) (-0.250000,-0.250000) (-0.235000,-0.235000) (-0.220000,-0.220000) (-0.205000,-0.205000)
(-0.190000,-0.190004) (-0.175000,-0.175036) (-0.160000,-0.160155) (-0.145000,-0.145462) (-0.130000,-0.131081) (-0.115000,-0.117144) (-0.100000,-0.103770) (-0.085000,-0.091061)
(-0.070000,-0.079099) (-0.055000,-0.067941) (-0.040000,-0.057629) (-0.025000,-0.048192) (-0.010000,-0.039644) (0.005000,-0.031994) (0.020000,-0.025243) (0.035000,-0.019385)
(0.050000,-0.014408) (0.065000,-0.010288) (0.080000,-0.006988) (0.095000,-0.004455) (0.110000,-0.002618) (0.125000,-0.001380) (0.140000,-0.000627) (0.155000,-0.000231)
(0.170000,-0.000062) (0.185000,-0.000010) (0.200000,-0.000001) (0.215000,0) (0.230000,0) (0.245000,0) (0.260000,0) (0.275000,0)
(0.290000,0) (0.305000,0) (0.320000,0) (0.335000,0) (0.350000,0) (0.365000,0) (0.380000,0) (0.395000,0)
(0.410000,0) (0.425000,0) (0.440000,0) (0.455000,0) (0.470000,0) (0.485000,0) (0.500000,0) (0.515000,0)
(0.530000,0) (0.545000,0) (0.560000,0) (0.575000,0) (0.590000,0) (0.605000,0) (0.620000,0) (0.635000,0)
(0.650000,0) (0.665000,0) (0.680000,0) (0.695000,0) (0.710000,0) (0.725000,0) (0.740000,0) (0.755000,0)
(0.770000,0) (0.785000,0) (0.800000,0.000001) (0.815000,0.000010) (0.830000,0.000062) (0.845000,0.000231) (0.860000,0.000627) (0.875000,0.001380)
(0.890000,0.002618) (0.905000,0.004455) (0.920000,0.006988) (0.935000,0.010288) (0.950000,0.014408) (0.965000,0.019385) (0.980000,0.025243) (0.995000,0.031994)
(1.010000,0.039644) (1.025000,0.048192) (1.040000,0.057629) (1.055000,0.067941) (1.070000,0.079099) (1.085000,0.091061) (1.100000,0.103770) (1.115000,0.117144)
(1.130000,0.131081) (1.145000,0.145462) (1.160000,0.160155) (1.175000,0.175036) (1.190000,0.190004) (1.205000,0.205000) (1.220000,0.220000) (1.235000,0.235000)
(1.250000,0.250000) (1.265000,0.265000) (1.280000,0.280000) (1.295000,0.295000) (1.310000,0.310000) (1.325000,0.325000) (1.340000,0.340000) (1.355000,0.355000)
(1.370000,0.370000) (1.385000,0.385000) (1.400000,0.400000) (1.415000,0.415000) (1.430000,0.430000) (1.445000,0.445000) (1.460000,0.460000) (1.475000,0.475000)
(1.490000,0.490000) (1.505000,0.505000) (1.520000,0.520000) (1.535000,0.535000) (1.550000,0.550000) (1.565000,0.565000) (1.580000,0.580000) (1.595000,0.595000)
(1.610000,0.610000) (1.625000,0.625000) (1.640000,0.640000) (1.655000,0.655000) (1.670000,0.670000) (1.685000,0.685000) (1.700000,0.700000) (1.715000,0.715000)
(1.730000,0.730000) (1.745000,0.745000) (1.760000,0.760000) (1.775000,0.775000) (1.790000,0.790000) (1.805000,0.805000) (1.820000,0.820000) (1.835000,0.835000)
(1.850000,0.850000) (1.865000,0.865000) (1.880000,0.880000) (1.895000,0.895000) (1.910000,0.910000) (1.925000,0.925000) (1.940000,0.940000) (1.955000,0.955000)
(1.970000,0.970000) (1.985000,0.985000) (2.000000,1.000000) (2.015000,1.015000) (2.030000,1.030000) (2.045000,1.045000) (2.060000,1.060000) (2.075000,1.075000)
(2.090000,1.090000) (2.105000,1.105000) (2.120000,1.120000) (2.135000,1.135000) (2.150000,1.150000) (2.165000,1.165000) (2.180000,1.180000) (2.195000,1.195000)
(2.210000,1.210000) (2.225000,1.225000) (2.240000,1.240000) (2.255000,1.255000) (2.270000,1.270000) (2.285000,1.285000) (2.300000,1.300000) (2.315000,1.315000)
(2.330000,1.330000) (2.345000,1.345000) (2.360000,1.360000) (2.375000,1.375000) (2.390000,1.390000) (2.405000,1.405000) (2.420000,1.420000) (2.435000,1.435000)
(2.450000,1.450000) (2.465000,1.465000) (2.480000,1.480000) (2.495000,1.495000) (2.510000,1.510000) (2.525000,1.525000) (2.540000,1.540000) (2.555000,1.555000)
(2.570000,1.570000) (2.585000,1.585000) (2.600000,1.600000) (2.615000,1.615000) (2.630000,1.630000) (2.645000,1.645000) (2.660000,1.660000) (2.675000,1.675000)
(2.690000,1.690000) (2.705000,1.705000) (2.720000,1.720000) (2.735000,1.735000) (2.750000,1.750000) (2.765000,1.765000) (2.780000,1.780000) (2.795000,1.795000)
(2.810000,1.809996) (2.825000,1.824964) (2.840000,1.839845) (2.855000,1.854538) (2.870000,1.868919) (2.885000,1.882856) (2.900000,1.896230) (2.915000,1.908939)
(2.930000,1.920901) (2.945000,1.932059) (2.960000,1.942371) (2.975000,1.951808) (2.990000,1.960356) (3.005000,1.968006) (3.020000,1.974757) (3.035000,1.980615)
(3.050000,1.985592) (3.065000,1.989712) (3.080000,1.993012) (3.095000,1.995545) (3.110000,1.997382) (3.125000,1.998620) (3.140000,1.999373) (3.155000,1.999769)
(3.170000,1.999938) (3.185000,1.999990) (3.200000,1.999999) (3.215000,2.000000) (3.230000,2.000000) (3.245000,2.000000) (3.260000,2.000000) (3.275000,2.000000)
(3.290000,2.000000) (3.305000,2.000000) (3.320000,2.000000) (3.335000,2.000000) (3.350000,2.000000) (3.365000,2.000000) (3.380000,2.000000) (3.395000,2.000000)
(3.410000,2.000000) (3.425000,2.000000) (3.440000,2.000000) (3.455000,2.000000) (3.470000,2.000000) (3.485000,2.000000) (3.500000,2.000000) (3.515000,2.000000)
(3.530000,2.000000) (3.545000,2.000000) (3.560000,2.000000) (3.575000,2.000000) (3.590000,2.000000) (3.605000,2.000000) (3.620000,2.000000) (3.635000,2.000000)
(3.650000,2.000000) (3.665000,2.000000) (3.680000,2.000000) (3.695000,2.000000) (3.710000,2.000000) (3.725000,2.000000) (3.740000,2.000000) (3.755000,2.000000)
(3.770000,2.000000) (3.785000,2.000000) (3.800000,2.000001) (3.815000,2.000010) (3.830000,2.000062) (3.845000,2.000231) (3.860000,2.000627) (3.875000,2.001380)
(3.890000,2.002618) (3.905000,2.004455) (3.920000,2.006988) (3.935000,2.010288) (3.950000,2.014408) (3.965000,2.019385) (3.980000,2.025243) (3.995000,2.031994)
(4.010000,2.039644) (4.025000,2.048192) (4.040000,2.057629) (4.055000,2.067941) (4.070000,2.079099) (4.085000,2.091061) (4.100000,2.103770) (4.115000,2.117144)
(4.130000,2.131081) (4.145000,2.145462) (4.160000,2.160155) (4.175000,2.175036) (4.190000,2.190004) (4.205000,2.205000) (4.220000,2.220000) (4.235000,2.235000)
(4.250000,2.250000) (4.265000,2.265000) (4.280000,2.280000) (4.295000,2.295000) (4.310000,2.310000) (4.325000,2.325000) (4.340000,2.340000) (4.355000,2.355000)
(4.370000,2.370000) (4.385000,2.385000) (4.400000,2.400000) (4.415000,2.415000) (4.430000,2.430000) (4.445000,2.445000) (4.460000,2.460000) (4.475000,2.475000)
(4.490000,2.490000) (4.505000,2.505000) (4.520000,2.520000) (4.535000,2.535000) (4.550000,2.550000) (4.565000,2.565000) (4.580000,2.580000) (4.595000,2.595000)
(4.610000,2.610000) (4.625000,2.625000) (4.640000,2.640000) (4.655000,2.655000) (4.670000,2.670000) (4.685000,2.685000) (4.700000,2.700000)
};
\draw[oldgray,line width=.85pt,dash pattern=on 4pt off 2pt]
(-.55,-.55)--(0,0)--(1,0)--(3,2)--(4,2)--(4.7,2.7);
\draw[black!65,line width=.5pt,->] (-.65,-.67)--(4.87,-.67) node[right] {$\ell$}; \draw[black!65,line width=.5pt,->] (-.65,-.67)--(-.65,2.93) node[above] {$t$};
\foreach \x/\lab in {0/{a_1=0},1/{b_1=1},3/{a_0=3},4/{b_0=4}} {
\draw[black!65] (\x,-.67)--(\x,-.72); \node[below,font=\scriptsize] at (\x,-.74) {$\lab$};
}
\draw[black!65] (-.65,0)--(-.7,0) node[left] {$\tau_1=0$}; \draw[black!65] (-.65,2)--(-.7,2) node[left] {$\tau_0=2$}; \node[smoothblue,above] at (.5,.05) {$t=\tau_1$};
\node[smoothblue,above] at (3.5,2.05) {$t=\tau_0$}; \node[rotate=42,above,text=black!65] at (2.0,1.15) {continuation};
\draw[smoothblue!65,line width=.55pt,densely dotted] (.7,-.09) rectangle (1.3,.36);
\end{scope}
\begin{scope}[shift={(10.75cm,0)},x=6.4cm,y=8.8cm]
\node at (.005,.486) {(b) Near $\ell=b_1=1$}; \fill[smoothblue!9] (-.25,-.035) rectangle (.25,.425);
\draw[black!20,densely dotted] (-.25,-.035)--(-.25,.425) (.25,-.035)--(.25,.425);
\draw[smoothblue,line width=1.2pt] plot coordinates {
(-0.360000,0) (-0.355889,0) (-0.351778,0) (-0.347667,0) (-0.343556,0) (-0.339444,0) (-0.335333,0) (-0.331222,0)
(-0.327111,0) (-0.323000,0) (-0.318889,0) (-0.314778,0) (-0.310667,0) (-0.306556,0) (-0.302444,0) (-0.298333,0)
(-0.294222,0) (-0.290111,0) (-0.286000,0) (-0.281889,0) (-0.277778,0) (-0.273667,0) (-0.269556,0) (-0.265444,0)
(-0.261333,0) (-0.257222,0) (-0.253111,0) (-0.249000,0) (-0.244889,0) (-0.240778,0) (-0.236667,0) (-0.232556,0)
(-0.228444,0) (-0.224333,0) (-0.220222,0) (-0.216111,0) (-0.212000,0) (-0.207889,0.000000) (-0.203778,0.000000) (-0.199667,0.000001)
(-0.195556,0.000002) (-0.191444,0.000003) (-0.187333,0.000007) (-0.183222,0.000013) (-0.179111,0.000022) (-0.175000,0.000036) (-0.170889,0.000056) (-0.166778,0.000085)
(-0.162667,0.000123) (-0.158556,0.000175) (-0.154444,0.000241) (-0.150333,0.000324) (-0.146222,0.000427) (-0.142111,0.000553) (-0.138000,0.000704) (-0.133889,0.000884)
(-0.129778,0.001094) (-0.125667,0.001337) (-0.121556,0.001617) (-0.117444,0.001935) (-0.113333,0.002295) (-0.109222,0.002697) (-0.105111,0.003146) (-0.101000,0.003642)
(-0.096889,0.004187) (-0.092778,0.004785) (-0.088667,0.005435) (-0.084556,0.006140) (-0.080444,0.006902) (-0.076333,0.007721) (-0.072222,0.008599) (-0.068111,0.009537)
(-0.064000,0.010536) (-0.059889,0.011597) (-0.055778,0.012721) (-0.051667,0.013909) (-0.047556,0.015160) (-0.043444,0.016476) (-0.039333,0.017858) (-0.035222,0.019305)
(-0.031111,0.020819) (-0.027000,0.022399) (-0.022889,0.024046) (-0.018778,0.025759) (-0.014667,0.027540) (-0.010556,0.029389) (-0.006444,0.031305) (-0.002333,0.033288)
(0.001778,0.035339) (0.005889,0.037458) (0.010000,0.039644) (0.014111,0.041897) (0.018222,0.044218) (0.022333,0.046607) (0.026444,0.049062) (0.030556,0.051584)
(0.034667,0.054173) (0.038778,0.056827) (0.042889,0.059548) (0.047000,0.062334) (0.051111,0.065185) (0.055222,0.068100) (0.059333,0.071079) (0.063444,0.074120)
(0.067556,0.077224) (0.071667,0.080389) (0.075778,0.083614) (0.079889,0.086898) (0.084000,0.090240) (0.088111,0.093638) (0.092222,0.097092) (0.096333,0.100598)
(0.100444,0.104157) (0.104556,0.107765) (0.108667,0.111422) (0.112778,0.115124) (0.116889,0.118870) (0.121000,0.122658) (0.125111,0.126484) (0.129222,0.130347)
(0.133333,0.134243) (0.137444,0.138171) (0.141556,0.142127) (0.145667,0.146109) (0.149778,0.150114) (0.153889,0.154140) (0.158000,0.158183) (0.162111,0.162241)
(0.166222,0.166311) (0.170333,0.170393) (0.174444,0.174483) (0.178556,0.178579) (0.182667,0.182680) (0.186778,0.186785) (0.190889,0.190893) (0.195000,0.195002)
(0.199111,0.199112) (0.203222,0.203222) (0.207333,0.207333) (0.211444,0.211444) (0.215556,0.215556) (0.219667,0.219667) (0.223778,0.223778) (0.227889,0.227889)
(0.232000,0.232000) (0.236111,0.236111) (0.240222,0.240222) (0.244333,0.244333) (0.248444,0.248444) (0.252556,0.252556) (0.256667,0.256667) (0.260778,0.260778)
(0.264889,0.264889) (0.269000,0.269000) (0.273111,0.273111) (0.277222,0.277222) (0.281333,0.281333) (0.285444,0.285444) (0.289556,0.289556) (0.293667,0.293667)
(0.297778,0.297778) (0.301889,0.301889) (0.306000,0.306000) (0.310111,0.310111) (0.314222,0.314222) (0.318333,0.318333) (0.322444,0.322444) (0.326556,0.326556)
(0.330667,0.330667) (0.334778,0.334778) (0.338889,0.338889) (0.343000,0.343000) (0.347111,0.347111) (0.351222,0.351222) (0.355333,0.355333) (0.359444,0.359444)
(0.363556,0.363556) (0.367667,0.367667) (0.371778,0.371778) (0.375889,0.375889) (0.380000,0.380000)
};
\draw[oldgray,line width=.85pt,dash pattern=on 4pt off 2pt] (-.36,0)--(0,0)--(.38,.38); \draw[black!65,line width=.5pt,->] (-.385,-.035)--(.415,-.035);
\node[below right] at (.365,-.06) {$\ell-b_1$}; \draw[black!65,line width=.5pt,->] (-.385,-.035)--(-.385,.446) node[above] {$t$};
\foreach \x/\lab in {-.25/{-\delta},0/0,.25/{\delta}} {
\draw[black!65] (\x,-.035)--(\x,-.042); \node[below] at (\x,-.045) {$\lab$};
}
\draw[black!65] (-.385,0)--(-.392,0) node[left] {$0$}; \draw[black!65] (-.385,.25)--(-.392,.25) node[left] {$\delta$};
\node[anchor=west,smoothblue] at (-.22,.34) {$T_I=H(\ell-b_1)$}; \draw[smoothblue!65,line width=.55pt,->] (-.10,.295)--(.005,.048);
\end{scope}
\begin{scope}[shift={(-.2cm,-1.85cm)}]
\draw[oldgray,line width=.85pt,dash pattern=on 4pt off 2pt] (0,0)--(.65,0); \node[anchor=west] at (.78,0) {Before smoothing}; \draw[smoothblue,line width=1.2pt] (4.8,0)--(5.45,0);
\node[anchor=west] at (5.58,0) {After smoothing}; \fill[smoothblue!9] (9.3,-.09) rectangle (9.8,.09); \node[anchor=west] at (9.94,0) {$\delta$-neighborhoods};
\end{scope}
\end{tikzpicture}
}
\caption{Smoothing the continuation parameter for $q_1=2$, $[a_1,b_1]=[0,1]$, $[a_0,b_0]=[3,4]$, and $\delta=\frac14$. Gray dashed: before smoothing; blue: $T_I(q,\ell)$. The curves agree outside the shaded neighborhoods. (b) enlarges the neighborhood of $b_1$.}\label{fig: smoothing continuation parameter}
\end{figure}

\bdefn[Smoothing collections]\label{defn: model function collections}
A \emph{smoothing collection} is a pair
\[
\bigl(\sigma,(T_I)_{I\in\cI_{l,k}}\bigr),
\]
where $\sigma\in(\PBI_a^{V,\aff})_{k,l}$ and each $T_I(q,\ell)$ is smooth for finite $q\in Q_I$ and $\ell\in\R$. Use the notation and product identifications of Construction \ref{cons: model metric family}. All functions and constructions are understood componentwise for disjoint unions. Require the following conditions.
\begin{itemize}
\item For each $\bW_{w_{j+1},w_j}(s_j)$, require
\begin{equation}\label{eq: T_I constant}
T_I(q,\ell)=\tau_j~\text{for }a_j+\delta\leq\ell\leq b_j-\delta
\end{equation}
for some $\delta>0$ such that the endpoint $\delta$-neighborhoods in the $\ell$-coordinate lie in the chosen product regions. On these neighborhoods, every product identification used must be defined at $\gamma_I(q,t)$ for all $t$ between $\tau_j$ and $T_I(q,\ell)$.
\item On the open ends, require \eqref{eq: T_I open end}. For sufficiently large $q_j$, require
\begin{equation}\label{eq: T_I long cylinder}
T_I(q,\ell)=\tau_j+\ell-b_j
\end{equation}
on $[b_j,a_{j-1}]$ outside neighborhoods of its two endpoints. The lengths of these neighborhoods are bounded independently of $q_j$.
\item At $q_j=0$, require
\begin{equation}\label{eq: T_I zero face}
T_I|_{q_j=0}=T_{I\backslash\{j\}}.
\end{equation}
At $q_j=\infty$, use
\begin{equation}\label{eq: mixed component coordinates}
(\ell,T_I)~\text{and}~(\ell-a_{j-1},T_I-\tau_{j-1})
\end{equation}
on the incoming and outgoing components, respectively. In these coordinates, require smooth extension to $T_{I_j^-}$ and $T_{I_j^+}$. All iterated restrictions agree.
\end{itemize}
For chains with fixed $w$, require $T_I=\ell$, recovering the continuation cobordism model. The neighborhood conditions and length bounds are locally uniform in the remaining parameters and are transported by the face identifications.

By \eqref{eq: mixed cutoff metric}, a smoothing collection defines decorated metric families
\[
\bigl(\bW_I(T_I),\bfg_I(T_I)\bigr)
\]
satisfying \eqref{eq: model metric faces}. Let $(\PBI_a^{V,\mathrm{met}})_{k,l}$ be the set of smoothing collections. With the notation and parameter identifications of Definition \ref{defn: abstract auxiliary model}, define the simplicial and spatial maps by
\begin{equation}\label{eq: pull back T}
    (\rho^*T)_J=T_{\rho J},~\rho=\alpha\times\beta.
\end{equation}
These maps preserve the above conditions and are compatible with the metric construction by $|\alpha|\circ\gamma_J=\gamma_{\rho J}$. The identity and composition laws follow from taking ordered image chains. Thus $\PBI_a^{V,\mathrm{met}}$ is a simplicial space with the forgetful morphism
\[
p_a^{V,\mathrm{met}}:\PBI_a^{V,\mathrm{met}}\longrightarrow\PBI_a^{V,\aff}.
\]
\edefn

\bprop\label{prop: trivial kan smoothing}
For every $k\geq0$, the map
\[
p_{a,k}^{V,\mathrm{met}}:
(\PBI_a^{V,\mathrm{met}})_{k,\bu}
\longrightarrow(\PBI_a^{V,\aff})_{k,\bu}
\]
is a trivial Kan fibration in the sense of Definition \ref{defn: trivial kan fibration}. Hence $p_a^{V,\mathrm{met}}$ is a levelwise weak homotopy equivalence.
\eprop
\bpf
Fix $\sigma\in(\PBI_a^{V,\aff})_{k,l}$ and a lift over $\partial\Delta^l$. Keep the given functions and construct the remaining $T_I$ by induction on $|I|$. The functions $T_{I\backslash\{j\}}$ and $T_{I_j^\pm}$ determine compatible values on $\partial Q_I$.

Let $T_I^0$ be the standard function of Construction \ref{cons: model metric family}. For each codimension-one face $F$ of $Q_I$, choose a neighborhood $U_F$ and a local extension $T_I^F$ agreeing with all given boundary functions on $U_F\cap\partial Q_I$. Extend the differences from $T_I^0$ in corner coordinates while preserving \eqref{eq: T_I constant}, \eqref{eq: T_I open end}, and \eqref{eq: T_I long cylinder}. At infinite faces, use the coordinates in \eqref{eq: mixed component coordinates} and $T_I^0$ on the middle of sufficiently long cylinders.

Choose a partition of unity subordinate to the neighborhoods $U_F$ and the interior, with $\chi_0=0$ on $\partial Q_I$ and $\chi_0=1$ outside a boundary neighborhood. Set
\[
\begin{aligned}
T_I&=\chi_0(q)T_I^0+\sum_F\chi_F(q)T_I^F,\\
\chi_0,\chi_F&\geq0,~\chi_0+\sum_F\chi_F=1.
\end{aligned}
\]
The convex combination preserves \eqref{eq: T_I constant}, \eqref{eq: T_I open end}, and \eqref{eq: T_I long cylinder}, as well as the allowed intervals for $T_I$ in each product identification. Locally, take common maxima of the widths $\delta_0,\delta_F$ and the endpoint bounds. Since all local extensions agree with the given function on each boundary face, so does their weighted sum. This gives \eqref{eq: T_I zero face} and \eqref{eq: mixed component coordinates}. For fixed $w$, use $T_I=\ell$.

Induction gives a lift over $\sigma$, also when $l=0$. Hence the map has the right lifting property for every $\partial\Delta^l\hookrightarrow\Delta^l$. The last assertion follows from Lemma \ref{lem: trivial Kan fibration}.
\epf

We now fix a smoothing collection and choose perturbations on the resulting decorated metric families.

\bdefn[Perturbations]\label{defn: broken perturbation}
An \emph{auxiliary collection} is a tuple
\begin{equation}\label{eq: aux col}
    \bigl(\sigma,(T_I)_{I\in\cI_{l,k}},(\pi_i)_{i\in P_{l,k}},(\Pi_I)_{I\in\cI_{l,k}}\bigr),
\end{equation}
where $(\sigma,(T_I)_{I\in \cI_{l,k}})\in(\PBI_a^{V,\met})_{k,l}$ and $\pi_i\in\cP(Y_i)$ is a primary perturbation. The secondary data $\Pi_I$ are chosen on the resulting metric families $\bW_I(T_I)$, suppressing $T_I$ from the notation. Require the following conditions.
\benum
\item For each $I\in\cI_{l,k}$, the datum $\Pi_I$ includes a smooth $Q_I$-family of external-end data in
\[
\cP(-Y_{i_-})\times\cP(Y_{i_+}).
\]
For every family, including constant chains, use \eqref{eq: perturbed ASD 2} at its two ends in outward collar coordinates. The primary boundary datum is $(-\pi_{i_-},\pi_{i_+})$.
\item Near $q_j=\infty$, include the secondary data at the adjoining ends of $\bW_{I_j^-}$ and $\bW_{I_j^+}$ required by \eqref{eq: abstract auxiliary face}. In the product coordinate $t$ increasing along the neck over $Y_{i_j}$, use
\begin{equation}\label{eq: perturbed ASD 3}
\begin{aligned}
F_A^{+_{\bfg_I}}+\widehat V_{\pi_{i_j},\bfg_I}(A)
&+\psi_-(t)\widehat V_{\pi'_{j,-},\bfg_I}(A)\\
&+\psi_+(t)\widehat V_{\pi'_{j,+},\bfg_I}(A)=0,
\end{aligned}
\end{equation}
where $\pi'_{j,\pm}$ and $\psi_\pm$ are obtained from the secondary data and the bump function in \eqref{eq: perturbed ASD 2} at the adjoining ends by expressing them in the common coordinate $t$.
\item The full data extend smoothly in the model component coordinates and satisfy the end and face conditions of Definition \ref{defn: cubical auxiliary data}.
\eenu
Let $(\PBI_a^{V,\col})_{k,l}$ be the set of these auxiliary collections. With the notation and parameter identifications of Definition \ref{defn: abstract auxiliary model}, pull back the smoothing functions as in \eqref{eq: pull back T} and set
\begin{equation}\label{eq: pull back pi}
    (\rho^*\pi)_j=\pi_{\rho(j)},~(\rho^*\Pi)_J=\Pi_{\rho J},~\rho=\alpha\times\beta.
\end{equation}
These maps preserve the conditions above. Their identity and composition laws, and the commutation of the simplicial and spatial maps, follow from taking ordered image chains. Thus $\PBI_a^{V,\col}$ is a simplicial space with a forgetful morphism
\[
\PBI_a^{V,\col}\longrightarrow\PBI_a^{V,\met}.
\]
\edefn

We impose the full boundary condition with superscript $\partial$ from Definition \ref{defn: PBIV k0}, so each connected cobordism component has a nonempty incoming and outgoing end. The variants with $\partial_+$ or $\partial_-$ are treated similarly.

\bdefn[Admissibility]\label{defn: admissible auxiliary simplex}
An auxiliary collection \eqref{eq: aux col} over $\sigma\in(\PBI_a^{\partial,V,\aff})_{k,l}$ is \emph{admissible} if the following conditions hold.
\begin{itemize}
\item For every $i\in P_{l,k}$, the primary perturbation $\pi_i$ belongs to $\cP^{\adm}(Y_i)$ as defined in \eqref{eq: adm Y}.
\item For every weak chain $I\in\cI_{l,k}$, the family $\Pi_I$ satisfies the parametrized admissibility condition of \eqref{eq: adm W}, using \eqref{eq: perturbed ASD 3} near the breaking faces. At a broken stratum, fix the matching critical points and use the product of the component moduli spaces, with independent parameters.
\end{itemize}
The simplicial and spatial pullbacks identify the equations over $Q_J=Q_{\rho J}$ and therefore preserve admissibility. Thus these collections form a simplicial subspace
\[
\PBI_a^{\partial,V,\adm}\subseteq\PBI_a^{\partial,V,\col},
\]
with the forgetful morphism
\[
p_a^{\partial,V,\adm}:\PBI_a^{\partial,V,\adm}\longrightarrow\PBI_a^{\partial,V,\aff}.
\]
For the corresponding $\Ga$-diagrams as in Definition \ref{defn: abstract auxiliary model}, require admissibility for every union of the component families over their common parameter space. Disjoint unions over independent parameter spaces are admissible by the product description of their moduli spaces.
\edefn
\bcons[Perturbations for comparison families]\label{cons: ga diagram}
Realize the comparison collections of Construction \ref{cons: comparison families} by secondary perturbations $\wti\Pi_{J,G}^\alpha$ as in Definition \ref{defn: broken perturbation}, keeping the metrics and primary perturbations fixed. On the remaining cell interiors, the perturbations may depend on every parameter and vary independently on each component.

The collection is \emph{admissible} if its original data satisfy Definition \ref{defn: admissible auxiliary simplex} and, on every stratum $G$ of every comparison cell, the component moduli spaces are regular and
\begin{equation}\label{eq: comparison admissibility}
\mathop{\times}_{G,\alpha\in S}\cM_G(\bW^\alpha,\wti\Pi^\alpha)
\end{equation}
is transverse for every nonempty $S\subseteq S_{\rho,\beta}$. Require the same after deleting unused parameters; at broken strata fix the matching critical points. These collections form the subdiagram
\[
\PBI_a^{\partial,V,\com,\adm}\subseteq\PBI_a^{\partial,V,\com},
\]
with the forgetful $\Ga$-morphism
\[
p_a^{\partial,V,\com,\adm}:\PBI_a^{\partial,V,\com,\adm}\longrightarrow\PBI_a^{\partial,V,\adm}.
\]
\econs

\bcons[Duality for collections]\label{cons: reflected model families}
For $\sigma\in(\PBI_a^{V,\aff})_{k,l}$, let $r_l(s_0,\ldots,s_l)=(s_l,\ldots,s_0)$ and write $1-[a,b]=[1-b,1-a]$. Define $\sigma^\vee(s)$ by reflecting the $\R$-coordinate of $M_{r_l(s)}\subset V\times\R$, transporting its decorations and product identifications, and replacing its interval list by
\[
\I^\vee(s)=\bigl(1-I_k(r_l(s))\leq\cdots\leq1-I_0(r_l(s))\bigr).
\]
For $i=(z,w)$, put $i^\vee=(l-z,k-w)$; reverse a weak chain $I$ and apply this operation to each entry to obtain $I^\vee$. For $p=|I|-2$, identify
\begin{equation}\label{eq: dual cube coordinates}
Q_I\longrightarrow Q_{I^\vee},~q=(q_1,\ldots,q_p)\longmapsto q^\vee=(q_p,\ldots,q_1).
\end{equation}
For finite parameters, the reflected metric and smoothing data are
\begin{equation}\label{eq: reflected continuation}
\begin{aligned}
\tau_j^\vee&=\tau_0-\tau_{p-j},&
\gamma_{I^\vee}(q^\vee,t)&=r_l\gamma_I(q,\tau_0-t),\\
a_j^\vee&=b_0-b_{p-j},&
b_j^\vee&=b_0-a_{p-j},
\end{aligned}
\end{equation}
\begin{equation}\label{eq: reflected smoothing}
(T^\vee)_{I^\vee}(q^\vee,\ell)=\tau_0-T_I(q,b_0-\ell).
\end{equation}
The standard smoothing formula follows from $H(-t)=H(t)-t$. These formulas preserve the conditions of Definition \ref{defn: model function collections} and extend over infinite faces by \eqref{eq: mixed component coordinates}, exchanging incoming and outgoing ends.

For $\bY=(Y,P,\phi,J)$ and $a=(g,\pi)$, set $\bY^\vee=(-Y,P,\phi,J)$ and $a^\vee=(g,-\pi)$. Reflect a family with fixed end data to
\[
(\bW_G^\vee,\bfa^\vee):(\bY_+^\vee,a_+^\vee)\longrightarrow(\bY_-^\vee,a_-^\vee).
\]
Pull back perturbations along \eqref{eq: dual cube coordinates} and exchange external secondary data in outward coordinates. In a centered directed neck coordinate, use
\begin{equation}\label{eq: reflected perturbations}
(\pi,\pi'_-,\pi'_+,t)\longmapsto(-\pi,-\pi'_+,-\pi'_-, -t),
\end{equation}
transporting cutoffs by $t\mapsto-t$. The resulting identifications of equations and linearizations preserve face conditions and admissibility. With $\rho^\vee$ as in \eqref{eq: alpha beta dual},
\begin{equation}\label{eq: reflected auxiliary operators}
(\rho^*x)^\vee=(\rho^\vee)^*(x^\vee),~x^{\vee\vee}=x.
\end{equation}

For comparison collections, use the chart reflection of Definition \ref{defn: diagrams with strict duality} and reverse $J$ while sending $(c,i)$ to $(u-c,i^\vee)$. Let $G^\vee$ be the image of $G$ under coordinate reversal. Keeping $\rho$ and $\alpha$ fixed, set
\[
(\wti{\bfa}^{\vee})_{J^\vee,G^\vee}^{\alpha}=(\wti{\bfa}_{J,G}^{\alpha})^\vee.
\]
Equation \eqref{eq: comparison family pullback} identifies the reflected families; reflection preserves the comparison conditions and admissibility.

Reflection commutes with disjoint union and the $\Ga$-maps. Thus the original $\Ga$-diagrams and $\PBI_a^{\partial,V,\com}$, $\PBI_a^{\partial,V,\com,\adm}$ take values in $\sSpaceo^{\sd}$ as in Construction \ref{cons: special dual}. The same reflection applies to $\PBI_a^{\partial,V,\aff}$, $\PBI^{\partial,V,\aff}$, and $\PBI^{\partial,V,\kih}$; vertex reversal in the Kihara model is smooth by \cite[Proposition 2.5]{Kihara2024smooth}. Stabilization and forgetful maps preserve these dualities, which exchange the $\partial_\pm$ variants.
\econs

\bprop\label{prop: prove kan perturb}
For every $k\geq0$, the morphisms of $\Ga$-diagrams
\[
\begin{aligned}
p_{a,k}^{\partial,V,\adm}&:(\PBI_a^{\partial,V,\adm})_{k,\bu}\longrightarrow(\PBI_a^{\partial,V,\aff})_{k,\bu},\\
p_{a,k}^{\partial,V,\com,\adm}&:(\PBI_a^{\partial,V,\com,\adm})_{k,\bu}\longrightarrow(\PBI_a^{\partial,V,\adm})_{k,\bu}
\end{aligned}
\]
are objectwise trivial Kan fibrations in the sense of Definition \ref{defn: trivial kan fibration}. Hence both forgetful maps are objectwise levelwise weak homotopy equivalences. The same assertions hold for the $\partial_\pm$ variants. These maps commute with the reflections of Construction \ref{cons: reflected model families}.
\eprop
\bpf
Fix $m,k,l\in\N$ and a lifting problem over $\partial\Delta^l\hookrightarrow\Delta^l$ in the $\langle m\rangle$-term. For $p_a^{\partial,V,\adm}$, extend the given smoothing collections componentwise by Proposition \ref{prop: trivial kan smoothing}. When $m=0$, there are no perturbation conditions. For $m>0$, choose admissible primary perturbations at the remaining vertices.

Construct the secondary perturbations $\Pi_I$ by induction on $|I|$, fixing the data on the spatial faces. Shorter chains determine compatible data on $\partial Q_I$ through \eqref{eq: abstract auxiliary face}. Extend these data smoothly on each component, with the required product form near infinite faces. For each nonempty $S\subseteq\{1,\ldots,m\}$, the relative Sard--Smale argument of Proposition \ref{prop: extension of perturbation} gives a residual regularity condition on the allowable componentwise perturbations for the subunion indexed by $S$ over the common $Q_I$. The finitely many conditions can be satisfied simultaneously.

On an infinite face, the subunion on each broken piece is regular by induction. Different broken pieces have independent parameters, so their product is regular. Repeating the construction for all weak chains gives the required lift for $p_a^{\partial,V,\adm}$.

For $p_a^{\partial,V,\com,\adm}$, fix the original collection and the comparison data over $\partial\Delta^l$. Apply the same relative argument to the comparison cells, treating faces and product factors first. Their data are admissible by the subunion and product conditions and admissibility after deleting unused parameters. Substitution reparametrizes whole factor cells as in Construction \ref{cons: comparison families}. On the remaining interiors, perturb the unbroken pieces independently to obtain the regularity and transversality of Construction \ref{cons: ga diagram}. The countable induction respects all identifications and pullbacks and gives the required lift.

Thus both maps have the right lifting property for every $\partial\Delta^l\hookrightarrow\Delta^l$, including $l=0$. The same argument applies to the $\partial_\pm$ variants. The weak homotopy equivalences follow from Lemma \ref{lem: trivial Kan fibration}, and compatibility with reflection follows from the definitions.
\epf

\bcor\label{cor: segal space PBIVa}
The simplicial spaces $(\PBI_a^{\partial,V,\adm})_{\bu,\bu}$ and $(\PBI_a^{\partial,V,\com,\adm})_{\bu,\bu}$ satisfy the Segal condition. The same holds for the $\partial_\pm$ variants.
\ecor
\bpf
Proposition \ref{prop: prove kan perturb} and \eqref{eq: Fg V} give levelwise weak homotopy equivalences
\[
(\PBI_a^{\partial,V,\com,\adm})_{\bu,\bu}\longrightarrow(\PBI_a^{\partial,V,\adm})_{\bu,\bu}\longrightarrow(\PBI^{\partial,V,\aff})_{\bu,\bu}.
\]
The result follows from Proposition \ref{prop: segal space PBIV} and homotopy invariance of the Segal condition by Lemmas \ref{lem: homotopy invariance fiber product} and \ref{lem: homotopy properties}.
\epf

\bdefn\label{defn: BIa}
Set
\[
\begin{aligned}
\PBI_a^{\partial,\adm}&=\hocolim_{V\subset\R^\infty}(\PBI_a^{\partial,V,\adm})_{\bu,\bu},\\
\PBI_a^{\partial,\com,\adm}&=\hocolim_{V\subset\R^\infty}(\PBI_a^{\partial,V,\com,\adm})_{\bu,\bu},
\end{aligned}
\]
taking the homotopy colimits objectwise in $\sSpace_c^{\sd}$. Using the replacement of Construction \ref{cons: special dual}, define, objectwise on the $\Ga$-diagram,
\[
\begin{aligned}
\BI^{\partial,\adm}_a&=\widehat R^{\sd}(\PBI_a^{\partial,\adm}),\\
\BIa^\partial=\BI_a^{\partial,\com,\adm}&=\widehat R^{\sd}(\PBI_a^{\partial,\com,\adm}).
\end{aligned}
\]
The latter is called the \emph{$\infty$-cobordism category for instantons with auxiliary data}. The $\Ga$-maps and strict dualities are induced functorially by $\widehat R^{\sd}$. Make the same construction for the $\partial_\pm$ variants and their inclusions, with reflection interchanging the signs.
\edefn

\bprop\label{prop: BIa=BI}
The completions in Definition \ref{defn: BIa} are $\Ga$-objects. The forgetful maps induce duality-preserving symmetric monoidal equivalences
\[
\BIa^\partial\xra{\simeq}\BI^{\partial,\adm}_a\xra{\simeq}\BI^\partial.
\]
The same holds for the $\partial_\pm$ variants, with duality switching the signs in the sense of Definition \ref{defn: dual}.
\eprop
\bpf
Proposition \ref{prop: prove kan perturb}, \eqref{eq: Fg V}, and Lemma \ref{lem: affine Kihara comparison}, after taking filtered homotopy colimits, give objectwise levelwise weak homotopy equivalences
\[
\PBI_a^{\partial,\com,\adm}\longrightarrow\PBI_a^{\partial,\adm}\longrightarrow\PBI^{\partial,\aff}\longrightarrow\PBI^{\partial,\kih}.
\]
These are morphisms of $\Ga$-diagrams preserving strict duality by Construction \ref{cons: reflected model families}. Thus $\PBI_a^{\partial,\adm}$ and $\PBI_a^{\partial,\com,\adm}$ satisfy the $\Ga$-condition by that for $\PBI^{\partial,\aff}$ from Definition \ref{defn: BI} and the displayed equivalences. Applying $\widehat R^{\sd}$ gives the asserted equivalences, computing the completion in Definition \ref{defn: BI} with the same replacement. Completion preserves the $\Ga$-condition by \cite[Lemma 3.7]{CS2019category}. The same argument applies simultaneously to the $\partial_\pm$ variants.
\epf

\section{Construction of the TQFT functor}\label{sec: Construction of the TQFT functor}
In this section, we construct the TQFT functor for instantons. We expect the reader to be familiar with the original construction in gauge theory, and only sketch it in \S \ref{subsect: Sketch of gauge theory construction}. The main theorem of this section (and this paper) is stated in \S \ref{subsec: Morphisms and higher morphisms}, where we mainly focus on the $\infty$-category construction.

\subsection{Gauge theory construction}\label{subsect: Sketch of gauge theory construction}
In this subsection, we sketch the gauge-theoretic construction of $CI$. We use the abbreviated vertex notation $(\bY=(Y,P,\phi,J),a)$ of Definition \ref{defn: BIa}, where $a=(g,\pi)$ consists of a metric and an admissible primary perturbation. A cobordism with auxiliary data, including a metric and admissible primary and secondary perturbations, is written as\[(\bW=(W,\bfP,\Phi,\bfJ),\bfa):(\bY_-=(Y_-,P_-,\phi_-,J_-),a_-)\to (\bY_+=(Y_+,P_+,\phi_+,J_+),a_+).\]

We start with the construction of $CI(\bY,a)$. Let $SO(3)=SU(2)/\{\pm1\}$ act on $SU(2)$ by inner automorphisms. Let $G(P)\to Y$ be the associated bundle of groups with fiber $SU(2)$, and let $\cG(P)$ be the group of smooth sections of $G(P)$, called \emph{determinant-$1$ gauge transformations} of $P$ in \cite[\S 2.5]{kronheimer2011khovanov}. Then there is a short exact sequence \[1\to\{\pm1\}\to G(P)\to\Aut(P)\to1\]of group bundles over $Y$. Taking global sections gives an exact sequence
$$1\to H^0(Y;\Z/2)\to\cG(P)\to\cAut(P)\to H^1(Y;\Z/2)\to 1.$$

Fix a sufficiently large integer $k$; indeed, $k\ge3$ suffices. Let $\cC_k(P)$ be the space of Sobolev $L^2_{k}$ connections on $P$, and let $\cG_{k+1}(P)$ be the group of Sobolev $L^2_{k+1}$ sections of $G(P)$. Let \[\cB_k(\bY)=\cC_k(P)/\cG_{k+1}(P),\] called the \emph{configuration space}. Then $\cAut(P)$ acts on $\cB_k(\bY)$ with $\cG(P)$ acting trivially, inducing an action of $H^1(Y;\Z/2)$ on $\cB_k(\bY)$. Let $\ov{\cB}_k(\bY)=\cB_k(\bY)/\phi$ be the quotient of $\cB_k(\bY)$ by the subgroup $\phi$ of $H^1(Y;\Z/2)$ under this action. We often omit the subscript $k$ and simply write $\cB(\bY)$ for the configuration space.

The admissible perturbation $\pi$ in $a$ determines an $\operatorname{ad}(P)$-valued $1$-form $V_\pi(B)$ for each connection $B$ representing a point of $\ov{\cB}(\bY)$ and is chosen so that there are only finitely many gauge-equivalence classes of perturbed flat connections, i.e.\ solutions of $*F_B+V_\pi(B)=0$. Then $CI(\bY,a)$ is the direct sum of the rank-one orientation modules associated to these gauge-equivalence classes, using the construction in \cite[\S 3.6 and Formula (84)]{kronheimer2011knot} and \cite[\S 5.1]{kronheimer2011khovanov} based on the almost complex structure $J$. After choosing a generator of each orientation module, it is freely generated over $\Z$ by these gauge-equivalence classes. There is also a canonical $\Z/2$-grading on generators from $J$.

More explicitly, for a critical point $x$, let $F_x$ be the Fredholm operator on $\R\times Y$ interpolating from the complex-linear reference operator determined by $J$ at $-\infty$ to the perturbed ASD operator at $+\infty$. Its determinant line $\det F_x$ defines the orientation module $\ov{\Lambda}(x)$ of $x$. Give this module degree
\begin{equation}\label{eq: J grading}
|x|_J=\operatorname{ind}_{\R}F_x\pmod2.
\end{equation}
We call the resulting $\Z/2$-grading the \emph{$J$-grading}. Let $|x|_{\mathrm{std}}$ denote the standard grading of \cite[\S 2.6]{kronheimer2010instanton}. Suppose that $(W,\bfJ)$ is an almost complex filling of $Y$, with $\bfJ$ agreeing with $J$ on the boundary collar, and that $P$ extends over $W$. Capping and index additivity give
\begin{equation}\label{eq: difference}
|x|_{\mathrm{std}}-|x|_J=\frac{\chi(W)+\sigma(W)+b_1(Y)-b_0(Y)}{2}\pmod2.
\end{equation}
Indeed, the standard grading is the index of the capping ASD operator plus the right-hand side. Gluing $F_x$ to the complex-linear reference cap shows that this index is congruent to $|x|_J$, since the reference cap has even real index. Thus the difference is independent of $x$ and of the compatible filling.

We will omit the subscript $J$ if it is clear. This construction applies to the quotient by $\phi$ used above. We regard the orientation modules as graded lines and use the sign conventions
\begin{equation}\label{eq: koszul convention}
\begin{aligned}
x\otimes y&\longmapsto(-1)^{|x||y|}y\otimes x,\\
(f\otimes g)(x\otimes y)&=(-1)^{|g||x|}f(x)\otimes g(y),
\end{aligned}
\end{equation}
for homogeneous elements $x,y$ and homogeneous maps $f,g$.

For the differential in $CI(\bY,a)$, we consider the perturbed anti-self-dual (ASD) equation $F_A^++\widehat{V}_\pi(A)=0$ on $\R\times Y$. Here $A=B(r)+c(r)\,dr$, $F_A^+=(F_A+*F_A)/2$, and $\widehat{V}_\pi(A)=P^+(dr\wedge V_\pi(B))$, where $*$ is the $4$-dimensional Hodge star for the product metric and $P^+$ is the projection onto self-dual $2$-forms.

For two generators $x$ and $y$ of $CI(\bY,a)$, we write $M(x,y)$ for the moduli space of the perturbed ASD equation on the configuration space $\ov{\cB}_{\loc}(\R\times Y)=\cB_{\loc}(\R\times Y)/\phi$ with asymptotic limits $x$ and $y$, which is a stratified space. We write $M_d(x,y)$ for the $d$-dimensional stratum of $M(x,y)$. Translation in the $r$-direction defines an $\R$-action on $M(x,y)$, free away from constant trajectories. Following \cite[\S 3.5]{kronheimer2011khovanov}, let $\breve M(x,y)$ denote the quotient of the nonconstant trajectories and $\breve M_d(x,y)$ its $d$-dimensional stratum. For the choice of the admissible perturbation $\pi$, these moduli spaces are regular and carry the determinant-line orientation data determined by $J$. Write $\breve M_d^+(x,y)$ for the completion of $\breve M_d(x,y)$ by broken trajectories, following the notation in \cite[\S 6.2]{kronheimer2011khovanov}. Since taking the quotient lowers the dimension by one, $\breve M_{d-1}^+(x,y)$ is compact when $d<8$ (or $d<4$ when the singular set is nonempty). In higher dimensions, bubbling strata may also occur. For all $d\in\N$, we write $\ov{M}_d(x,y)$ for the full compactification of $\breve M_d(x,y)$, including bubbles when necessary. Let the coefficient $\langle \partial x,y\rangle$ be the signed count of the $0$-dimensional compact manifold $\ov{M}_0(x,y)$, where the orientation uses these determinant-line data and the standard orientation of $\R$. The result $\partial^2=0$ comes from the counting of the boundaries of the $1$-dimensional compact manifold $\ov{M}_1(x,y)$.

For the construction of $CI(\bW,\bfa)$, we consider the more general setting of a cobordism with family auxiliary data, called a \emph{family cobordism} and written as \[(\bW_G,\bfa):(\bY_-,a_-)\to (\bY_+,a_+)\]when the family is parametrized by a compact oriented smooth manifold with corners $G$. In the case of higher morphisms in $\BIa$, the auxiliary data $\bfa$ consist of the model families of metrics and perturbations over the cubes $Q_I$ from \S \ref{subsec: perturbation}. The end data are fixed throughout the family.

We extend the structures to the open manifold $W^*$ in \eqref{eq: open end} and write $\bW^*$ for the resulting tuple. Let $M_d(x_-,\bW_G,x_+)$ be the $d$-dimensional stratum of the moduli space of pairs $(q,A)$ consisting of $q\in G$ and a solution $A$ of the perturbed ASD equation on the corresponding fiber $\bW_q^*$, with respect to the restriction of $\bfa$ at $q$ and with boundary conditions given by generators $x_\pm$ in $CI(\bY_\pm,a_\pm)$. Here $d$ includes the parameter directions. In a local trivialization over $U\subseteq G$, a pair $(q,A)$ with $q\in U$ lives in $U\times\ov{\cB}_{\loc}(\bW^*)=U\times\cB_{\loc}(\bW^*)/\Phi$. We also write $M(x_-,\bW_G,x_+)$ as the union of all moduli spaces, which is again a stratified space. The choices of $\bfJ$ and admissible perturbations in $\bfa$ make the unbroken strata regular and oriented. They can be compactified by adding broken trajectories when $d<8$ (or $d<4$ when the singular set is nonempty). In higher dimensions, bubbling strata may also occur. We write $\ov M_d(x_-,\bW_G,x_+)$ for the resulting compactified moduli space for all $d\in\N$. The map $CI(\bW_G,\bfa)$ is obtained by the signed count of $\ov M_0(x_-,\bW_G,x_+)$. To show it is a chain map, a homotopy, or a higher homotopy, we use the signed count of the boundaries of $\ov M_1(x_-,\bW_G,x_+)$.

We also allow metrics broken on parameter faces or throughout $G$. With matching critical points fixed, the parametrized moduli spaces are fiber products of those on the unbroken pieces over the parameter stratum. We require regularity of the component moduli spaces and transversality of these fiber products on every stratum. On faces with the product data of Definition \ref{defn: cubical auxiliary data}, this gives the corresponding identifications, summed over matching critical points. Near a breaking face, use \eqref{eq: perturbed ASD 3} with the component equations at the face. The subgroup condition in \S \ref{subsec: basic setup} gives the hypothesis of \cite[Proposition 5.4]{kronheimer2011khovanov}. The oriented boundary argument then applies, with contributions at matched internal ends cancelling.

Write $m_G=CI(\bW_G,\bfa)$ for the family count defined above. The boundary and broken-product formulas in \cite[\S 3.9, pp.\ 138--139]{kronheimer2011khovanov} give
\begin{equation}\label{eq: J instanton family identities}
\begin{aligned}
\partial_+m_G-(-1)^{\dim G}m_G\partial_-&=m_{\partial G},\\
m_{G_+\times G_-}&=m_{G_+}\circ m_{G_-}.
\end{aligned}
\end{equation}
The second formula uses the outgoing-to-incoming order of the parameter factors.

We now determine the $\Z/2$-degree of $m_G$. In the standard grading of \cite[\S 2.6]{kronheimer2010instanton}, cobordism maps can carry a topological grading shift. Here we instead use the $J$-grading of \eqref{eq: J grading}.

Let $D_A$ be the linearized operator at a solution with limits $x_-$ and $x_+$. Gluing the operators defining \eqref{eq: J grading} at its ends replaces the asymptotic operators by their complex-linear references. The compatible almost complex structure on the cobordism gives a Fredholm homotopy from the resulting operator to a complex-linear operator, whose real index is even. Index additivity therefore gives
\begin{equation}\label{eq: J index parity}
\operatorname{ind}_{\R}D_A=|x_-|_{J_-}-|x_+|_{J_+}\pmod2.
\end{equation}
A solution counted by $m_G$ satisfies $\operatorname{ind}_{\R}D_A+\dim G=0$, so
\begin{equation}\label{eq: J family degree}
\deg m_G=\dim G\pmod2.
\end{equation}
In particular, cobordism maps preserve the $J$-grading. The differential is odd because its trajectories have index $1$ before quotienting by translation. Together with \eqref{eq: J instanton family identities}, this gives the degree and boundary identities required in Proposition \ref{prop: abstract auxiliary theory}.

Finally, we summarize the results for disjoint unions and duals in the following propositions.
\bprop[Disjoint unions]\label{prop: J disjoint union complexes}
There are canonical isomorphisms of $2$-periodic chain complexes
\begin{equation}\label{eq: J disjoint union complexes}
\iota:CI\Bigl(\bigsqcup_{\alpha=1}^m(\bY^\alpha,a^\alpha)\Bigr)\xra{\cong}\bigotimes_{\alpha=1}^m CI(\bY^\alpha,a^\alpha),
\end{equation}
where the subgroup on the disjoint union is the direct sum of the component subgroups. Set $CI(\varnothing)=\Z_{(0)}$ in even degree. These isomorphisms are associative and unital, with symmetry given by \eqref{eq: koszul convention}.

For $\alpha\in[1,m]$, let $\bfa^\alpha$ be admissible families over independent compact oriented parameter manifolds with corners $G_\alpha$, possibly broken throughout $G_\alpha$ as above, and let $\iota_\pm$ be the identifications at their ends. With the product orientation ordered by $\alpha$ and the corresponding determinant-line orientations,
\begin{equation}\label{eq: J independent parameter product}
\iota_+\circ m_{\prod_\alpha G_\alpha}\Bigl(\bigsqcup_\alpha\pr_\alpha^*\bfa^\alpha\Bigr)=\Bigl(\bigotimes_\alpha m_{G_\alpha}(\bfa^\alpha)\Bigr)\circ\iota_-.
\end{equation}
This identity is compatible with parameter faces and broken-family composition.
\eprop
\bpf
The product decomposition of the quotient configuration spaces identifies generators and orientation modules with graded tensor products. The trajectory count gives the tensor-product differential, as in \cite[\S 5.2]{kronheimer2011khovanov}.

For independent parameters, the parametrized moduli spaces are products, also for broken families after fixing matching critical points and taking the transverse fiber products within each factor. Their ordered determinant-line and parameter orientations give \eqref{eq: J independent parameter product} with the sign in \eqref{eq: koszul convention}. Compatibility with parameter faces follows from the boundary Leibniz rule, and compatibility with broken-family composition follows from graded interchange.
\epf

\bprop[Duality]\label{prop: J duality}
Use $a^\vee$ and $\bfa^\vee$ from Construction \ref{cons: reflected model families}. There is a canonical chain isomorphism
\begin{equation}\label{eq: J dual complex}
\xi_{(\bY,a)}:CI(\bY^\vee,a^\vee)\xra{\cong}CI(\bY,a)^\vee.
\end{equation}
For an admissible family $(\bW_G,\bfa)$, with the same oriented parameter space $G$ for its dual, we have
\begin{equation}\label{eq: J dual family}
\xi_{\bY_-,a_-}\circ m_G(\bW_G^\vee,\bfa^\vee)
=m_G(\bW_G,\bfa)^\vee\circ\xi_{\bY_+,a_+},
\end{equation}
where we use the graded dual of the family map. From the identification $(\bY^{\vee\vee},a^{\vee\vee})=(\bY,a)$, we have
\begin{equation}\label{eq: J bidual compatibility}
\xi_{(\bY^\vee,a^\vee)}=\xi_{(\bY,a)}^\vee\circ\eta_{CI(\bY,a)},
\end{equation}where $\eta_{CI(\bY,a)}$ is the graded bidual isomorphism from Definition \ref{defn: dual on D}. These identities are natural under family identifications and compatible with faces, broken products, and the tensor identifications of Proposition \ref{prop: J disjoint union complexes}, including all iterated operations.
\eprop
\bpf
Exchanging the ends identifies the perturbed equations and their deformation operators, as in \eqref{eq: moduli dual}. For each critical point $x$, identify its orientation module on $\bY^\vee$ with the dual of that on $\bY$ by reversing $F_x$ and then reflecting the cylinder. Repeating the graded determinant-line interchange of \cite[Lemma 22.5.4]{kronheimer2007monopoles} with the complex-oriented $J$-reference gives the sign $(-1)^{|y|(1+|x|)}$ before quotienting a trajectory $x\to y$ by translation. Reversal of the translation orientation changes this to $(-1)^{|y|(1+|x|)+1}$, as in \cite[Corollary 22.5.5]{kronheimer2007monopoles}. For a trajectory contributing to the differential, $|y|=|x|+1\pmod2$, so the sign is $(-1)^{|x|}$. This is the graded dual differential and proves \eqref{eq: J dual complex}.

For a solution counted by a $G$-family, no translation quotient is taken. By \eqref{eq: J index parity}, $|x_-|+|x_+|=\dim G\pmod2$, so the same interchange sign is $(-1)^{|x_+|(1+|x_-|)}=(-1)^{(\dim G)|x_+|}$. Keeping the orientation of $G$, we obtain
\[
m_G^\vee(\lambda)=(-1)^{(\dim G)|\lambda|}\lambda\circ m_G,
\]
proving \eqref{eq: J dual family}.

Compatibility with parameter faces follows from $d(f^\vee)=(df)^\vee$. For broken products, the sign in
\[
(g\circ f)^\vee=(-1)^{|f||g|}f^\vee\circ g^\vee
\]
agrees with the product-orientation sign from reversing the component parameter spaces, since $\deg m_G=\dim G$ modulo $2$.

The canonical graded determinant-line identifications respect tensor products and give \eqref{eq: J bidual compatibility} upon dualizing twice. Their naturality, associativity, and graded symmetry give all iterated compatibilities.
\epf

The following trace formulas are not needed for the functor construction. We include them because self-gluing is used in Floer excision; see \cite[Proposition 5.5]{kronheimer2011khovanov} and the proof of Theorem 5.6 in this reference.
\bprop[Trace, cotrace, and self-gluing]\label{prop: J trace self gluing}
Suppose $\phi=0$ and put $C=CI(\bY,a)$. For product auxiliary data, the trace and cotrace cylinders with zero subgroup induce the alternating evaluation and coevaluation of Definition \ref{defn: dual on D}
\begin{equation}\label{eq: J evaluation coevaluation}
\begin{aligned}
\ev &: C\otimes C^\vee\longrightarrow\Z,
& \ev(x\otimes y^\vee)&=(-1)^{|x|}\delta_{xy},\\
\coev &: \Z\longrightarrow C^\vee\otimes C,
& \coev(1)&=\sum_x(-1)^{|x|}x^\vee\otimes x.
\end{aligned}
\end{equation}
Here $x,y$ are homogeneous basis elements. These maps satisfy \eqref{eq: dual comp}.

For self-gluing, allow $\phi\ne0$. Let $\bW^{\bY}$ be obtained from $\bW:\bY_-\sqcup\bY\to\bY_+\sqcup\bY$ by identifying the two copies of a connected $\bY$, with matching additional data. Write $q:W\to W^{\bY}$ for the quotient map. Suppose the boundary subgroups are $\phi_\pm\oplus\phi$, respectively. Choose
\begin{equation}\label{eq: self gluing subgroup}
\Phi^{\bY}\subset H^1(W^{\bY};\Z/2),~
\operatorname{PD}[Y]\in\Phi^{\bY},~
\Phi=q^*\Phi^{\bY}.
\end{equation}
Choose admissible auxiliary data $\bfa$ on $\bW$ inducing the same datum $a$ on the two $\bY$. Let $\bfa^{\bY}$ be admissible auxiliary data obtained by gluing $\bfa$ along the two copies of $\bY$. Whenever $\bW^{\bY}$ satisfies the positive, negative, or the full boundary condition, there is a chain homotopy
\begin{equation}\label{eq: J self gluing trace}
CI(\bW^{\bY},\bfa^{\bY})\simeq\tr_C CI(\bW,\bfa),
\end{equation}
where $\tr_C$ is the partial alternating trace over $C=CI(\bY,a)$. The same statement holds for several disjoint pairs if the boundary subgroups split into the corresponding factors and $\Phi^{\bY}$ contains every $\operatorname{PD}[Y_\alpha]$.
\eprop
\bpf
For product cylinders, only constant trajectories contribute. The graded determinant-line identifications give \eqref{eq: J evaluation coevaluation}, whose signs cancel in \eqref{eq: dual comp}.

For self-gluing, stretch the neck along $Y$ as in Remark \ref{rem: generalization of broken}, starting at $\bfa^{\bY}$ and with broken limit $\bfa$. Choose admissible perturbations for this family $G=[0,\infty]$ relative to its endpoint data, using the argument of Proposition \ref{prop: extension of perturbation}. By the gluing argument of \cite[Proposition 5.5]{kronheimer2011khovanov} and the subgroup condition, the contribution at $\infty$ is $\tr_C CI(\bW,\bfa)$, with the alternating sign explained in \cite[Lemma 2.4]{kronheimer2010knots}. Set $H=m_G$. The same oriented-boundary argument as in \eqref{eq: J instanton family identities} gives
\[
\tr_C CI(\bW,\bfa)-CI(\bW^{\bY},\bfa^{\bY})
=\partial_+H+H\partial_-,
\]
which proves \eqref{eq: J self gluing trace}. For several pairs, apply the same argument successively; the pullbacks of $\Phi^{\bY}$ satisfy the subgroup hypothesis at each step.
\epf
\subsection{Construction of the functor}\label{subsec: Morphisms and higher morphisms}
In this section, we prove the main theorem of the paper.
\bthm\label{thm: TQFT functor}
The instanton family counts determine a duality-preserving symmetric monoidal functor
\[
CI_a:\BIa^{\partial}\longrightarrow\Psi(\D),
\]where $\D=\D_2(\Mod(\Z)^\fin)$ from Example \ref{exmp: finite module} and $\Psi(-)$ from \eqref{eq: Reedy fiber replacement}. The same construction applies to the $\partial_\pm$ variants, compatibly with the inclusions from $\BIa^{\partial}$, with duality switching the signs in the sense of Definition \ref{defn: dual}.
\ethm
\bpf
We verify the hypotheses of Propositions \ref{prop: abstract auxiliary theory} and \ref{prop: symmetric monoidal auxiliary theory} for the instanton family counts.

\textbf{Objects and family maps.}
For $(\sigma,\bfa)\in(\PBI_a^{\partial,V,\adm})_{k,l}$, set
\[
X_i=CI(\bY_i,a_i),~f_I=m_{Q_I}(\bfa_I),~I\in\cI_{l,k}.
\]
\S \ref{subsect: Sketch of gauge theory construction} gives finite free $2$-periodic complexes $X_i$ with the $J$-grading. By \eqref{eq: J family degree},
\[
\deg f_I=\dim Q_I=|I|-2\pmod2.
\]
Together with \eqref{eq: J instanton family identities} and the orientation and invariance conventions of \S\ref{subsect: Sketch of gauge theory construction}, this verifies Proposition \ref{prop: abstract auxiliary theory}(1).

\textbf{Simplicial and spatial maps.}
For $\rho=\alpha\times\beta$, the identifications in \eqref{eq: weak chain pullback} give
\[
m_{Q_J}\bigl((\rho^*\bfa)_J\bigr)=m_{Q_{\rho J}}(\bfa_{\rho J}).
\]
This holds also for noninjective operators, since repetitions in $\rho J$ are kept and the parameters of $Q_J=Q_{\rho J}$ are identified by position.

For $I=((z,w),(z',w))$ with $z\leq z'$, the underlying cobordism is a continuation cylinder, with the perturbations on the incoming and outgoing ends as in \eqref{eq: perturbed ASD 2}. Its map is a chain homotopy equivalence by \cite[\S3.8]{kronheimer2011khovanov}, hence an isomorphism in $\D$, with the subgroup and orientation conventions of \cite[\S5.1--5.2]{kronheimer2011khovanov}. This verifies Proposition \ref{prop: abstract auxiliary theory}(2), including when $I=(i,i)$. 

\textbf{Stabilization and completion.}
Proposition \ref{prop: abstract auxiliary theory} now gives morphisms of simplicial spaces
\begin{equation}\label{eq: final CI}
CI_a^V:\PBI_a^{\partial,V,\adm}\longrightarrow\Psi(\D).
\end{equation}
The family assignments are unchanged by enlarging $V$, so Proposition \ref{prop: abstract auxiliary theory} makes the maps $CI_a^V$ coherently compatible with these inclusions. Taking the filtered homotopy colimit and complete Segal space replacement as in Definition \ref{defn: BIa} and Remark \ref{rem: affine construction} gives
\[
\BIa^{\partial,\adm}\longrightarrow\Psi(\D).
\]
Precomposing with the forgetful equivalence $\BIa^\partial\to\BIa^{\partial,\adm}$ of Proposition \ref{prop: BIa=BI} gives the underlying functor in the statement. The construction applies simultaneously to the boundary variants and their inclusions.

\textbf{Symmetric monoidal structure.}
Put $\sfC=\Ch_2((\Mod(\Z)^\fin)_{\proj})$. Use the comparison collections of Construction \ref{cons: ga diagram}. Their admissibility allows the same counts and boundary argument on every comparison cell and its faces, including broken products. Regularity on the subdivision strata ensures that the solutions counted by $m_G$ lie off the internal subdivision faces, so the counts are additive under the prescribed subdivisions. Proposition \ref{prop: J disjoint union complexes} gives the vertex tensor identifications and the independent-product formula. If an interval factor is unused, the moduli space is its product with the regular reduced moduli space, so its zero-dimensional count vanishes. For a two-entry comparison chain with constant $[k]$-coordinate, the family is a continuation cylinder, whose map is a chain homotopy equivalence by the argument above. These verify Proposition \ref{prop: symmetric monoidal auxiliary theory}(1)--(4).

Consequently, the maps induced by $CI_a^V$ at each $\langle m\rangle\in\Ga$, pulled back along $p_a^{\partial,V,\com,\adm}$, extend to the homotopy coherent natural transformation of that proposition, with target $\Psi(\fB)$ from Construction \ref{cons: symm on D}.

\textbf{Duality.}
For each vertex $v=(\bY,a)$, set
\[
\widetilde X_v=(CI(\bY,a),CI(\bY^\vee,a^\vee),\xi_{(\bY,a)})\in\sfC^{\sd}.
\]
Equation \eqref{eq: J bidual compatibility} gives $\widetilde X_{v^\vee}=\widetilde X_v^{\vee_{\sd}}$. The maps $f_I$ also define morphisms in $\sfC^{\sd}$. For $x=(\sigma,\bfa)$ and $p=|I|-2$, \eqref{eq: J dual family} and the orientation of the coordinate reversal \eqref{eq: dual cube coordinates} give
\[
f_{I^\vee}(x^\vee)=(-1)^{p(p-1)/2}f_I(x)^{\vee_{\sd}}.
\]
This is the sign in \eqref{eq: sign reversal fI}; compatibility with both simplicial and spatial maps follows from \eqref{eq: reflected auxiliary operators}. The same calculation on every top-dimensional comparison cell gives \eqref{eq: comparison family duality}. The tensor identifications lift to $\sfC^{\sd}$ compatibly with duality by Proposition \ref{prop: J duality}.

The duality assertion of Proposition \ref{prop: symmetric monoidal auxiliary theory} therefore applies, with $B^{\sd}$ from Construction \ref{cons: strict tensor diagram}. Its coherent naturality under enlarging $V$ permits stabilization and completion as in Definition \ref{defn: BIa}. Proposition \ref{prop: BIa=BI} and Construction \ref{cons: strict tensor diagram} give the $\Ga$-conditions on the source and target. The duality-preserving symmetric monoidal equivalence of Example \ref{exmp: finite module} identifies the target at $\langle1\rangle$ with $\Psi(\D)$. This gives the asserted functor, also for the $\partial_\pm$ variants and their inclusions.
\epf

\section{Cube complexes and exact triangles}\label{sec: Cube complexes and exact triangles}
In this section, we reinterpret the results of \cite[\S 5--6]{kronheimer2011khovanov} and \cite[\S 6--7]{scaduto2015instantons} in the $\infty$-categorical setting. As in \S \ref{sec: Construction of the TQFT functor}, we write $\BIa$ and $\BI$ for the $\infty$-cobordism category for instantons with and without auxiliary data, respectively, and write \[\D=\D_2(\Mod(\Z)^\fin)\]for the corresponding $\infty$-derived category. Let \[CI:\BI^{\pa}\simeq \BIa^{\pa}\to \D\]be the instanton TQFT functor constructed in \S \ref{sec: Construction of the TQFT functor}, where the equivalence is given by Proposition \ref{prop: BIa=BI}.
\subsection{Cube complexes and hypercubes}\label{subsec: Cubes of chain complexes}
In this subsection, we introduce cubes of chain complexes and study their relation to hypercubes of chain complexes. Useful references are \cite[\S 1]{Goodwillie1992cube} and \cite[\S 2.1]{Varolgunes2021cube}.

\bdefn\label{defn: subcube}
    Given a finite set $S$, let $\bbox^S$ denote the power category $(\Delta^1)^S$. We identify its objects with subsets of $S$ using the bijection
    $$\{0,1\}^S\to\cP(S),$$
    $$(x_i)_{i\in S}\mapsto\{i\in S\mid x_i=1\}$$for the power set $\cP(S)$ of $S$.

    For an $\infty$-category $\cC$, an \emph{$S$-cube in $\cC$} is a functor \[\sX:\bbox^S\to\cC.\] An \emph{$S$-cube complex} is an $S$-cube in $\D$.    For $n\in\N$, an \emph{$n$-cube} is an $\un{n}$-cube where $\un{n}=\{1,2,\ldots,n\}$. As
    \[\Fun(\bbox^S,\Fun(\bbox^T,\cC))=\Fun(\bbox^{S\sqcup T},\cC),\]
    we can regard an $(S\sqcup T)$-cube in $\cC$ as an $S$-cube of $T$-cubes in $\cC$; in particular, we can regard an $(S\sqcup T)$-cube complex as an $S$-cube of $T$-cube complexes.

    We define the evident inclusion $\iota_S:\bbox^S\to \bbox^n$ associated to the map sending $S'\in\cP(S)$ to $S'\in\cP(\un{n})$. More generally, for any subset $T\subset \un{n}\bs S$, we define the inclusion $\iota_{S,T}:\bbox^S\to \bbox^n$ associated to the map sending $S'\in\cP(S)$ to $S'\cup T\in\cP(\un{n})$. Note that $\iota_{S,\varnothing}=\iota_S$ and there are $2^{|\un{n}\bs S|}$ choices of $T$.

    Then we define
    \begin{equation}\label{eq: subcube}
        \sX_S=\sX\circ\iota_S\aand \sX_{S,T}=\sX\circ\iota_{S,T}:\bbox^S\to\D
    \end{equation}and call them \emph{$S$-subcubes} of $\sX$.
\edefn

We use the following subscript notation for the values of an $n$-cube complex on objects.

\bdefn\label{defn: ep S}
    For $S\subseteq\un n$, let $\ep(S)\in\{0,1\}^n$ be its characteristic vector. Under the coordinatewise order, $S\subseteq S'$ if and only if $\ep(S)\le\ep(S')$, and strict inclusion corresponds to strict inequality.

    From \cite[\texttt{002Q}]{kerodon}, a $k$-simplex (in the quasicategory model) of $\bbox^n$ is a sequence of subsets \[S_0\subset\cdots\subset S_k\subset \un{n}.\]This $k$-simplex is nondegenerate if and only if $S_{g-1}\subsetneqq S_g$ for each $g\in\un{k}$.

\edefn

\bexmp\label{exmp: cube complex}
    A $0$-cube complex is just an object $\sX(\varnothing)$ in $\D$, i.e.\ a chain complex. A $1$-cube complex encodes two objects
    \[\sX(\varnothing)=(C_0,\pa_0),\,\sX(\{1\})=(C_1,\pa_1),\]
    and a $1$-morphism
    \[f_{0,1}=\sX(\varnothing\subset \{1\})\colon(C_0,\pa_0)\to (C_1,\pa_1)\]
    in $\D$, where the $1$-morphism is just a chain map. We refer to the images of the functor on objects as \emph{vertex complexes}, and to the images on $1$-morphisms as \emph{edge maps}.

    A $2$-cube complex consists of four vertex complexes $C_\ep$ for $\ep\in\{00,01,10,11\}=\{0,1\}^2$, five morphisms $f_{\ep,\epp}$ indexed by $\ep<\epp$ and two homotopies (i.e.\ $2$-morphisms) $h_{00,10,11}$ and $h_{00,01,11}$ relating $f_{00,11}$ to the two composites in the following diagram.
    \[\xymatrix@C=3pc@R=3.5pc{C_{01}\ar@(lu,ld)_{\pa_{01}}\ar[r]^{f_{01,11}}&C_{11}\ar@(rd,ru)_{\pa_{11}}\\
    C_{00}\ar@(lu,ld)_{\pa_{00}}\ar[r]_{f_{00,10}}\ar[u]^{f_{00,01}}_>>>>{h_{00,01,11}}\ar[ur]|-{f_{00,11}}&C_{10}\ar@(rd,ru)_{\pa_{10}}\ar[u]_{f_{10,11}}^<<<<{h_{00,10,11}}}.\]
\eexmp

We generalize the cofiber operation in $\D$ (i.e.\ the mapping cone) to cube complexes using the following definition; see \cite[Definition 1.4]{Goodwillie1992cube} for equivalent formulations.

\bdefn\label{defn: total cofiber}

 Let $\sX:\bbox^n\to\D$ be a cube complex. We define the \emph{total cofiber} $\Tcofib\sX\in\D$ as follows; the same definition applies with $\D$ replaced by any stable $\infty$-category. We also call $\Tcofib\sX$ the \emph{total complex} of $\sX$.

 Fix $p\in\un{n}$. From $\bbox^n= \bbox^{\{p\}}\times \bbox^{\un{n}\bs \{p\}}$, we regard $\sX$ as a $1$-cube of $(n-1)$-cube complexes
    \begin{equation}\label{eq: 1 cube n-1 cube}
        \sX^{\{p\}}:\bbox^{\{p\}}\to \Fun(\bbox^{\un{n}\bs \{p\}},\D)
    \end{equation}
    Then $\cofib(\sX^{\{p\}})$ is an $(n-1)$-cube complex. Given an order of elements in $\un{n}$, we apply the above procedure iteratively on all elements and finally obtain a $0$-cube complex, i.e.\ an object in $\D$. The choice of the order is not important. From this procedure, we obtain a functor
    \begin{equation*}\label{eq: total functor}
        \Tcofib:\Fun(\bbox^n,\D)\to\D.
    \end{equation*}

\edefn

We construct the cube filtration on $\Tcofib\sX$ of a cube complex $\sX$ as follows.

\bdefn\label{defn: cube filtration}
    Recall from \cite[Definition 1.2.2.9]{HAlurie2017}, a \emph{filtered object} of an $\infty$-category $\cC$ is a functor
    \[X\colon(\Z,\le)\to \cC,\]
    where we view the poset $(\Z,\le)$ as a category. The category of filtered objects of $\cC$ is denoted $\Fil(\cC)$. The finite cube filtration below can be modeled by a $\Z/2$-graded chain complex with a finite filtration by subcomplexes, using iterated mapping cones. It therefore has an associated spectral sequence $\{E_r,d_r\}_{r\ge0}$ (cf.\ \cite[\S 2.1]{BHL2019cancel}).

    Given an $S$-cube complex $\sX:\bbox^S\to\D$, we construct the \emph{cube filtration} on the total complex $\Fil_\bullet\Tcofib\sX\in\Fil(\D)$ as follows. For $k\in\Z$, let $\bbox^S_k$ be the full subcategory of $\bbox^S$ consisting of subsets of $S$ whose complement has at most $k$ elements. Then we define
    \begin{equation}\label{eq: def of fil}
        \Fil_k\Tcofib\sX=\begin{cases}
            0& \rmif k<0;\\
            \Tcofib (\sX|_{\bbox^S_k})& \rmif  k\ge0,
        \end{cases}
    \end{equation}where $\sX|_{\bbox^S_k}$ is obtained from $\sX$ by replacing all images of $\bbox^S\backslash \bbox^S_k$ by zeros. It is obvious that the spectral sequence associated to this cube filtration collapses at the $(|S|+1)$-st page.
\edefn
Finally, we describe the relation of the cube complex and the hypercube complexes in Definition \ref{defn: hypercube}.
\bcons\label{cons: cube complex to hypercube}
Let $\bK$ be a commutative Noetherian ring with $2=0$ and let $m\in\N_+$. An $n$-cube complex
\[
\sX:\bbox^n\longrightarrow\D_m(\Mod(\bK)^\fin)
\]
determines a hypercube complex over $\bK$ as follows.

For $\ep\in\{0,1\}^n$, let $C_\ep$ be the underlying graded module of $\sX(\ep)$ and set $D_{\ep,\ep}=\partial_{\sX(\ep)}$. For $\ep<\epp$ with $k=|\epp-\ep|$, set
\begin{equation}\label{eq: defn of D ep epp}
D_{\ep,\epp}=\sum_{\ep=\ep_0<\cdots<\ep_k=\epp}\sX(\ep_0<\cdots<\ep_k),
\end{equation}
where the sum runs over maximal chains in the face from $\ep$ to $\epp$, and each summand denotes the degree-$(k-1)$ component of the corresponding simplex in Definition \ref{defn: dg nerve}.

Summing the dg nerve identities cancels the common internal faces and gives
\[
\sum_{\ep\leq\de\leq\epp}D_{\de,\epp}\circ D_{\ep,\de}=0.
\]
Hence
\[
(C,D)=\left(\bigoplus_{\ep\in\{0,1\}^n}C_\ep[|\ep|-n],\sum_{\ep\leq\epp}D_{\ep,\epp}\right)
\]
is a hypercube complex in Definition \ref{defn: hypercube}. The iterated cone construction gives
\[
\Tcofib\sX\simeq(C,D).
\]
\econs
The signs for general coefficients are specified in the following Construction.
\bcons\label{cons: sign on hypercube}
When $2\neq0\in\bK$, one needs to be careful about the signs in \eqref{eq: compatibility condition}; in the periodic case, assume that $m$ is even. There are two sign conventions, discussed in \cite[\S 2.1]{Varolgunes2021cube}. Fix the standard order on $\un{n}$. For $\ep\leq\epp$, put $k=|\epp-\ep|$. When $k>0$, define the oriented version of \eqref{eq: defn of D ep epp} by
\[
D^{\mathrm{sign}}_{\ep,\epp}=\sum_{\ep=\ep_0<\cdots<\ep_k=\epp}\pm\sX(\ep_0<\cdots<\ep_k),
\]
where each sign compares the vertex orientation of the simplex with the standard orientation of the corresponding cube face, equivalently the sign of the permutation recording the order in which the coordinates change. Thus common internal faces cancel. Set $D^{\mathrm{sign}}_{\ep,\ep}=\partial_{\sX(\ep)}$.

The first convention is used in \cite[Definition 2.1.1]{Varolgunes2021cube}, where the corresponding hypercubes are called \emph{$n$-cubes}. The dg nerve identities imply that the maps $D^{\mathrm{sign}}_{\ep,\epp}$ form such an $n$-cube. For each $\ep\leq\epp$, put
\begin{equation}\label{eq: a ep}
  a(\ep,\epp)=\sum_{i<j}(1-\ep^{\prime}_i)(\ep^{\prime}_j-\ep_j).  
\end{equation}
The iterated cone construction of \cite[\S 2.1.4]{Varolgunes2021cube} gives
\[
D^{\mathrm{tot}}_{\ep,\epp}=(-1)^{a(\ep,\epp)+n-|\epp|}D^{\mathrm{sign}}_{\ep,\epp},~
D^{\mathrm{tot}}=\sum_{\ep\leq\epp}D^{\mathrm{tot}}_{\ep,\epp}.
\]
With the grading convention of Definition \ref{defn: hypercube},
\[
C=\bigoplus_{\ep\in\{0,1\}^n}C_\ep[|\ep|-n],~
\Tcofib\sX\simeq(C,D^{\mathrm{tot}}).
\]
In particular, $D^{\mathrm{tot}}_{\ep,\ep}=(-1)^{n-|\ep|}\partial_{\sX(\ep)}$, as required by the vertex shifts.

The second convention appears in \cite[Definition 6.6]{kronheimer2011khovanov} and \cite[\S 6.2]{scaduto2015instantons}, and corresponds to \emph{$n$-cubes with positive signs} in \cite[\S 2.1.3]{Varolgunes2021cube}. The sign change from the first convention is multiplication by $(-1)^{a(\ep,\epp)+n-|\epp|}$; its exponent is denoted by $n(F)$ there, where $F$ is the face from $\ep$ to $\epp$. Hence
\[
\sum_{\ep\leq\de\leq\epp}D^{\mathrm{tot}}_{\de,\epp}\circ D^{\mathrm{tot}}_{\ep,\de}=0.
\]
The inverse sign change uses the same signs.

For comparison with the signs in the instanton construction of \cite{kronheimer2011khovanov,scaduto2015instantons}, see Construction \ref{cons: hypercube to cube complex}.

\econs

\subsection{Cobordisms supported in disjoint regions}

In this subsection, we consider several cobordisms with the same source but supported in disjoint regions and use them to construct a cube of cobordisms. We start with the following definition.

\bdefn\label{defn: supported region}
From \S \ref{subsec: basic setup}, we write an object in $\BI$ by $\bY=(Y,P,\phi,J)$ and a morphism by \[\bW=(W,\bfP,\Phi,\bfJ):\bY\to \bY_1.\]We also call such a morphism a \emph{cobordism}. Suppose $\sR$ is an open $3$-dimensional submanifold of $Y$ such that the pullback of $\phi$ on $\sR$ is trivial, which is called a \emph{cobordism region}, or simply \emph{region}. A cobordism $\bW$ is called \emph{supported} in $\sR$ if it satisfies the following.

\begin{itemize}
    \item $W$ is obtained from $I\times Y$ by attaching $4$-dimensional handles supported in $\{1\}\times \sR$;
    \item $\bfP$ coincides with $I\times P$ outside $\sR$;
    \item The pullback of $\Phi$ on $\sR$ is trivial, which is induced from $\phi$ on the product $I\times (Y\backslash\sR)$;
    \item $\bfJ$ coincides with $J$ outside $I\times \sR$.
    \item There exists a closed oriented surface $\Sigma$ in $Y\backslash\sR$ satisfying \eqref{eq: witness}, which is a witness surface for the non-integral condition.
\end{itemize}
If for $i\in \un{n}=\{1,\dots,n\}$, there are cobordisms\[\bW_i:\bY\to \bY_i\]supported in $\sR_i$ such that all regions are disjoint, i.e.\ $\sR_i\cap \sR_j=\varnothing$ for $i\neq j\in \un{n}$, and there exists a witness surface for non-integral condition in $Y\backslash\bigcup_{i=1}^n\sR_i$, then we write\[(\bY,\{\bW_i\}_{i\in\un{n}})\]and call them \emph{cobordisms supported in disjoint regions}. We require every $3$-manifold obtained from a subcollection of the handle attachments to satisfy the non-integral condition with its induced data. All induced cobordisms must have full boundary.

\edefn

\bexmp\label{exmp: surgery cube}
Suppose $\bY$ is an object in $\BI$. Suppose the almost complex structure $J$ in $\bY$ is determined by an overtwisted contact structure $\xi$ (i.e.\ the unique almost complex structure up to isotopy compatible with the induced symplectic structure on $\R\times Y$ as in \cite[Remark 4.1.19]{OzbS2004contactsurgery}). Note that the overtwisted contact always exists by Eliashberg's classification \cite{Eliashberg1989overtwisted}.

Suppose $L=\bigcup_{i=1}^nK_i$ is a framed link. Choose disjoint neighborhoods $N(K_i)$ containing a meridian $m_i$ of $K_i$ and a meridian $m_i'$ of $m_i$, both with surgery coefficient $-1$ relative to their disk framings in $Y$. Choose an overtwisted disk in each component of $Y$, isotope the entire link $\bigcup_i(K_i\cup m_i\cup m_i')$ away from these disks, and take a Legendrian approximation disjoint from them. Adjust the contact framings componentwise as in the proof of \cite[Lemma 11]{GayKirby2004}, so that $0$-surgery on each framed $K_i$ and the indicated $(-1)$-surgeries are all contact $(-1)$-surgeries. Transport the neighborhoods along the ambient isotopy and set $\sR_i=\inr N(K_i)$. Suppose $\phi$ restricts trivially to each $\sR_i$ and $Y\backslash\bigcup_i\sR_i$ contains a witness surface for the non-integral condition.

For $i\in\un{n}$, suppose $\bW_i:\bY\to \bY_i$ is obtained from $\bY$ by $0$-surgery on the framed knot $K_i$, with bundle $\bfP_i$ from \cite[\S 3.3]{scaduto2015instantons}, subgroup $\Phi_i$ induced from $\phi$, and the almost complex structure $\bfJ$ from the Legendrian surgery on $K_i$ as in \cite[Theorem 7.2.4]{OzbS2004contactsurgery}. These 2-handle attachments preserve the connected components, so the same witness surface establishes non-integrality at every vertex. The induced cobordisms have full boundary. Thus the $\bW_i$ are supported in disjoint regions. We call \[(\bY,\{\bW_i\}_{i\in\un{n}})\]\emph{surgery cobordisms} induced by $(\bY,L)$, which is the special case of cobordisms supported in disjoint regions. Sometimes we will simply use $(\bY,L)$ to denote the surgery cobordisms.

With these choices, Legendrian surgeries along $K_i$ and $m_i$ give an almost complex cobordism from $Y$ to $Y_1(K_i)$, where $Y_1(K_i)$ is obtained by $1$-surgery on $K_i$ with its original framing. The contact framing of $m_i'$ is unchanged by these surgeries. Further Legendrian surgery along $m_i'$ gives a $3$-manifold canonically diffeomorphic to $Y$, but the induced almost complex structure differs from the original one.

Note that the above construction of almost contact structure from the overtwisted contact structure and Legendrian surgery can be generalized to any cobordism obtained by attaching $4$-dimensional $2$-handles.
\eexmp

\bcons\label{cons: omega construction}
Suppose $Y$ is a closed, oriented $3$-manifold. Given an unoriented closed $1$-submanifold $\omega\subset Y$, there is a standard way in \cite[\S 4.1]{kronheimer2011khovanov} to construct an $SO(3)$-bundle $P$ over $Y$, unique up to isomorphism with $w_2(P)$ dual to $\omega$. Fix one such bundle $P$. See also \cite[\S 3.7]{scaduto2015instantons} for a Dehn surgery understanding of the construction of the bundle from $\omega$.

Similarly, given an unoriented $2$-submanifold $\nu$ in an oriented $4$-manifold $W$ which restricts to $\omega_i$ on the boundary $Y_i$, the construction in \cite[\S 4.1]{kronheimer2011khovanov} gives an $SO(3)$-bundle over $W$ extending the boundary bundles. Fix such an extension $\bfP$ together with its boundary identifications, which are used for gluing.

Suppressing these choices and the almost complex structures from Example \ref{exmp: surgery cube}, we write $\bY=(Y,\omega,\phi)$ and $\bW=(W,\nu,\Phi)$, omitting $\phi$ and $\Phi$ when they are trivial. We call $\omega$ and $\nu$ \emph{bundle sets}.
\econs
To construct a cube complex in the $\infty$-derived category, it suffices to construct a cube in the instanton $\infty$-cobordism category and apply $CI$. The construction of $CI$ in \S \ref{sec: Construction of the TQFT functor} already incorporates the families of metrics and perturbations and the associated signs, so we need not repeat the constructions in \cite[\S 6.2]{kronheimer2011khovanov} and \cite[\S 6.1-6.2]{scaduto2015instantons}.

\bcons\label{cons: cube of cobordisms}
Given cobordisms supported in disjoint regions $(\bY,\{\bW_i\}_{i\in\un{n}})$, we construct an $n$-cube
\[\sX=\sX[\bY,\{\bW_i\}_{i\in\un{n}}]:\bbox^n\to\BI^\pa.\]
Fix collared representatives of the cobordisms $\bW_i$, together with their bundle gluing maps and product identifications outside the corresponding regions $\sR_i$.

\begin{itemize}
\item For $S\subset\un{n}$, let
\[\sX(S)=\bY_S=\bY_{\ep(S)}\]
be obtained from $\bY$ by attaching the handles in $\{\bW_i\}_{i\in S}$, with the induced additional structures.

\item For $S\subset S'$, let
\[\bW_{S,S'}=\bW_{\ep(S),\ep(S')}:\bY_S\to\bY_{S'}\]
be obtained from $I\times\bY_S$ by attaching the handles in $\{\bW_j\}_{j\in S'\backslash S}$. In particular, $\bW_{j,S}=\bW_{S,S\cup\{j\}}$ represents the edge $\sX(S\subset S\cup\{j\})$.

\item For a strict chain $S_0\subset\cdots\subset S_k$ with $k\geq1$, the consecutive edges are represented by $(\bW_{S_{p-1},S_p})_{p\in\un{k}}$. We construct the coherence data over $[0,1]^{k-1}$, the geometric realization of the corresponding cubical chart of $\fC[\bbox^n]$ in Definition \ref{defn: homotopy coherent realization}. The face $t_i=0$ omits $S_i$, while $t_i=1$ composes the data for the two chains meeting at $S_i$.
\end{itemize}

Since the regions are disjoint, we compare all compositions by identifying their handle attachments and reparametrizing the product collars. The increasing interval reparametrizations, fixed near the endpoints, form convex spaces, so compatible comparison families extend over parameter cubes relative to their faces. These comparisons transport the bundles and almost complex structures. The subgroup data glue compatibly with composition and restriction by \eqref{eq: gluing subgroup} and Proposition \ref{prop: segal space PBIV}.

Choose the embedding families inductively on chain length. On each chart boundary, the face rules and fixed gluing identifications give a family of embeddings of one fixed collared cobordism. Extend the coordinate functions relative to the parameter boundary and the outer collars, then apply the relative parametric Whitney embedding argument in additional Euclidean coordinates vanishing on these fixed subsets. Compare \cite[proofs of Propositions 5.19 and 8.17]{CS2019category}. Transporting the additional structures gives the extension without changing the non-integrality or full-boundary conditions.

These families, with repeated vertices handled by the degeneracy rules, define a homotopy coherent diagram. By Theorems \ref{thm: coherent nerve Quillen equivalence} and \ref{thm: quillen equ}, represent it by a simplicial functor $\fC[\bbox^n]\to\cD$, where $\cD$ is locally Kan with a fixed categorical equivalence $\Nhc_\bu(\cD)\simeq(\BI^\pa)_{\bu,0}$. The homotopy coherent nerve adjunction gives
\[\bbox^n\longrightarrow\Nhc_\bu(\cD)\xrightarrow{\simeq}(\BI^\pa)_{\bu,0},\]
which corresponds under $p_1^*\dashv i_1^*$ to a map $\disc\bbox^n\to\BI^\pa$ representing the required functor $\sX$. All choices on subcubes are obtained by restriction.

If $(\bY,\{\bW_i\}_{i\in\un{n}})$ are surgery cobordisms induced by $(\bY,L=\bigcup_{i=1}^nK_i)$ as in Example \ref{exmp: surgery cube}, then we also write
\[\sX[\bY,L]=\sX[\bY,\{K_i\}_{i\in\un{n}}]=\sX[\bY,\{\bW_i\}_{i\in\un{n}}]\]
and call it a \emph{surgery $n$-cube}.
\econs
\brem
Construction \ref{cons: cube of cobordisms} is obviously compatible with the subcube, i.e.\ for any subset $S\subset \un{n}$, the cube \begin{equation*}\label{eq: cob S}
        \sX[\bY,\{\bW_i\}_{i\in S}]:\bbox^S\to \BI^\pa
    \end{equation*}coincides with the subcube $\sX_S$ from \eqref{eq: subcube}. Moreover, for any subset $T\subset \un{n}\bs S$, the functor \begin{equation*}\label{eq: cob S T}
        \sX[\bY_{\ep(T)},\{\bW_{i,T}\}_{i\in S}]:\bbox^S\to \BI^\pa
    \end{equation*}coincides with the subcube $\sX_{S,T}$ from \eqref{eq: subcube}.
\erem
\brem
If we replace the cobordism $\bW_i$ by $d_i$ composable cobordisms $\{\bW_{i,j}\}_{j\in \un{d_i}}$, then Construction \ref{cons: cube of cobordisms} extends to the hyperbox described in Remark \ref{rem: hyperbox}.
\erem
\bcons\label{cons: hypercube to cube complex}
The constructions in \cite{kronheimer2011khovanov,scaduto2015instantons} yield hypercube complexes whose maps count solutions over families of metrics on cobordisms. For $\ep\leq\epp$, put $k=|\epp-\ep|$. When $k>0$, let $G_{\ep,\epp}\cong\R^k$ be the family with translation parameters $\tau_i$ as in \cite[\S 6.1, Figure 8]{kronheimer2011khovanov}. Write $\breve{G}_{\ep,\epp}$ for its quotient by adding a common constant to all translation parameters, and $\breve{G}_{\ep,\epp}^{+}$ for its compactification.

For a maximal simplex in the subcube from $\ep$ to $\epp$, let $I$ be its vertex chain and $i_1,\dots,i_k$ the coordinates changed in this order. The cube $Q_I$ in Definition \ref{defn: cubical auxiliary data} parametrizes the closure in $\breve{G}_{\ep,\epp}^{+}$ of the region $\tau_{i_1}\leq\cdots\leq\tau_{i_k}$, with coordinates
\[
q_{k-j}=\tau_{i_{j+1}}-\tau_{i_j},~1\leq j<k.
\]
The cubes $Q_I$ therefore subdivide $\breve{G}_{\ep,\epp}^{+}$, with zero faces shared by adjacent regions and infinite faces describing the broken families in \cite[(33)]{kronheimer2011khovanov}. See also \cite[\S 2]{Bloom2011link} and \cite[\S 4.4]{scaduto2015instantons}. The $\infty$-categorical construction records the individual maps associated to these cubes, whose sum appears in \eqref{eq: defn of D ep epp}.

For the sign comparison, choose common almost complex structures with their induced homology orientations, and common metric and perturbation families regular on the subdivision strata. Write $D^{\mathrm{KM}}_{\ep,\epp}$ for $f_{1-\ep,1-\epp}$ in \cite[Definition 6.6]{kronheimer2011khovanov}. For $k>0$, the parameter orientation is the product of the edge orientations along the increasing coordinate order $i_1<\cdots<i_k$ by \cite[Lemma 6.2]{kronheimer2011khovanov}. Formula (41) there gives the total exponent
\[
\begin{aligned}
\sum_{\ell=1}^k\left(\sum_{i<i_\ell}(1-\ep_i)-(\ell-1)\right)
&=\sum_{\ell=1}^k\sum_{i<i_\ell}(1-\ep^{\prime}_i)\\
&=\sum_{i<j}(1-\ep^{\prime}_i)(\ep^{\prime}_j-\ep_j)=a(\ep,\epp),
\end{aligned}
\]
as defined in \eqref{eq: a ep}. With simultaneous translation placed first in the quotient orientation, as in \cite[Lemma 6.3]{kronheimer2011khovanov}, the reversed order of the successive differences in $(q_1,\ldots,q_{k-1})$ contributes $(-1)^{\binom{k-1}{2}}$. The coordinate-permutation signs are already included in $D^{\mathrm{sign}}_{\ep,\epp}$. The additional signs in \cite[(36) and (38)]{kronheimer2011khovanov} therefore give
\[
\begin{aligned}
D^{\mathrm{KM}}_{\ep,\epp}
&=(-1)^{a(\ep,\epp)+\binom{k-1}{2}}(-1)^{k-1}(-1)^{\binom{k}{2}+n-|\ep|}D^{\mathrm{sign}}_{\ep,\epp}\\
&=(-1)^{a(\ep,\epp)+n-|\ep|}D^{\mathrm{sign}}_{\ep,\epp}\\
&=(-1)^kD^{\mathrm{tot}}_{\ep,\epp},
\end{aligned}
\]
since $\binom{k-1}{2}+k-1+\binom{k}{2}=k(k-1)$ is even. The last equality also holds on the diagonal, where both entries are $(-1)^{n-|\ep|}\partial_{\sX(\ep)}$. Thus multiplication by $(-1)^{|\ep|}$ at each vertex identifies the two $n$-cubes with positive signs.

With the same parameter orientations as above, let $D^{\mathrm{Sca}}_{\ep,\epp}$ denote the coefficients in \cite[\S 6.2]{scaduto2015instantons}, with increasing vertex indices. For $k>0$, Scaduto's family counts differ from the quotient counts of Kronheimer--Mrowka by $(-1)^{k-1}$ \cite[(36)]{kronheimer2011khovanov}; for $k=0$, the same factor comes from reversing the trajectory-quotient orientations \cite[\S 4.4]{scaduto2015instantons}. With $(v,u)=(1-\ep,1-\epp)$ in \cite[(38)]{kronheimer2011khovanov} and $(v,w)=(\ep,\epp)$ in \cite[\S 6.2, immediately after (6.2)]{scaduto2015instantons}, the sign exponents become $\binom{k}{2}+n-|\ep|$ and $\binom{k}{2}+|\ep|$, respectively, and hence differ by $n$ modulo $2$. Hence
\[
D^{\mathrm{Sca}}_{\ep,\epp}=(-1)^{n+k-1}D^{\mathrm{KM}}_{\ep,\epp}=(-1)^{n-1}D^{\mathrm{tot}}_{\ep,\epp}.
\]
Hence $x\mapsto(-1)^{(n-1)|x|}x$, where $|x|$ is the degree in $C_\ep[|\ep|-n]$, identifies this hypercube complex with $(C,D^{\mathrm{tot}})$.

The inverse sign change of Varolgunes in Construction \ref{cons: sign on hypercube} also identifies the corresponding $n$-cubes. With these common choices, Constructions \ref{cons: cube complex to hypercube} and \ref{cons: sign on hypercube} therefore recover, up to isomorphism, the hypercube complexes in \cite{kronheimer2011khovanov,scaduto2015instantons} from Construction \ref{cons: cube of cobordisms}.
\econs

\subsection{Fiber squares for surgery cobordisms}\label{subsec: Fiber squares for surgery cobordisms}
For surgery cobordisms in Example \ref{exmp: surgery cube}, we can further translate the surgery exact triangle and link spectral sequence in \cite[\S 5-6]{scaduto2015instantons} to the $\infty$-categorical language. We start with the following construction.

\bcons\label{cons: surgery term}
Let\[\sX=\sX[\bY,L=\bigcup_{i=1}^nK_i]\]be a surgery $n$-cube from Construction \ref{cons: cube of cobordisms}. For fixed $q\in\un{n}$, let $\bY'$ be obtained from $\bY$ by doing $1$-surgery on the knot $K_q$ with the meridian $m_q$ of $K_q$ as the extra bundle set (cf.\ the second row in \cite[Figure 1]{scaduto2015instantons}). Let $L'=\bigcup_{i\in\un{n}\backslash\{q\}} K_i$ be the induced knot in $\bY'$ and let\[\sX'=\sX[\bY',L']\]be the surgery $(n-1)$-cube. From \eqref{eq: 1 cube n-1 cube}, we regard the original surgery $n$-cube $\sX$ as a $1$-cube of $(n-1)$-cube
\begin{equation*}
    \sX^{\{q\}}:\bbox^{\{q\}}=\Delta^1\to \Fun(\bbox^{\un{n}\bs \{q\}}=\bbox^{n-1},\BI^\pa).
\end{equation*}Let \[\sY_i=CI( \sX^{\{q\}}(i))\text{ for }i=0,1,~ \sY_2=CI(\sX')\in \Fun(\bbox^{\un{n}\bs \{q\}},\D)\]\[\aand \sW_{0}=CI(\sX^{\{q\}}(0\to 1)):\sY_0\to \sY_1.\]

From Example \ref{exmp: surgery cube}, there are two further $1$-cubes of $(n-1)$-cube from surgery cobordisms on $m_q$ and the meridian $m_q'$ of $m_q$ \[\sX^{\{q\}}(1)\to \sX'\aand\sX'\to \sX^{\{q\}}(0)',\]where the prime on $\sX^{\{q\}}(0)$ indicates that the almost complex structures are different. From the discussion before \cite[Proposition 5.2]{kronheimer2011khovanov} and the remark after \cite[Proposition 7.6]{kronheimer2011khovanov}, when applying $CI$, we have\[CI(\sX^{\{q\}}(0)')=CI(\sX^{\{q\}}(0))[1],\]where $[1]$ denotes the grading shift operator of the $\Z/2$-grading. Note that this perspective is useful because it makes anti-chain maps in \cite[\S 7]{kronheimer2011khovanov} become chain maps because the shift operator introduces a minus sign on the differential.

From those $1$-cubes of $(n-1)$-cube, we define\[\sW_1:\sY_1\to \sY_2\aand \sW_2:\sY_2\to \sY_0[1]\]by applying $CI$ to the morphisms. For $i\in\{0,1,2\}$ and $k\in\Z$, we define \[\sY_{i+3k}=\sY_i[k]\aand\sW_{i+3k}=\sW_i[k].\]Note that \[\sY_{i+6k}=\sY_i\aand\sW_{i+6k}=\sW_i\]because we only have $\Z/2$-homological grading.
\econs
Then we state the $\infty$-categorical lift of the surgery exact triangle and its iteration for link surgeries.
\bthm\label{thm: surgery fiber sequence}
Suppose $\sY_j\in \Fun(\bbox^{n-1}=\bbox^{\un{n}\backslash\{q\}},\D)$ and $\sW_j:\sY_j\to \sY_{j+1}$ for $j\in \Z$ are from Construction \ref{cons: surgery term}. Then, for any $j\in\Z$, the following diagram is a fiber square in $\Fun(\bbox^{n-1},\D)$.
\begin{equation}\label{eq: fiber sequence}
\xymatrix{\sY_j\ar[r]^{\sW_j}\ar[d]&\sY_{j+1}\ar[d]^{\sW_{j+1}}\\0\ar[r]&\sY_{j+2}}
\end{equation}
\ethm
\bpf
Following \cite[\S 5, \S 6.3]{scaduto2015instantons} (see also \cite[\S 6-7]{Bloom2011link}), we apply the triangle detection lemma in Appendix \S \ref{sec: Triangle detection lemma}, namely Lemma \ref{six-lemma}, to the stable $\infty$-category $\Fun(\bbox^{n-1},\D)$. Apply the cubical construction of Construction \ref{cons: hypercube to cube complex} to the surgery families in \cite[\S 5.2, \S 6.1]{scaduto2015instantons}. Choose metric and perturbation families extending those defining $\sY_j$ and $\sW_j$. Their restrictions agree on identified parameter faces, and on broken faces they are the products of the families on the pieces. Besides the usual breaking faces, these families have faces obtained by stretching along local copies of $S^3$. Their counts vanish by \cite[\S 5.2]{scaduto2015instantons}. The remaining boundary identities give successive nullhomotopy squares, forming the diagram of Lemma \ref{six-lemma} with $X_i=\sY_{j+i}$ for $0\le i\le5$. At each vertex, the induced maps $X_i[1]\to X_{i+3}$ for $i=0,1,2$ are the maps studied in \cite[\S 5.3]{scaduto2015instantons} and are isomorphisms. Since isomorphisms in the functor category are detected pointwise, Lemma \ref{six-lemma} applies and gives \eqref{eq: fiber sequence}.
\epf
\brem
The nullhomotopies above use $S^3$ cuts, while the pointwise verification of these isomorphisms in \cite[\S 5.3]{scaduto2015instantons} also uses an $S^1\times S^2$ cut. These stretching families are of the type described in Remark \ref{rem: generalization of broken}, but the bundle restrictions to these hypersurfaces do not satisfy the non-integral condition. Thus the corresponding families do not define higher morphisms in $\BI^\partial$ and are not included in the construction of $CI$.
\erem

\brem\label{rem: iterating}
We iterate Construction \ref{cons: cube of cobordisms} and Theorem \ref{thm: surgery fiber sequence} in the decreasing order $n,n-1,\dots,1$ to obtain an isomorphism in $\D$ (i.e.\ a quasi-isomorphism in the usual sense)\begin{equation}\label{eq: quasi iso}
        \Tcofib (CI(\sX))\to CI(\bY_{\vec{2}})
    \end{equation}where $\bY_{\vec{2}}$ is obtained from $\bY$ by doing $1$-surgeries on  $\{K_i\}_{i\in\un{n}}$ with the union of meridians $\omega=\bigcup_{i=1}^nm_i$ of $K_i$ as the extra bundle set. The constructions in \cite[\S 5, especially the discussion after Remark 5.4]{Bloom2011link} and \cite[\S 6.1]{scaduto2015instantons} may provide a comparison between the isomorphisms obtained from different orders. Here we only need the existence of \eqref{eq: quasi iso} and use the decreasing order throughout.

   The isomorphism \eqref{eq: quasi iso} preserves the $\Z/2$-grading in $\D$. For the corresponding grading shifts in the standard instanton grading, see \cite[\S 6.4]{scaduto2015instantons}.
\erem
\brem
    From \cite[Lemma 2.5]{LY2025dimension}, if $\omega=\bigcup_{i=1}^nm_i$ represents trivial homology class in $H_1(Y_{\vec{2}};\Z/2)$ and we pick a possibly unorientable embedded surface $S\subset Y_{\vec{2}}$ with $\partial S=\omega$, then there are isomorphisms\[\bI_S:CI(\bY_{\vec{2}})\xra{\cong} CI(\bY_{\vec{2}}')\aand \bI_S':CI(\bY_{\vec{2}}')\xra{\cong} CI(\bY_{\vec{2}})\]obtained from applying $CI$ to some morphisms in $\BI$ induced by $S$, where $\bY_{\vec{2}}'$ is obtained from $\bY_{\vec{2}}$ by removing the extra bundle set.
\erem

\section{Enhanced cobordisms and \texorpdfstring{$\mu$}{mu}-operators}\label{sec: Generalized cap product}
In this section, we provide the construction of the $\mu$-operators at the chain level and study the interaction with cobordisms and cubes of chain complexes. We mainly follow \cite[\S II.(ii)]{KM1995structure} for the detailed construction. Throughout this section, we assume that all subgroup data $\phi$ and $\Phi$ are trivial, including those on the composite cobordisms. We use determinant-$1$ gauge transformations and omit the trivial subgroup data from the notation.
\subsection{Idea of the construction}\label{subsec: enhanced cobordisms}
In this subsection, we first review the idea of the construction for $\mu$-operators from references in \S \ref{subsec: enhanced cob}, so that we can focus on the construction of enhanced data in the next subsection.

Following \S \ref{subsect: Sketch of gauge theory construction}, we fix a morphism in $\BIa^\pa$ (cf.\ Definition \ref{defn: BIa}), i.e.\ a cobordism with auxiliary data including a Riemannian metric and admissible primary and secondary perturbations, written as\[(\bW=(W,\bfP,\Phi,\bfJ),\bfa):(\bY_0=(Y_0,P_0,\phi_0,J_0),a_0)\to (\bY_1=(Y_1,P_1,\phi_1,J_1),a_1).\]Note that \cite[\S II.(ii)]{KM1995structure} focuses on the case when $\Phi$ and $\phi_i$ are trivial, so we only consider the configuration spaces $\cB$ instead of the quotient $\ov{\cB}=\cB/\Phi$\footnote{In the first arXiv version of this paper, we work on nontrivial subgroups and the quotient, which indeed needs extra care. We remove that general case for simplicity and only consider $\cB$.}.

Following \S \ref{subsec: enhanced cob}, all chain complexes and maps in this section are over $\Z[1/2]$. Let \begin{equation}\label{eq: A(W)}
    \bA(W;1/2)=\operatorname{Sym}^*H_{\operatorname{even}}(W;\Z[1/2])\ot \Lambda^*H_{\operatorname{odd}}(W;\Z[1/2]),
\end{equation}be the graded-commutative algebra generated by $H_*(W;\Z[1/2])$, where the subscripts denote the even and odd degree parts, respectively. Recall $W^*$ is defined in \eqref{eq: open end} by attaching half-infinite cylinders to the ends of $W$. Note that $\bA(W^*;1/2)\cong \bA(W;1/2)$ by a canonical isomorphism from the inclusion. Hence we do not distinguish them. We also write $\bW^*$ for the tuple obtained from $\bW$ by extending the rest structures to the cylinders. The rational algebra $\bA(W)$ of \S \ref{subsec: enhanced cob} is recovered by extension of scalars to $\Q$. The following rational calculation motivates the construction of the $\mu$-operators but is not needed to define them.

\bprop[{\cite[Proposition (5.1.15)]{DK1992instanton}}]\label{prop: homotopy type}
Suppose $X$ is a closed, oriented, connected, and simply connected $4$-manifold and suppose $\bfP$ is a principal $SU(2)$ bundle over $X$. Let $\cB^*(X,\bfP)$ be the configuration space of irreducible connections. Then we have \[H^*(\cB^*(X,\bfP);\Q)\cong \operatorname{Sym}^*\left(H_{0}(X;\Q)\op H_{2}(X;\Q)\right).\] The odd homology $H_{\operatorname{odd}}(X;\Q)$ vanishes because $X$ is simply connected.
\eprop
\brem
When $X$ is not simply-connected, the cohomology $H^*(\cB^*(X,\bfP);\Q)$ is more complicated, but also closely related to $\bA(X)$ by the discussion after \cite[Proposition (5.1.15)]{DK1992instanton}. The cohomology over $\Z$ is even more complicated as there could be torsion classes.
\erem

For the concrete construction, fix
\begin{equation}\label{eq: enhanced monomial}
z=c\prod_{p=1}^r\be_p\in\bA(W;1/2),~c\in\Z[1/2],
\end{equation}
where the $\be_p$ are integral homology classes of dimensions $0,1,2,3$. Note that we remove the assumption that $\be_p$ is divisible by $2$ in dimension $2$ as in \cite[\S II.(ii)]{KM1995structure} by using the construction of complex line bundle in \cite[\S 2.1]{sasahira2013floer} instead of \cite[\S II]{donaldson1986connection} and \cite[\S 5.2.1]{DK1992instanton}; see Construction \ref{cons: sasahira} for more details.

Given $z$ as in \eqref{eq: enhanced monomial} and some auxiliary choice $\tau$ about $z$, called \emph{enhanced data}, we define a map \begin{equation}\label{eq: enhanced cob}
    CI(\bW,\bfa, z,\tau):CI(\bY_0,a_0)\to CI(\bY_1,a_1),
\end{equation}which has $\Z/2$-degree $\sum_{p=1}^r(4-\deg\be_p)$ and satisfies
\begin{equation}\label{eq: chain relation for z}
    \partial\,CI(\bW,\bfa,z,\tau)=(-1)^{\sum_{p=1}^r(4-\deg\be_p)}CI(\bW,\bfa,z,\tau)\,\partial.
\end{equation}
For $z=1$, this is $CI(\bW,\bfa)$. The coefficient $c$ multiplies the entire map. For the even-dimensional homology classes considered below, a family of enhanced data joining $\tau^0$ and $\tau^1$ determines a chain homotopy \begin{equation}\label{eq: enhanced cob family}
    CI(\bW,\bfa,z,\wti{\tau}(I)):CI(\bY_0,a_0)\to CI(\bY_1,a_1),
\end{equation}where $\wti{\tau}(I)$ is some auxiliary choice parametrized by $I=[0,1]$ and restricting to $\tau^i$ at $i\in\{0,1\}$ (also called \emph{enhanced data over $I$}, or \emph{family enhanced data}). Here we use the superscript to reserve the subscript for factors of $z$. Higher families similarly determine higher homotopies, which we also construct explicitly.

\brem
To adapt results in the references directly and avoid detailed computations, we will only focus on an \emph{ad hoc} setup to construct homotopies and higher homotopies and leave the full general construction and the proof of independence to the future.
\erem

Recall from \S \ref{subsect: Sketch of gauge theory construction}, the generators in $CI(\bY_i,a_i)$ correspond to perturbed flat connections on $(Y_i,P_i)$ and the construction of the cobordism map $CI(\bW,\bfa)$ uses the $d$-dimensional compactified moduli spaces $\ov{M}(x_0,\bW,x_1)$ of perturbed ASD connections in the irreducible configuration space $\cB_{\loc}(\bW^*)$ modulo determinant-$1$ gauge transformations, for $d=0,1$. To construct the map $CI(\bW,\bfa, z,\tau)$, we need to use moduli spaces of higher dimensions. Even though the bubble terms do appear in the compactification for $d\ge 8$, the choice of $\tau$ will make the construction disjoint from them.

The central stabilizers act trivially on the $SO(3)$ bundle, so there is a universal $SO(3)$ bundle over $W^*\times\cB_{\loc}(\bW^*)$. As in \cite[\S II.(i)--(ii)]{KM1995structure}, slant product with $-1/4$ times its first Pontryagin class defines
\begin{equation}\label{eq: mu maps}
\ov{\mu}:H_*(W^*;\Z[1/2])\to H^{4-*}(\cB_{\loc}(\bW^*);\Z[1/2]).
\end{equation}

For $z$ as in \eqref{eq: enhanced monomial}, the intuition behind the map $CI(\bW,\bfa, z,\tau)$ is the formal pairing\begin{equation}\label{eq: pairing cob}
    \begin{aligned}
        \langle CI(\bW,\bfa, z,\tau)(x_0),x_1\rangle=&c\cdot \langle \ov{\mu}(\be_1)\cup\cdots\cup \ov{\mu}(\be_r),[\ov{M}_d(x_0,\bW,x_1)]\rangle \\&\rmfor d=\sum_{p=1}^r (4-\deg\be_p),
    \end{aligned}
\end{equation}while the existence of bubbles makes the genuine construction more complicated. Indeed, one needs to pick geometric representatives $V(\be_p)=V(\be_p,\tau)\subset \ov{M}_d(x_0,\bW,x_1)$ (not in $\cB_{\loc}(\bW^*)$ for some technical reason) of $\ov{\mu}(\be_p)$, with coefficients $b_p$. These representatives are obtained by factoring through the configuration space of the connections on some \emph{suitable} neighborhood $U(\be_p)$ of an embedded manifold $Z_p\subset W$ representing the fixed integral class $\be_p$. Then \eqref{eq: pairing cob} becomes
\begin{equation}\label{eq: pairing 2}
    \langle CI(\bW,\bfa, z,\tau)(x_0),x_1\rangle=c\cdot\Bigl(\prod_{p=1}^r b_p\Bigr)\#\Bigl(\ov{M}_d(x_0,\bW,x_1)\cap V(\be_1)\cap\cdots\cap V(\be_r)\Bigr),
\end{equation}where $\#$ is the signed count and $b_p$ is the coefficient in Table \ref{tab: representative} for $V(\be_p)$. The choices of $V(\be_p)$ make all intersections transverse and disjoint from bubbles. For two fixed generators $x_0$ and $x_1$, the equation showing that $CI(\bW,\bfa, z,\tau)$ is a chain map comes from the pairing
\begin{equation}\label{eq: pairing 3}
    0=\# \partial \left(\ov{M}_{d+1}(x_0,\bW,x_1)\cap V(\be_1)\cap\cdots\cap V(\be_r)\right),
\end{equation}since the choices of $V(\be_p)$ make the boundary all come from the broken trajectories.
\subsection{Enhanced data}
In this subsection, we will only sketch the construction of $V(\be_p)$ from \cite[\S II.(ii)]{KM1995structure}, aiming to point out the enhanced data $\tau$ and its relation to $V(\be_p)=V(\be_p,\tau)$. Note that the reference uses $X$ to denote the closed $4$-manifold and we use $W$ to denote the compact $4$-manifold with nonempty boundary, and $W^*$ for the corresponding open manifold from \eqref{eq: open end}.

We call a subset $U\subset W$ \emph{suitable} if it is a smooth submanifold with boundary and $U\neq W$, and the map induced by inclusion $H_1(U;\Z/2)\to H_1(W;\Z/2)$ is surjective. Note that a suitable subset always exists by taking the regular neighborhood of a collection of embedded curves that generate $H_1(W;\Z/2)$ and have the same basepoint $x_0$. The \emph{suitable neighborhood} of a submanifold $Z$ in $W$ is the union of the above regular neighborhood of curves and the regular neighborhood of $Z$. The suitable condition implies that twisted reducible connections on $W$ remain irreducible when restricted to $U$, where the twisted reducible connections are irreducible but involve $H_1(W;\Z/2)$, discussed in \cite[\S II.(i)]{KM1995structure}.

For each factor $\be_p$ in \eqref{eq: enhanced monomial}, choose an embedded oriented submanifold $Z_p$ with $[Z_p]=\be_p$ in integral homology. For the collection of homology classes in \eqref{eq: enhanced monomial}, we choose a collection of submanifolds $\{Z_p\}_{p\in\un{r}}$ which intersect each other transversely. Note that $Z_p$ can be disconnected. Here $\un{r}=\{1,\dots,r\}$, and the underline on an integer $n$ denotes the set $\{1,\dots,n\}$ later.

Moreover, we choose a basepoint $x_p$ in $W$ and a collection of arcs $\al_p$ connecting $x_p$ to each component of $Z_p$. For each $x_p$, we choose a collection of curves $C_p=\{C_{p,s}\}_{s\in\un{t_p}}$ based at $x_p$ that generates $H_1(W;\Z/2)$. We assume all $\al_p$ and $C_p$ are disjoint from each other and transverse to $Z_{q}$ for $q\in\un{r}\bs \{p\}$.

For $\ep>0$, let \begin{equation}\label{eq: regular U}
    U_p=U(Z_p,x_p,\al_p,C_p,\ep)
\end{equation}be the regular neighborhood of $Z_p\cup \al_p\cup C_p$ with radius $\ep$ with respect to the metric on $W$. We assume $\ep$ is small enough so that the $U_p$ intersect pairwise transversely. We call $\{U_p\}_{p\in \un{r}}$ a collection of \emph{suitable neighborhoods}.

By the proof of \cite[Proposition B.1]{Froyshov2024mod2}, which also applies with the cutoff functions, solutions of \eqref{eq: perturbed ASD 2} on $W^*$ with irreducible limiting critical points restrict irreducibly to the unperturbed compact part $W$.

Suitable neighborhoods are used to ensure that restriction from $W$ to $U_p$ preserves irreducibility. For irreducible non-flat ASD connections, this follows from the surjectivity of $H_1(U_p;\Z/2)\to H_1(W;\Z/2)$ and unique continuation, as in \cite[\S II.(ii)]{KM1995structure}. In that construction, flat connections are excluded by passing to a blow-up and changing the bundle on the exceptional sphere. Here we include additional curves to obtain the same restriction property for flat connections.

\bcons\label{cons: flat restrictions}
For each connected component $Y_\al$ of $\partial W$, choose a connected witness surface $\Sigma_\al\subset Y_\al$. Choose loops representing generators of $\pi_1(\Sigma_\al)$, based at a common point of $Y_\al$ by a common connecting path. We call these loops \emph{witness curves}. Enlarge this collection by adding finitely many based loops so that their homology classes generate $H_1(Y_\al;\Z/2)$. We call the enlarged collection \emph{boundary curves}.

On each connected component of $W$, push its boundary curves into the collars and include them in each $C_p$ on that component, using one connecting path from $x_p$ to the basepoint of each collection. Choose the curves and paths in general position. For a neighborhood $U_p$ on a cylinder over $Y_\al$, use the corresponding boundary curves in a common slice and translate them with the neighborhood. The projection then induces a surjection $H_1(U_p;\Z/2)\to H_1(Y_\al;\Z/2)$, which is preserved by translation.

If a flat connection had reducible restriction to $U_p$, the holonomies of the witness curves would fix a common nonzero vector. Its restriction to the corresponding $\Sigma_\al$ would therefore be a reducible flat connection, contradicting \eqref{eq: witness}, since a flat $SO(2)$-bundle over $\Sigma_\al$ has zero Euler class.
\econs

For each fixed suitable neighborhood $U_p$, we omit the subscript $p$. Let $\cB^*(U)\subset \cB(U)$ be the configuration space of irreducible connections. Let $\cB^{o,*}(U)$ be the configuration space of framed irreducible connections at $x_0=x_p$ (i.e.\ irreducible connections modulo determinant-$1$ gauge transformations that fix the fiber over the basepoint $x_p\in U_p$). From \cite[\S II.(ii)]{KM1995structure} and Construction \ref{cons: flat restrictions}, there is a well-defined restriction map for all $d'\in\N$ \begin{equation}\label{eq: restriction}
    r_{U,d'}:M_{d'}(x_0,\bW,x_1)\to \cB^*(U),
\end{equation}which extends over broken trajectories by neglecting the trajectories on the cylinders and is compatible with Uhlenbeck limits when no bubbling occurs in $U$. Depending on $\deg\be_p$, we consider some bundle $\wti{\cB}$ over $\cB^*(U)$ in Table \ref{tab: representative}. When $\ga$ is disconnected, fix an ordering of its components and define $\operatorname{hol}_{\ga,x_0}(A)$ as the ordered product of their holonomies, transported to $x_0$ along the chosen arcs $\al$. Include this ordering in $\tau$.

\begin{table}
\begin{center}
\small
\begin{tabular}{ |c|c|c|c|c| }
 \hline
 $\deg\be_p$ & $Z_p$ &Bundle $\wti{\cB}$ & Sections & Coefficient $b_p$\\
 \hline
$0$ & $x$ & \makecell{$\Xi\otimes_{\R}\C$, where\\ $\Xi=\cB^{o,*}(U)\times_\rho\R^3$} & \makecell{a generic pair of complex\\ sections $(s_1,s_2)$} & $|x|/4$\\
 \hline
 $1$ & $\ga$ & \makecell{$\operatorname{Ad}_{SO(3)}\cB^{o,*}(U)$\\$=\cB^{o,*}(U)\times_{\operatorname{Ad}}SO(3)$} & \makecell{$s_0([A])=[A,\id]$ and\\a generic $s'$ homotopic to \\$s([A])=[A,\operatorname{hol}_{\ga, x_0}(A)]$}& $-1/2$\\
 \hline
 $2$ & $\Sigma$ &\makecell{The complex line bundle\\ $\cL_U^{
 \ot 2
 }$ from \cite[\S 2.1]{sasahira2013floer}} & \makecell{the zero section $s_0$\\and a generic $s'$} & $1/2$\\
 \hline
 $3$ & $Y$ & $\cB^*(U)\times\R/\frac{1}{4}\Z$ & \makecell{the zero section $s_0$ and\\a generic $s'$ homotopic to \\$s([A])=C.S.(A|_Y)$}&$1/4$\\
 \hline
\end{tabular}

\caption{\label{tab: representative}Bundles, sections, and coefficients for different degrees, where $\rho$ is the standard representation of $SO(3)$ on $\R^3$, and $C.S.$ is the Chern-Simons functional. In degree $0$, $|x|$ denotes the signed number of points in $x$. For disconnected $W$, decompose the class over its connected components and extend the construction linearly.}
\end{center}
\end{table}
\bcons\label{cons: sasahira}
In degree $2$, the construction in \cite[\S II.(ii)]{KM1995structure} uses suitable neighborhoods $U_p$ with $w_2(\bfP)|_{U_p}=0$ and assumes $\be_p\in2H_2(W;\Z)$ when choosing the representatives. If a connected embedded oriented surface $\Sigma$ has $\Sigma\cdot\Sigma=0$, joining two parallel copies by a tube gives a connected embedded representative of $2[\Sigma]$. However, for $\Sigma$ in a $3$-dimensional submanifold $Y$, this replacement need not remain in $Y$, making the product condition become complicated in Definition \ref{defn: special enhanced cobordism}.

We use the squared determinant-line construction in \cite[\S 2.1]{sasahira2013floer}. Choose a $U(2)$ lift of $\bfP|_\Sigma$, a fixed determinant connection, and a spin structure on $\Sigma$ that induces a square root $K_\Sigma^{1/2}$ of the canonical line bundle of $\Sigma$. Restricting connections on $U$ to $\Sigma$ gives a family of spin-twisted $\bar\partial^*$-operators and its determinant line over $\cB^{o,*}(U)$. The central elements $\{\pm1\}$ act trivially on its square, which therefore descends to the line bundle $\cL_U^{\ot 2}$ over $\cB^*(U)$ in Table \ref{tab: representative}. Despite this notation, $\cL_U^{\ot 2}$ need not be the square of a genuine line bundle on $\cB^*(U)$. This local construction does not require $W$ to be simply connected or $\langle w_2(\bfP),[\Sigma]\rangle$ to vanish.

The isomorphism class of $\cL_U^{\ot 2}$ is independent of these choices. Changes of spin structure or determinant connection give homotopic operator families after identifying the underlying complex bundles. Two $U(2)$ lifts differ by tensoring with a line bundle on $\Sigma$. By excision in integral $K$-theory, their family index classes differ by a sum of integer multiples of evaluation bundles, whose determinants are equivariantly trivial under determinant-$1$ gauge transformations.

By the families index theorem, the first Chern class of $\cL_U^{\ot 2}$, after pullback by \eqref{eq: restriction}, is twice the restriction of $\ov{\mu}([\Sigma])$. Thus its generic zero set, with coefficient $1/2$, represents the usual Donaldson $\mu$-class. We may therefore use the original surface $\Sigma\subset Y$, including the case in Remark \ref{rem: special mu case}, and take $\wti{Z}=I\times\Sigma$ in the product cobordism without doubling $\Sigma$.
\econs

Then we choose a generic pair of sections $(s_1,s_2)$ (or $(s_0,s')$) of $\wti{\cB}$ as in Table \ref{tab: representative} such that they are related to $\ov{\mu}(\be_p)$. The transversality condition can be achieved by \cite[Lemma (5.2.9)]{DK1992instanton} for all dimensions $0\le d'<D$, where $D>d+1$ and $d$ is given by \eqref{eq: pairing cob}. Let \begin{equation}\label{eq: construction of V}
    V(\be_p)=V(U_p,D,s_1,s_2)\subset \bigcup_{0\le d'<D}\ov{M}_{d'}(x_0,\bW,x_1)
\end{equation}be the locus where the two complex sections are linearly dependent when $\deg\be_p=0$, and where $s_0=s'$ otherwise, after pullback by the restriction maps. In degree $2$ this is the zero set of $s'$. Use the complex normal orientation in degree $0$, so that $\frac14c_2(\Xi\otimes_{\R}\C)=-\frac14p_1(\Xi)$, as in \cite[\S II.(ii)]{KM1995structure}. Such a construction makes the pairing in \eqref{eq: pairing 2} well-defined by a dimension argument about bubbles in \cite[\S II.(ii)]{KM1995structure}. The signed count in \eqref{eq: pairing 2} is multiplied by $c\prod_{p=1}^r b_p$, with $b_p$ as in Table \ref{tab: representative}.

\bcons\label{cons: tau}
Assume that $\bW$ has trivial $\Phi$. The auxiliary choice $\tau$ in \eqref{eq: enhanced cob} for $z$ in \eqref{eq: enhanced monomial} consists of the following data for each $p\in\un{r}$, satisfying the transversality conditions.
\benum
    \item an embedded oriented submanifold $Z_p\subset W$ such that $[Z_p]=\be_p$ in integral homology;
    \item a basepoint $x_p\in W$, a collection of arcs $\al_p$ connecting $x_p$ to each component of $Z_p$, a collection of curves $C_p$ based at $x_p$ that generates $H_1(W;\Z/2)$ and includes the boundary curves of Construction \ref{cons: flat restrictions}, and a small number $\ep>0$ for the choice of the regular neighborhood; the first two terms determine a collection of suitable neighborhoods $\{U_p\}_{p\in \un{r}}$;
    \item a fixed integer $D>d+1$ (for all $p$) and two sections $s_{p,1}$ and $s_{p,2}$ of some bundle $\wti{\cB}_p$ over $\cB^*(U_p)$.
\eenu
\econs
\brem\label{rem: abstract homotopy}
Over a field containing $\Z[1/2]$, if we only want to show the chain map $CI(\bW,\bfa,z,\tau)$ is independent of $\tau$ up to an \emph{abstract} homotopy, i.e.\ we only need its existence and not an explicit construction, then the strategy in \cite[\S II.(ii)]{KM1995structure} can be adapted directly. Roughly, if the data for all but one $\be_q$ in $\tau^0$ and $\tau^1$ are the same, then we consider the intersection\begin{equation}\label{eq: the intersection}
    V(z\backslash \be_{q})=\ov{M}_d(x_0,\bW,x_1)\cap V(\be_1)\cap\cdots\cap\widehat{V}(\be_{q})\cap\cdots \cap V(\be_r),
\end{equation}where $\widehat{V}(\be_{q})$ denotes the absence of the space in the intersection. By a dimension argument, the space \eqref{eq: the intersection} is already disjoint from the bubbles (when the degree of some $\be_p$ for $p\in \un{r}\backslash\{q\}$ is $3$, the argument works with a little modification). Hence the further intersection with $V(\be_{q})$ is equivalent to \begin{equation}\label{eq: cohomology pairing}
    \langle \ov{\mu}(\be_{q}),[V(z\backslash \be_{q})]\rangle.
\end{equation}Since different choices of $\tau$ correspond to the same cohomology class $\ov{\mu}(\be_{q})$, the chain maps for $\tau^0$ and $\tau^1$ are abstractly chain homotopic. If the data for more $\be_p$ are different, we modify the data for one $\be_p$ at a time and then show different $\tau$ produce chain homotopic maps. Using a similar argument and the linearity of $\ov{\mu}$, we can show the chain map depends on the homology class $\be_p$ linearly up to abstract homotopy. The graded commutative up to abstract chain homotopy follows from the cup product in the construction \eqref{eq: pairing cob}.
\erem
\bcons\label{cons: mu maps}
For a product cobordism $\bW:\bY\to\bY$ with trivial $\Phi$ and $z\in\bA(Y;1/2)=\bA(W;1/2)$ as in \eqref{eq: enhanced monomial}, the construction above gives
\begin{equation}\label{eq: mu z}
H_*(CI(\bW,\bfa,z,\tau))=\mu(z):I(\bY;\Z[1/2])\to I(\bY;\Z[1/2])
\end{equation}
where $\mu(z)$ is the standard operator recalled in \S\ref{subsec: enhanced cob}.
\econs
\subsection{Family enhanced data}

Beyond Remark \ref{rem: abstract homotopy}, we want to construct the homotopy explicitly via the enhanced data $\wti{\tau}(I)$ in \eqref{eq: enhanced cob family}. For our concrete purpose, we only focus on the cases when $\deg\be_p\in\{0,2\}$. After identifying the bundles, the sections in these cases can be joined by linear interpolation. For $\deg\be_p\in\{1,3\}$, one must additionally check that the endpoint sections are homotopic; this holds when both are homotopic to the same section $s$ in Table \ref{tab: representative}. The idea of this subsection comes from the proof of \cite[Lemma (3.3)(i)]{Donaldson1990polynomial}.

\bcons\label{cons: tau tilde}
For $z=c\prod_{p=1}^r\be_p\in\bA(W;1/2)$ as in \eqref{eq: enhanced monomial} with $\deg\be_p\in\{0,2\}$, the family enhanced data in $\wti{\tau}(I)$ is similar to $\tau$ in Construction \ref{cons: tau} with mild modifications. We will use the superscript $i$ to denote the data in $\wti{\tau}(I)$ for $i=0,1$ and use the notation with tilde to denote the data for $\wti{\tau}$. We list the choices in $\wti{\tau}$ as follows.

\benum
    \item Assume $Z_p^0\cap Z_p^1=\varnothing$. Since $[Z_p^0]=[Z_p^1]=\be_p$, there is a submanifold $\wti{Z}_p$ with $\partial\wti{Z}_p=-Z_p^0\sqcup Z_p^1$. We also pick a path $\wti{x}_p$ connecting $x_p^0$ to $x_p^1$. In degree $2$, choose a spin structure on $\wti{Z}_p$ and a $U(2)$ lift of $\bfP|_{\wti{Z}_p}$ with a fixed determinant connection, all of product form near the boundary. Use their restrictions as endpoint data. By Construction \ref{cons: sasahira}, choose isomorphisms from the original squared determinant lines to those obtained from these restrictions, and transport the endpoint sections through these isomorphisms.
    \item Let $\wti{\ep}>0$ be a small number for the choice of the regular neighborhood $\wti{U}_p$ of \[\wti{Z}_p\cup\wti{x}_p\cup C_p^0\cup C_p^1\cup\al_p^0\cup\al_p^1\] which contains $U_p^0$ and $U_p^1$.
    \item For $i=0,1$, the restriction \eqref{eq: restriction} onto $\cB^*(U_p^i)$ factors through the open subset $\cB^{*,*}(\wti{U}_p)\subset\cB^*(\wti{U}_p)$ where both endpoint restrictions are irreducible, as in \cite[Lemma 2.6]{KM1995structure}. Here and in subsequent pullbacks, we suppress these open subsets from the notation. With these compatible choices, fix an isomorphism between the two pullback bundles
\begin{equation}\label{eq: isomorphism iota}
\iota_p:\wti{\cB}_p^1\xra{\cong}\wti{\cB}_p^0,
\end{equation}
where we abuse the notation to denote the pullback bundles. More explicitly, in degree $0$, use parallel transport along $\wti{x}_p$ to obtain the isomorphism. In degree $2$, the endpoint families in Construction \ref{cons: sasahira} have equal gauge-equivariant index classes by bordism invariance. Taking determinants, squaring, and descending gives $\iota_p$.
    \item Consider the product $I\times \cB^*(\wti{U}_p)$ and the pullback bundles (identified by $\iota_p$ from \eqref{eq: isomorphism iota}) and sections on the product. Since $\deg\be_p\in\{0,2\}$, the bundles $\wti{\cB}_p^i$ are complex vector bundles with fibers $\C^3$ or $\C$, hence any two sections can be joined by linear interpolation. In degree $2$, keep the zero section fixed. Choose $\wti{D}$ with $d+1<\wti{D}\le\min\{D^0,D^1\}$. Keeping the endpoint sections fixed, perturb the interpolation relative to endpoint collars to obtain sections $\wti{s}_{p,1},\wti{s}_{p,2}$ satisfying the transversality conditions on the strata of the parametrized moduli spaces of total dimension less than $\wti{D}$. Here the total dimension includes the parameter in $I$.
\eenu
\econs
\brem\label{rem: translation}
The sections $\wti{s}_{p,1},\wti{s}_{p,2}$ in Construction \ref{cons: tau tilde} are equivalent to two (smooth) families of sections $s_{p,1}^t,s_{p,2}^t$ of the two identified pullback bundles over $\cB^*(\wti{U}_p)$ for $t\in I$ that restrict to the prescribed sections at $t\in\{0,1\}$. Moreover, if $\wti{Z}_p$ is diffeomorphic to a product $I\times Z_p$, we use the slices $\{t\}\times Z_p$ as the underlying submanifolds.
\erem

Given the enhanced data $\wti{\tau}(I)$ as in Construction \ref{cons: tau tilde}, we can define similarly as in \eqref{eq: construction of V} \begin{equation}\label{eq: construction of V family}
V(\wti{U}_p)=V(\wti{U}_p,\iota_p,\wti{D},\wti{s}_{p,1},\wti{s}_{p,2})\subset\bigcup_{0\le d'<\wti{D}}\ov{M}_{d'}(x_0,\bW(I),x_1)
\end{equation}using the complex linear-dependence condition for degree $0$ and the zero-set condition for degree $2$, after pullback by the restriction maps, with the coorientations specified in \eqref{eq: construction of V}.

\brem\label{rem: moduli space diff}
The moduli space $\ov{M}_d(x_0,\bW(I),x_1)$ is the compactification of the solution space in $I\times \cB_{\loc}(\bW^*)$, which is different from $\ov{M}_d(x_0,\bW,x_1)$. Other than broken trajectories and bubbles, the boundary of $\ov{M}_d(x_0,\bW(I),x_1)$ also contains some terms from $\partial I$, which is related to $\ov{M}_{d-1}(x_0,\bW,x_1)$.
\erem

One may expect to define the map \eqref{eq: enhanced cob family} similarly as in \eqref{eq: pairing 2} and verify it is a homotopy by a similar argument in \eqref{eq: pairing 3}. However, just as in \cite[\S II.(ii)]{KM1995structure} (more originally in the proofs of \cite[Lemma (3.3)(i)]{Donaldson1990polynomial} and \cite[Theorem (9.2.12)]{DK1992instanton}), we need to deal with $V(\wti{U}_p)$ one by one, and hence a choice of order of the indices in $\un{r}$ is needed. More precisely, we first suppose the data in $\wti{\tau}(I)$ only differ in the index $q$, and suppose the data in $\wti{\tau}$ are trivial for $p\in\un{r}\backslash \{q\}$. Then we define
\begin{equation}\label{eq: pairing 2 family}
    \langle CI(\bW,\bfa, z,\wti{\tau}(I))(x_0),x_1\rangle=\# \ov{M}_d(x_0,\bW(I),x_1)\cap \bigcap_{p\in\un{r}\backslash \{q\}} \big(I\times V(\be_{p})\big)\cap V(\wti{U}_{q}),
\end{equation}up to multiplication by the overall coefficient $c$ and the coefficients $b_p$ in Table \ref{tab: representative}, where $d$ is the number in \eqref{eq: pairing cob}.

Again by the construction of $V(\wti{U}_p)$, no bubble terms appear in the intersection. It is a homotopy between $CI(\bW,\bfa,z,\tau^0)$ and $CI(\bW,\bfa,z,\tau^1)$ by the equation
\begin{equation}\label{eq: pairing 3 family}
    0=\# \partial \left(\ov{M}_{d+1}(x_0,\bW(I),x_1)\cap \bigcap_{p\in\un{r}\backslash \{q\}} \big(I\times V(\be_{p})\big)\cap V(\wti{U}_{q})\right).
\end{equation}

In the general case, we use the order of $\un{r}$ to construct the homotopies for each index as in \eqref{eq: pairing 2 family} and take the sum. To study the dependence of the order, we consider the following generalization of Construction \ref{cons: ordered cylinder model}.

\bdefn\label{defn: family of tau}
For $r\in\N$ and a subset $S\subset \un{r}$, an enhanced cobordism $(\bW,\bfa,z,\wti{\tau}(I^S))$ consists of the following.

\benum
\item $(\bW,\bfa)$ is a morphism in $\BIa$ with trivial $\Phi$, i.e.\ a cobordism with auxiliary data.
\item $z=c\prod_{p=1}^r\be_p\in \bA(W;1/2)$ with $\deg\be_p\in\{0,2\}$.
\item We consider the enhanced data $\tau^i$ for $i=0,1$ and $\wti{\tau}$ from Constructions \ref{cons: tau} and \ref{cons: tau tilde}, respectively, and use $\tau^0_p,\tau^1_p$, and $\wti{\tau}_p$ to denote the data involving $\be_p$. We assume that $\tau^0_p$ is different from $\tau^1_p$ and $\wti{\tau}_p$ is nontrivial only for $p\in S$. For $p\in \un{r}\bs S$, we write $\tau_p$ for either $\tau^0_p$ or $\tau^1_p$.
\item For each $p\in S$, use $\iota_p$ from Construction \ref{cons: tau tilde} to identify the two pullback bundles over $\cB^*(\wti{U}_p)$. Choose sections $\wti{s}_{p,1},\wti{s}_{p,2}$ over $I\times\cB^*(\wti{U}_p)$, where $I$ has coordinate $t_p$, restricting to the prescribed sections $s_{p,1}^i,s_{p,2}^i$ at $t_p=i$, and pull them back along $I^S\to I$. In degree $2$, keep the zero section fixed. Pull these data back to
\begin{equation}\label{eq: underlying config}
I^S\times\cB^*(\bigcup_{p\in S}\wti{U}_p).
\end{equation}
\eenu
For a family cobordism $(\bW(G),\bfa)$, choose $Z_p^0,Z_p^1,\wti{Z}_p$ and the auxiliary paths, arcs and curves independently of $G$ on each cobordism piece, jointly transverse across distinct factors, with the submanifolds disjoint from the cuts. Choose the boundary curves of Construction \ref{cons: flat restrictions} on each cobordism piece and translate them with the neighborhoods on cylindrical ends. For cobordisms supported in disjoint regions, use a common witness surface and choose the witness curves and their connecting paths in the common product complement. 

Take fixed regular neighborhoods $U_p,\wti{U}_p$ as in Constructions \ref{cons: tau} and \ref{cons: tau tilde}, sufficiently small to preserve the intersection pattern of the corresponding unions. The metric and perturbation data
are independent of $I^S$. The $p$-th sections depend only on $G$ and, for $p\in S$, on $t_p$. These data define an \emph{enhanced family cobordism} $(\bW(G),\bfa,z,\wti{\tau}(I^S))$, including the preceding case $G=\pt$.

Fix a common integer $D>1+\sum_{p=1}^r(4-\deg\be_p)$ before choosing the sections. Choose the endpoint sections for all factors jointly, then choose the factor families simultaneously and obtain the data on faces of $I^S$ by restriction. On faces of $G$, all data agree with those already chosen. On broken product strata, the bundles and sections satisfy the pullback conditions of \eqref{eq: abstract auxiliary face}, and $\iota_p$ is induced by the identifications on the pieces, compatibly with further breaking. For every $T\subseteq\un{r}$, the defining sections for $p\in T$ are jointly transverse to the relevant rank strata or zero sections on all moduli strata of total dimension less than $D$ over $G\times I^{T\cap S}$, including lower-charge spaces, parameter faces, and broken strata. Here total dimension includes all parameter directions.
\edefn

For the unbroken cobordisms in a neck-stretching family, use Construction \ref{cons: flat restrictions} for the whole cobordism. Restriction to $U_p$ is irreducible by applying the argument preceding \eqref{eq: restriction} to the boundary curves in an unperturbed external collar. Near broken strata, use pullbacks from the neighborhoods on the pieces as in Construction \ref{cons: tau tilde}, where irreducibility follows from the same argument and compactness.

\bdefn\label{defn: special enhanced cobordism}
An enhanced cobordism $(\bW,\bfa,z,\wti{\tau}(I^S))$ is called \emph{special} if the following conditions hold.
 \benum
 \item The subset $S=\un{r}$.
 \item The cobordism $\bW$ is a product cobordism of an object $\bY$ in $\BI$ with trivial $\phi$.
 \item The submanifold $Z_p^i$ in $\tau_p^i$ for $i=0,1$ is in $\{i\}\times Y$ and the submanifold $\wti{Z}_p$ in $\wti{\tau}_p$ is a product. Hence we use $Z_p\subset Y$ to denote either $Z_p^i$.
 \eenu
\edefn
\brem\label{rem: translation family}
For a special enhanced cobordism, the underlying submanifolds at $\vec t=(t_1,\ldots,t_r)\in I^r$ are $\{t_p\}\times Z_p$ for $p\in\un{r}$. The sections form smooth families satisfying the transversality conditions in Definition \ref{defn: family of tau}. This observation will be useful in the construction of cube complexes for $\mu$ maps.
\erem

For each $p\in S$, define $V(\wti{U}_p)$ using the sections $\wti{s}_{p,1},\wti{s}_{p,2}$ over $I^S$ as in \eqref{eq: construction of V family}, and put
\[
V(\{\wti{U}_p\}_{p\in S})=\bigcap_{p\in S}V(\wti{U}_p).
\]
Replacing the corresponding factors in \eqref{eq: pairing 2 family} and \eqref{eq: pairing 3 family}, we obtain\begin{equation}\label{eq: pairing 2 family mor}
   \begin{aligned}
       \langle CI(\bW,\bfa, z,\wti{\tau}(I^S))(x_0),x_1\rangle=&\# \ov{M}_d(x_0,\bW(I^S),x_1)\cap \\&\bigcap_{p\in \un{r}\bs S} \big(I^S\times V(\be_{p})\big)\cap V(\{\wti{U}_p\}_{p\in S}),
   \end{aligned}
\end{equation}up to multiplication by the overall coefficient $c$ and the coefficients $b_p$ in Table \ref{tab: representative}, where $d$ is the number in \eqref{eq: pairing cob}, and $\ov{M}_d(x,\bW(I^S),y)$ is the compactification of the solution space in $I^S\times \cB_{\loc}(\bW^*)$. Similarly, we have an equation \begin{equation}\label{eq: pairing 3 family mor}
    0=\# \partial\left(\ov{M}_{d+1}(x_0,\bW(I^S),x_1)\cap \bigcap_{p\in \un{r}\bs S} \big(I^S\times V(\be_{p})\big)\cap V(\{\wti{U}_p\}_{p\in S})\right)
\end{equation}

\brem\label{rem: no bubble}
One may worry that bubble terms appear in the counting of \eqref{eq: pairing 2 family mor} and \eqref{eq: pairing 3 family mor}. However, this does not happen by Remark \ref{rem: moduli space diff} and the choices in Definition \ref{defn: family of tau}. Indeed, consider a zero- or one-dimensional intersection and suppose a limit has $\ell\ge1$ bubble points. Let $T$ index the factors whose neighborhoods contain no bubble point. Adapting the compactness argument following \cite[\S II.(ii), (2.8)]{KM1995structure}, the intersection defined by the remaining conditions over $G\times I^{T\cap S}$ has expected dimension at most
\[
1-8\ell+\sum_{p\notin T}(4-\deg\be_p)-|S\setminus T|
\le1-8\ell+5\ell=1-3\ell<0.
\]
Here $1$ bounds the original intersection dimension, and bubbling lowers the expected dimension by at least $8\ell$. The sum restores the codimensions of the conditions that need not survive at the bubble points. The term $-|S\setminus T|$ comes from forgetting the parameters $t_p$ for $p\in S\setminus T$, since the remaining conditions are independent of them.

To obtain the bound $5\ell$, note that each omitted factor contributes $4-\deg\be_p$ when $p\notin S$ and $3-\deg\be_p$ when $p\in S$, the codimensions of $Z_p$ and $\wti{Z}_p$, respectively. As in the neighborhood argument for \cite[\S II.(ii), (2.7)]{KM1995structure}, the neighborhood condition reduces the estimate to common intersections of the corresponding submanifolds and auxiliary curves. General position bounds the sum of geometric codimensions by $4$. An auxiliary curve of a fixed degree $0$ factor has codimension $3$ but contributes $4$, adding at most $1$ since auxiliary curves of distinct factors are disjoint. Thus each bubble point contributes at most $5$. By the transversality, the remaining intersection is empty. The same argument applies on parameter faces and broken strata.
\erem

We also provide the construction related to the additive relation in \eqref{eq: additive z} at the chain level, also only for $\deg\be_p\in\{0,2\}$.
\bcons\label{cons: add z}
Fix $q\in\un{r}$ and choose enhanced data $(z_i,\tau_i)$ for $i=1,2$, with
\[
z_i=c\cdot \be_{q,i}\prod_{p\ne q}\be_p,~c\in\Z[1/2],
\]
identical at every $p\ne q$. Assume $\deg\be_{q,1}=\deg\be_{q,2}\in\{0,2\}$ and $Z_{q,1}\cap Z_{q,2}=\varnothing$, with the usual transversality among the representatives. Use a common $D$, basepoint $x_q$, and collection of curves $C_q$. Set $Z_q=Z_{q,1}\sqcup Z_{q,2}$, use the connecting arcs from both collections, and choose a good suitable neighborhood $U_q$ containing both original neighborhoods. Pull the $q$-th bundles and sections back to $\cB^*(U_q)$.

Keeping the data at $p\ne q$ fixed, define
\begin{equation}\label{eq: sum of z}
(z,\tau)=(z_1,\tau_1)+(z_2,\tau_2),~z=c\cdot (\be_{q,1}+\be_{q,2})\prod_{p\ne q}\be_p=z_1+z_2,
\end{equation}
with the following section choices at $q$.

If $\deg\be_{q,i}=0$, use the same generic pair of sections for both summands and their sum. The degree $0$ coefficient in Table \ref{tab: representative}, interpreted using signed multiplicities, is additive under $Z_q=Z_{q,1}\sqcup Z_{q,2}$.

If $\deg\be_{q,i}=2$, use $\cL=\cL_1\ot\cL_2$ for $Z_q$, where $\cL_i$ is the pullback of the corresponding line bundle in Table \ref{tab: representative}. Choose $s'_1,s'_2$ jointly transverse together with the common factors. After imposing the conditions at $p\ne q$, the spaces relevant to the zero- and one-dimensional counts have dimensions $2$ and $3$. Thus $s'_1$ and $s'_2$ have no common zero there, since their simultaneous zero locus has codimension $4$. The section $s'_1\ot s'_2$ is therefore transverse there, with zero set the oriented disjoint union of the two zero sets. Choose the section for $\tau$ to agree with this tensor product near these loci and to satisfy the remaining transversality conditions.

For families over $G\times I^S$, make the same choices with the data at $p\ne q$ fixed and with the given restrictions on parameter faces and broken strata, as in Definition \ref{defn: family of tau}. The signed counts add, including the coefficients in Table \ref{tab: representative} and the common coefficient $c$.
\econs
With the definition in \eqref{eq: pairing 2 family mor} and the equation in \eqref{eq: pairing 3 family mor}, we have the following theorem for enhanced family cobordisms, which is the generalization of the usual cobordism maps in \cite[\S 3.8-3.9]{kronheimer2011khovanov}

\bthm\label{thm: instanton construction 2}
Suppose $(\bW(G),\bfa, z,\wti{\tau}(I^S)): (\bY_0,a_0)\to (\bY_1,a_1)$ is an enhanced family cobordism as in Definition \ref{defn: family of tau}. Then the following hold over $\Z[1/2]$.

\benum
    \item\label{enh1} There exists a map (called an \emph{enhanced family cobordism map} or an \emph{enhanced cobordism map} if $G=\pt$)\begin{equation}\label{eq: enhanced cobordism map}
    CI(\bW(G),\bfa,z,\wti{\tau}(I^S)): CI(\bY_0,a_0)\to CI(\bY_1,a_1).
\end{equation}with the relation \begin{equation}\label{eq: relation of cobordism map mu}
\begin{aligned}
    CI((\bW(G),\bfa,z,\wti{\tau}(I^S))|_{\partial (G\times I^S)})&+(-1)^{\dim (G\times I^S)}CI(\bW(G),\bfa,z,\wti{\tau}(I^S))\circ \partial \\&=\partial \circ CI(\bW(G),\bfa,z,\wti{\tau}(I^S)),
\end{aligned}
\end{equation}where the enhanced family cobordism map associated with \[\partial (G\times I^S)=\partial G\times I^S\cup G\times \partial I^S\]is the sum of enhanced family cobordism maps as $\partial I^S$ is the union of copies of $I^{|S|-1}$ along boundaries. If $z=1$, then $S=\varnothing$ and $CI(\bW(G),\bfa,z,\wti{\tau}(I^S))$ reduces to $CI(\bW(G),\bfa)$ in \S \ref{subsect: Sketch of gauge theory construction}.
    \item\label{enh2} The $\Z/2$-degree of \eqref{eq: enhanced cobordism map} is
\begin{equation*}
    \begin{aligned}
        \deg(CI(\bW(G),\bfa,z,\wti{\tau}(I^S)))&=\deg(CI(\bW(G),\bfa))+\sum_{p=1}^r(4-\deg\be_p)+|S| ,\\&\equiv \dim G+|S|\pmod 2
    \end{aligned}
\end{equation*}where last equation follows from \eqref{eq: J family degree} and the fact $\deg\be_p\in\{0,2\}$.
    
    \item\label{enh3} If $(z,\tau)=(z_1,\tau_1)+(z_2,\tau_2)$ as in Construction \ref{cons: add z}, then we have\begin{equation*}\label{eq: sum of z cob}
        CI(\bW(G),\bfa,z,\wti{\tau}(I^S))=CI(\bW(G),\bfa,z_1,\wti{\tau}_1(I^S))+CI(\bW(G),\bfa,z_2,\wti{\tau}_2(I^S)).
    \end{equation*}
\item\label{enh4} On broken product families, the enhanced maps satisfy the composition formula in \eqref{eq: J instanton family identities}, under the pullback conditions in Definition \ref{defn: family of tau}. For disjoint unions, choose the enhanced data componentwise and expand each homology factor over the components. For each resulting term, write $\bW(G)=\bigsqcup_{\alpha=1}^m\bW_\alpha(G_\alpha)$, where $G=\prod_\alpha G_\alpha$, $I^S=\prod_\alpha I^{S_\alpha}$, and $z=\prod_\alpha z_\alpha$. Pull all data back from the corresponding factors. With the product orientation ordered by the factors $G_\alpha\times I^{S_\alpha}$ and the identifications \eqref{eq: J disjoint union complexes}, we have
\[
\iota_+\circ CI(\bW(G),\bfa,z,\wti{\tau}(I^S))
=\Bigl(\bigotimes_{\alpha=1}^m CI(\bW_\alpha(G_\alpha),\bfa_\alpha,z_\alpha,\wti{\tau}_\alpha(I^{S_\alpha}))\Bigr)\circ\iota_-.
\]
Both formulas are compatible with parameter faces and further breaking.
\eenu
\ethm
\subsection{Cube complexes from enhanced cobordisms}\label{subsec: Cube complexes for enhanced cobordisms}
In this subsection, we construct the cube complex related to enhanced cobordism \[(\bW,\bfa,z,\wti{\tau}(I^r))\]from Definition \ref{defn: family of tau}, where $I^r$ denotes $I^{\un{r}}$ for $\un{r}=\{1,\dots,r\}$. Recall that we restrict the choice of $z$ such that\[z\in\bA_{\e}(W;1/2)=\sym^*(H_{\mathrm{even}}(W;\Z[1/2]))=\sym^*(H_0(W;\Z[1/2])\op H_2(W;\Z[1/2])).\]

\bdefn\label{defn: enhanced supported}
An enhanced cobordism $(\bW,\bfa,z,\wti{\tau}(I^r))$ is called \emph{supported} in $\sR$ if it satisfies the following.
\benum
    \item $\bW$ is supported in $\sR$ in the sense of Definition \ref{defn: supported region};
    \item Recall that the enhanced data $\wti{\tau}(I^r)$ contains the information that for any $p\in \un{r}$, a $1$- or $3$-submanifold $\wti{Z}_p$  with $\partial \wti{Z}_p=- Z_p^0\sqcup  Z_p^1$, and also a suitable neighborhood $\wti{U}_p$ of $\wti{Z}_p$. The remaining data are omitted. Note that\[z=c\prod_{p=1}^r\be_p=c\prod_{p=1}^r[Z_p^0]=c\prod_{p=1}^r[Z_p^1]\in\sym^*(H_{\mathrm{even}}(W;\Z[1/2]))\] is determined by $c$ and $\wti{\tau}(I^r)$. Furthermore, for any $p\in\un{r}$, we need that $Z_p^0$ and $Z_p^1$ lie in the incoming and outgoing ends of $W$, respectively, and $\wti{Z}_p$ is a product outside $I\times \sR$.
\eenu
In particular, the special enhanced cobordism in Definition \ref{defn: special enhanced cobordism} is supported in the empty set.
\edefn
\brem
    Since later we will only consider an enhanced cobordism $(\bW,\bfa,z,\wti{\tau}(I^r))$ supported in a region $\sR$, we simply write it as \begin{equation}\label{eq: enhanced cob short}
    (\bW,\{\wti{Z}_p\}_{p\in\un{r}}):(\bY,\{Z_p^0\}_{p\in\un{r}})\to (\bY_1,\{Z_p^1\}_{p\in\un{r}})
\end{equation}and fix certain choices of the rest information. The coefficient $c$ is omitted from this notation, since the cube constructions below concern the individual operators $\mu(\be_p)$. Indeed, we leave the study of the dependence of the rest information to the future.
\erem

We start with the construction for a special enhanced cobordism and then combine it with Construction \ref{cons: cube of cobordisms} to obtain the most general case. Recall from Definition \ref{defn: BNr} the inner anodyne inclusion \[(\Delta^1/\partial\Delta^1)^r\hookrightarrow\BN^r\]and the notation \[\D[1/2]=\D_2(\Mod(\Z[1/2])^\fin).\]

The following theorem is a generalization of \eqref{eq: commute} at the chain level for multiple $\mu$-operators.

\bthm\label{thm: mu Y}
 For a special enhanced cobordism $(\bW,\bfa,z,\wti{\tau}(I^r))$ from Definition \ref{defn: special enhanced cobordism}, there is a functor\begin{equation}\label{eq: mu Y}
    \mu[\bY,\{Z_p\}_{p\in\un{r}}] =(CI(\bY),\mu(Z_p)_{p\in\un{r}})\in\Fun(\BN^r,\D[1/2]),
\end{equation}where the notation is described below.
\benum
\item The notation in \eqref{eq: mu Y} suppresses the higher morphisms encoding the homotopy coherent commutativity of the endomorphisms.

\item The notation $\mu(Z_p)$ denotes a lift of $\mu([Z_p])$ from Construction \ref{cons: mu maps} at the chain level.
\eenu
\ethm
\bpf
By Definition \ref{defn: BNr}, it suffices first to construct a functor $(\Delta^1/\partial\Delta^1)^r\to\D[1/2]$, equivalently a functor $\sX:\bbox^r\to\D[1/2]$ whose pairs of opposite $(r-1)$-faces agree, and then extend it along the inner anodyne inclusion $(\Delta^1/\partial\Delta^1)^r\hookrightarrow\BN^r$. The construction of the cube complex follows from similar construction as in \cite[\S 4-7]{Bloom2011link}, \cite[\S 5-6]{scaduto2015instantons}, and \cite[\S 6-7]{kronheimer2011khovanov}. For further generalization to cobordisms supported in disjoint regions, we follow the construction in \cite[\S 6]{kronheimer2011khovanov}, but now we consider families of sections in $\wti{\tau}(I^r)$ instead of families of metrics.

Note that the cobordism $\bW$ is a product cobordism in the case of special cobordism. From Remark \ref{rem: translation family}, the submanifolds $\{t_p\}\times Z_p$ vary by translation, while the sections form smooth parametrized families. See Figure \ref{fig:cube2} without the small squares about $\bW_i$. Note that we do not need the assumption that $Z_p$ are disjoint in $Y$, because the map for the enhanced family cobordism is still well-defined from Remark \ref{rem: no bubble} when all $Z_p$ are in the generic positions (i.e.\ no bubble appears in the moduli spaces). Choose the auxiliary arcs and curves used to construct the suitable neighborhoods jointly with the $Z_p$ in general position in $Y$.

Furthermore, we consider $\bW^*$ with the underlying manifold $\R\times Y$ and extend the family of sections over $I^r$ to $\R^r$ by considering the translations into the attached cylindrical ends. Similar to the case of the family of metrics in \cite[\S 6]{kronheimer2011khovanov}, when the difference of coordinates in $\R^r$ goes to $\infty$, we can compactify the family of sections by broken cobordisms, which correspond to compositions of $\mu(Z_p)$ in some order.

For fixed choices of $u,v\in\{0,1\}^r$ with $v<u$, write \[i_1,\dots,i_d\]for the coordinates where $u$ and $v$ differ. Consider the normalized families $\breve{G}_{vu}$ of sections parametrized by
\[
(\tau_{i_1},\dots,\tau_{i_d})\in\R^d\subset\R^r~\mathrm{with}~\sum_{j=1}^d\tau_{i_j}=0,
\]
and their compactifications $\breve{G}_{vu}^{+}$ by broken families as in \cite[(33)]{kronheimer2011khovanov}. Use the cubical subdivision of $\breve{G}_{vu}^{+}$ into cubes $Q_I$ indexed by the vertex chains of maximal simplices in the subcube from $v$ to $u$ from Construction \ref{cons: hypercube to cube complex}.

Translate the suitable neighborhoods together with finite unperturbed product cylinders containing them. Choose the cutoff functions relative to the fixed data and faces, using \eqref{eq: perturbed ASD 3} between these cylinders near breaking. The remaining end collars support the secondary perturbations, which are chosen by the proof of Proposition \ref{prop: extension of perturbation}, keeping the primary perturbations and the data without insertions fixed. The restriction argument preceding \eqref{eq: restriction} applies on these unperturbed cylinders.

Using the bundle identifications induced by translation, choose the sections for the individual factors jointly by parametric transversality as in Definition \ref{defn: family of tau} for every subcollection on each $Q_I$ and its faces. Use the same data on identified faces and product data on broken faces.

We apply the argument of Remark \ref{rem: no bubble} on each $Q_I$ and its faces. Choose the suitable neighborhoods sufficiently small so that their projections to $Y$ preserve the intersection pattern of the submanifolds, arcs and curves. Consider a zero- or one-dimensional intersection and suppose a limit has $\ell\ge1$ bubble points. The defining conditions whose neighborhoods contain no bubble point remain satisfied by the limiting connection. General position in $Y$ bounds the total codimension of the omitted conditions by $6\ell$, since each bubble point can affect at most three degree $2$ conditions, or one degree $0$ condition and one degree $2$ condition. The remaining intersection therefore has expected dimension at most
\[
1-8\ell+6\ell=1-2\ell<0.
\]
By transversality, this intersection is empty, so bubbles do not occur.

For the factors $\be_{i_1},\dots,\be_{i_d}$, the counting argument of Theorem \ref{thm: instanton construction 2} therefore gives the counts and boundary identities on each $Q_I$. With Constructions \ref{cons: cube complex to hypercube}, \ref{cons: sign on hypercube}, and \ref{cons: hypercube to cube complex}, their sums give the hypercube maps, while the individual counts give the morphisms and higher morphisms of a cube complex. Use the same choices whenever the same factors occur, so that opposite cube faces agree. This gives the desired cube complex $\mu[\bY,\{Z_p\}_{p\in\un{r}}]$.

For homogeneous $\beta_p\in H_{\mathrm{even}}(Y;\Z[1/2])$, choose $N\geq0$ such that $2^N\beta_p=[Z_p]$ for every $p$. In the cube constructed above, multiply each dg nerve component from $v$ to $u$ by $2^{-N|u-v|}$. This preserves the dg nerve identities and the agreement of opposite faces, and gives the required operators $\mu(\beta_p)$.
\epf

Finally, we obtain the most general result that combines Construction \ref{cons: cube of cobordisms}, Theorem \ref{thm: surgery fiber sequence}, and Theorem \ref{thm: mu Y}.
\bthm[]\label{thm: mu cube}
Suppose $(\bY=(Y,P,\phi,J),a)$ is an object in $\BIa$ with trivial $\phi$ and there are enhanced cobordisms \[(\bW_i,\{\wti{Z}_{p,i}\}_{p\in\un{r}}):(\bY,\{Z_{p,i}^0\}_{p\in\un{r}})\to (\bY_i,\{Z_{p,i}^1\}_{p\in\un{r}})\] for $i\in \un{n}$ supported in disjoint regions $\sR_i\subset Y$ as in Definition \ref{defn: enhanced supported}. We suppose $Z_{p,i}^0$ is independent of $i$ and written as $Z_p^0$ and suppose that there are no quadruple intersections among $3$-dimensional members of $\{\wti{Z}_{p,i}\}_{p\in\un{r}}$ for each $i\in\un{n}$. Assume that all cobordisms $\bW_{S,S'}$ in Construction \ref{cons: cube of cobordisms} have trivial $\Phi$. Then there is a functor\begin{equation}\label{eq: end}
    \ov{\sX}=\ov{\sX}[\bY,\{\bW_i\}_{i\in\un{n}},\{\wti{Z}_{p,i}\}_{p\in\un{r},i\in\un{n}}] \in\Fun(\BN^r\times \bbox^n,\D[1/2])
\end{equation}satisfying the following properties.

\benum
\item\label{mu1} As an object in $\Fun(\BN^r,\Fun(\bbox^n,\D[1/2]))$, the functor $\ov{\sX}$ sends the unique object in $\BN^r$ to $CI\circ\sX[\bY,\{\bW_i\}_{i\in\un{n}}]$ from Construction \ref{cons: cube of cobordisms} and Theorem \ref{thm: TQFT functor}.
\item\label{mu2} As an object of $\Fun(\bbox^n,\Fun(\BN^r,\D[1/2]))$, we have\[\ov{\sX}(\varnothing)=\mu[\bY,\{Z_{p}^0\}_{p\in\un{r}}]\aand \ov{\sX}(\{i\})=\mu[\bY_i,\{Z_{p,i}^1\}_{p\in\un{r}}]\]for notation from \eqref{eq: mu Y}.
\item\label{mu3} If $\sX=\sX[\bY,\{\bW_i\}_{i\in\un{n}}]=\sX(\bY,\{K_i\}_{i\in\un{n}})$ is a surgery $n$-cube as in Example \ref{exmp: surgery cube}, with $Z_p\subset Y\backslash\bigcup_i\sR_i$ and $\wti{Z}_{p,i}=I\times Z_p$ in the common product region for every $p,i$, then for a fixed $q\in\un{n}$, the $(n-1)$-cube complexes \[\sY_j\in \Fun(\bbox^{\un{n}\backslash\{q\}},\D)\aand \sW_j:\sY_j\to \sY_{j+1}\]for $j\in \Z$ from Construction \ref{cons: surgery term} all have lifts \[\ov{\sY}_j\in\Fun(\BN^r,\Fun(\bbox^{\un{n}\backslash\{q\}},\D[1/2]))\aand \ov{\sW}_j:\ov{\sY}_j\to \ov{\sY}_{j+1},\]respectively. Moreover, the fiber square \eqref{eq: fiber sequence} lifts to \begin{equation}\label{eq: fiber sequence mu}
\xymatrix{\ov{\sY}_j\ar[r]^{\ov{\sW}_j}\ar[d]&\ov{\sY}_{j+1}\ar[d]^{\ov{\sW}_{j+1}}\\0\ar[r]&\ov{\sY}_{j+2}}
\end{equation}for any $j\in \Z$.
\eenu
\ethm
\bpf
For item (\ref{mu3}), combine the metric and perturbation families, including the nullhomotopies, from the proof of Theorem \ref{thm: surgery fiber sequence} with the section families of Theorem \ref{thm: mu Y}. Fix suitable neighborhoods and bundle identifications in the common product complement of the surgery regions. Choose sections relative to parameter faces, agreeing on identified faces and satisfying the pullback conditions of Definition \ref{defn: family of tau} on broken faces. Bubbling is excluded as in Theorem \ref{thm: mu Y}. For $r=1$, compare Xie's $\mu$-action on the total resolution complex in singular instanton homology \cite[\S 5.1]{Yi2021spectralsequence}.

At each additional $S^3$ face, fix the metric and perturbation on the local punctured $\ov{\mathbb{CP}^2}$, with all section conditions imposed on the complementary cobordism and all extra parameters varying only its data. Nonempty cut-down moduli spaces on the complement have nonnegative total dimension. The index and gluing formulae in \cite[\S 5.2]{scaduto2015instantons} therefore give dimension at least $0-1+3=2$, excluding these boundary terms. The remaining identities give nullhomotopies $\ov{\sW}_{j+1}\circ\ov{\sW}_j\simeq0$ over $(\Delta^1/\partial\Delta^1)^r$, compatible with all remaining surgery directions. Extend these diagrams together to $\BN^r$ as in Theorem \ref{thm: mu Y}, identifying common objects and maps. Their cofiber comparisons are isomorphisms pointwise by Theorem \ref{thm: surgery fiber sequence}, so \eqref{eq: fiber sequence mu} is a cofiber square, which is also a fiber square by stability.

The general case of items (\ref{mu1}) and (\ref{mu2}), not necessarily a surgery cube, follows from a similar strategy. We again consider both families of metrics and families of sections over $\R^n\times \R^r$. To be clear, we first consider the case $r=0$.

Let \[W:Y_{\vec{0}}=Y\to Y_{\vec{1}}\]be obtained from $I\times Y$ by attaching handles in all $\bW_i$, called the \emph{total cobordism}. Let $W^*$ be obtained from $W$ by attaching infinite half cylinders to the ends as in \eqref{eq: open end}. Following the discussion near \cite[Fig. 8]{kronheimer2011khovanov}, we start with a metric $\bfg$ on $W^*$ that is the product metric outside $\bigcup_{i=1}^n\R\times \sR_i$ and construct a family of metrics over $\R^n$ from $\bfg$ by modifying the metrics in $\bigcup_{i=1}^n\R\times \sR_i$ via translation, where the $i$-th coordinate in $\R^n$ corresponds to the translation in $\R\times \sR_i$. The original metric $\bfg$ corresponds to the origin in $\R^n$. See Figure \ref{fig:cube1} for an illustration similar to \cite[Fig. 8]{kronheimer2011khovanov}.

\begin{figure}[htbp]
\centering
\resizebox{0.65\textwidth}{!}{%
\begin{circuitikz}
\tikzstyle{every node}=[font=\LARGE]
\draw  (6.25,15.25) rectangle (19.25,7);
\draw (6.25,14) to[short] (19.25,14);
\draw (6.25,12.75) to[short] (19.25,12.75);
\draw (6.25,11.5) to[short] (19.25,11.5);
\draw (6.25,10.25) to[short] (19.25,10.25);
\draw (6.25,9) to[short] (19.25,9);
\draw (6.25,7.75) to[short] (19.25,7.75);
\draw (8.25,14) to[short] (8.25,12.75);
\draw (9.75,14) to[short] (9.75,12.75);
\draw (15.5,11.5) to[short] (15.5,10.25);
\draw (17,11.5) to[short] (17,10.25);
\draw [dashed] (19.25,15.25) -- (20.75,15.25);
\draw [dashed] (19.25,7) -- (20.75,7);
\draw [dashed] (6.25,15.25) -- (4.75,15.25);
\draw [dashed] (6.25,7) -- (4.75,7);
\draw (11.75,9) to[short] (11.75,7.75);
\draw (13.25,9) to[short] (13.25,7.75);
\node [font=\LARGE] at (12.5,8.25) {$\bW_3$};
\node [font=\LARGE] at (5.25,13.25) {$\sR_1$};
\node [font=\LARGE] at (9,13.25) {$\bW_1$};
\node [font=\LARGE] at (16.25,10.75) {$\bW_2$};
\node [font=\LARGE] at (5.25,10.75) {$\sR_2$};
\node [font=\LARGE] at (5.25,8.25) {$\sR_3$};
\node [font=\LARGE] at (12.5,6.25) {$[0,1]$};
\node [font=\LARGE] at (6.25,15.75) {$\bY$};
\node [font=\LARGE] at (5.25,6.25) {$(-\infty,0]$};
\node [font=\LARGE] at (19.75,6.25) {$[1,\infty)$};
\end{circuitikz}
}%
\caption{A schematic picture for $n=3$ and some $\tau\in\R^3$. We use small squares with $\bW_i$ to denote the handles in $\bW_i$. Moving the small squares corresponds to modifying the metrics. \label{fig:cube1}}
\end{figure}

\begin{figure}[htbp]
\centering
\resizebox{0.65\textwidth}{!}{%
\begin{circuitikz}
\tikzstyle{every node}=[font=\LARGE]
\draw  (6.25,15.25) rectangle (19.25,7);
\draw (6.25,14) to[short] (19.25,14);
\draw (6.25,12.75) to[short] (19.25,12.75);
\draw (6.25,11.5) to[short] (19.25,11.5);
\draw (6.25,10.25) to[short] (19.25,10.25);
\draw (6.25,9) to[short] (19.25,9);
\draw (6.25,7.75) to[short] (19.25,7.75);
\draw (8.25,14) to[short] (8.25,12.75);
\draw (9.75,14) to[short] (9.75,12.75);
\draw (15.5,11.5) to[short] (15.5,10.25);
\draw (17,11.5) to[short] (17,10.25);
\draw [dashed] (19.25,15.25) -- (20.75,15.25);
\draw [dashed] (19.25,7) -- (20.75,7);
\draw [dashed] (6.25,15.25) -- (4.75,15.25);
\draw [dashed] (6.25,7) -- (4.75,7);
\draw (11.75,9) to[short] (11.75,7.75);
\draw (13.25,9) to[short] (13.25,7.75);
\draw (10.25,14.75) to[short] (10.25,9.5);
\draw (13.75,12) to[short] (13.75,7.25);
\draw (17.25,14.5) to[short] (17.25,12);
\node [font=\LARGE] at (12.5,8.25) {$\bW_3$};
\node [font=\LARGE] at (5.25,13.25) {$\sR_1$};
\node [font=\LARGE] at (9,13.25) {$\bW_1$};
\node [font=\LARGE] at (16.25,10.75) {$\bW_2$};
\node [font=\LARGE] at (5.25,10.75) {$\sR_2$};
\node [font=\LARGE] at (5.25,8.25) {$\sR_3$};
\node [font=\LARGE] at (12.5,6.25) {$[0,1]$};
\node [font=\LARGE] at (6.25,15.75) {$\bY$};
\node [font=\LARGE] at (5.25,6.25) {$(-\infty,0]$};
\node [font=\LARGE] at (19.75,6.25) {$[1,\infty)$};
\node [font=\LARGE] at (11,14.5) {$Z_{1}^{1/3}$};
\node [font=\LARGE] at (14.5,12) {$Z_{3}^{2/3}$};
\node [font=\LARGE] at (18,14.5) {$Z_{2}^{1}$};
\end{circuitikz}
}%
\caption{A schematic picture for $n=r=3$ and some $\tau\in\R^3\times\R^3$. We use small squares with $\bW_i$ to denote the handles and use vertical lines to denote submanifolds. Moving the small squares corresponds to modifying the metrics and moving the vertical lines corresponds to modifying the sections. \label{fig:cube2}}
\end{figure}
Then we consider the most general case. See Figure \ref{fig:cube2}. We point out the main modifications as follows.

\benum
\item Choose the auxiliary paths, arcs and curves used to construct the suitable neighborhoods jointly with the $\wti{Z}_{p,i}$ in general position on each cobordism piece, as in Construction \ref{cons: tau tilde} and Definition \ref{defn: family of tau}.
\item Once we fix the coordinates in $\R^n$, the total cobordism $W$ is obtained by attaching handles in $\{\bW_i\}_{i\in\un{n}}$ in some order of $\un{n}$. Then we consider the total cobordism \begin{equation*}\label{eq: total Z}
    \wti{Z}_p: Z_p^{\vec{0}}=Z_p^0\to Z_{p}^{\vec{1}}
\end{equation*}obtained by attaching handles in $\{\wti{Z}_{p,i}\}_{i\in\un{n}}$ in the same order of $\un{n}$. For $j\in\un{n}$, let\[\wti{Z}_p^{j/n}: Z_p^{(j-1)/n}\to Z_p^{j/n}\]be the part of $\wti{Z}_p$ corresponding to the $j$-th element in the given order of $\un{n}$. For any choice of $(j_p)_{p\in\un{r}}\in \un{n}^r$, we consider a family of sections over $I^r$ from Definition \ref{defn: family of tau} corresponding to $\{\wti{Z}_p^{j_p/n}\}_{p\in\un{r}}$. We glue those families together along boundaries to obtain a family of sections over $[0,1]^r$. More precisely, these sections are pulled back to the corresponding bundles over\[\cB^*(\bigcup_{p=1}^r\wti{U}_p),\]where $\wti{U}_p$ is the suitable neighborhood of $\wti{Z}_p$ (rather than just $\wti{Z}_p^{j/n}$). The isomorphism $\iota_p$ in \eqref{eq: isomorphism iota} of pullback bundles is chosen to be the composition of isomorphisms related to $\wti{Z}_p^{j/n}$ for all $j\in\un{n}$.
\item Note that the family of sections over $[0,1]^r$ cannot be interpreted as translations of $\{Z_p^{\vec{0}}\}_{p\in \un{r}}$ because $\wti{Z}_p$ is a general cobordism of submanifolds. However, we can extend the family of sections over $\R^r$ by considering translations of $\{Z_p^{\vec{0}},Z_p^{\vec{1}}\}_{p\in\un{r}}$ in the ends of $W^*$.
\eenu
As in Theorem \ref{thm: mu Y}, normalize the combined parameters and subdivide the compactification into the cubes $Q_I$ of Construction \ref{cons: hypercube to cube complex}. Choose the cutoff functions and perturbations as in that proof, and choose the section families jointly, satisfying the transversality and compatibility conditions of Definition \ref{defn: family of tau} on each $Q_I$ and its faces, with the same choices for the same factors.

Choose sufficiently small suitable neighborhoods on each cobordism piece and use their translates on the ends. The assumption on the $\wti{Z}_{p,i}$ excludes quadruple intersections among the degree $2$ neighborhoods. At a triple intersection, three degree $2$ conditions may be lost, with total codimension $6$. General choices for suitable neighborhoods also allow at most one degree $0$ condition and one degree $2$ condition to be lost at the same point, again with total codimension $6$. This gives the bound $6\ell$ for $\ell$ bubble points. As in Theorem \ref{thm: mu Y}, the remaining intersection for a zero- or one-dimensional count over any $Q_I$ or its faces has expected dimension at most $1-8\ell+6\ell=1-2\ell<0$ for $\ell\ge1$, so it is empty by transversality. The counting argument of Theorem \ref{thm: instanton construction 2} therefore gives the counts and boundary identities on these cubes. These define the required cube complex, whose opposite faces in the section directions agree. Extend to $\BN^r$ as in Theorem \ref{thm: mu Y}. Then items (\ref{mu1}) and (\ref{mu2}) follow by restriction.
\epf

\bexmp\label{exmp: excision surface intersections}
An example relevant to Theorem \ref{thm: mu cube} comes from comparing sutured closures by Floer excision, as in \cite{kronheimer2010knots,BS2015naturality,li2019decomposition}. Let
\[
W:Y\sqcup Y_\varphi\longrightarrow Y'
\]
be such a cobordism, where $Y_\varphi$ is the mapping torus of a diffeomorphism $\varphi$ of the cutting surface $R$. Suppose the surfaces $Z_p\subset Y$ are jointly in general position and simultaneously in product form near $R$, which is transverse to all their intersection strata. Assume that $\varphi$ preserves each collection of oriented curves $Z_p\cap R$. These curves determine vertical tori $F_p\subset Y_\varphi$, and excision produces surfaces $Z_p'\subset Y'$ and cobordisms
\[
\wti Z_p:Z_p\sqcup F_p\longrightarrow Z_p'.
\]

Using the excision model of \cite[\S 3.2]{kronheimer2010knots}, each $\wti Z_p$ is a product over the cut-open incoming manifolds and equals $(Z_p\cap R)\times P$ in the piece $R\times P$, where $P$ is the saddle surface. General position in the incoming $3$-manifolds excludes quadruple intersections in the first pieces, while general position of the curves on $R$ excludes triple intersections in the second. Hence the $\wti Z_p$ have no quadruple intersections. Furthermore, choose small parallel copies of $R$ inside the cut-open pieces of $Y$ and $Y_\varphi$, near boundary components identified in $Y'$. Their product cylinders give two disjoint product cobordisms to parallel copies of $R$ in $Y'$, entirely contained in the product regions over the cut-open manifolds. Adjoining these two cylinders introduces no quadruple intersections, since sufficiently small parallel copies of $R$ avoid all triple intersections of the incoming surfaces.
\eexmp
\newpage
\appendix
\section{Proof of the Segal condition}\label{sec: segal}
In this appendix, we prove Proposition \ref{prop: segal space}. The affine case is mainly an adaptation of the proof of \cite[Proposition 5.19]{CS2019category} to the $(\infty,1)$-categorical case, in which the globularity condition is vacuous. As mentioned in Remark \ref{rem: special proof}, we must first show that the affine path model computes the relevant homotopy pullbacks, for which we use the following Rezk's criterion.

\blem[Rezk's pullback criterion]\label{lem: Rezk pullback criterion}
Let $p:E\to B$ be a map of simplicial sets. For a simplex $\sigma:\Delta^q\to B$, set $E_\sigma=\Delta^q\times_BE$. Suppose that for every simplex $\sigma:\Delta^q\to B$ and every simplicial operator $\alpha:\Delta^r\to\Delta^q$, the induced map $E_{\sigma\alpha}\to E_\sigma$ is a weak homotopy equivalence. Then, for every $Z\to B$, the canonical map
\[
Z\times_BE\longrightarrow Z\times_B^hE
\]
is a weak homotopy equivalence.
\elem
\bpf
The hypothesis is condition (3) of \cite[Remark 4.2]{Rezk1998fibrations}, and hence is equivalent to $p$ being sharp in the sense of \cite[\S 2]{Rezk1998fibrations}. Since the Kan--Quillen model category of simplicial sets is right proper, \cite[Proposition 2.7]{Rezk1998fibrations} implies that, for every map $Z\to B$, the strict pullback square
\[
\xymatrix{Z\times_BE\ar[r]\ar[d]&E\ar[d]^p\\Z\ar[r]&B}
\]
is homotopy Cartesian. By Definition \ref{defn: homotopy fiber product}, the canonical map
\[
Z\times_BE\longrightarrow Z\times_B^hE
\]
is therefore a weak homotopy equivalence.
\epf

Choose a smooth nondecreasing function $\rho:\R\to[0,1]$ which is $0$ for $t\leq1/3$ and $1$ for $t\geq2/3$. It gives a smooth concatenation $a\star b$ of composable smooth paths by
\[
(a\star b)(t)=\begin{cases}a(\rho(2t)),&0\leq t\leq1/2,\\b(\rho(2t-1)),&1/2\leq t\leq1.
\end{cases}
\]
The two pieces are stationary near $t=1/2$, so this construction is smooth in all parameters. Reversal and the usual contractions of $a\star\bar a$ and $\bar a\star a$ are smooth as well.

\blem\label{lem: affine endpoint evaluation}
Let $X$ be a diffeological space, let $\mathcal PX=C^\infty([0,1],X)$ carry the functional diffeology, and set
\[
S=S^{\aff}_\bu,~ X_\bu=S(X),~ P_\bu=S(\mathcal PX),~ e_\bu=S(e_0,e_1):P_\bu\longrightarrow X_\bu^2.
\]
Then $e_\bu$ satisfies the hypothesis of Lemma \ref{lem: Rezk pullback criterion}. More precisely, for every simplex $\sigma:\Delta^q\to X_\bu^2$ and every simplicial operator $\delta:\Delta^r\to\Delta^q$, the induced morphism of simplicial sets
\[
E_{\sigma\delta}\longrightarrow E_\sigma,~ E_\sigma=\Delta^q\times_{X_\bu^2}P_\bu,
\]
is a weak homotopy equivalence.
\elem
\bpf
Let $\sigma$ be represented by $(\sigma_0,\sigma_1):|\Delta^q|_{\aff}\to X^2$, let $\delta:\Delta^r\to\Delta^q$ be a simplicial operator, and write $|\delta|$ for its affine realization. Set
\[
\mathcal P_\sigma X=\{h\in\mathcal PX\mid h(i)=\sigma_i(v_0^q),\ i=0,1\},
\]
where $\mathcal P_\sigma X$ carries the sub-diffeology. Define
\[
d_i^\delta(t)=\sigma_i\bigl((1-\rho(t))v_0^q+\rho(t)|\delta|(v_0^r)\bigr),~
\tau_\delta(h)=d_0^\delta\star h\star\bar d_1^\delta.
\]
The path $d_i^\delta$ transports the $i$-th endpoint from $\sigma_i(v_0^q)$ to $\sigma_i(|\delta|(v_0^r))$ along the image of the affine segment joining the two vertices. Thus $\tau_\delta$ transports both endpoints to the reference vertex $v_0^q$. It is a smooth homotopy equivalence, with homotopy inverse $h\mapsto\bar d_0^\delta\star h\star d_1^\delta$. Lemma \ref{lem: affine Kihara comparison} implies that $S(\tau_\delta)$ is a simplicial homotopy equivalence in the sense of Definition \ref{defn: homotopy equivalence}.

We now trivialize the varying endpoints over each simplex. If $\alpha:\Delta^s\to\Delta^q$ is a simplicial operator, write $|\alpha|$ for its affine realization and define
\[
p_{i,\alpha}^\sigma(z,t)=\sigma_i\bigl((1-\rho(t))v_0^q+\rho(t)|\alpha|(z)\bigr),~ i=0,1.
\]
For fixed $z$, this transports the $i$-th endpoint from $\sigma_i(v_0^q)$ to $\sigma_i(|\alpha|(z))$. An $s$-simplex of $E_\sigma$ is a pair $(\alpha,\gamma)$, where $\gamma$ is a smooth $|\Delta^s|_{\aff}$-family of paths from $\sigma_0\alpha$ to $\sigma_1\alpha$. The formulas
\[
T_\sigma(\alpha,\gamma)=(\alpha,p_{0,\alpha}^\sigma\star\gamma\star\bar p_{1,\alpha}^\sigma),~
R_\sigma(\alpha,\eta)=(\alpha,\bar p_{0,\alpha}^\sigma\star\eta\star p_{1,\alpha}^\sigma)
\]
define morphisms of simplicial sets over $\Delta^q$
\[
T_\sigma:E_\sigma\longrightarrow\Delta^q\times S(\mathcal P_\sigma X),~
R_\sigma:\Delta^q\times S(\mathcal P_\sigma X)\longrightarrow E_\sigma.
\]
The map $T_\sigma$ replaces the varying endpoints by fixed endpoints, while $R_\sigma$ restores them. Their composites differ from the respective identities only by smooth reparametrizations and inserted backtracking paths $p\star\bar p$ or $\bar p\star p$. Contracting each backtracking path gives contractions smooth and natural in all parameters. Applied pointwise to the $|\Delta^s|_{\aff}$-families above, they define simplicial homotopies $R_\sigma T_\sigma\simeq\id$ and $T_\sigma R_\sigma\simeq\id$. Thus $T_\sigma$ and $R_\sigma$ are simplicial homotopy inverses.

Consider the diagram
\[
\xymatrix{
E_{\sigma\delta}\ar[r]\ar[d]_{T_{\sigma\delta}}&E_\sigma\ar[d]^{T_\sigma}\\
\Delta^r\times S(\mathcal P_{\sigma\delta}X)\ar[r]_{\delta\times S(\tau_\delta)}&\Delta^q\times S(\mathcal P_\sigma X).
}
\]
The vertical morphisms are simplicial homotopy equivalences. Since the standard simplices are simplicially contractible \cite[\texttt{00XA}]{kerodon}, every simplicial operator $\delta:\Delta^r\to\Delta^q$ is a simplicial homotopy equivalence. Hence the lower horizontal morphism is also a simplicial homotopy equivalence by Lemma \ref{lem: homotopy properties}. Finally, convex interpolation in $|\Delta^q|_{\aff}$ gives, naturally in $\alpha:\Delta^s\to\Delta^r$, an endpoint-preserving smooth homotopy
\[
p_{i,\delta\alpha}^\sigma\simeq d_i^\delta\star p_{i,\alpha}^{\sigma\delta}.
\]
Hence the diagram commutes after passing to simplicial homotopy classes. Since the other three morphisms are invertible there, the two-out-of-three property from Lemma \ref{lem: homotopy properties} shows that its upper horizontal morphism is a simplicial homotopy equivalence, and hence a weak homotopy equivalence.
\epf

We next show that for affine complex, the path model from Definition \ref{defn: homotopy fiber product} computes the corresponding homotopy fiber product, even without the Kan assumption. We call it the \emph{affine path model}.

\blem\label{lem: smooth mapping path}
Let $A\to X\leftarrow B$ be smooth maps of diffeological spaces, let $\mathcal PX=C^\infty([0,1],X)$ carry the functional diffeology, and let $e_i:\mathcal PX\to X$ be endpoint evaluation. Then the natural comparison map
\[
S^{\aff}_\bu\bigl(A\times_{X,e_0}\mathcal PX\times_{e_1,X}B\bigr)
\longrightarrow
S^{\aff}_\bu(A)\times^h_{S^{\aff}_\bu(X)}S^{\aff}_\bu(B)
\]
is a weak homotopy equivalence.
\elem
\bpf
Write $S=S^{\aff}_\bu$ and set
\[
A_\bu=S(A),~ B_\bu=S(B),~ X_\bu=S(X),~ P_\bu=S(\mathcal PX),~ e_\bu=S(e_0,e_1):P_\bu\longrightarrow X_\bu^2.
\]
By Lemma \ref{lem: affine Kihara comparison}, $S$ preserves products and pullbacks. Hence
\[
S\bigl(A\times_{X,e_0}\mathcal PX\times_{e_1,X}B\bigr)=(A_\bu\times B_\bu)\times_{X_\bu^2}P_\bu.
\]
Let $c:X\to\mathcal PX$ be the constant-path map and choose a smooth map $\lambda:\R\to[0,1]$ which vanishes on $(-\infty,0]$ and equals $1$ on $[1,\infty)$. By the functional diffeology in Definition \ref{defn: diffeological space}, the formula
\[
H:\mathcal PX\times\R\longrightarrow\mathcal PX,~H(h,s)(t)=h(\lambda(s)t),
\]
defines a smooth map. Its restriction to $\mathcal PX\times[0,1]$ is a smooth homotopy from $c\circ e_0$ to $\id_{\mathcal PX}$, while $e_0\circ c=\id_X$. Thus $c$ is a smooth homotopy equivalence in the sense of Definition \ref{defn: smooth homotopy equivalence}. Lemma \ref{lem: affine Kihara comparison} implies that
\[
c_\bu=S(c):X_\bu\longrightarrow P_\bu
\]
is a simplicial homotopy equivalence. 

Suppose $e_\bu=S(e_0,e_1):P_\bu\to X_\bu^2$ is from Lemma \ref{lem: affine endpoint evaluation}. Since $e_\bu c_\bu$ is the diagonal $X_\bu\to X_\bu^2$, the morphism $c_\bu$ is compatible with the maps to $X_\bu^2$.

By Lemma \ref{lem: affine endpoint evaluation}, $e_\bu$ satisfies the hypothesis of Lemma \ref{lem: Rezk pullback criterion}. Hence
\[
(A_\bu\times B_\bu)\times_{X_\bu^2}P_\bu\simeq(A_\bu\times B_\bu)\times_{X_\bu^2}^hP_\bu.
\]
Lemma \ref{lem: homotopy invariance fiber product}, applied to the weak homotopy equivalence $c_\bu$ over $X_\bu^2$, gives
\[
(A_\bu\times B_\bu)\times_{X_\bu^2}^hP_\bu\simeq(A_\bu\times B_\bu)\times_{X_\bu^2}^hX_\bu\simeq A_\bu\times_{X_\bu}^hB_\bu.
\]
This proves the assertion.
\epf

\brem\label{rem: affine path limitation}
We do not claim that $e_\bu=S^\aff_\bu(e_0,e_1)$ is a Kan fibration, which would be the Christensen--Wu fibration condition \cite[Definition 4.8]{CW2014diff}; compare \cite[Remark 5.6]{Kuribayashi2020simplicial}. Lemma \ref{lem: affine endpoint evaluation} proves only the sharpness needed for the affine path models above to compute homotopy fiber products.
\erem
For the reader's convenience, we repeat the statement of Proposition \ref{prop: segal space}.
\begin{repproposition}{prop: segal space}
Both $(\PB_n^{V,\aff})_{\bu,\bu}$ and $(\PB_n^{V,\kih})_{\bu,\bu}$ satisfy the Segal condition. Thus $(\PB_n^{V,\aff})_{\bu,\bu}$ is a Segal simplicial set in the sense of Definition \ref{defn: segal space}, while $(\PB_n^{V,\kih})_{\bu,\bu}$ is a Segal space.
\end{repproposition}
\bpf
The proof follows the proof of \cite[Proposition 5.19]{CS2019category}, specialized to its single composition direction. For simplicity, write
\[
X_k^{\mmod}=(\PB_n^{V,\mmod})_k=S^{\mmod}_\bu((\fPB_n^V)_k),~\mmod\in\{\aff,\kih\}.
\]
We first prove the Segal condition for the affine model. For $m,l\in\N$ and $k=m+l$, consider the Segal map
\begin{equation}\label{eq: affine Segal map}
\ga_{m,l}^{\aff}:X_k^{\aff}\longrightarrow X_m^{\aff}\times_{X_0^{\aff}}^hX_l^{\aff}.
\end{equation}
If $m=0$ or $l=0$, one of the two maps to $X_0^{\aff}$ is the identity and hence a Kan fibration. Definition \ref{defn: homotopy fiber product} therefore identifies the homotopy fiber product with the strict pullback. We henceforth assume $m,l>0$ and fix $k=m+l$.

By Lemma \ref{lem: smooth mapping path}, it suffices to replace the homotopy fiber product in \eqref{eq: affine Segal map} by its affine path model, denoted by
\[
X_m^{\aff}\times_{X_0^{\aff}}^{\mathrm{path}}X_l^{\aff}.
\]
A simplex of this model consists of $(M,\un I)$ and $(N,\un J)$ together with a smooth path $h$ from the target of the former to the source of the latter, all varying smoothly over an affine simplex. Cutting defines
\[
\ga_{m,l}^{\aff,\mathrm{path}}:X_k^{\aff}\longrightarrow X_m^{\aff}\times_{X_0^{\aff}}^{\mathrm{path}}X_l^{\aff},
\]
whose composite with the comparison
\begin{equation}\label{eq: affine path comparison}
X_m^{\aff}\times_{X_0^{\aff}}^{\mathrm{path}}X_l^{\aff}\longrightarrow X_m^{\aff}\times_{X_0^{\aff}}^hX_l^{\aff}
\end{equation}
is the Segal map.

As in the proof of \cite[Proposition 5.19]{CS2019category}, the canonical map
\[
\Int_m^{\aff}\times_{\Int_0^{\aff}}\Int_l^{\aff}\longrightarrow\Int_m^{\aff}\times_{\Int_0^{\aff}}^{\mathrm{path}}\Int_l^{\aff}
\]
is a weak homotopy equivalence. Indeed, its source is the affine complex of the nonempty convex diffeological space $\fInt_m\times_{\fInt_0}\fInt_l$, while its target is weakly contractible by Proposition \ref{prop: Int constant} and Lemmas \ref{lem: smooth mapping path} and \ref{lem: homotopy invariance fiber product}. Let $(X^{\aff})^{m,l}_\bu$ be the strict pullback
\begin{equation}\label{eq: pullback}
\xymatrix{
(X^{\aff})^{m,l}_\bu\ar[r]\ar[d]&\Int_m^{\aff}\times_{\Int_0^{\aff}}\Int_l^{\aff}\ar[d]^{\simeq}\\
X_m^{\aff}\times_{X_0^{\aff}}^{\mathrm{path}}X_l^{\aff}\ar[r]&\Int_m^{\aff}\times_{\Int_0^{\aff}}^{\mathrm{path}}\Int_l^{\aff}.
}
\end{equation}
Its simplices are precisely the triples for which $I_m=J_0$ and the interval component of $h$ is constant. The map \eqref{eq: affine Segal map} is represented by the composite
\begin{equation}\label{eq: factor through}
\xymatrix{X_k^{\aff}\ar[rr]^-{\ga_{m,l}^{\aff,\mathrm{path}}}\ar[d]_{\ga'_{m,l}}&&X_m^{\aff}\times_{X_0^{\aff}}^{\mathrm{path}}X_l^{\aff}\\(X^{\aff})^{m,l}_\bu\ar@{^{(}->}[urr]}
\end{equation}
where $\ga'_{m,l}$ is cutting. We show that the vertical and diagonal maps are simplicial homotopy equivalences by gluing and rescaling, respectively. The additional intermediate object used in the proof of \cite[Proposition 5.19]{CS2019category} for the other composition directions is unnecessary here.

\emph{Gluing.} Fix a smoothed diagonal $D\subset[0,1]_{\mathrm{path}}\times[0,1]_{\mathrm{time}}$, given by the graph of a smooth map $\zeta:[0,1]\to[0,1]$ which is $0$ near $0$, is $1$ near $1$, and is a diffeomorphism between the remaining open intervals. For $\tau\in[0,1]$, set
\[
\zeta_\tau^s=\tau\zeta,~\zeta_\tau^t=1+\tau(\zeta-1).
\]
Given $((M,\un I),(N,\un J),h)$, let $P$ be the smooth family of embedded manifolds corresponding to $h$. Pulling $P$ back to the rescaled diagonal and gluing along the collars gives an embedded bordism
\[
M\cup_{P_{\mathrm{diag}}}N
\]
with interval tuple $(I_0,\ldots,I_m=J_0,J_1,\ldots,J_l)$. The union of the two horizontal collar segments with the diagonal projects diffeomorphically to the common time interval, so projection to $V\times\R$ remains an embedding. The construction is smooth in affine families and commutes with the face and degeneracy maps, and hence defines
\[
\glue:(X^{\aff})^{m,l}_\bu\longrightarrow X_k^{\aff}.
\]
The composite $\glue\circ\ga'_{m,l}$ is the identity. Pulling $P$ back along the graphs of $\zeta_\tau^s$ and $\zeta_\tau^t$, gluing the two pullbacks to $M$ and $N$, and retaining the part of $P$ between the two graphs gives a smooth homotopy from the identity to $\ga'_{m,l}\circ\glue$. By Lemma \ref{lem: affine Kihara comparison}, this smooth homotopy exhibits $\ga'_{m,l}$ and $\glue$ as simplicial homotopy inverses.

\emph{Rescaling.} For a $q$-simplex $((M,\un I),(N,\un J),h)$, write the interval component of $h(s,u)$ as $[c(s,u),d(s,u)]$, where $(s,u)\in|\Delta^q|_{\aff}\times[0,1]$.
For $\tau\in[0,1]$, set
\[
c_\tau(s,u)=(1-\tau)c(s,u)+\tau c(s,0),~d_\tau(s,u)=(1-\tau)d(s,u)+\tau d(s,0),
\]
and
\[
f_{\tau,s,u}(x)=\frac{d_\tau(s,u)-c_\tau(s,u)}{d(s,u)-c(s,u)}(x-c(s,u))+c_\tau(s,u).
\]
Apply $\id_V\times f_{\tau,s,u}$ to $h(s,u)$ and $\id_V\times f_{\tau,s,1}$ to $(N,\un J)$. This gives a smooth deformation retraction onto the subcomplex for which the interval component of $h$ is constant and the endpoint intervals agree. It commutes with the face and degeneracy maps, so Lemma \ref{lem: affine Kihara comparison} makes the inclusion in \eqref{eq: factor through} a simplicial homotopy equivalence.

The gluing and rescaling constructions show that $\ga_{m,l}^{\aff,\mathrm{path}}$ is a simplicial homotopy equivalence. The comparison \eqref{eq: affine path comparison} is a weak homotopy equivalence by Lemma \ref{lem: smooth mapping path}. The two-out-of-three property from Lemma \ref{lem: homotopy properties} therefore shows that the Segal map \eqref{eq: affine Segal map} is a weak homotopy equivalence. Thus $X_\bu^{\aff}$ satisfies the Segal condition.

Lemma \ref{lem: homotopy invariance fiber product}, the two-out-of-three property from Lemma \ref{lem: homotopy properties}, and Lemma \ref{lem: affine Kihara comparison} transfer the Segal condition to $X_\bu^{\kih}$. Its levels are Kan by the latter, so it is a Segal space.
\epf
\brem
For the general $d$-fold construction, fix one simplicial direction and all remaining indices. Lemma \ref{lem: smooth mapping path} replaces the use of \cite[Proposition 5.11]{CS2019category} in the proof of \cite[Proposition 5.19]{CS2019category}; the remaining cutting, gluing, rescaling, and essential-constancy arguments apply without change. Lemma \ref{lem: affine Kihara comparison} then transfers the result to the Kihara model.
\erem
\newpage
\section{Proof of the symmetric monoidal condition}\label{sec: symmetric monoidal}
Following \cite[Proposition 7.2]{CS2019category}, this appendix gives the proof of Proposition \ref{prop: symmetric condition} and, for the reader's convenience, repeats its statement.
\begin{repproposition}{prop: symmetric condition}
The stabilized Kihara $\Ga$-object $A^\kih$ satisfies the $\Ga$-condition in the sense of Definition \ref{defn: sym str Ga object}. Its objectwise completion endows $\Bord_n$ with a symmetric monoidal structure whose unit is represented by $(\varnothing,\I)$.
\end{repproposition}
\bpf
We first verify the $\Ga$-condition for the affine model. For $m=0$, Proposition \ref{prop: Int constant} shows that $A^\aff\langle0\rangle$ is levelwise weakly contractible. Let $m>0$ and fix a simplicial level $k$. The canonical map
\[
\prod_{\beta=1}^m\ga_\beta:A^\aff\langle m\rangle_k\longrightarrow A^\aff\langle1\rangle_k^m
\]
forgets the common-interval and disjointness conditions. Let $P_k^{m,V}$ be the affine complex of $m$-tuples in a fixed ambient space $V$ with common interval data, and set $P_k^m=\hocolim_VP_k^{m,V}$. The common-interval comparison in the proof of \cite[Proposition 7.2]{CS2019category}, using the interval adjustments from the proof of \cite[Proposition 5.19]{CS2019category}, applies to these smooth families and gives a weak homotopy equivalence
\[
P_k^m\longrightarrow A^\aff\langle1\rangle_k^m.
\]
Here finite products commute with filtered homotopy colimits of spaces \cite[Proposition 5.3.3.3]{HTTlurie2009}.

It remains to separate the embedded manifolds. For a tuple represented in $V$, consider
\[
\iota_{\alpha,s}:V\longrightarrow\R\oplus V,~v\longmapsto(s\alpha,v),~1\leq\alpha\leq m.
\]
For $s>0$, the images of the $m$ manifolds are pairwise disjoint, while $s=0$ gives the stabilized original tuple. If the original tuple is already disjoint, it remains so throughout the deformation. After stabilization, this gives a simplicial homotopy inverse to
\[
A^\aff\langle m\rangle_k\longrightarrow P_k^m.
\]
Thus $A^\aff$ satisfies the $\Ga$-condition.

Lemma \ref{lem: affine Kihara comparison}, homotopy invariance of filtered homotopy colimits, and the finite-product and two-out-of-three properties from Lemma \ref{lem: homotopy properties} transfer the $\Ga$-condition to $A^\kih$. Finally, objectwise completion preserves the $\Ga$-condition by \cite[Lemma 3.7]{CS2019category}. Hence it endows $\Bord_n$ with the asserted symmetric monoidal structure, whose unit is represented by $(\varnothing,\I)$.
\epf

\newpage
\section{Triangle detection lemma in stable \texorpdfstring{$\infty$}{infty}-categories}\label{sec: Triangle detection lemma}
\centerline{Longke Tang}

\medskip
In this appendix, we provide a version of the triangle detection lemma in stable $\infty$-categories. The original version in the triangulated derived category was developed by Ozsv\'{a}th--Szab\'{o} \cite[Lemma 4.2]{ozsvath2005double}; see also \cite[Lemma 7.1]{kronheimer2011khovanov} and \cite[Lemma 5.1]{scaduto2015instantons} for the version over $\Z$. In a stable $\infty$-category, let $X[1]$ denote the suspension of an object $X$.


\blem\label{five-lemma}
    Let $\cC$ be a stable $\infty$-category and let
    $$\begin{tikzcd}
        X_1\ar[r]\ar[d]&0\ar[d]&\\
        X_2\ar[r]\ar[d]&X_3\ar[r]\ar[d]&0\ar[d]\\
        0\ar[r]&X_4\ar[r]&X_5
    \end{tikzcd}$$
    be a diagram in $\cC$. Let $C=\cofib(X_1\to X_2)$ and $F=\fib(X_4\to X_5)$. Define two maps from $C$ to $F$ as follows:
    \begin{itemize}
        \item By definition of cofiber as the pushout along $0$, the upper left square gives a map $C\to X_3$. Similarly, the lower right square gives a map $X_3\to F$. They compose to a map $C\to F$.
        \item Ignore $X_3$. The big left rectangle gives a map $X_1[1]\to X_4$ and the big lower rectangle gives a map $X_2[1]\to X_5$. Note that these maps are compatible with the maps $X_1\to X_2$ and $X_4\to X_5$ in the diagram, so taking fibers gives a map $C=\fib(X_1[1]\to X_2[1])\to\fib(X_4\to X_5)=F$.
    \end{itemize}
    The two maps $C\to F$ are canonically equal.
\elem
\bpf
    By definition of pushouts, among those diagrams of the above shape with $X_1\to X_2$ as the upper left arrow,
    $$\begin{tikzcd}
        X_1\ar[r]\ar[d]&0\ar[d]&\\
        X_2\ar[r]\ar[d]&C\ar[r]\ar[d]&0\ar[d]\\
        0\ar[r]&X_1[1]\ar[r]&X_2[1]
    \end{tikzcd}$$
    is the initial one, as all its three squares are pushouts. Similarly,
    $$\begin{tikzcd}
        X_4[-1]\ar[r]\ar[d]&0\ar[d]&\\
        X_5[-1]\ar[r]\ar[d]&F\ar[r]\ar[d]&0\ar[d]\\
        0\ar[r]&X_4\ar[r]&X_5
    \end{tikzcd}$$
    is the final one with $X_4\to X_5$ as the lower right arrow. Therefore, we have maps of diagrams
    $$\begin{tikzcd}
        X_1\ar[r]\ar[d,equal]&X_2\ar[r]\ar[d,equal]&C\ar[r]\ar[d]&X_1[1]\ar[r]\ar[d]&X_2[1]\ar[d]\\
        X_1\ar[r]\ar[d]&X_2\ar[r]\ar[d]&X_3\ar[r]\ar[d]&X_4\ar[r]\ar[d,equal]&X_5\ar[d,equal]\\
        X_4[-1]\ar[r]&X_5[-1]\ar[r]&F\ar[r]&X_4\ar[r]&X_5
    \end{tikzcd}$$
    where we omit all the zeros. Now the middle vertical map $C\to F$ is on the one hand the composite $C\to X_3\to F$, and on the other hand the fiber of $X_1[1]\to X_4$ mapping to $X_2[1]\to X_5$, which proves the lemma.
\epf

\blem[Triangle detection lemma]\label{six-lemma}
    Let $\cC$ be a stable $\infty$-category and let
    $$\begin{tikzcd}
        X_0\ar[r]\ar[d]&X_1\ar[r]\ar[d]&0\ar[d]&\\
        0\ar[r]&X_2\ar[r]\ar[d]&X_3\ar[r]\ar[d]&0\ar[d]\\
        &0\ar[r]&X_4\ar[r]&X_5
    \end{tikzcd}$$
    be a diagram in $\cC$. Assume that the maps $X_0[1]\to X_3$, $X_1[1]\to X_4$, and $X_2[1]\to X_5$ obtained from the diagram are isomorphisms. Then all the squares in the diagram are fiber squares.
\elem
\bpf
    Let $C=\cofib(X_1\to X_2)$ and $F=\fib(X_4\to X_5)$. Then we obtain maps
    $$X_0[1]\to C\to X_3\to F$$
    from the diagram. By assumption, the composite $X_0[1]\to X_3$ is an isomorphism. By Lemma \ref{five-lemma} and the assumption that $X_1[1]\to X_4$ and $X_2[1]\to X_5$ are isomorphisms, the composite $C\to F$ is also an isomorphism. Therefore, $C\to X_3$ is an isomorphism, being a map with both left and right inverses. Then $X_0[1]\to C$ is also an isomorphism. From these, we see that the two upper squares are fiber squares. The same reasoning in $\cC^{\opp}$ with the diagram flipped shows that the two lower squares are also fiber squares.
\epf


\bibliographystyle{alpha}

\end{document}